\documentclass[nonblindrev]{informs3noheader}

\OneAndAHalfSpacedXI

\usepackage{mathtools}
\usepackage[inline]{enumitem}
\usepackage{subfig}
\usepackage{booktabs}
\usepackage{array}
\usepackage{multirow}
\usepackage{booktabs,caption,fixltx2e}
\usepackage[flushleft]{threeparttable}
\usepackage{lscape}
\usepackage{makecell}
\usepackage{amssymb}
\usepackage[mathscr]{euscript}
\usepackage{algorithm}
\usepackage{algpseudocode}
\usepackage{pifont}
\usepackage{caption}
\usepackage{algorithmicx}
\usepackage{hyperref}
\usepackage{mathrsfs}
\usepackage{tabularx}
\usepackage{bbm}
\usepackage{soul}
\usepackage{blindtext}
\usepackage{longtable}
\hypersetup{
colorlinks=true,
linkcolor=black,
citecolor=black,
filecolor=black,
urlcolor=black,
}

\usepackage{comment}

\usepackage{tikz}
\usetikzlibrary{shapes.geometric}
\usetikzlibrary{positioning}
\usetikzlibrary{decorations.pathreplacing}
\definecolor{mygreen}    {RGB}{0,90,0}
\definecolor{myblue}     {RGB}{0,51,140}
\definecolor{myorange}   {RGB}{238,118,0}
\definecolor{myred}      {RGB}{126,0,0}
\definecolor{mygray}     {RGB}{100,100,105}
\definecolor{mygrayblue} {RGB}{0,128,128}
\definecolor{mygraygreen}{RGB}{128,128,0}
\definecolor{DarkPurple}     {RGB}{142, 36, 170}
\definecolor{LightPurple}    {RGB}{57, 130, 7}

\usepackage{natbib}
\bibpunct[, ]{(}{)}{,}{a}{}{,}%
\def\bibfont{\small}%

\usepackage[utf8]{inputenc}
\usepackage[english]{babel}
\usepackage{graphicx}

\newcolumntype{F}[1]{%
    >{\raggedright\arraybackslash\hspace{0pt}}p{#1}}%
\newcolumntype{T}[1]{%
    >{\centering\arraybackslash\hspace{0pt}}p{#1}}%
\newcolumntype{H}{>{\setbox0=\hbox\bgroup}c<{\egroup}@{}}

\TheoremsNumberedThrough

\EquationsNumberedThrough
{\theoremstyle{EXkey}}

\MANUSCRIPTNO{}

\def\R{\mathbb{R}}
\def\Z{\mathbb{Z}}

\def\calA{\mathcal{A}}
\def\calB{\mathcal{B}}
\def\calC{\mathcal{C}}
\def\calD{\mathcal{D}}
\def\calE{\mathcal{E}}
\def\calF{\mathcal{F}}

\def\calH{\mathcal{H}}
\def\calI{\mathcal{I}}
\def\calJ{\mathcal{J}}
\def\calK{\mathcal{K}}
\def\calL{\mathcal{L}}
\def\calM{\mathcal{M}}
\def\calN{\mathcal{N}}
\def\calO{\mathcal{O}}
\def\calP{\mathcal{P}}
\def\calQ{\mathcal{Q}}
\def\calR{\mathcal{R}}
\def\calS{\mathcal{S}}
\def\calT{\mathcal{T}}
\def\calU{\mathcal{U}}
\def\calV{\mathcal{V}}
\def\calW{\mathcal{W}}
\def\calX{\mathcal{X}}
\def\calY{\mathcal{Y}}

\def\bA{\boldsymbol{A}}

\def\bD{\boldsymbol{D}}
\def\bE{\boldsymbol{E}}

\def\bb{\boldsymbol{b}}
\def\bc{\boldsymbol{c}}

\def\bff{\boldsymbol{f}}
\def\bg{\boldsymbol{g}}
\def\bh{\boldsymbol{h}}
\def\bi{\boldsymbol{i}}

\def\bq{\boldsymbol{q}}

\def\bt{\boldsymbol{t}}

\def\bv{\boldsymbol{v}}
\def\bw{\boldsymbol{w}}
\def\bx{\boldsymbol{x}}
\def\by{\boldsymbol{y}}
\def\bz{\boldsymbol{z}}

\def\bo{\boldsymbol{0}}

\def\bh{\boldsymbol{h}}

\def\bbeta{\boldsymbol{\beta}}
\def\bgamma{\boldsymbol{\gamma}}

\def\bnu{\boldsymbol{\nu}}

\def\btheta{\boldsymbol{\theta}}

\def\bxi{\boldsymbol{\xi}}

\def\bpsi{\boldsymbol{\psi}}
\def\bzeta{\boldsymbol{\zeta}}

\def\bxi{\boldsymbol{\xi}}

\def\st{\text{s.t.}}

\definecolor{mygreen}    {RGB}{0,90,0}
\definecolor{myblue}     {RGB}{0,51,140}
\definecolor{myorange}   {RGB}{238,118,0}
\definecolor{myred}      {RGB}{126,0,0}
\definecolor{mygray}     {RGB}{100,100,105}
\definecolor{mygrayblue} {RGB}{0,128,128}
\definecolor{mygraygreen}{RGB}{128,128,0}
\definecolor{DarkPurple}     {RGB}{142, 36, 170}
\definecolor{LightPurple}    {RGB}{57, 130, 7}

\begin{document}

\ARTICLEAUTHORS{
	\AUTHOR{Bernardo Martin-Iradi, Alexandria Schmid, Kayla Cummings, Alexandre Jacquillat}
}	

\ARTICLEAUTHORS{
        \AUTHOR{Bernardo Martin-Iradi}
	\AFF{Institute for Transport Planning and Systems, ETH Zurich, Zurich,
Switzerland}
	\AUTHOR{Alexandria Schmid, Kayla Cummings, Alexandre Jacquillat}
	\AFF{Sloan School of Management and Operations Research Center, MIT, Cambridge, MA}
}	

\RUNAUTHOR{Martin-Iradi, Schmid, Cummings, Jacquillat}
\RUNTITLE{}

\TITLE{\Large A Double Decomposition Algorithm for Network Planning and Operations in Deviated Fixed-route Microtransit}

\ABSTRACT{
Microtransit offers opportunities to enhance urban mobility by combining the reliability of public transit and the flexibility of ride-sharing. This paper optimizes the design and operations of a deviated fixed-route microtransit system that relies on reference lines but can deviate on demand in response to passenger requests. We formulate a \textit{Microtransit Network Design (MiND)} model via two-stage stochastic integer optimization, with a first-stage network design and service scheduling structure and a second-stage vehicle routing structure. We derive a tight second-stage relaxation using a subpath-based representation of microtransit operations in a load-expanded network. We develop a double-decomposition algorithm combining Benders decomposition and subpath-based column generation. We prove that the algorithm maintains a valid optimality gap and converges to an optimal solution in a finite number of iterations. Results obtained with real-world data from Manhattan show that the methodology scales to large and otherwise-intractable instances, with up to 10-100 candidate lines and hundreds of stops. Comparisons with transit and ride-sharing suggest that microtransit can provide win-win outcomes toward efficient mobility (high demand coverage, low costs, high level of service), equitable mobility (broad geographic reach) and sustainable mobility (limited environmental footprint). We provide an open-source implementation to enable replication.
}


\maketitle

\vspace{-12pt}
\section{Introduction}

Major cities face critical challenges to meet mobility needs in the midst of rising congestion, greenhouse gas emissions and socioeconomic inequalities. Static transit infrastructure offers limited flexibility to respond to ever-changing mobility needs, resulting in a ridership decline \citep{economist} and transit deserts \citep{allen2017lost}. Simultaneously, ride-sharing provides flexible, on-demand mobility services, but low-occupancy vehicles still lead to high fares, congestion, and emissions. This context identifies opportunities to leverage emerging \textit{microtransit} services toward efficient, equitable, and sustainable mobility. Broadly defined by the \cite{dotmicrotransit}  as ``privately owned and operated shared transportation system(s) that can offer fixed routes and schedules, as well as flexible routes and on-demand scheduling,'' microtransit shepherds the digital capabilities and operating flexibility of ride-sharing into the realm of public transit. Yet, microtransit raises critical questions about \textit{how} to combine transit and ride-sharing components into low-cost, high-quality services and how to develop dedicated analytics capabilities \citep{mck_analytics}.

This paper develops a two-stage stochastic integer optimization methodology to support the design and operations of \textit{deviated fixed-route microtransit}---a hybrid microtransit system based on transit lines to consolidate passenger demand into high-capacity vehicles and on-demand deviations in response to passenger requests. Our model optimizes network design and service scheduling under demand uncertainty in the first stage, and on-demand operations in the second stage. We propose a scalable methodology, relying on (i) a network-based second-stage formulation with a tight linear relaxation but an exponential number of subpath-based variables; and (ii) a double-decomposition algorithm combining Benders decomposition and subpath-based column generation. The objectives of the paper are to establish its scalability in large and practical problem instances, and to assess the performance of deviated fixed-route microtransit in the urban mobility ecosystem.

Our first contribution is to formulate a \textit{Microtransit Network Design (MiND)} model via two-stage stochastic optimization (Section~\ref{sec:model}). The first-stage problem optimizes network design and service scheduling by selecting line-based \textit{reference trips}. The second-stage problem reflects on-demand routing operations in response to passenger requests, under vehicle capacity and time window constraints; in particular, microtransit operations must visit checkpoints on the reference lines at designated times. The objective combines minimizing planning costs, maximizing ridership and maximizing passenger level of service. For simplicity, we focus primarily on a MiND-VRP problem, corresponding to a vehicle routing setting in which all passengers have the same origin or the same destination---motivated by use cases such as airport or university shuttles. In~\ref{app:DAR}, we extend the methodology and main results to a MiND-DAR problem, corresponding to a dial-a-ride setting in which passengers request transportation from origin to destination; and in~\ref{app:Tr}, we extend them to a MiND-Tr problem that allows transfers between microtransit lines.

One of the main complexities of the problem lies in its discrete second-stage structure---a capacitated vehicle routing problem with time windows. To retain a tight recourse formulation, we propose a subpath-based representation of second-stage microtransit operations in a load-expanded network. Each node encodes a checkpoint on the reference line and a vehicle load, and each arc characterizes on-demand operations between checkpoints. Subpaths encapsulate time window constraints, and load expansion accommodates vehicle capacities without big-$M$ constraints, enabling a continuous recourse function approximation. We show that our subpath-based variables enable a more effective formulation than a compact formulation (tighter second-stage formulation), than a segment-based benchmark with variables connecting consecutive stops in a granular time-load-expanded network (sparser network) and than a path-based benchmark with variables connecting the start to the end of each line (much slower rate of exponential growth in the number of variables).

Our second contribution is a scalable double-decomposition (DD) algorithm combining Benders decomposition and subpath-based column generation (Section~\ref{sec:alg}). The Benders decomposition scheme iterates between a first-stage network design problem and second-stage routing problems, exploiting the nested block-diagonal structure to decompose on-demand operations for each reference trip in each scenario. The column generation scheme adds subpath-based variables iteratively in the Benders subproblem. We develop exact and heuristic label-setting algorithms to generate subpaths of negative reduced cost while keeping track of vehicle load and level of service. Our methodology induces a double-decomposition structure: the column generation pricing problem adds subpath-based arcs into load-expanded networks (i.e., local on-demand deviations between checkpoints); the Benders subproblem combines them into full second-stage paths (i.e., full microtransit trips for each reference trip in each scenario), and the Benders master problem optimizes first-stage planning decisions accordingly (i.e., network design and service scheduling). This algorithm converges to the partial relaxation of the (tight) subpath-based formulation with mixed-integer first-stage variables and continuous second-stage variables. To guarantee convergence and second-stage integrality, we augment the DD scheme with integer L-shaped cuts from \cite{laporte1993integer} (DD\&ILS) and the unified branch-and-Benders-cut (UB\&BC) algorithm from \cite{maheo2024unified} (UB\&DD). Ultimately, the algorithm yields a finite and exact double decomposition methodology for the two-stage stochastic integer optimization problem.

Our third contribution is to demonstrate the scalability of our methodology to large and practical MiND instances (Section~\ref{sec:comp}). We develop a real-world setup in Manhattan using data from the \cite{taxi}. Results show the benefits of our subpath-based formulation and our double decomposition methodology. Specifically, the DD methodology yields certifiably optimal, or near-optimal solutions to the problem, and the broader DD\&ILS and UB\&DD algorithms converge to an exact solution. Altogether, our methodology scales to large and otherwise intractable instances of the MiND-VRP, of the size of the full Manhattan network with up to 100 candidate lines, hundreds of stations, thousands of passenger requests, and 5-20 demand scenarios; and it scales to instances with 10--40 candidate lines, over a hundred stations, and thousands of passenger requests for the MiND-DAR and MiND-Tr extensions. Our results also show the practical benefits of our stochastic optimization methodology, with a value of the stochastic solution of 5--7\% against a deterministic benchmark. Ultimately, this paper contributes the first integrated methodology to optimize microtransit design and operations under uncertainty, and a methodology that scales to much larger instances than previous approaches in microtransit operations. 

Our final contribution is to derive evidence that deviated fixed-route microtransit can provide win-win mobility outcomes (Section~\ref{sec:practical}). As compared to ride-sharing, microtransit consolidates demand into high-capacity vehicles along reference lines. As compared to fixed-route transit, it increases demand coverage and improve passenger level of service by leveraging on-demand flexibility. In turn, the optimized microtransit network has a broader catchment area than its fixed-route counterpart, especially in otherwise unserved regions. Finally, demand consolidation and high coverage result in a significant decrease in distance traveled per passenger, with cost benefits and environmental benefits. Since Manhattan represents a high-density region, these results can be seen as conservative estimates of the impact of microtransit in lower-density areas with fewer transit alternatives. Altogether, deviated fixed-route microtransit can contribute to efficient  mobility (high demand coverage, low costs per passenger, high service levels), equitable mobility (broad geographic reach), and sustainable mobility (limited environmental footprint). These results have inspired ongoing collaborations toward the pilot deployment of new microtransit solutions.
\section{Background, motivation and literature review}

\textit{Practical.} Microtransit seeks to combine the efficiencies of public transit with the flexibility of ride-sharing. Our first-stage problem relates to transit planning \citep{desaulniers2007public,ortega2018line,wei2022transit,sun2022distributionally}. These problems have been solved with heuristics \citep{ceder1986bus,walteros2015hybrid} and exact methods in small instances with 10-25 stops \citep{marin2009urban}. \cite{bertsimas2021data} developed a column generation methodology for transit network design that scales to large instances with hundreds of stops.

In microtransit, one possible operating model is to design a joint system combining fixed-route transit and ride-sharing, which \cite{chopra2023mobility} framed via dual sourcing. Another model is to provide on-demand transportation with high-capacity vehicles, which \cite{alonso2017demand} optimized in a request-trip-vehicle network. However, on-demand high-capacity operations may induce detours and delays. \cite{blanchard2022probabilistic} showed that the optimal latency in the traveling repairman problem grows with the size and dispersion of the geographic area, and grows at a supra-linear rate of $\Theta(n\sqrt{n})$ where $n$ is the number of customers. This convex function reflects negative spatial externalities across customers induced by on-demand operations with high-capacity vehicles---that is, on-demand deviations become more costly with more passengers onboard.

This theoretical result outlines two approaches to alleviate spatiotemporal externalities in on-demand mobility: restricting vehicle occupancy---as in ride-sharing---or operating in small or concentrated areas. In practice, high-capacity microtransit has been successful in small municipalities\footnote{See, e.g., city.ridewithvia.com/salem-skipper, city.ridewithvia.com/newmo-newton} and university campuses.\footnote{See, e.g., ridewithvia.com/news/northeastern-university-taps-via-to-power-new-on-demand-safety-shuttle} In larger regions, zone-based microtransit operates in limited geographic locales; for instance, MetroConnect operates in 12 areas of Miami, and Metro Micro operates in eight areas of Los Angeles. It also acts as a first- and last-mile feeder into fixed-route transit \citep{steiner2020strategic,banerjee2021real,silva2022demand,guan2023path,cummings2023multimodal}.\footnote{city.ridewithvia.com/go-connect-miami, micro.metro.net, www.dart.org/guide/transit-and-use} In practice, multimodal microtransit introduces complexity to establish first- and last-mile zones; for instance, DART's GoLink service in Dallas partitions the service region into 34 zones. Such partitioning raises similar trade-offs as in door-to-door microtransit, between high costs with low-occupancy vehicles in small zones vs. detours and delays with larger vehicles in larger zones.

Deviated fixed-route microtransit, in contrast, consolidates demand into high-occupancy vehicles along transit routes while allowing on-demand deviations in response to passenger requests.\footnote{www.nationalrtap.org/Toolkits/ADA-Toolkit/Service-Type-Requirements/Route-Deviation-Requirements; https://www.rideuta.com/Services/Flex-Routes; https://www.tricountytransit.org/flex-routes.html} This model leverages \textit{virtual bus stops} to consolidate pickups and dropoffs in central locations. On-demand operations with virtual bus stops induce challenging routing problems \citep{zhang2023routing}. Viewed through this lens, deviated fixed-route microtransit leverages transit lines as a natural regularization, while allowing on-demand deviations with virtual bus stops.

Deviated fixed-route microtransit has been subject to limited research. \cite{quadrifoglio2007insertion,quadrifoglio2008mobility} optimized on-demand deviations with a single vehicle. \cite{quadrifoglio2006performance} and \cite{zhao2008service} quantified trade-offs between frequencies, deviations, and service levels. \cite{galarza2022column} optimized operations in a related system in which transit vehicles can skip stops. \cite{lqm2021} formulated a mixed-integer linear optimization model to optimize on-demand deviations with autonomous vehicles. All these methods focus on the operations alone (our second-stage problem), and scale to small instances with 1-5 vehicles and 10-50 stops.

\textit{Stochastic programming.}
Our problem combines a network design structure and a capacitated vehicle routing structure with time windows, under uncertainty. It is cast as a two-stage stochastic program with discrete recourse, a challenging class of problems \citep{caroe1999dual,sen2006decomposition,gade2014decomposition,zhang2014finitely,kim2015two,bodur2017strengthened,wang2020stochastic}. \cite{maheo2019benders} developed a unified branch-and-Benders-cut (UB\&BC) algorithm as a general-purpose single-tree Benders decomposition algorithm with tailored branching rules, for stochastic mixed-integer programming with any mixed-integer structure in the first stage and in the second stage, and uncertainty in any parameters. In this paper, we leverage a network-based reformulation to retain a tight second-stage representation of routing operations. This approach applies extended formulations principles \citep{conforti2010extended,conforti2014integer}, which have been used in lot sizing \citep{eppen1987solving,ahuja2008solving}, machine scheduling \citep{sousa1992time,pessoa2010exact}, bin packing \citep{valerio1999exact,delorme2020enhanced}, network design \citep{frangioni20090}, unit commitment \citep{queyranne2017tight}, etc. In vehicle routing, a common approach involves modeling routing problems with temporal coordination requirements in time-expanded networks \citep{crainic2016service,lee2020dynamic,agarwal2008ship,liebchen2008first}. In our setting, this representation leads to a very large formulation in a granular time-space-load network (Section~\ref{sec:problem_structure}).

Instead, we develop a subpath-based formulation of microtransit operations in a load-expanded network. This approach avoids time expansion by capturing time window requirements within subpaths, and enables a coarse spatial discretization between checkpoints. \cite{macedo2011solving} proposed a pseudo-polynomial network flow model based on timed routes; \cite{vazifeh2018addressing,bertsimas2019online} relied on vehicle-sharing networks for fleet sizing and ride-sharing; in pickup-and-delivery and ride-pooling, \cite{alyasiry2019exact,zhang2023routing} defined subpaths from a point where the vehicle is empty to another; \cite{rist2021new,hasan2021benefits} considered subpaths consisting of a sequence of consecutive pickups and dropoffs. \cite{schulz2024using} developed a fixed-path procedure to accelerate branch-and-cut algorithms for routing problems with precedence constraints, using a compact formulation. All these models focus on single-stage routing optimization. Our paper contributes a new load-expanded subpath-based network representation of microtransit operations; and it embeds it into a two-stage stochastic optimization framework to jointly optimize microtransit design and operations.

The problem is formulated as a two-stage stochastic integer program with an exponential number of second-stage variables, which we solve via Benders decomposition and column generation. Column generation has been applied to network-based extended formulations \citep{sadykov2013column,delorme2020enhanced,hasan2021benefits,jacquillat2024relay}; our methodology embeds it into a Benders decomposition scheme to solve the two-stage stochastic optimization problem. Combinations of column generation and Benders decomposition encompass simultaneous column-and-row generation \citep{muter2013simultaneous}, as well as path-based column generation in the Benders master problem or subproblem \citep{karsten2018simultaneous,zeighami2019a,wu2022vessel}. In contrast, our algorithm adds subpath-based variables to the Benders subproblem, which gives rise to a novel double-decomposition structure (Section~\ref{sec:alg}).

In summary, our methodology contributes: (i) a subpath-based extended formulation of microtransit operations; (ii) an integrated two-stage stochastic optimization formulation with a tight subpath-based representation of second-stage vehicle routing; and (iii) a double-decomposition solution algorithm combining Benders decomposition and subpath-based column generation.
\section{Microtransit Network Design (MiND) Model}\label{sec:model}

The MiND optimizes the design and operations of a deviated fixed-route microtransit system, under demand uncertainty. The first-stage problem defines reference lines and service schedules (Section~\ref{sec:schedule}). The second-stage problem defines on-demand deviations in response to passenger requests, using a subpath-based representation in a load-expanded network (Section~\ref{sec:subpath_formulation}). We formulate the MiND-VRP in Section~\ref{sec:model_formulation} and discuss its two-stage stochastic discrete optimization structure in Section~\ref{sec:problem_structure}. The extensions to MiND-DAR and Mind-Tr are in~\ref{subsec:MiND-DAR} and~\ref{subsec:MiND-Tr}.

\subsection{First-stage Problem: Network Design and Frequency Planning}\label{sec:schedule}

The first-stage problem defines \textit{reference trips}, each characterized by a reference line and a departure time. Each reference line is defined as an ordered set of checkpoints, and each reference trip determines the scheduled time at each checkpoint. Vehicles are required to visit some checkpoints at the scheduled times, but will also be allowed to visit other locations in-between (Section~\ref{sec:subpath_formulation}).

Operations occur over a roadway network. Let $\calN^S$ denote the set of stations, including all candidate checkpoints and possible stopping locations. We represent demand as a set of passenger requests $\calP$. In the MiND-VRP, each request $p\in\calP$ is characterized by an origin $o(p)\in\calN^S$ and a requested drop-off time $t^{\text{req}}_p$. Demand uncertainty is modeled via a set of scenarios $\calS$; each scenario $s \in \calS$ has probability $\pi_s$ and comprises $D_{ps}$ passengers from request $p \in \calP$.

\paragraph{Network design.} We pre-process candidate reference lines in a set $\calL$. Let $h_{\ell}$ denote the cost to operate one trip of line $\ell\in\calL$. Let $\calT_{\ell}$ store time periods when a vehicle can depart from the first checkpoint in line $\ell \in \calL$. We introduce the following decision variables to define reference trips:
\begin{equation*}
    x_{\ell t} = \begin{cases}
    1 &\text{reference trip $(\ell, t)\in\calL\times\calT_\ell$ is selected,}\\
    0 &\text{otherwise.}
    \end{cases}
\end{equation*}

Let $\calI_{\ell}\subseteq\calN^S$ index the checkpoints in reference line $\ell$, of cardinality $I_{\ell} = |\calI_{\ell}|$. Let $\calI_{\ell}^{(i)}$ refer to the $i^{th}$ checkpoint in the line, for $i\in\{1,\cdots,I_l\}$. All reference lines share the same final checkpoint $\calI^{\text{end}}=\calI_{\ell}^{(I_{\ell})}$. Reference trip $(\ell, t) \in \calL \times \calT_{\ell}$ is scheduled in checkpoint $\calI_{\ell}^{(i)}$ at time $T_{\ell t}(\calI_\ell^{(i)})$.

We impose a fleet budget constraint by limiting the number of active trips at any time $t$:
\begin{align}\label{eq:budget}
    &\sum_{\ell \in \calL}\ \sum_{t' \in \calT_\ell \, : \, t' \leq t \leq t' + T_{\ell t}\left(\calI^{\text{end}}\right)-T_{\ell t}\left(\calI_{\ell}^{(1)}\right)} x_{\ell t} \leq F,\qquad\forall t \in \bigcup_{\ell \in \calL} \calT_\ell 
\end{align}

\paragraph{Internal passenger assignments.} Let $\calM_p \subseteq \calL \times \calT_{\ell}$ denote the subset of reference trips that can serve request $p\in\calP$ within a tolerance $\alpha$ of their requested drop-off time:
$$\calM_p=\left\{(\ell,t)\in\calL\times\calT_\ell:\left|T_{\ell t}\left(\calI^{\text{end}}\right)-t^{\text{req}}_p\right|\leq\alpha\right\},\ \forall p\in\calP$$

We define assignment variables to identify a candidate reference trip for each passenger:
\begin{equation*}
    z_{\ell pst} = \begin{cases}
    1 &\text{if passenger request $p\in \calP$ is internally assigned to trip $(\ell, t)\in \calM_p$ in scenario $s \in \calS$,} \\
    0 &\text{otherwise.} 
    \end{cases}
\end{equation*}

We impose packing constraints so that each passenger is assigned to at most one reference trip.
\begin{align}\label{eq:assign}
    &\sum_{(\ell, t) \in \calM_p} z_{\ell pst} \leq 1,\qquad \forall p \in \calP, \forall s \in \calS
\end{align}

The assignments $z_{\ell pst}$ link first-stage and second-stage decisions but are not executed in practice, since passenger service is optimized at the operational level. These variables are scenario-dependent but will be treated as first-stage variables in the algorithm; this choice enables separability across reference trips and leads to a slightly tighter second-stage relaxation. A similar methodology and similar computational results could be obtained by treating $z_{\ell pst}$ as second-stage variables.

\paragraph{Vehicle load.} We assume that vehicles are homogeneous within each reference line $\ell \in \calL$, with capacity $C_\ell$. We impose a {\it target load factor} $\kappa \in (0,1)$ to induce high vehicle utilization. We also allow first-stage assignments to exceed vehicle capacities by a factor $\kappa$ to create operating flexibility, but the second-stage passenger service decisions will strictly comply with vehicle capacities.
\vspace{-6pt}
\begin{align}
    &\sum_{p \in \calP \, : \, (\ell, t) \in \calM_p} D_{ps}z_{\ell pst} \geq (1-\kappa) C_\ell x_{\ell t} && \forall (\ell, t) \in \calL \times \calT_\ell, \forall s \in \calS \label{eq:load}\\
    &\sum_{p \in \calP \, : \, (\ell, t) \in \calM_p} D_{ps}z_{\ell pst} \leq (1+\kappa)C_\ell x_{\ell t} && \forall (\ell, t) \in \calL \times \calT_\ell, \forall s \in \calS\label{eq:load2}
\end{align}

\subsection{Second-stage Problem: On-demand Deviations} \label{sec:subpath_formulation}

Second-stage deviations must stay within a distance $\Delta$ of the reference line and must respect scheduled times at the checkpoints. The reference schedule includes buffers between checkpoints to allow for deviations. Moreover, vehicles may skip up to $K$ checkpoints in a row: $K=0$ induces closer adherence to the reference trip, whereas $K\geq1$ provides more flexibility. We denote by $\Gamma_\ell$ the checkpoint pairs separated by up to $K$ checkpoints on line $\ell\in\calL$.

The second-stage problem involves capacitated vehicle routing with time windows for each reference trip and in each scenario. We formulate it in a load-expanded network with subpath variables characterizing on-demand deviations between checkpoints, leveraging the structure of microtransit.

\subsubsection*{Subpaths.} For reference trip $(\ell, t) \in \calL \times \calT_{\ell}$ and demand scenario $s \in \calS$, a subpath $r \in \calR_{\ell st}$ is identified by its starting checkpoint $u_r\in\calI_{\ell}$, its ending checkpoint $v_r\in\calI_{\ell}$, and the passenger requests $\calP_r \subseteq \calP$ served in between. The set $\calR_{\ell st}$ includes all subpaths such that the distance to the reference line never exceeds $\Delta$; the load satisfies $\sum_{p \in \calP_r} D_{ps} \leq C_\ell$; the travel time does not exceed  $T_{\ell t}(v_r) - T_{\ell t}(u_r)$; and up to $K$ checkpoints are skipped. The second-stage problem selects a sequence of subpaths that (i) starts at the origin of the reference line and ends at its destination while maintaining flow balance; and (ii) serves up to $C_\ell$ passengers overall.

\subsubsection*{Load-expanded subpath network.} 

\begin{figure}[htbp!]
    \begin{center}
    \includegraphics[width=0.2\textwidth]{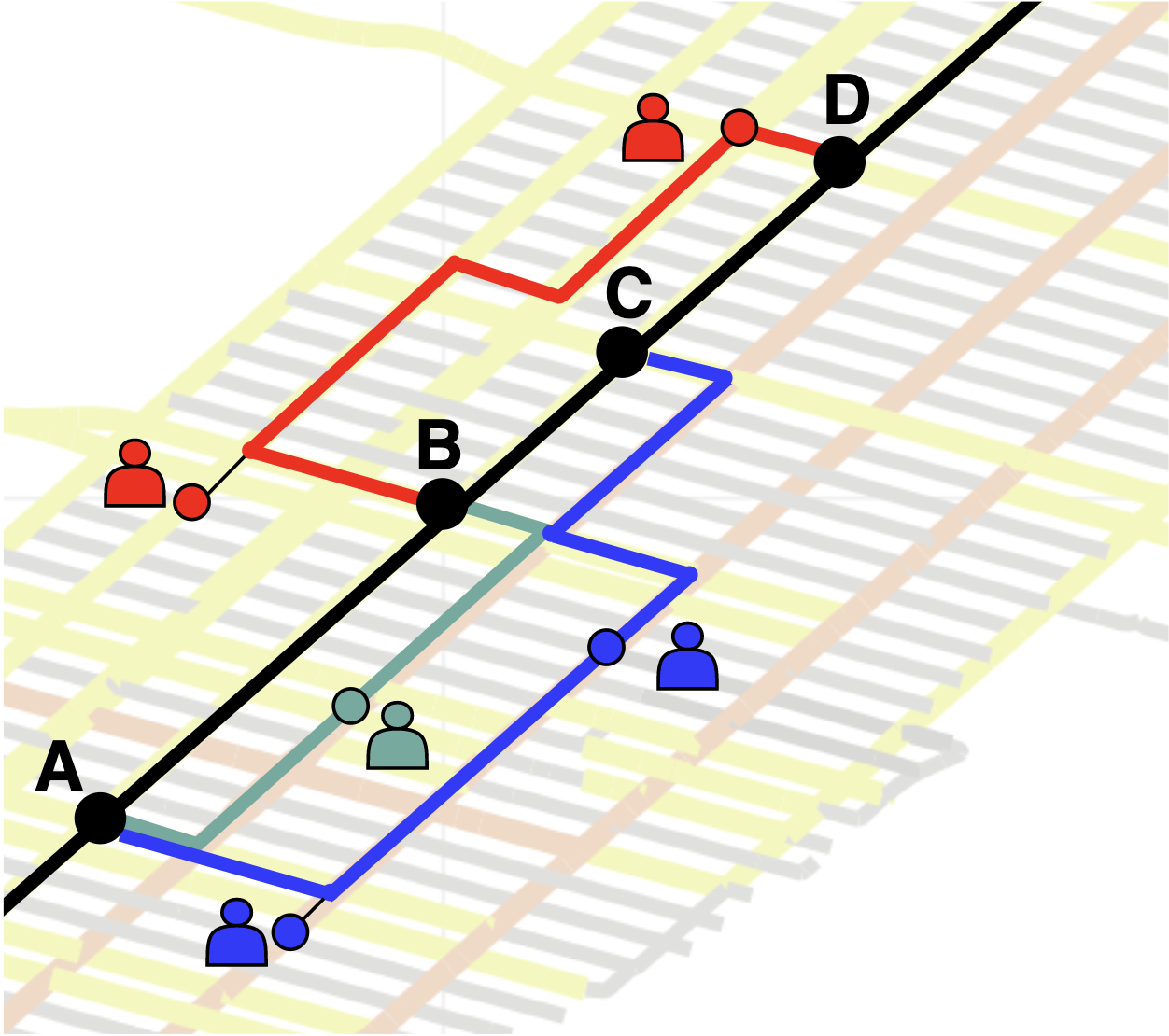}
    \hspace{1 cm}
    \includegraphics[width=0.5\textwidth]{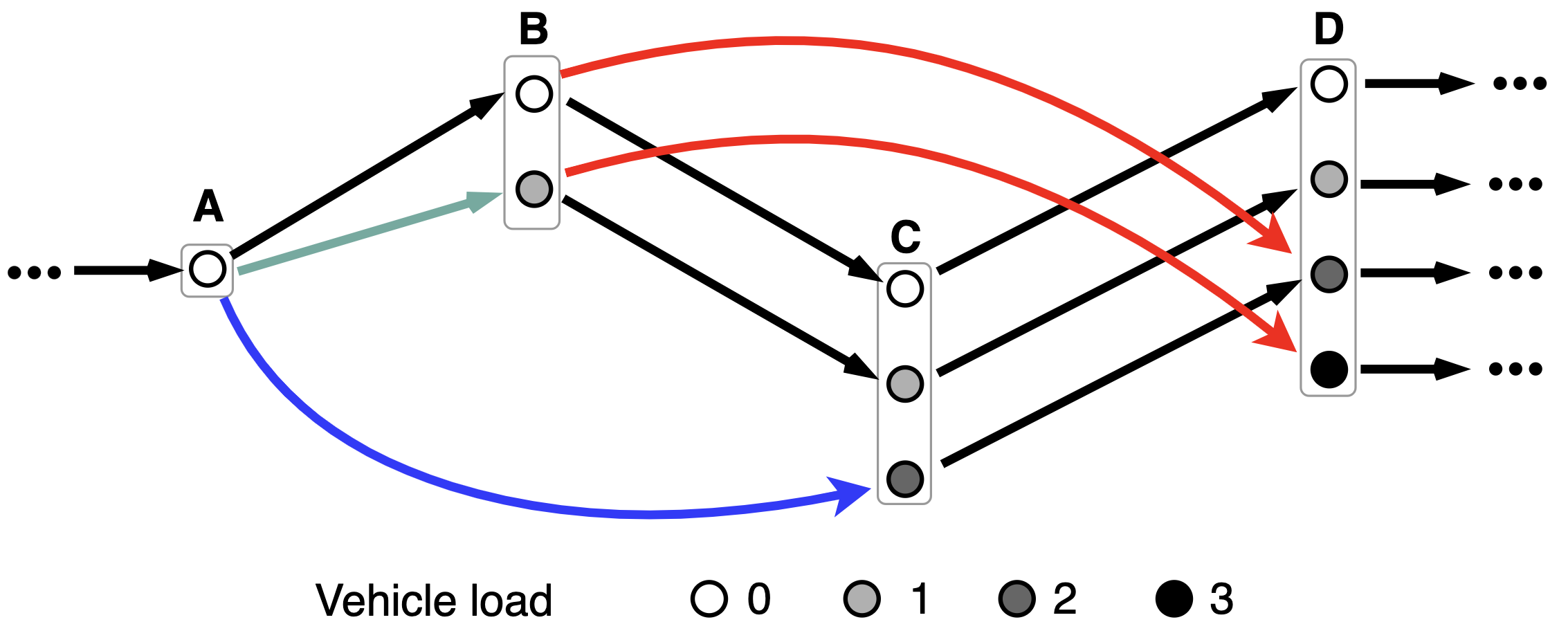}
    \end{center}
    \caption{{\it Left:} Physical network with three candidate deviations. {\it Right:} Load-expanded subpath network.}
    \label{F:subpath_comparisons}
\end{figure}

We represent routing operations in a load-expanded network (Figure~\ref{F:subpath_comparisons}). Each node tracks the checkpoint and the vehicle load, and each arc encapsulates a subpath. Flow balance constraints capture physical flows and vehicle capacity constraints.

Let $\calC_\ell = \{0, 1, \cdots, C_\ell\}$ store vehicle loads. For reference trip $(\ell, t) \in \calL \times \calT_{\ell}$ and scenario $s \in \calS$, we denote the load-expanded network by $(\calV_{\ell st}, \calA_{\ell st})$. Node $v_{\ell st}$ represents the end of a trip. Each other node $n\in \calV_{\ell st}$ corresponds to a tuple $(k_n, c_n)$ consisting of checkpoint $k_n \in \calI_{\ell}$ and load $c_n \in \calC_\ell$; node $u_{\ell st}$ encodes the start at stop $\calI_{\ell}^{(1)}$ and time $T_{\ell t}(\calI_\ell^{(1)})$. Each arc $a\in\calA_{\ell st}$ connects nodes $start(a)\in\calV_{\ell st}$ and $end(a)\in\calV_{\ell st}$; vice versa, $r(a) \in \calR_{\ell st}$ is the subpath corresponding to arc $a \in \bigcup_{r \in \calR_{\ell st}} \calA_r$.  We partition $\calA_{\ell st}= \bigcup_{r \in \calR_{\ell st}}\calA_r \cup \calA_{\ell st}^v$ into traveling arcs $\calA_r$ and terminating arcs $\calA_{\ell st}^v$:
\begin{align}
\calA_r &= \left\{(n,m) \in \calV_{\ell st}\times \calV_{\ell st} \, : \, k_n = u_r, k_m = v_r, c_m - c_n = \sum_{p \in \calP_{r}} D_{ps}\right\}, \quad \forall r \in \calR_{\ell st}\label{eq:ASP},\\
\calA^v_{\ell st} &= \{(n,m) \in \calV_{\ell st}\times\calV_{\ell st} \, : \, k_n = \calI^{\text{end}}, m = v_{\ell st} \}.\label{eq:AT}
\end{align}
Our second-stage decisions select subpaths in the load-expanded networks via the following variables. These define on-demand deviations and pickups for each reference trip and each scenario.
\begin{equation*}
    y_a = \begin{cases}
    1 &\text{if arc $a\in \calA_{\ell st}$ is selected, for $(\ell, t) \in \calL \times \calT_{\ell},\ s \in \calS$,} \\
    0 &\text{otherwise.}
    \end{cases}
\end{equation*}

\subsubsection*{Passenger service.}

The second-stage formulation aims to maximize demand coverage and passenger level of service. Coverage is formalized via a large reward $M$ incurred for each successful pickup. Level of service is formalized via the following objectives:
\begin{enumerate}
    \item $\tau^{walk}_{rp}$: walking time from passenger $p$'s origin to the pickup location via subpath $r \in \calR_{\ell st}$;
    \item $\tau^{wait}_{rp}$: waiting time of passenger $p$ prior to pickup via subpath $r \in \calR_{\ell st}$; 
    \item $\frac{\tau^{travel}_{rp}}{\tau_p^{dir}}$: relative detour, defined as the in-vehicle travel time of passenger $p$ via subpath $r \in \calR_{\ell st}$ normalized with respect to the direct trip time (e.g., a taxi trip); and
    \item $\frac{\tau^{late}_{\ell tp}}{\tau_p^{dir}}$, $\frac{\tau^{early}_{\ell tp}}{\tau_p^{dir}}$: relative delay and earliness of passenger $p$ at the destination via trip $(\ell, t) \in \calM_p$.
\end{enumerate}

The model can be extended to incorporate additional practical considerations in arc cost parameters. For example, $g_a$ could penalize subpaths that skip checkpoints to promote closer adherence to the reference line, even when $K>0$ (although shorter subpaths arise naturally in our solutions).

In addition, we impose three restrictions to guarantee convenient service. First, passengers must not leave before their planned departure times $t_{\ell pt}^0$, corresponding to the time at which they would start walking to the nearest checkpoint without deviations. Second, they must not walk more than a limit $\Omega$. Third, they must not wait more than a limit $\Psi$. Thus, pickup $p\in\calP$ is only acceptable in location $i\in\calN^S$ at time $\bar{t}$ via subpath $r\in\calR_{\ell st}$ if $t^0_{\ell p t} + \tau^{walk}_{rp} \leq \bar{t}\leq t^0_{lpt} + \tau^{walk}_{rp}+\Psi$ and $\tau^{walk}_{rp} \leq \Omega$.

We define non-negative hyperparameters $\lambda$, $\mu$, $\sigma,$ and $\delta$ to weigh the level of service cost components. The arc costs in the load-expanded network are defined as follows for all $(\ell, t) \in \calL \times \calT_{\ell},\ s \in \calS$.
\begin{align}\label{eq:arccost}
g_{a} = &\begin{cases}
    \sum_{p \in \calP_{r(a)}}D_{ps} \left(\lambda \tau^{walk}_{r(a)p} + \mu \tau^{wait}_{r(a)p} + \sigma \frac{\tau_{r(a)p}^{travel}}{\tau_p^{dir}} + \delta \frac{\tau^{late}_{\ell tp}}{\tau_p^{dir}} + \frac{\delta}{2} \frac{\tau^{early}_{\ell tp}}{\tau_p^{dir}} - M  \right) &\forall a \in \bigcup_{r \in \calR_{\ell st}}\calA_r,\\
    0 &\forall a  \in\calA^v_{\ell st}.
\end{cases}
\end{align}

\subsection{Two-stage Stochastic Optimization Formulation (MiND-VRP)}\label{sec:model_formulation}

The MiND-VRP minimizes planning costs, maximizes demand coverage, and maximizes level of service (Equation~\eqref{obj}). The constraints apply fleet size, target load factors, and packing constraints (Equations~\eqref{eq:budget}--\eqref{eq:load2}); enforce flow balance over load-expanded networks (Constraint~\eqref{eq:2Sflow}); and link first-stage assignments to second-stage operations (Constraint~\eqref{eq:2Sy2z}). Notation is summarized in~\ref{app:notation}. The MiND-DAR and MiND-Tr are formulated similarly, albeit with extra constraints to ensure consistency between pickups, transfers and dropoffs.
\begin{align}
    \min   \quad&\sum_{\ell \in \calL} \sum_{t \in \calT_\ell} \left(h_\ell x_{\ell t} +  \sum_{s \in \calS} \pi_s \sum_{a \in \calA_{\ell st}} g_a y_a\right)\label{obj}\\
    \st\quad & \text{First-stage constraints: Equations~\eqref{eq:budget}--\eqref{eq:load2}}\nonumber\\
    &\sum_{m:(n,m) \in \calA_{\ell st}} 
    y_{(n,m)} - \sum_{m:(m,n) \in \calA_{\ell st}} 
    y_{(m,n)} = \begin{cases}
    x_{\ell t} &\text{if } n = u_{\ell st} \\
    -x_{\ell t} &\text{if } n = v_{\ell st} \\
    0 &\text{ otherwise}
    \end{cases} \ \forall \ell\in \calL, t \in \calT_{\ell}, s \in \calS, n \in \calV_{\ell st}\label{eq:2Sflow} \\
    &\sum_{a \in \calA_{\ell st} \, : \, p \in \calP_{r(a)}} y_a \leq z_{\ell pst} \quad \forall
    s \in \calS, p \in \calP, (\ell,t) \in \calM_p \label{eq:2Sy2z} \\
    &\bx, \by, \bz \text{ binary} \label{domain}
\end{align}

\subsection{General structure, and comparison to benchmarks}\label{sec:problem_structure}

We provide a general-purpose description of our two-stage stochastic discrete optimization setting to compare our subpath-based network model to alternative approaches based on compact formulations or other network-extended formulations. This description will also be used in Section~\ref{sec:alg} to develop our double decomposition algorithm in a general-purpose environment.

\paragraph{Compact formulation.}
Let $\bx\in\calX^{\text{MIO}}$ and $\by^C_s\in\calY^{\text{MIO-C}}_s$ denote first-stage and second-stage variables, respectively. We use $\calX^{\text{MIO}}$ and $\calY^{\text{MIO-C}}_s$ to represent mixed-integer optimization regions; for instance, in a setting with $n_Z$ integer variables and $n_R$ continuous variables, they encode the mixed-integer set $\Z^{n_Z}\times\R^{n_R}$. The optimization problem is formulated as follows, where $\calK$ stores indices of the first-stage decisions (e.g., line-time pairs in the MiND), $\bA\bx\geq\bb$ encode all first-stage constraints (e.g., fleet budget), and $\bD_s\bx+\bE_s\by^C_s\geq\widetilde{\bh}_s$ encode all second-stage and linking constraints (e.g., checkpoint-to-checkpoint operations for selected microtransit lines).
\begin{align}
    \min\quad   &   \sum_{j \in \calJ}\sum_{k\in\calK_j}c_kx_k+\sum_{s\in\calS}\pi_s\widetilde{\bg}_s^\top\by^C_s\label{prob:xz}\tag{$C$}\\
    \st\quad    &   \bA\bx\geq\bb\nonumber\\
                &   \bD_s\bx+\bE_s\by^C_s\geq\widetilde{\bh}_s, \forall s\in\calS\nonumber\\
                &   \bx\in\calX^{\text{MIO}};\quad\by^C_s\in\calY^{\text{MIO-C}}_s,\ \forall s\in\calS\nonumber
\end{align}

\paragraph{Extended reformulations.}
Our subpath-based model lifts the second-stage problem in a higher-dimensional space. Let $\calK_j,$ for $j\in\calJ$, denote a partition of the first-stage indices such that the second-stage problem is independent across subsets, with block-diagonal matrices $\bD_s$ across $\calK_j$. In the MiND, this corresponds to the partition across line-time pairs $(\ell,t)$. We then derive a network representation induced by scenario $s$ and subset $\calK_j$, for $j\in\calJ$, with node set $\calN_{sj}$ and arc set $\calA_{sj}$. We denote the new second-stage variables by $y_{a}\in\calY^{\text{MIO}}_a$ for each arc $a\in\calA_{sj}$ in each scenario $s\in\calS$ and for each partition element $j\in\calJ$; again, the sets $\calY^{\text{MIO}}_a$ encode the discrete requirements of variables $y_{a}$. The formulation encompass flow balance constraints in the network $(\calN_{sj},\calA_{sj})$ (e.g., combining subpaths into paths, in Equation~\eqref{eq:2Sflow}) as well as general-purpose linear constraints $\sum_{a\in\calA_{sj}}\bff_{asj} y_{a}\geq\bh_{sj}$ (e.g., demand constraints in Equation~\eqref{eq:2Sy2z}). It is given in Problem~\eqref{OPT}:
\begin{align}
    \min\quad   &   \sum_{j \in \calJ}\sum_{k\in\calK_j}c_kx_k+\sum_{s\in\calS}\sum_{j \in \calJ}\pi_s\left(\sum_{a\in\calA_{sj}}g_{a}y_{a}\right)\label{OPT}\tag{$\star$}\\
    \st\quad    &   \bA\bx\geq\bb\nonumber\\
                &   \sum_{m:(n,m)\in\calA_{sj}}y_{(n,m)}-\sum_{m:(m,n)\in\calA_{sj}}y_{(m,n)}=\sum_{k\in\calK_j}b_{nsk}x_k,\ \forall s\in\calS,\ \forall j\in\calJ,\ \forall n\in\calN_{sj}\nonumber\\
                &   \sum_{a\in\calA_{sj}}\bff_{asj} y_{a}\geq\bh_{sj},\ \forall s\in\calS,\ \forall j\in\calJ\nonumber\\
                &   \bx\in\calX^{\text{MIO}};\quad y_{a}\in\calY^{\text{MIO}}_a,\ \forall a\in\calA_{sj},\ \forall s\in\calS\ \forall j\in\calJ\nonumber
\end{align}

\paragraph{Discussion.}

Proposition~\ref{prop:compact} shows that the extended network-based formulation defines the same optimization problem as the compact formulation in the MiND. Moreover, it defines a tighter second-stage linear relaxation by embedding time windows into the network representation itself and capturing vehicle capacities via flow balance constraints, as compared to relying on explicit big-M constraints. This reformulation, however, involves more subpath-based variables.

\begin{proposition}\label{prop:compact}
    The subpath-based second-stage formulation defines an equivalent mixed-integer model as the compact formulation in the MiND, with a tighter linear relaxation.
\end{proposition}

This discussion invites comparisons to other possible network-based representations of routing operations in our second-stage problem, detailed in~\ref{app:segment_path_benchmarks} and illustrated in Figure~\ref{fig:formulations}:
\begin{itemize}
    \item[--] A \textit{segment-based} model with nodes encoding time-station-load tuples and arcs encoding segments between consecutive stops. This network is much more granular along the time dimension (each station in Figure~\ref{subfig:segment}, versus each checkpoint in Figure~\ref{subfig:subpath}) and is also expanded in a granular temporal discretization (to capture passengers' and vehicles' time windows). Moreover, the model is further complicated by two multi-commodity flow structures with linking constraints from checkpoint to checkpoint (so the vehicle does not skip more than $K$ checkpoints in a row) and from station to station (to maintain continuity in time and space).
    \item[--] A \textit{path-based} model in a simple origin-destination network. with arcs encoding full paths that start at the line's origin, end at its destination, and serve at most $C_\ell$ passengers.
\end{itemize}

\begin{figure}[h!]
    \centering
    \subfloat[Path-based formulation]{\label{subfig:path}
        \begin{tikzpicture}[scale=0.6,transform shape]\node[] at (-2,0){};
            \node[] at (0,1){Start};
            \node[] at (15,1){End};
            \node[ultra thick,draw=black,fill=black,minimum size=10pt] (O) at (0,0){};
            \node[ultra thick,draw=black,fill=black,minimum size=10pt] (D) at (15,0){};
            \draw[->,draw=black,thick,dotted] (O) to (D);
            \foreach \i in {1,...,20}
            {
                    \draw[->,draw=black,bend left=\i,thick,dotted] (O) to (D);
                    \draw[->,draw=black,bend right=\i,thick,dotted] (O) to (D);
            }
            \foreach \i in {21,...,37}
            {
                    \draw[->,draw=black,bend left=\i,thick,dotted] (O) to (D);
                    \draw[->,draw=black,bend right=\i,thick,dotted] (O) to (D);
            }
            \draw[->,draw=myred,bend left=20,ultra thick] (O) to (D);
            \draw[->,draw=myblue,ultra thick] (O) to (D);
            \draw[->,draw=mygreen,bend right=20,ultra thick] (O) to (D);
        \end{tikzpicture}}\hspace{0.1cm}
    \subfloat[Subpath-based formulation (squares: checkpoints)]{\label{subfig:subpath}
        \begin{tikzpicture}[scale=0.6,transform shape]\node[] at (-2,0){};
            \node[] at (0,4){Start};
            \node[] at (-1.5,0){load: 0};
            \node[] at (-1.5,1){load: 1};
            \node[] at (-1.5,2){load: 2};
            \node[] at (-1.5,3){load: 3};
            \node[] at (15,4){End};
            \node[thick,draw=black,fill=black,minimum size=10pt] (O) at (0,0){};
            \node[thick,draw=black,fill=white,minimum size=10pt] (c1l0) at (3,0){};
            \node[thick,draw=black,fill=white,minimum size=10pt] (c1l1) at (3,1){};
            \node[thick,draw=black,fill=white,minimum size=10pt] (c1l2) at (3,2){};
            \node[thick,draw=black,fill=white,minimum size=10pt] (c1l3) at (3,3){};
            \node[thick,draw=black,fill=white,minimum size=10pt] (c2l0) at (6,0){};
            \node[thick,draw=black,fill=white,minimum size=10pt] (c2l1) at (6,1){};
            \node[thick,draw=black,fill=white,minimum size=10pt] (c2l2) at (6,2){};
            \node[thick,draw=black,fill=white,minimum size=10pt] (c2l3) at (6,3){};
            \node[thick,draw=black,fill=white,minimum size=10pt] (c3l0) at (9,0){};
            \node[thick,draw=black,fill=white,minimum size=10pt] (c3l1) at (9,1){};
            \node[thick,draw=black,fill=white,minimum size=10pt] (c3l2) at (9,2){};
            \node[thick,draw=black,fill=white,minimum size=10pt] (c3l3) at (9,3){};
            \node[thick,draw=black,fill=white,minimum size=10pt] (c4l0) at (12,0){};
            \node[thick,draw=black,fill=white,minimum size=10pt] (c4l1) at (12,1){};
            \node[thick,draw=black,fill=white,minimum size=10pt] (c4l2) at (12,2){};
            \node[thick,draw=black,fill=white,minimum size=10pt] (c4l3) at (12,3){};
            \node[thick,draw=black,fill=black,minimum size=10pt] (Dl0) at (15,0){};
            \node[thick,draw=black,fill=black,minimum size=10pt] (Dl1) at (15,1){};
            \node[thick,draw=black,fill=black,minimum size=10pt] (Dl2) at (15,2){};
            \node[thick,draw=black,fill=black,minimum size=10pt] (Dl3) at (15,3){};
            \draw[->,black,thick,dotted] (O) to (c1l0);
            \draw[->,black,thick,dotted] (c1l0) to (c2l0);
            \draw[->,black,thick,dotted] (c2l0) to (c3l0);
            \draw[->,black,thick,dotted] (c3l0) to (c4l0);
            \draw[->,black,thick,dotted] (c4l0) to (Dl0);
            \draw[->,black,thick,dotted] (c3l1) to (c4l1);
            \draw[->,black,thick,dotted] (c4l1) to (Dl1);
            \draw[->,black,thick,dotted] (c1l2) to (c2l2);
            \draw[->,black,thick,dotted] (c2l2) to (c3l2);
            \draw[->,black,thick,dotted] (c3l2) to (c4l2);
            \draw[->,black,thick,dotted] (c1l3) to (c2l3);
            \draw[->,black,thick,dotted] (c2l3) to (c3l3);
            \draw[->,black,thick,dotted] (c3l3) to (c4l3);
            \draw[->,black,thick,dotted] (c4l2) to (Dl2);
            \draw[->,black,thick,dotted] (O) to (c1l1);
            \draw[->,black,bend right=10,thick,dotted] (O) to (c1l1);
            \draw[->,black,thick,dotted] (O) to (c1l2);
            \draw[->,black,bend left=10,thick,dotted] (O) to (c1l2);
            \draw[->,black,bend right=10,thick,dotted] (O) to (c1l2);
            \draw[->,black,thick,dotted] (O) to (c1l3);
            \draw[->,black,bend left=10,thick,dotted] (O) to (c1l3);
            \draw[->,black,bend right=10,thick,dotted] (O) to (c1l3);
            \draw[->,black,thick,dotted] (c1l0) to (c2l1);
            \draw[->,black,bend left=10,thick,dotted] (c1l0) to (c2l1);
            \draw[->,black,bend right=10,thick,dotted] (c1l0) to (c2l1);
            \draw[->,black,thick,dotted] (c1l0) to (c2l2);
            \draw[->,black,bend left=10,thick,dotted] (c1l0) to (c2l2);
            \draw[->,black,bend right=10,thick,dotted] (c1l0) to (c2l2);
            \draw[->,black,thick,dotted] (c1l0) to (c2l3);
            \draw[->,black,bend left=10,thick,dotted] (c1l0) to (c2l3);
            \draw[->,black,bend right=10,thick,dotted] (c1l0) to (c2l3);
            \draw[->,black,thick,dotted] (c1l1) to (c2l2);
            \draw[->,black,bend left=10,thick,dotted] (c1l1) to (c2l2);
            \draw[->,black,bend right=10,thick,dotted] (c1l1) to (c2l2);
            \draw[->,black,thick,dotted] (c1l1) to (c2l3);
            \draw[->,black,bend left=10,thick,dotted] (c1l1) to (c2l3);
            \draw[->,black,bend right=10,thick,dotted] (c1l1) to (c2l3);
            \draw[->,black,thick,dotted] (c1l2) to (c2l3);
            \draw[->,black,bend left=10,thick,dotted] (c1l2) to (c2l3);
            \draw[->,black,bend right=10,thick,dotted] (c1l2) to (c2l3);
            \draw[->,black,thick,dotted] (c2l0) to (c3l1);
            \draw[->,black,bend left=10,thick,dotted] (c2l0) to (c3l1);
            \draw[->,black,bend right=10,thick,dotted] (c2l0) to (c3l1);
            \draw[->,black,thick,dotted] (c2l0) to (c3l2);
            \draw[->,black,bend left=10,thick,dotted] (c2l0) to (c3l2);
            \draw[->,black,bend right=10,thick,dotted] (c2l0) to (c3l2);
            \draw[->,black,thick,dotted] (c2l0) to (c3l3);
            \draw[->,black,bend left=10,thick,dotted] (c2l0) to (c3l3);
            \draw[->,black,bend right=10,thick,dotted] (c2l0) to (c3l3);
            \draw[->,black,thick,dotted] (c2l1) to (c3l2);
            \draw[->,black,bend left=10,thick,dotted] (c2l1) to (c3l2);
            \draw[->,black,bend right=10,thick,dotted] (c2l1) to (c3l2);
            \draw[->,black,thick,dotted] (c2l1) to (c3l3);
            \draw[->,black,bend left=10,thick,dotted] (c2l1) to (c3l3);
            \draw[->,black,bend right=10,thick,dotted] (c2l1) to (c3l3);
            \draw[->,black,thick,dotted] (c2l2) to (c3l3);
            \draw[->,black,bend left=10,thick,dotted] (c2l2) to (c3l3);
            \draw[->,black,bend right=10,thick,dotted] (c2l2) to (c3l3);
            \draw[->,black,thick,dotted] (c3l0) to (c4l1);
            \draw[->,black,bend left=10,thick,dotted] (c3l0) to (c4l1);
            \draw[->,black,bend right=10,thick,dotted] (c3l0) to (c4l1);
            \draw[->,black,thick,dotted] (c3l0) to (c4l2);
            \draw[->,black,bend left=10,thick,dotted] (c3l0) to (c4l2);
            \draw[->,black,bend right=10,thick,dotted] (c3l0) to (c4l2);
            \draw[->,black,thick,dotted] (c3l0) to (c4l3);
            \draw[->,black,bend left=10,thick,dotted] (c3l0) to (c4l3);
            \draw[->,black,bend right=10,thick,dotted] (c3l0) to (c4l3);
            \draw[->,black,thick,dotted] (c3l1) to (c4l2);
            \draw[->,black,bend left=10,thick,dotted] (c3l1) to (c4l2);
            \draw[->,black,bend right=10,thick,dotted] (c3l1) to (c4l2);
            \draw[->,black,thick,dotted] (c3l1) to (c4l3);
            \draw[->,black,bend left=10,thick,dotted] (c3l1) to (c4l3);
            \draw[->,black,bend right=10,thick,dotted] (c3l1) to (c4l3);
            \draw[->,black,thick,dotted] (c3l2) to (c4l3);
            \draw[->,black,bend left=10,thick,dotted] (c3l2) to (c4l3);
            \draw[->,black,bend right=10,thick,dotted] (c3l2) to (c4l3);
            \draw[->,black,thick,dotted] (c4l0) to (Dl1);
            \draw[->,black,bend left=10,thick,dotted] (c4l0) to (Dl1);
            \draw[->,black,bend right=10,thick,dotted] (c4l0) to (Dl1);
            \draw[->,black,thick,dotted] (c4l0) to (Dl2);
            \draw[->,black,bend left=10,thick,dotted] (c4l0) to (Dl2);
            \draw[->,black,bend right=10,thick,dotted] (c4l0) to (Dl2);
            \draw[->,black,thick,dotted] (c4l0) to (Dl3);
            \draw[->,black,bend left=10,thick,dotted] (c4l0) to (Dl3);
            \draw[->,black,bend right=10,thick,dotted] (c4l0) to (Dl3);
            \draw[->,black,thick,dotted] (c4l1) to (Dl2);
            \draw[->,black,bend left=10,thick,dotted] (c4l1) to (Dl2);
            \draw[->,black,bend right=10,thick,dotted] (c4l1) to (Dl2);
            \draw[->,black,thick,dotted] (c4l1) to (Dl3);
            \draw[->,black,bend left=10,thick,dotted] (c4l1) to (Dl3);
            \draw[->,black,bend right=10,thick,dotted] (c4l1) to (Dl3);
            \draw[->,black,thick,dotted] (c4l2) to (Dl3);
            \draw[->,black,bend left=10,thick,dotted] (c4l2) to (Dl3);
            \draw[->,black,bend right=10,thick,dotted] (c4l2) to (Dl3);
            \draw[->,myred,bend left=10,ultra thick] (O) to (c1l1);
            \draw[->,myred,ultra thick] (c1l1) to (c2l1);
            \draw[->,myred,ultra thick] (c2l1) to (c3l1);
            \draw[->,myred,bend right=10,ultra thick] (c3l1) to (c4l3);
            \draw[->,myred,ultra thick] (c4l3) to (Dl3);
            \draw[->,mygreen,ultra thick] (O) to (c1l0);
            \draw[->,mygreen,ultra thick] (c1l0) to (c2l0);
            \draw[->,mygreen,ultra thick] (c2l0) to (c3l0);
            \draw[->,mygreen,bend right=10,ultra thick] (c3l0) to (c4l1);
            \draw[->,mygreen,bend left=10,ultra thick] (c4l1) to (Dl3);
            \draw[->,myblue,bend left=10,ultra thick] (O) to (c1l2);
            \draw[->,myblue,ultra thick] (c1l2) to (c2l2);
            \draw[->,myblue,ultra thick] (c2l2) to (c3l2);
            \draw[->,myblue,ultra thick] (c3l2) to (c4l2);
            \draw[->,myblue,bend left=10,ultra thick] (c4l2) to (Dl3);
        \end{tikzpicture}}\hspace{0.1cm}
    \subfloat[Segment-based formulation (squares: checkpoints; circles: other stations)]{\label{subfig:segment}
        \begin{tikzpicture}[scale=0.6,transform shape]\node[] at (-2,0){};
            \node[] at (0,4){Start};
            \node[] at (-1.5,0){load: 0};
            \node[] at (-1.5,1){load: 1};
            \node[] at (-1.5,2){load: 2};
            \node[] at (-1.5,3){load: 3};
            \draw[->] (0.9,-0.5) to (1.4,-1);
            \node[] at (1.15,-1.2){time};
            \node[] at (15,4){End};
            \foreach \i in {0,...,3}{\foreach \j in {0,...,3}{
                    \draw[->,draw=black,dotted] (0,0) to (1+0.1*\i,\j-0.1*\i);
            }}
            \foreach \k in {1,...,14}{\foreach \i in {0,...,3}{\foreach \j in {0,...,3}{\foreach \ii in {0,...,3}{\foreach \jj in {\j,...,3}{
                    \draw[->,draw=black,dotted] (\k+0.1*\i,\j-0.1*\i) to (\k+1+0.1*\ii,\jj-0.1*\ii);
            }}}}}
            \node[thick,draw=black,fill=black,minimum size=10pt] (Ot1) at (0,0){};
            \node[thick,draw=black,fill=white,minimum size=10pt] (c1l0t1) at (3,0){};
            \node[thick,draw=black,fill=white,minimum size=10pt] (c1l0t2) at (3.1,-.1){};
            \node[thick,draw=black,fill=white,minimum size=10pt] (c1l0t3) at (3.2,-.2){};
            \node[thick,draw=black,fill=white,minimum size=10pt] (c1l0t4) at (3.3,-.3){};
            \node[thick,draw=black,fill=white,minimum size=10pt] (c1l1t1) at (3,1){};
            \node[thick,draw=black,fill=white,minimum size=10pt] (c1l1t2) at (3.1,.9){};
            \node[thick,draw=black,fill=white,minimum size=10pt] (c1l1t3) at (3.2,.8){};
            \node[thick,draw=black,fill=white,minimum size=10pt] (c1l1t4) at (3.3,.8){};
            \node[thick,draw=black,fill=white,minimum size=10pt] (c1l2t1) at (3,2){};
            \node[thick,draw=black,fill=white,minimum size=10pt] (c1l2t2) at (3.1,1.9){};
            \node[thick,draw=black,fill=white,minimum size=10pt] (c1l2t3) at (3.2,1.8){};
            \node[thick,draw=black,fill=white,minimum size=10pt] (c1l2t4) at (3.3,1.8){};
            \node[thick,draw=black,fill=white,minimum size=10pt] (c1l3t1) at (3,3){};
            \node[thick,draw=black,fill=white,minimum size=10pt] (c1l3t2) at (3.1,2.9){};
            \node[thick,draw=black,fill=white,minimum size=10pt] (c1l3t3) at (3.2,2.8){};
            \node[thick,draw=black,fill=white,minimum size=10pt] (c1l3t4) at (3.3,2.8){};
            \node[thick,draw=black,fill=white,minimum size=10pt] (c2l0t1) at (6,0){};
            \node[thick,draw=black,fill=white,minimum size=10pt] (c2l0t2) at (6.1,-.1){};
            \node[thick,draw=black,fill=white,minimum size=10pt] (c2l0t3) at (6.2,-.2){};
            \node[thick,draw=black,fill=white,minimum size=10pt] (c2l0t4) at (6.3,-.3){};
            \node[thick,draw=black,fill=white,minimum size=10pt] (c2l1t1) at (6,1){};
            \node[thick,draw=black,fill=white,minimum size=10pt] (c2l1t2) at (6.1,.9){};
            \node[thick,draw=black,fill=white,minimum size=10pt] (c2l1t3) at (6.2,.8){};
            \node[thick,draw=black,fill=white,minimum size=10pt] (c2l1t4) at (6.3,.8){};
            \node[thick,draw=black,fill=white,minimum size=10pt] (c2l2t1) at (6,2){};
            \node[thick,draw=black,fill=white,minimum size=10pt] (c2l2t2) at (6.1,1.9){};
            \node[thick,draw=black,fill=white,minimum size=10pt] (c2l2t3) at (6.2,1.8){};
            \node[thick,draw=black,fill=white,minimum size=10pt] (c2l2t4) at (6.3,1.8){};
            \node[thick,draw=black,fill=white,minimum size=10pt] (c2l3t1) at (6,3){};
            \node[thick,draw=black,fill=white,minimum size=10pt] (c2l3t2) at (6.1,2.9){};
            \node[thick,draw=black,fill=white,minimum size=10pt] (c2l3t3) at (6.2,2.8){};
            \node[thick,draw=black,fill=white,minimum size=10pt] (c2l3t4) at (6.3,2.8){};
            \node[thick,draw=black,fill=white,minimum size=10pt] (c3l0t1) at (9,0){};
            \node[thick,draw=black,fill=white,minimum size=10pt] (c3l0t2) at (9.1,-.1){};
            \node[thick,draw=black,fill=white,minimum size=10pt] (c3l0t3) at (9.2,-.2){};
            \node[thick,draw=black,fill=white,minimum size=10pt] (c3l0t4) at (9.3,-.3){};
            \node[thick,draw=black,fill=white,minimum size=10pt] (c3l1t1) at (9,1){};
            \node[thick,draw=black,fill=white,minimum size=10pt] (c3l1t2) at (9.1,.9){};
            \node[thick,draw=black,fill=white,minimum size=10pt] (c3l1t3) at (9.2,.8){};
            \node[thick,draw=black,fill=white,minimum size=10pt] (c3l1t4) at (9.3,.8){};
            \node[thick,draw=black,fill=white,minimum size=10pt] (c3l2t1) at (9,2){};
            \node[thick,draw=black,fill=white,minimum size=10pt] (c3l2t2) at (9.1,1.9){};
            \node[thick,draw=black,fill=white,minimum size=10pt] (c3l2t3) at (9.2,1.8){};
            \node[thick,draw=black,fill=white,minimum size=10pt] (c3l2t4) at (9.3,1.8){};
            \node[thick,draw=black,fill=white,minimum size=10pt] (c3l3t1) at (9,3){};
            \node[thick,draw=black,fill=white,minimum size=10pt] (c3l3t2) at (9.1,2.9){};
            \node[thick,draw=black,fill=white,minimum size=10pt] (c3l3t3) at (9.2,2.8){};
            \node[thick,draw=black,fill=white,minimum size=10pt] (c3l3t4) at (9.3,2.8){};
            \node[thick,draw=black,fill=white,minimum size=10pt] (c4l0t1) at (12,0){};
            \node[thick,draw=black,fill=white,minimum size=10pt] (c4l0t2) at (12.1,-.1){};
            \node[thick,draw=black,fill=white,minimum size=10pt] (c4l0t3) at (12.2,-.2){};
            \node[thick,draw=black,fill=white,minimum size=10pt] (c4l0t4) at (12.3,-.3){};
            \node[thick,draw=black,fill=white,minimum size=10pt] (c4l1t1) at (12,1){};
            \node[thick,draw=black,fill=white,minimum size=10pt] (c4l1t2) at (12.1,.9){};
            \node[thick,draw=black,fill=white,minimum size=10pt] (c4l1t3) at (12.2,.8){};
            \node[thick,draw=black,fill=white,minimum size=10pt] (c4l1t4) at (12.3,.8){};
            \node[thick,draw=black,fill=white,minimum size=10pt] (c4l2t1) at (12,2){};
            \node[thick,draw=black,fill=white,minimum size=10pt] (c4l2t2) at (12.1,1.9){};
            \node[thick,draw=black,fill=white,minimum size=10pt] (c4l2t3) at (12.2,1.8){};
            \node[thick,draw=black,fill=white,minimum size=10pt] (c4l2t4) at (12.3,1.8){};
            \node[thick,draw=black,fill=white,minimum size=10pt] (c4l3t1) at (12,3){};
            \node[thick,draw=black,fill=white,minimum size=10pt] (c4l3t2) at (12.1,2.9){};
            \node[thick,draw=black,fill=white,minimum size=10pt] (c4l3t3) at (12.2,2.8){};
            \node[thick,draw=black,fill=white,minimum size=10pt] (c4l3t4) at (12.3,2.8){};
            \node[thick,draw=black,fill=black,minimum size=10pt] (Dl0t1) at (15,0){};
            \node[thick,draw=black,fill=black,minimum size=10pt] (Dl0t2) at (15.1,-.1){};
            \node[thick,draw=black,fill=black,minimum size=10pt] (Dl0t3) at (15.2,-.2){};
            \node[thick,draw=black,fill=black,minimum size=10pt] (Dl0t4) at (15.3,-.3){};
            \node[thick,draw=black,fill=black,minimum size=10pt] (Dl1t1) at (15,1){};
            \node[thick,draw=black,fill=black,minimum size=10pt] (Dl1t2) at (15.1,.9){};
            \node[thick,draw=black,fill=black,minimum size=10pt] (Dl1t3) at (15.2,.8){};
            \node[thick,draw=black,fill=black,minimum size=10pt] (Dl1t4) at (15.3,.8){};
            \node[thick,draw=black,fill=black,minimum size=10pt] (Dl2t1) at (15,2){};
            \node[thick,draw=black,fill=black,minimum size=10pt] (Dl2t2) at (15.1,1.9){};
            \node[thick,draw=black,fill=black,minimum size=10pt] (Dl2t3) at (15.2,1.8){};
            \node[thick,draw=black,fill=black,minimum size=10pt] (Dl2t4) at (15.3,1.8){};
            \node[thick,draw=black,fill=black,minimum size=10pt] (Dl3t1) at (15,3){};
            \node[thick,draw=black,fill=black,minimum size=10pt] (Dl3t2) at (15.1,2.9){};
            \node[thick,draw=black,fill=black,minimum size=10pt] (Dl3t3) at (15.2,2.8){};
            \node[thick,draw=black,fill=black,minimum size=10pt] (Dl3t4) at (15.3,2.8){};
            \node[draw=black,fill=white,minimum size=2pt,circle] (s1l0t1) at (1,0){};
            \node[draw=black,fill=white,minimum size=2pt,circle] (s1l0t2) at (1.1,-.1){};
            \node[draw=black,fill=white,minimum size=2pt,circle] (s1l0t3) at (1.2,-.2){};
            \node[draw=black,fill=white,minimum size=2pt,circle] (s1l0t4) at (1.3,-.3){};
            \node[draw=black,fill=white,minimum size=2pt,circle] (s1l1t1) at (1,1){};
            \node[draw=black,fill=white,minimum size=2pt,circle] (s1l1t2) at (1.1,.9){};
            \node[draw=black,fill=white,minimum size=2pt,circle] (s1l1t3) at (1.2,.8){};
            \node[draw=black,fill=white,minimum size=2pt,circle] (s1l1t4) at (1.3,.7){};
            \node[draw=black,fill=white,minimum size=2pt,circle] (s1l2t1) at (1,2){};
            \node[draw=black,fill=white,minimum size=2pt,circle] (s1l2t2) at (1.1,1.9){};
            \node[draw=black,fill=white,minimum size=2pt,circle] (s1l2t3) at (1.2,1.8){};
            \node[draw=black,fill=white,minimum size=2pt,circle] (s1l2t4) at (1.3,1.7){};
            \node[draw=black,fill=white,minimum size=2pt,circle] (s1l3t1) at (1,3){};
            \node[draw=black,fill=white,minimum size=2pt,circle] (s1l3t2) at (1.1,2.9){};
            \node[draw=black,fill=white,minimum size=2pt,circle] (s1l3t3) at (1.2,2.8){};
            \node[draw=black,fill=white,minimum size=2pt,circle] (s1l3t4) at (1.3,2.7){};
            \node[draw=black,fill=white,minimum size=2pt,circle] (s2l0t1) at (2,0){};
            \node[draw=black,fill=white,minimum size=2pt,circle] (s2l0t2) at (2.1,-.1){};
            \node[draw=black,fill=white,minimum size=2pt,circle] (s2l0t3) at (2.2,-.2){};
            \node[draw=black,fill=white,minimum size=2pt,circle] (s2l0t4) at (2.3,-.3){};
            \node[draw=black,fill=white,minimum size=2pt,circle] (s2l1t1) at (2,1){};
            \node[draw=black,fill=white,minimum size=2pt,circle] (s2l1t2) at (2.1,.9){};
            \node[draw=black,fill=white,minimum size=2pt,circle] (s2l1t3) at (2.2,.8){};
            \node[draw=black,fill=white,minimum size=2pt,circle] (s2l1t4) at (2.3,.7){};
            \node[draw=black,fill=white,minimum size=2pt,circle] (s2l2t1) at (2,2){};
            \node[draw=black,fill=white,minimum size=2pt,circle] (s2l2t2) at (2.1,1.9){};
            \node[draw=black,fill=white,minimum size=2pt,circle] (s2l2t3) at (2.2,1.8){};
            \node[draw=black,fill=white,minimum size=2pt,circle] (s2l2t4) at (2.3,1.7){};
            \node[draw=black,fill=white,minimum size=2pt,circle] (s2l3t1) at (2,3){};
            \node[draw=black,fill=white,minimum size=2pt,circle] (s2l3t2) at (2.1,2.9){};
            \node[draw=black,fill=white,minimum size=2pt,circle] (s2l3t3) at (2.2,2.8){};
            \node[draw=black,fill=white,minimum size=2pt,circle] (s2l3t4) at (2.3,2.7){};
            \node[draw=black,fill=white,minimum size=2pt,circle] (s4l0t1) at (4,0){};
            \node[draw=black,fill=white,minimum size=2pt,circle] (s4l0t2) at (4.1,-.1){};
            \node[draw=black,fill=white,minimum size=2pt,circle] (s4l0t3) at (4.2,-.2){};
            \node[draw=black,fill=white,minimum size=2pt,circle] (s4l0t4) at (4.3,-.3){};
            \node[draw=black,fill=white,minimum size=2pt,circle] (s4l1t1) at (4,1){};
            \node[draw=black,fill=white,minimum size=2pt,circle] (s4l1t2) at (4.1,.9){};
            \node[draw=black,fill=white,minimum size=2pt,circle] (s4l1t3) at (4.2,.8){};
            \node[draw=black,fill=white,minimum size=2pt,circle] (s4l1t4) at (4.3,.7){};
            \node[draw=black,fill=white,minimum size=2pt,circle] (s4l2t1) at (4,2){};
            \node[draw=black,fill=white,minimum size=2pt,circle] (s4l2t2) at (4.1,1.9){};
            \node[draw=black,fill=white,minimum size=2pt,circle] (s4l2t3) at (4.2,1.8){};
            \node[draw=black,fill=white,minimum size=2pt,circle] (s4l2t4) at (4.3,1.7){};
            \node[draw=black,fill=white,minimum size=2pt,circle] (s4l3t1) at (4,3){};
            \node[draw=black,fill=white,minimum size=2pt,circle] (s4l3t2) at (4.1,2.9){};
            \node[draw=black,fill=white,minimum size=2pt,circle] (s4l3t3) at (4.2,2.8){};
            \node[draw=black,fill=white,minimum size=2pt,circle] (s4l3t4) at (4.3,2.7){};
            \node[draw=black,fill=white,minimum size=2pt,circle] (s5l0t1) at (5,0){};
            \node[draw=black,fill=white,minimum size=2pt,circle] (s5l0t2) at (5.1,-.1){};
            \node[draw=black,fill=white,minimum size=2pt,circle] (s5l0t3) at (5.2,-.2){};
            \node[draw=black,fill=white,minimum size=2pt,circle] (s5l0t4) at (5.3,-.3){};
            \node[draw=black,fill=white,minimum size=2pt,circle] (s5l1t1) at (5,1){};
            \node[draw=black,fill=white,minimum size=2pt,circle] (s5l1t2) at (5.1,.9){};
            \node[draw=black,fill=white,minimum size=2pt,circle] (s5l1t3) at (5.2,.8){};
            \node[draw=black,fill=white,minimum size=2pt,circle] (s5l1t4) at (5.3,.7){};
            \node[draw=black,fill=white,minimum size=2pt,circle] (s5l2t1) at (5,2){};
            \node[draw=black,fill=white,minimum size=2pt,circle] (s5l2t2) at (5.1,1.9){};
            \node[draw=black,fill=white,minimum size=2pt,circle] (s5l2t3) at (5.2,1.8){};
            \node[draw=black,fill=white,minimum size=2pt,circle] (s5l2t4) at (5.3,1.7){};
            \node[draw=black,fill=white,minimum size=2pt,circle] (s5l3t1) at (5,3){};
            \node[draw=black,fill=white,minimum size=2pt,circle] (s5l3t2) at (5.1,2.9){};
            \node[draw=black,fill=white,minimum size=2pt,circle] (s5l3t3) at (5.2,2.8){};
            \node[draw=black,fill=white,minimum size=2pt,circle] (s5l3t4) at (5.3,2.7){};
            \node[draw=black,fill=white,minimum size=2pt,circle] (s7l0t1) at (7,0){};
            \node[draw=black,fill=white,minimum size=2pt,circle] (s7l0t2) at (7.1,-.1){};
            \node[draw=black,fill=white,minimum size=2pt,circle] (s7l0t3) at (7.2,-.2){};
            \node[draw=black,fill=white,minimum size=2pt,circle] (s7l0t4) at (7.3,-.3){};
            \node[draw=black,fill=white,minimum size=2pt,circle] (s7l1t1) at (7,1){};
            \node[draw=black,fill=white,minimum size=2pt,circle] (s7l1t2) at (7.1,.9){};
            \node[draw=black,fill=white,minimum size=2pt,circle] (s7l1t3) at (7.2,.8){};
            \node[draw=black,fill=white,minimum size=2pt,circle] (s7l1t4) at (7.3,.7){};
            \node[draw=black,fill=white,minimum size=2pt,circle] (s7l2t1) at (7,2){};
            \node[draw=black,fill=white,minimum size=2pt,circle] (s7l2t2) at (7.1,1.9){};
            \node[draw=black,fill=white,minimum size=2pt,circle] (s7l2t3) at (7.2,1.8){};
            \node[draw=black,fill=white,minimum size=2pt,circle] (s7l2t4) at (7.3,1.7){};
            \node[draw=black,fill=white,minimum size=2pt,circle] (s7l3t1) at (7,3){};
            \node[draw=black,fill=white,minimum size=2pt,circle] (s7l3t2) at (7.1,2.9){};
            \node[draw=black,fill=white,minimum size=2pt,circle] (s7l3t3) at (7.2,2.8){};
            \node[draw=black,fill=white,minimum size=2pt,circle] (s7l3t4) at (7.3,2.7){};
            \node[draw=black,fill=white,minimum size=2pt,circle] (s8l0t1) at (8,0){};
            \node[draw=black,fill=white,minimum size=2pt,circle] (s8l0t2) at (8.1,-.1){};
            \node[draw=black,fill=white,minimum size=2pt,circle] (s8l0t3) at (8.2,-.2){};
            \node[draw=black,fill=white,minimum size=2pt,circle] (s8l0t4) at (8.3,-.3){};
            \node[draw=black,fill=white,minimum size=2pt,circle] (s8l1t1) at (8,1){};
            \node[draw=black,fill=white,minimum size=2pt,circle] (s8l1t2) at (8.1,.9){};
            \node[draw=black,fill=white,minimum size=2pt,circle] (s8l1t3) at (8.2,.8){};
            \node[draw=black,fill=white,minimum size=2pt,circle] (s8l1t4) at (8.3,.7){};
            \node[draw=black,fill=white,minimum size=2pt,circle] (s8l2t1) at (8,2){};
            \node[draw=black,fill=white,minimum size=2pt,circle] (s8l2t2) at (8.1,1.9){};
            \node[draw=black,fill=white,minimum size=2pt,circle] (s8l2t3) at (8.2,1.8){};
            \node[draw=black,fill=white,minimum size=2pt,circle] (s8l2t4) at (8.3,1.7){};
            \node[draw=black,fill=white,minimum size=2pt,circle] (s8l3t1) at (8,3){};
            \node[draw=black,fill=white,minimum size=2pt,circle] (s8l3t2) at (8.1,2.9){};
            \node[draw=black,fill=white,minimum size=2pt,circle] (s8l3t3) at (8.2,2.8){};
            \node[draw=black,fill=white,minimum size=2pt,circle] (s8l3t4) at (8.3,2.7){};
            \node[draw=black,fill=white,minimum size=2pt,circle] (s10l0t1) at (10,0){};
            \node[draw=black,fill=white,minimum size=2pt,circle] (s10l0t2) at (10.1,-.1){};
            \node[draw=black,fill=white,minimum size=2pt,circle] (s10l0t3) at (10.2,-.2){};
            \node[draw=black,fill=white,minimum size=2pt,circle] (s10l0t4) at (10.3,-.3){};
            \node[draw=black,fill=white,minimum size=2pt,circle] (s10l1t1) at (10,1){};
            \node[draw=black,fill=white,minimum size=2pt,circle] (s10l1t2) at (10.1,.9){};
            \node[draw=black,fill=white,minimum size=2pt,circle] (s10l1t3) at (10.2,.8){};
            \node[draw=black,fill=white,minimum size=2pt,circle] (s10l1t4) at (10.3,.7){};
            \node[draw=black,fill=white,minimum size=2pt,circle] (s10l2t1) at (10,2){};
            \node[draw=black,fill=white,minimum size=2pt,circle] (s10l2t2) at (10.1,1.9){};
            \node[draw=black,fill=white,minimum size=2pt,circle] (s10l2t3) at (10.2,1.8){};
            \node[draw=black,fill=white,minimum size=2pt,circle] (s10l2t4) at (10.3,1.7){};
            \node[draw=black,fill=white,minimum size=2pt,circle] (s10l3t1) at (10,3){};
            \node[draw=black,fill=white,minimum size=2pt,circle] (s10l3t2) at (10.1,2.9){};
            \node[draw=black,fill=white,minimum size=2pt,circle] (s10l3t3) at (10.2,2.8){};
            \node[draw=black,fill=white,minimum size=2pt,circle] (s10l3t4) at (10.3,2.7){};
            \node[draw=black,fill=white,minimum size=2pt,circle] (s11l0t1) at (11,0){};
            \node[draw=black,fill=white,minimum size=2pt,circle] (s11l0t2) at (11.1,-.1){};
            \node[draw=black,fill=white,minimum size=2pt,circle] (s11l0t3) at (11.2,-.2){};
            \node[draw=black,fill=white,minimum size=2pt,circle] (s11l0t4) at (11.3,-.3){};
            \node[draw=black,fill=white,minimum size=2pt,circle] (s11l1t1) at (11,1){};
            \node[draw=black,fill=white,minimum size=2pt,circle] (s11l1t2) at (11.1,.9){};
            \node[draw=black,fill=white,minimum size=2pt,circle] (s11l1t3) at (11.2,.8){};
            \node[draw=black,fill=white,minimum size=2pt,circle] (s11l1t4) at (11.3,.7){};
            \node[draw=black,fill=white,minimum size=2pt,circle] (s11l2t1) at (11,2){};
            \node[draw=black,fill=white,minimum size=2pt,circle] (s11l2t2) at (11.1,1.9){};
            \node[draw=black,fill=white,minimum size=2pt,circle] (s11l2t3) at (11.2,1.8){};
            \node[draw=black,fill=white,minimum size=2pt,circle] (s11l2t4) at (11.3,1.7){};
            \node[draw=black,fill=white,minimum size=2pt,circle] (s11l3t1) at (11,3){};
            \node[draw=black,fill=white,minimum size=2pt,circle] (s11l3t2) at (11.1,2.9){};
            \node[draw=black,fill=white,minimum size=2pt,circle] (s11l3t3) at (11.2,2.8){};
            \node[draw=black,fill=white,minimum size=2pt,circle] (s11l3t4) at (11.3,2.7){};
            \node[draw=black,fill=white,minimum size=2pt,circle] (s13l0t1) at (13,0){};
            \node[draw=black,fill=white,minimum size=2pt,circle] (s13l0t2) at (13.1,-.1){};
            \node[draw=black,fill=white,minimum size=2pt,circle] (s13l0t3) at (13.2,-.2){};
            \node[draw=black,fill=white,minimum size=2pt,circle] (s13l0t4) at (13.3,-.3){};
            \node[draw=black,fill=white,minimum size=2pt,circle] (s13l1t1) at (13,1){};
            \node[draw=black,fill=white,minimum size=2pt,circle] (s13l1t2) at (13.1,.9){};
            \node[draw=black,fill=white,minimum size=2pt,circle] (s13l1t3) at (13.2,.8){};
            \node[draw=black,fill=white,minimum size=2pt,circle] (s13l1t4) at (13.3,.7){};
            \node[draw=black,fill=white,minimum size=2pt,circle] (s13l2t1) at (13,2){};
            \node[draw=black,fill=white,minimum size=2pt,circle] (s13l2t2) at (13.1,1.9){};
            \node[draw=black,fill=white,minimum size=2pt,circle] (s13l2t3) at (13.2,1.8){};
            \node[draw=black,fill=white,minimum size=2pt,circle] (s13l2t4) at (13.3,1.7){};
            \node[draw=black,fill=white,minimum size=2pt,circle] (s13l3t1) at (13,3){};
            \node[draw=black,fill=white,minimum size=2pt,circle] (s13l3t2) at (13.1,2.9){};
            \node[draw=black,fill=white,minimum size=2pt,circle] (s13l3t3) at (13.2,2.8){};
            \node[draw=black,fill=white,minimum size=2pt,circle] (s13l3t4) at (13.3,2.7){};
            \node[draw=black,fill=white,minimum size=2pt,circle] (s14l0t1) at (14,0){};
            \node[draw=black,fill=white,minimum size=2pt,circle] (s14l0t2) at (14.1,-.1){};
            \node[draw=black,fill=white,minimum size=2pt,circle] (s14l0t3) at (14.2,-.2){};
            \node[draw=black,fill=white,minimum size=2pt,circle] (s14l0t4) at (14.3,-.3){};
            \node[draw=black,fill=white,minimum size=2pt,circle] (s14l1t1) at (14,1){};
            \node[draw=black,fill=white,minimum size=2pt,circle] (s14l1t2) at (14.1,.9){};
            \node[draw=black,fill=white,minimum size=2pt,circle] (s14l1t3) at (14.2,.8){};
            \node[draw=black,fill=white,minimum size=2pt,circle] (s14l1t4) at (14.3,.7){};
            \node[draw=black,fill=white,minimum size=2pt,circle] (s14l2t1) at (14,2){};
            \node[draw=black,fill=white,minimum size=2pt,circle] (s14l2t2) at (14.1,1.9){};
            \node[draw=black,fill=white,minimum size=2pt,circle] (s14l2t3) at (14.2,1.8){};
            \node[draw=black,fill=white,minimum size=2pt,circle] (s14l2t4) at (14.3,1.7){};
            \node[draw=black,fill=white,minimum size=2pt,circle] (s14l3t1) at (14,3){};
            \node[draw=black,fill=white,minimum size=2pt,circle] (s14l3t2) at (14.1,2.9){};
            \node[draw=black,fill=white,minimum size=2pt,circle] (s14l3t3) at (14.2,2.8){};
            \node[draw=black,fill=white,minimum size=2pt,circle] (s14l3t4) at (14.3,2.7){};
            \draw[->,myred,ultra thick] (O) to (s1l0t1);
            \draw[->,myred,ultra thick] (s1l0t1) to (s2l1t1);
            \draw[->,myred,ultra thick] (s2l1t1) to (c1l1t1);
            \draw[->,myred,ultra thick] (c1l1t1) to (s4l1t1);
            \draw[->,myred,ultra thick] (s4l1t1) to (s5l1t2);
            \draw[->,myred,ultra thick] (s5l1t2) to (c2l1t2);
            \draw[->,myred,ultra thick] (c2l1t2) to (s7l1t2);
            \draw[->,myred,ultra thick] (s7l1t2) to (s8l1t2);
            \draw[->,myred,ultra thick] (s8l1t2) to (c3l1t3);
            \draw[->,myred,ultra thick] (c3l1t3) to (s10l2t3);
            \draw[->,myred,ultra thick] (s10l2t3) to (s11l3t3);
            \draw[->,myred,ultra thick] (s11l3t3) to (c4l3t3);
            \draw[->,myred,ultra thick] (c4l3t3) to (s13l3t4);
            \draw[->,myred,ultra thick] (s13l3t4) to (s14l3t4);
            \draw[->,myred,ultra thick] (s14l3t4) to (Dl3t4);
            \draw[->,mygreen,ultra thick] (O) to (s1l0t1);
            \draw[->,mygreen,ultra thick] (s1l0t1) to (s2l0t1);
            \draw[->,mygreen,ultra thick] (s2l0t1) to (c1l0t1);
            \draw[->,mygreen,ultra thick] (c1l0t1) to (s4l0t1);
            \draw[->,mygreen,ultra thick] (s4l0t1) to (s5l0t2);
            \draw[->,mygreen,ultra thick] (s5l0t2) to (c2l0t2);
            \draw[->,mygreen,ultra thick] (c2l0t2) to (s7l0t2);
            \draw[->,mygreen,ultra thick] (s7l0t2) to (s8l0t2);
            \draw[->,mygreen,ultra thick] (s8l0t2) to (c3l0t3);
            \draw[->,mygreen,ultra thick] (c3l0t3) to (s10l0t3);
            \draw[->,mygreen,ultra thick] (s10l0t3) to (s11l1t3);
            \draw[->,mygreen,ultra thick] (s11l1t3) to (c4l1t3);
            \draw[->,mygreen,ultra thick] (c4l1t3) to (s13l2t4);
            \draw[->,mygreen,ultra thick] (s13l2t4) to (s14l2t4);
            \draw[->,mygreen,ultra thick] (s14l2t4) to (Dl3t4);
            \draw[->,myblue,ultra thick] (O) to (s1l1t1);
            \draw[->,myblue,ultra thick] (s1l1t1) to (s2l2t1);
            \draw[->,myblue,ultra thick] (s2l2t1) to (c1l2t1);
            \draw[->,myblue,ultra thick] (c1l2t1) to (s4l2t1);
            \draw[->,myblue,ultra thick] (s4l2t1) to (s5l2t2);
            \draw[->,myblue,ultra thick] (s5l2t2) to (c2l2t2);
            \draw[->,myblue,ultra thick] (c2l2t2) to (s7l2t2);
            \draw[->,myblue,ultra thick] (s7l2t2) to (s8l2t2);
            \draw[->,myblue,ultra thick] (s8l2t2) to (c3l2t3);
            \draw[->,myblue,ultra thick] (c3l2t3) to (s10l2t3);
            \draw[->,myblue,ultra thick] (s10l2t3) to (s11l2t3);
            \draw[->,myblue,ultra thick] (s11l2t3) to (c4l2t3);
            \draw[->,myblue,ultra thick] (c4l2t3) to (s13l2t4);
            \draw[->,myblue,ultra thick] (s13l2t4) to (s14l3t4);
            \draw[->,myblue,ultra thick] (s14l3t4) to (Dl3t4);
        \end{tikzpicture}}
        \caption{Visualization of the path-based, subpath-based and segment-based formulations.}
        \label{fig:formulations}
        \vspace{-12pt}
\end{figure}

Proposition~\ref{subpath_form} shows that the three formulations are equivalent, as long as time discretization is sufficiently granular in the segment-based benchmark (we formalize this condition in~\ref{A:subpath_structure}). The segment-based benchmark induces a weaker relaxation due to the double flow structure with linking constraints. In constrast, the subpath-based formulation achieves an equally strong relaxation as the path-based benchmark thanks to the flow balance structure on the load-expanded network. Most importantly, Proposition~\ref{form_complexity} shows the size benefits of the subpath-based formulation. The segment-based formulation features a polynomial number of variables in a dense time-station-load network (Figure~\ref{subfig:segment}). In the path-based formulation, the number of variables scales exponentially with the total number of stations along the reference line (Figure~\ref{subfig:path}). In the subpath-based formulation, the number of variables scales exponentially with the number of stations between checkpoints, and relies on a coarser checkpoint-load network representation (Figure~\ref{subfig:subpath}). The proposition also shows that the model becomes increasingly complex as more checkpoints can be skipped (higher $K$).

\begin{proposition}\label{subpath_form} 
    The path-based and subpath-based formulations are equivalent and define identical linear relaxations. If all subpath travel times are strictly less than the elapsed time between the scheduled arrival times at the checkpoints, there exists a time discretization such that the segment-based formulation is also equivalent but its linear relaxation is at most as strong.
\end{proposition}

\begin{proposition}\label{form_complexity}
    Consider the second-stage problem for reference trip $(\ell, t) \in \calL \times \calT_{\ell}$ in scenario $s \in \calS$. Let $\Xi$ be the maximum number of stations between any pair of checkpoints in $\Gamma_\ell$. The segment-based formulation has $\calO(T_S \cdot C_{\ell}^2 \cdot I_\ell\cdot \Xi^2)$ variables and $\calO(|\calP| + T_S \cdot C_{\ell} \cdot |\calN^S| + T_S\cdot C_{\ell}^2 \cdot I_\ell\cdot \Xi^2)$ constraints. The subpath-based formulation has $\calO(I_{\ell} \cdot C_{\ell} \cdot 2^\Xi)$ variables and $\calO(|\calP| + C_{\ell} \cdot I_\ell)$ constraints. The path-based formulation has $\calO(2^{\Xi\cdot I_{\ell}})$ variables and $\calO\left(|\calP|\right)$ constraints.
\end{proposition}
\section{Double-Decomposition (DD) Algorithm}\label{sec:alg}

The MiND-VRP exhibits a two-stage stochastic optimization structure with a tight recourse function and exponentially many second-stage variables. We propose a double decomposition algorithm to solve its partial relaxation with discrete first-stage variables and continuous second-stage variables. The methodology relies on Benders decomposition to exploit the nested block-angular structure (Section~\ref{sec:4Benders}), and on subpath-based column generation in the Benders subproblem (Section~\ref{sec:4CG}). We formalize the algorithm and establish its exactness in Section~\ref{sec:4overview}, and augment it with integer L-shaped cuts and UB\&BC to retrieve an exact algorithm for Problem~\eqref{OPT} in Section~\ref{subsec:integer}. For generalizability, we describe the algorithm using the general-purpose notation from Problem~\eqref{OPT}; we develop it for the MiND-VRP in~\ref{app:VRPalg}, for the MiND-DAR in~\ref{app:DARalg}, and for the MiND-Tr in~\ref{app:TRalg}. Section~\ref{subsec:label} describes the label-setting algorithm for the pricing problem. 

\subsection{Multi-cut Benders Decomposition}\label{sec:4Benders}

Let $\texttt{OPT}$ denote the optimal value of~\eqref{OPT}. We denote the second-stage problem for a first-stage solution $\bx$ by $\texttt{SP}(\bx)$, and its optimal solution by $\varPhi(\bx)$. Due to its nested block-angular structure, $\texttt{SP}(\bx)$ is decomposable across $s\in\calS$ and $j\in\calJ$. Namely, $\varPhi(\bx)=\sum_{s\in\calS}\sum_{j\in\calJ}\pi_s\varphi_{sj}(\bx)$ where:
\begin{align}
    \varphi_{sj}(\bx)=\min\quad&\sum_{a\in\calA_{sj}}g_{a}y_{a}\label{eq:SP_obj}\\
    \st\quad    &   \sum_{m:(n,m)\in\calA_{sj}}y_{(n,m)}-\sum_{m:(m,n)\in\calA_{sj}}y_{(m,n)}=\sum_{k\in\calK_j}b_{nsk}x_k,\ \forall n\in\calN_{sj}\label{eq:SP_FB}\\
                &   \sum_{a\in\calA_{sj}}\bff_{asj} y_{a}\geq\bh_{sj}\label{eq:SP_side}\\
                &   y_{a}\in\calY^{\text{MIO}}_a,\ \forall a\in\calA_{sj}\label{eq:SP_domain}
\end{align}

We define and evaluate the total cost associated with each first-stage solution as follows:
$$\texttt{OPT}(\bx)=\bc^\top\bx+\varPhi(\bx)$$

We refer to the partial relaxation as $(\texttt{MIO}-\texttt{LO})$. We define the Benders subproblem (BSP) as its second-stage relaxation $\overline{\texttt{SP}}(\bx)$, with optimal value $\overline\varPhi(\bx)=\sum_{s\in\calS}\sum_{j\in\calJ}\pi_s\overline\varphi_{sj}(\bx)$ where:
\begin{align*}
    \overline\varphi_{sj}(\bx)=\min\quad\left\{\sum_{a\in\calA_{sj}}g_{a}y_{a}\ :\ \text{Equations~\eqref{eq:SP_FB}--\eqref{eq:SP_side}},\ \by \geq \bo\right\}
\end{align*}

Its dual formulation in scenario $s\in\calS$ and for partition element $j\in\calJ$ is given as follows, using $\bpsi$ and $\bgamma$ to denote the dual variables of the flow balance constraints and the side constraints.
\begin{align}
    & \max\quad
    \left\{
    \sum_{n\in\calN_{sj}}\sum_{k\in\calK_j}b_{nsk}x_k\psi_{sjn}+\bh_{sj}^\top\bgamma_{sj}
    :
    (\bpsi_{sj},\bgamma_{sj})\in\calP_{sj}
    \right\},\quad\text{where:} \label{cg:dualsp}\\
    & \calP_{sj}=\left\{(\bpsi_{sj},\bgamma_{sj})\ :\ \bgamma_{sj}\geq\bo;\ \psi_{sjm}-\psi_{sjn}+\bgamma_{sj}^\top\bff_{(m,n),s,j}\leq g_{(m,n)},\ \forall (m,n)\in\calA_{sj}\right\}\nonumber
\end{align}

We index the extreme points of the dual polyhedron $\calP_{sj}$ by $\{(\bpsi^u_{sj},\bgamma^u_{sj}):u\in\calU_{sj}\}$ and its extreme rays by $\{(\bpsi^v_{sj},\bgamma^v_{sj}):v\in\calV_{sj}\}$. The Benders master problem, denoted by~\ref{MP}, is given by
\begin{align}
    \min\quad   &   \sum_{j \in \calJ}\sum_{k\in\calK_j}c_kx_k+\sum_{s\in\calS}\sum_{j\in\calJ}\pi_s\theta_{sj}\label{MP}\tag{$\texttt{MP}(\calU^0,\calV^0)$}\\
    \st\quad    &   \bA\bx\geq\bb\nonumber\\
                &   \theta_{sj}\geq\sum_{n\in\calN_{sj}}\sum_{k\in\calK_j}b_{nsk}x_k\psi^u_{sjn}+\bh_{sj}^\top\bgamma^u_{sj},\ \forall s\in\calS,\ \forall j\in\calJ,\ \forall u\in\calU_{sj}^0 \label{bendersoptcut}\\
                &   0\geq\sum_{n\in\calN_{sj}}\sum_{k\in\calK_j}b_{nsk}x_k\psi^v_{sjn}+\bh_{sj}^\top\bgamma^v_{sj},\ \forall s\in\calS,\ \forall j\in\calJ,\ \forall v\in\calV_{sj}^0 \label{bendersfeascut}\\
                &   \bx\in\calX^{\text{MIO}}\nonumber
\end{align}

The Benders master problem solves a relaxation $\texttt{MP}(\calU^0,\calV^0)$ containing a subset of constraints indexed by $\calU^0_{sj}\subseteq\calU_{sj}$ and $\calV^0_{sj}\subseteq\calV_{sj}$. By design, it yields a lower bound of $(\texttt{MIO}-\texttt{LO})$ and its combination with the $\overline{\texttt{SP}}(\bx)$ yield an upper bound of $(\texttt{MIO}-\texttt{LO})$. If the gap lies within a given tolerance, Benders decomposition stops. Otherwise, we add optimality cuts by augmenting $\calU^0_{sj}$ with an optimal extreme point $(\bpsi^u_{sj},\bgamma^u_{sj})$ for all scenario-partition combinations $s\in\calS,\ j\in\calJ$ such that the Benders subproblem admits an optimal solution, and by augmenting $\calV^0_{sj}$ with an extreme ray defining a direction of unboundeness $(\bpsi^v_{sj},\bgamma^v_{sj})$ for all others.

\subsection{Subpath-based Column Generation for Benders Subproblem}\label{sec:4CG}

The main challenge with the Benders decomposition scheme lies in the large number of second-stage variables $y_{a}$ in the expanded network representation. In the MiND, the arc set $\mathcal{A}_{\ell st}$ grows exponentially with the number of stations between checkpoints. To prevent the subproblem $\texttt{SP}(\bx)$ from becoming too computationally intensive at each iteration, we propose a column generation procedure iterating between a restricted Benders subproblem and a pricing problem.

The restricted Benders subproblem (RBSP) is formulated as $\overline{\texttt{SP}}(\bx)$, except that the arc sets $\calA_{sj}$ are replaced by restricted arc sets $\calA'_{sj}\subseteq\calA_{sj}$. We refer to it as $\texttt{RSP}(\bx,\calA'_{sj})$ and to its optimal value as $\varPhi'(\bx,\calA'_{sj})$. Namely, $\varPhi'(\bx,\calA'_{sj})=\sum_{s\in\calS}\sum_{j\in\calJ}\pi_s\varphi'_{s,j}(\bx,\calA'_{sj})$ where:
\begin{align*}
    \varphi'_{s,j}(\bx,\calA'_{sj})=
    \min_{\by \geq \boldsymbol{0}}\quad   &   \sum_{a\in\calA'_{sj}}g_{a}y_{a}\\
    \st\quad    &   \sum_{m:(n,m)\in\calA'_{sj}}y_{(n,m)}-\sum_{m:(m,n)\in\calA'_{sj}}y_{(m,n)}=\sum_{k\in\calK_j}b_{nsk}x_k,\ \forall n\in\calN_{sj}\\
                &   \sum_{a\in\calA'_{sj}}\bff_{asj} y_{a}\geq\bh_{sj}
\end{align*}

The pricing problem, denoted $\texttt{PP}(\bpsi,\bgamma)$, serves as dual separation by seeking a variable of negative reduced cost or proving that none exists. It is given as follows for $s\in\calS$ and $j\in\calJ$:
\begin{equation}
    RC=\min\ \left\{g_{(m,n)}-\left(\psi_{sjm}-\psi_{sjn}\right)-\bgamma_{sj}^\top\bff_{(m,n),s,j}:(m,n)\in\calA_{sj}\right\}\label{eq:RC}
\end{equation}

\subsection{Combining Benders Decomposition and Subpath-based Column Generation}\label{sec:4overview}

Our solution algorithm, summarized in Figure~\ref{fig:multitree_DD}, involves two interconnected decomposition structures. In an outer Benders decomposition loop, the master problem generates a feasible first-stage solution and a lower bound; the Benders subproblem generates a second-stage fractional solution and an upper bound to the partial relaxation $(\texttt{MIO}-\texttt{LO})$. At each outer iteration, the algorithm adds optimality and feasibility cuts, or certifies the optimality of the $(\texttt{MIO}-\texttt{LO})$ solution. The inner loop solves the Benders subproblem via subpath-based column generation: the restricted Benders subproblem generates a feasible solution; the pricing problem identifies variables with negative reduced cost for each scenario $s\in\calS$ and partition $j\in\calJ$, or yields a certificate of optimality in the restricted Benders subproblem. The algorithm continues until convergence in both loops.\footnote{Figure~\ref{fig:multitree_DD} shows a multi-tree implementation of the DD algorithm, based on cutting-plane version of Benders decomposition; for completeness, we describe the single-tree implementation of the algorithm in Figure~\ref{fig:singletree_DD}, based on lazy constraints within the branch-and-cut algorithm \citep{fortz2009improved,bodur2017mixed}.}

\begin{figure}[h!]
    \centering
    \begin{tikzpicture}[scale=0.4, every node/.style={scale=0.75},every text node part/.style={align=center}]
        \tikzset{test/.style={diamond,fill=mygray!20,draw,aspect=2,minimum width = 50pt, minimum height=30pt,text badly centered}}
        \tikzset{reg/.style={rectangle,draw, minimum width = 60pt, minimum height=20pt}}

        \node[reg] (LO) at (-10.5,13){Choose node, solve $\texttt{MP}(\calU^0,\calV^0)$\\Solution $(\bx,\by)$, objective $MP$};
        
        \node[test] (bad) at (-10.5,8.5){Infeasible};
        \node[reg,anchor=north,rectangle,draw] (terminateF) at (-10.5,4.5){Terminate\\Problem infeasible};
        \draw[->,thick] (bad.south) -- node[right,midway]{YES} (terminateF.north);
        
        \node[reg,fill=myred!20] (BSP) at (2,8.5){Solve $\texttt{RSP}(\bx,\calA'_{sj})$};
        \node[reg,fill=myred!20] (PP) at (2,4){Solve $\texttt{PP}(\bpsi,\bgamma)$};
        \node[test] (RCneg) at (2,0){$RC<0$};
        \node[test] (violation) at (12,8.5){Violated cuts};
        
        \draw[->,thick,color=myred] (BSP.south) -- (PP.north);
        \draw[->,thick,color=myred] (PP.south) -- (RCneg);
        \draw[->,thick,color=myred] (RCneg.west) -- node[above,midway]{YES} +(-2,0) |- (BSP.south west);
        \draw[->,thick,color=myred] (RCneg.east) -- node[above,midway]{NO} +(2,0) |- (violation.west);
        
        \node[reg,anchor=north,rectangle,draw] (terminateO) at (12,4.5){Terminate\\Return optimal solution};
        \draw[->,thick] (violation.south) -- node[right,very near start]{NO} (terminateO.north);
        
        \node[reg,draw] (Benders) at (2,13){Add Benders cuts\\Expand $\calU^0,\calV^0$};
        
        \draw[->,thick] (LO.south) -- (bad.north);
        \draw[->,thick] (bad.east) -- node[above,midway]{NO}(BSP.west);
        \draw[->,thick] (violation.north) -- node[right,very near start]{YES} (12,13) -- (Benders.east);
        \draw[->,thick] (Benders.west) -- (LO.east);

        \draw[color=myred,dotted,ultra thick] (-3.5,-2) rectangle (7.5,10);
    \end{tikzpicture}
    \caption{Multi-tree DD algorithm to solve Problem~$(\texttt{MIO}-\texttt{LO})$.}
    \label{fig:multitree_DD}
    \vspace{-12pt}
\end{figure}

Theorem~\ref{thm:exact} shows the exactness of the DD algorithm toward solving $(\texttt{MIO}-\texttt{LO})$. Upon termination, we solve the second-stage problem with integrality constraints to derive a feasible solution to Problem~\eqref{OPT}. Thus, the DD algorithm returns a valid optimality gap for Problem~\eqref{OPT}.

\begin{theorem}\label{thm:exact}
The double decomposition algorithm returns an optimal solution to Problem $(\texttt{MIO}-\texttt{LO})$ in a finite number of iterations along with an optimality gap for Problem~\eqref{OPT}.
\end{theorem}

This integrated Benders decomposition and column generation algorithm gives rise to our double-decomposition structure, illustrated in Figure~\ref{fig:DD} for the MiND. The master problem solves a network design and service scheduling problem with all reference lines in all scenarios (Figure~\ref{subfig:BMP}). The Benders subproblem is then separable across scenarios $s\in\calS$ and partition elements $j\in\calJ$ (i.e., across reference trips, in Figure~\ref{subfig:BSP}). Column generation exploits the network-based second-stage formulation to further decompose the Benders subproblem across arcs (i.e., across subpaths between checkpoints, in Figure~\ref{subfig:PP}). The other way around, the pricing problem adds subpaths locally between checkpoints; the restricted Benders subproblem combines them into full paths to optimize operations along each reference trip in each scenario; and the Benders master problem brings all reference trips together to optimize first-stage network design and service scheduling decisions.

\begin{figure}[h!]
    \centering
    \subfloat[Benders master problem]{\label{subfig:BMP}
        \begin{tikzpicture}[scale=0.55,transform shape]
            \foreach \i in {1,...,5}{\foreach \j in {1,...,7}{
                    \node[circle,draw=black,fill=white,minimum size=.5pt] at (1.3*\j,1.3*\i){};
            }}
            \node[circle,draw=myred,fill=myred,minimum size=.5pt] at (1.3,1.3){};
            \node[circle,draw=myred,fill=myred,minimum size=.5pt] at (3.9,1.3){};
            \node[circle,draw=myred,fill=myred,minimum size=.5pt] at (6.5,2.6){};
            \node[circle,draw=myred,fill=myred,minimum size=.5pt] at (6.5,5.2){};
            \node[circle,draw=myred,fill=myred,minimum size=.5pt] at (9.1,5.2){};
            \draw[ultra thick,color=myred] (1.3,1.3) -- (3.9,1.3);
            \draw[ultra thick,color=myred] (3.9,1.3) -- (3.9,2.6) -- (6.5,2.6);
            \draw[ultra thick,color=myred] (6.5,2.6) -- (6.5,5.2);
            \draw[->,ultra thick,color=myred] (6.5,5.2) -- (9.1,5.2);
            \node[circle,draw=myblue,fill=myblue,minimum size=.5pt] at (1.3,6.5){};
            \node[circle,draw=myblue,fill=myblue,minimum size=.5pt] at (3.9,6.5){};
            \node[circle,draw=myblue,fill=myblue,minimum size=.5pt] at (5.2,5.2){};
            \node[circle,draw=myblue,fill=myblue,minimum size=.5pt] at (5.2,1.3){};
            \node[circle,draw=myblue,fill=myblue,minimum size=.5pt] at (9.1,1.3){};
            \draw[ultra thick,color=myblue] (1.3,6.5) -- (3.9,6.5);
            \draw[ultra thick,color=myblue] (3.9,6.5) -- (5.2,6.5) -- (5.2,5.2);
            \draw[ultra thick,color=myblue] (5.2,5.2) -- (5.2,1.3);
            \draw[->,ultra thick,color=myblue] (5.2,1.3) -- (9.1,1.3);
            \node[circle,draw=mygreen,fill=mygreen,minimum size=.5pt] at (1.3,3.9){};
            \node[circle,draw=mygreen,fill=mygreen,minimum size=.5pt] at (3.9,3.9){};
            \node[circle,draw=mygreen,fill=mygreen,minimum size=.5pt] at (6.5,3.9){};
            \node[circle,draw=mygreen,fill=mygreen,minimum size=.5pt] at (9.1,3.9){};
            \draw[ultra thick,color=mygreen] (1.3,3.9) -- (3.9,3.9);
            \draw[ultra thick,color=mygreen] (3.9,3.9) -- (6.5,3.9);
            \draw[->,ultra thick,color=mygreen] (6.5,3.9) -- (9.1,3.9);
            \node[very thick,draw=myblue,fill=myblue,minimum size=10pt] at (1.9,6.0){};
            \node[very thick,draw=mygreen,fill=mygreen,minimum size=10pt] at (3.3,4.5){};
            \node[very thick,draw=myblue,fill=myblue,minimum size=10pt] at (4.5,5.4){};
            \node[very thick,draw=myred,fill=myred,minimum size=10pt] at (5.5,3.2){};
            \node[very thick,draw=myred,fill=myred,minimum size=10pt] at (3.2,2.2){};
            \node[very thick,draw=mygreen,fill=mygreen,minimum size=10pt] at (1.8,3.2){};
            \node[very thick,draw=myblue,fill=myblue,minimum size=10pt] at (7.2,1.8){};
            \node[very thick,draw=myblue,fill=myblue,minimum size=10pt] at (8.7,2.2){};
            \node[very thick,draw=myred,fill=myred,minimum size=10pt] at (8.3,6.0){};
            \node[diamond,very thick,draw=myblue,fill=myblue,minimum size=6pt] at (6,5.8){};
            \node[diamond,very thick,draw=myred,fill=myred,minimum size=6pt] at (1.9,2.6){};
            \node[diamond,very thick,draw=myred,fill=myred,minimum size=6pt] at (7.4,4.4){};
            \node[diamond,very thick,draw=myred,fill=myred,minimum size=6pt] at (4.3,3.3){};
            \node[diamond,very thick,draw=mygreen,fill=mygreen,minimum size=6pt] at (1.9,4.4){};
            \node[diamond,very thick,draw=myred,fill=myred,minimum size=6pt] at (5.6,1.7){};
            \node[diamond,very thick,draw=mygreen,fill=mygreen,minimum size=6pt] at (8.8,3.1){};
            \node[diamond,very thick,draw=myred,fill=myred,minimum size=6pt] at (8.5,4.4){};
        \end{tikzpicture}}\hspace{1 cm}
    \subfloat[Benders subproblem]{\label{subfig:BSP}
        \begin{tikzpicture}[scale=0.55,transform shape]
            \foreach \i in {1,...,5}{\foreach \j in {1,...,7}{
                    \node[circle,draw=black,fill=white,minimum size=.5pt] at (1.3*\j,1.3*\i){};
            }}
            \node[circle,draw=myred,fill=myred,minimum size=.5pt] at (1.3,1.3){};
            \node[circle,draw=myred,fill=myred,minimum size=.5pt] at (3.9,1.3){};
            \node[circle,draw=myred,fill=myred,minimum size=.5pt] at (6.5,2.6){};
            \node[circle,draw=myred,fill=myred,minimum size=.5pt] at (6.5,5.2){};
            \node[circle,draw=myred,fill=myred,minimum size=.5pt] at (9.1,5.2){};
            \draw[->,ultra thick,dotted,color=myred] (1.3,1.3) -- (3.9,1.3);
            \draw[->,draw=myred,ultra thick] (1.3,1.3) -- (1.3,2.0) -- (3.9,2.0) -- (3.9,1.3);
            \draw[ultra thick,color=myred,dotted] (3.9,1.3) -- (3.9,2.6) -- (6.5,2.6);
            \draw[->,draw=myred,ultra thick] (3.9,1.3) -- (5.6,1.3) -- (5.6,2.6) -- (6.5,2.6);
            \draw[->,draw=myred,ultra thick,dotted] (3.9,1.3) -- (4.3,1.3) -- (4.3,3.3) -- (5.2,3.3) -- (5.2,2.6) -- (6.5,2.6);
            \draw[ultra thick,color=myred] (6.5,2.6) -- (6.5,5.2);
            \draw[->,draw=myred,ultra thick,dotted] (6.5,2.6) -- (6.5,3.9) -- (6,3.9) -- (6,5.2) -- (6.5,5.2);
            \draw[->,draw=myred,ultra thick,dotted] (6.5,2.6) -- (7.8,2.6) -- (7.8,4.4) -- (6.5,4.4) -- (6.5,5.2);
            \draw[->,ultra thick,color=myred] (6.5,5.2) -- (9.1,5.2);
            \draw[->,draw=myred,ultra thick,dotted] (6.5,5.2) -- (6.5,5.8) -- (9.1,5.8) -- (9.1,5.2);
            \node[diamond,draw=myred,fill=myred,minimum size=6pt] at (1.9,2.6){};
            \node[diamond,draw=myred,fill=white,ultra thick,minimum size=6pt] at (7.4,4.4){};
            \node[diamond,draw=myred,fill=white,ultra thick,minimum size=6pt] at (4.3,3.3){};
            \node[diamond,draw=myred,fill=myred,minimum size=6pt] at (5.6,1.7){};
            \node[diamond,draw=myred,fill=myred,minimum size=6pt] at (8.5,4.4){};
        \end{tikzpicture}}\hspace{1 cm}
    \subfloat[Pricing problem]{\label{subfig:PP}
        \begin{tikzpicture}[scale=0.55,transform shape]
            \foreach \i in {1,...,5}{\foreach \j in {1,...,7}{
                    \node[circle,draw=black,fill=white,minimum size=.5pt] at (1.3*\j,1.3*\i){};
            }}
            \node[circle,draw=myred,fill=myred,minimum size=.5pt] at (6.5,5.2){};
            \node[circle,draw=myred,fill=myred,minimum size=.5pt] at (9.1,5.2){};
            \draw[->,ultra thick,color=myred,dotted] (6.5,5.2) -- (9.1,5.2);
            \draw[->,draw=myred,ultra thick,dotted] (6.5,5.2) -- (6.5,5.8) -- (9.1,5.8) -- (9.1,5.2);
            \draw[->,draw=myred,ultra thick] (6.5,5.2) -- (6.5,5.5) -- (7.8,5.5) -- (7.8,4.8) -- (9.1,4.8) -- (9.1,5.2);
            \node[diamond,draw=myred,fill=myred,minimum size=6pt] at (1.9,2.6){};
            \node[diamond,draw=myred,fill=myred,minimum size=6pt] at (7.4,4.4){};
            \node[diamond,draw=myred,fill=white,ultra thick,minimum size=6pt] at (4.3,3.3){};
            \node[diamond,draw=myred,fill=myred,minimum size=6pt] at (5.6,1.7){};
            \node[diamond,draw=myred,fill=myred,minimum size=6pt] at (8.5,4.4){};
        \end{tikzpicture}}
    \caption{Double-decomposition algorithm. {\it Left:}: BMP with three reference lines (blue, red, green); passenger requests in two scenarios (squares, diamonds) with their first-stage assignments (colors). {\it Middle:} BSP for one reference trip and one scenario; full diamonds encode served passengers; solid lines characterize selected subpaths in RBSP. {\it Right:} PP to generate new subpath between checkpoints (solid line).}
    \label{fig:DD}
\end{figure}

This algorithm raises two final questions. First, the DD structure solves the partial relaxation $(\texttt{MIO}-\texttt{LO})$ to optimality, so it may still leave an optimality gap in the full problem~\eqref{OPT}. As we shall see experimentally, the gap is very small in the MiND due to the tight second-stage formulation. Still, we augment the DD methodology toward an exact and finitely convergent algorithm for Problem~\eqref{OPT} in Section~\ref{subsec:integer}. Second, the scalability of the DD algorithm hinges on the efficiency of the pricing algorithm. In the MiND, we propose a label-setting algorithm that exploits an additional decomposition of the pricing problem into routing and load components, in Section~\ref{subsec:label}.

\subsection{Exact double-decomposition algorithms for Problem~\eqref{OPT}.}
\label{subsec:integer}

\paragraph{Double decomposition with integer L-shaped cuts (DD\&ILS).} This approach assumes that the first-stage variables are binary, i.e., $\calX^{\text{MIO}}=\{0,1\}^{n_Z}$, and that the problem has relatively complete recourse; both conditions are satisfied in the MiND. We adopt a multi-cut variant of the integer L-shaped method from \cite{laporte1993integer} by adding the following optimality cut to the Benders master problem, where $\underline{\varPhi}_{sj}$ denotes a global lower bound of the second-stage cost $\varphi_{sj}(\bx)$ and $\{\widehat\bx^{w}:w\in\calW^0\}$ indexes the (binary) first-stage variables visited through the algorithm.
\begin{align}
   \theta_{sj}\geq\underline{\varPhi}_{sj}+(\varphi_{sj}(\widehat\bx^{w})-\underline{\varPhi}_{sj})\left(1-\sum_{k\in\calK_j:\widehat{x}^{w}_k=1}(1-x_k)-\sum_{k\in\calK_j:\widehat{x}^{w}_k=0}x_k\right),\ \forall s\in\calS,j\in\calJ,w\in\calW^0 \label{ils:cut}
\end{align}

The DD\&ILS algorithm (Figure~\ref{fig:DD_ILS}) involves three interconnected loops. The two upper loops solve the partial relaxation $(\texttt{MIO}-\texttt{LO})$ with integer L-shaped cuts, via DD (Figure~\ref{fig:multitree_DD}). The lower loop adds new integer L-shaped cuts if the second-stage solution violates integrality requirements.

\begin{theorem}\label{thm:DDILS}
The DD\&ILS algorithm converges in a finite number of iterations to an optimal solution of Problem~\eqref{OPT}, if $\calX^{\text{MIO}}=\{0,1\}^{n_Z}$ and if the problem has relatively complete recourse.
\end{theorem}

\paragraph{Unified branch-and-double-decomposition algorithm (UB\&DD).}

The UB\&BC algorithm from \cite{maheo2024unified} develops a single-tree algorithm in two-stage stochastic mixed-integer programming. It relies on tailored branching rules whenever a solution satisfies first-stage integrality requirements and Benders cuts but violates second-stage integrality constraints. Throughout, it maintains lower bounds from the linear and Benders relaxations, and upper bounds from second-stage heuristics. A post-processing procedure solves second-stage mixed-integer problems to find the optimum out of all candidates first-stage solutions. In contrast, the DD methodology derives a tight network-based reformulation of the second-stage problem and solves the partial relaxation (with mixed-integer first-stage solutions and continuous second-stage solutions) via Benders decomposition and column generation. We propose a UB\&DD algorithm that combines these two elements---namely, the DD methodology to circumvent the exponential number of second-stage variables, and the UB\&BC methodology to restore second-stage integrality (Figure~\ref{fig:UBDD}).

\begin{theorem}\label{thm:UBDD}
The UB\&DD algorithm converges finitely to an optimal solution of Problem~\eqref{OPT}.
\end{theorem}

\subsection{Solving the pricing problem in the MiND-VRP.}
\label{subsec:label}

Consider two nodes in the load-expanded network $(u,c_1),(v,c_2)\in\calV_{\ell st}$. The pricing problem seeks a subpath that starts in checkpoint $u\in\calN^S$ at time $T_{\ell t}(u)$ with vehicle load $c_1$, and ends in checkpoint $v \in \calN^S$ at time $T_{\ell t}(v)$ with load $c_2\ge c_1$, while satisfying the maximum deviation.

We characterize subpaths in a time-expanded network $(\calU^{uv}_{\ell st}, \calH^{uv}_{\ell st})$. Let $\calT^{uv}_{\ell t}$ be a set of discretized time intervals between $T_{\ell t}(u)$ and $T_{\ell t}(v)$. As in the segment-based benchmark (Section~\ref{sec:problem_structure}), the sets $\calT^{uv}_{\ell t}$ need to be much more granular than the first-stage sets $\calT_{\ell}$ (30 seconds vs. 15 minutes, in our experiments). Of course, this discretization is only applied locally in each pricing problem, thus retaining a much more manageable structure than in the segment-based benchmark. Each node $m\in\calU^{uv}_{\ell st}$ is represented by a tuple $(k_{m}, t_{m})\in \calN^S_{uv}\times\calT_{\ell t}^{uv}$; $(u,T_{\ell t}(u)) \in \calU^{uv}_{\ell st}$ is the source node and $(v,T_{\ell t}(v)) \in \calU^{uv}_{\ell st}$ is the sink node. The arc set $\calH^{uv}_{\ell st}$ comprises traveling arcs connecting any node pair $(i,t) \to (j,t+tt_{ij})$ with travel time $tt_{ij}$, and idling arcs connecting $(i,t) \to (i,t+1)$. In particular, the graph $(\calU^{uv}_{\ell st}, \calH^{uv}_{\ell st})$ is acyclic due to time moving forward in the time-space network. Each node $m \in \calU^{uv}_{\ell st}$ also defines passengers' waiting, walking and travel times, as well as arrival delays and earliness, which we store in parameters $\tau^{\text{walk}}_{mp}$, $\tau^{\text{wait}}_{mp}$, $\tau^{\text{travel}}_{mp}$, $\tau^{\text{late}}_{mp}$, and $\tau^{\text{early}}_{mp}$. We denote by $\calP_{m}\subset\calP$ the set of passengers that can be picked up at node $m \in \calU^{uv}_{\ell st}$ given the walking and waiting restrictions (Section \ref{sec:subpath_formulation}). Table~\ref{T:notation} summarizes notation. We define the following variables:
\begin{align*}
    f_{mq}&=\begin{cases}1&\text{if arc $(m,q) \in \calH^{uv}_{\ell st}$ is traversed in the time-expanded road segment network,}\\0&\text{otherwise.}\end{cases}\\
    w_{mp}&=\begin{cases}1&\text{if passenger $p \in \calP_{m}$ is picked up in node $m\in\calU^{uv}_{\ell st}$,}\\0&\text{otherwise.}\end{cases}\\
    \xi_m&=\ \text{vehicle load in node $m \in \calU^{uv}_{\ell st}$}
\end{align*}

Let $\widehat{g}_{a}$ denote the reduced cost of arc-based variable $a=((u,c_1),(v,c_2)) \in \calA_{\ell st}$. In the MiND-VRP, each partition element $j\in\calJ$ corresponds to a reference trip $(\ell,t)\in\calL\times\calT_\ell$ and each node $n\in\calN_{sj}$ to a checkpoint-load pair $(u,c)\in\calV_{\ell st}$; the dual variables can be written as $\psi_{\ell, s, t, (u,c)}$ and $\gamma_{\ell s t p}$ (because the second-stage problem admits one side constraint per passenger, in Equation~\eqref{eq:2Sy2z}). From Equation~\eqref{eq:RC} (or Equation~\eqref{Benders:dualArc}), the reduced cost can be separated into a routing component and a load component. The routing component comprises (i) the level-of-service penalty for served passengers, and (ii) the value of serving a passenger, captured by the reward $M$ and the dual price $\gamma_{\ell s t p}$. The load component reflects the dual cost differential $\psi_{\ell, s, t, (v,c_2)} - \psi_{\ell, s, t, (u,c_1)}$
\begin{align}\label{eq:RCarc}
    &\widehat{g}_{a} = \underbrace{\sum_{m \in \calU^{uv}_{\ell st}}
\sum_{p\in\calP_{m}}d_{mp}w_{mp}}_{\text{routing component}} + \underbrace{\psi_{\ell, s, t, (v,c_2)} - \psi_{\ell, s, t, (u,c_1)}}_{\text{load component}}\\
&\qquad\qquad\text{with}\quad d_{mp}= D_{ps} \left(\frac{\delta\tau^{\text{late}}_{mp}+\frac{\delta}{2}\tau^{\text{early}}_{mp} + \sigma\tau^{\text{travel}}_{mp}}{\tau^{dir}_p}+ \lambda \tau^{\text{walk}}_{mp} + \mu \tau^{\text{wait}}_{mp}-M  \right) + \gamma_{\ell s t p}.\nonumber
\end{align}

The pricing problem seeks a subpath with minimum reduced cost (Equation~\eqref{eq:PPobj}). Constraints~\eqref{eq:PPse}--\eqref{eq:PPdiff2} define the load at each node, starting from load $c_1$ and ending with load $c_2$. Constraints~\eqref{eq:PPw2f} and~\eqref{eq:PPOneP} ensure that passenger pickups occur only in visited nodes, and at most once. Constraints~\eqref{eq:PPflow} apply flow balance in the time-expanded network.
\begin{align}
    \text{PP}_{\ell st}^{u,v,c_1,c_2} \quad \min  \quad
    &\sum_{m \in \calU^{uv}_{\ell st}} \sum_{p \in \calP_{m}} d_{mp} w_{mp}  + \psi_{\ell, s, t, (v,c_2)} - \psi_{\ell, s, t, (u,c_1)}\label{eq:PPobj}\\
    \st\quad
    &\xi_{(u,T_{\ell t}(u))}=c_1,\ \xi_{(v,T_{\ell t}(v))}=c_2 \label{eq:PPse}\\
    & \xi_q-\xi_m \leq \sum_{p\in\calP_m}D_{mp}w_{mp} + C_\ell (1 - f_{mq}),\quad\forall(m,q)\in\calH^{uv}_{\ell st} \label{eq:PPdiff}\\
    & \xi_q-\xi_m \geq \sum_{p\in\calP_m}D_{mp}w_{mp} - C_\ell (1 - f_{mq}),\quad\forall(m,q)\in\calH^{uv}_{\ell st} \label{eq:PPdiff2}\\
    &w_{mp}\leq \sum_{q:(m,q) \in \calH^{uv}_{\ell st}} f_{mq} \quad \forall m \in \calU^{uv}_{\ell st},\ \forall p \in \calP_{m} \label{eq:PPw2f}\\
    &\sum_{m \in \calU^{uv}_{\ell st} \, : \, p \in \calP_m} w_{mp} \leq 1 \quad \forall p \in \calP \, : \, (\ell, t) \in \calM_p \label{eq:PPOneP} \\
    &\sum_{q:(m,q) \in \calH^{uv}_{\ell st}} f_{mq} - \sum_{q:(q,m) \in \calH^{uv}_{\ell st}} f_{qm} = \begin{cases}
    1 &\text{if } m = (u,T_{\ell t}(u)), \\
    -1 &\text{if } m = (v,T_{\ell t}(v)), \\
    0 &\text{ otherwise.}
    \end{cases}\quad \forall m \in \calU^{uv}_{\ell st} \label{eq:PPflow}\\
    &\bff, \bw \text{ binary},\ \bxi\ \text{non-negative integer} \label{eq:PPdomain}
\end{align}

Whenever the solution of the pricing problem is negative, we add the corresponding subpath-based arc $a\in \calA_{\ell st}$ to the load-expanded network, by defining its cost parameter $g_a$ as:
\begin{align}\label{eq:RCsubpath}
    g_a = \sum_{m \in \calU^{uv}_{\ell st}}
\sum_{p\in\calP_{m}}D_{ps} \left(\frac{\delta\tau^{late}_{mp}+\frac{\delta}{2}\tau^{early}_{mp} + \sigma\tau^{travel}_{mp}}{\tau^{dir}_p}+ \lambda \tau^{\text{walk}}_{mp} + \mu \tau^{\text{wait}}_{mp}-M\right) w_{mp}
\end{align}

Remark~\ref{R:ppRHS} notes that the the pricing problem searches over all subpaths, including those from non-selected reference lines (Equation~\eqref{eq:PPflow}) and non-assigned passengers (Equation~\eqref{eq:PPOneP}). Such subpaths will be primal infeasible in the restricted Benders subproblem. However, the corresponding constraints take the form ``$0\leq0$'' and cannot be assumed to have zero duals. This is essential to certify the validity of Benders decomposition. In general terms, the dual polyhedron $\calP_{sj}$ is independent on the incumbent first-stage variables $\bx$, and so is the pricing problem.

\begin{remark}\label{R:ppRHS}
    The right-hand side of Equation~\eqref{eq:PPflow} (resp, Equation~\eqref{eq:PPOneP}) must be 1 rather than $x_{\ell t}$ (resp., $z_{\ell pst}$) to certify optimality of the RBSP solution and guarantee the algorithm's exactness.
\end{remark}

For each reference trip $(\ell,t) \in \mathcal{L} \times \mathcal{T}_{\ell}$ and scenario $s \in \mathcal{S}$, the pricing problem is defined for each node pair $((u,c_1),(v,c_2))\in\calV_{\ell st}\times\calV_{\ell st}$ in the load-expanded network. In fact, we can reduce the number of pricing problems by exploiting the decomposition of the reduced cost into routing and load component. We first maximize the load component for each differential $\nu \in \calC_\ell$:
\begin{equation*}
    \Delta\psi^{u,v,\nu}_{\ell st}=\max\left\{\psi_{\ell s t n} - \psi_{\ell s t m}:(m,n)\in\calA_{\ell st},k_{m}=u,k_{n}=v,c_{n}-c_{m}=\nu\right\}
\end{equation*}
We then seek a subpath that serves $\nu$ passengers and minimizes the routing component:
\begin{align*}
    Z^{u,v,\nu}_{\ell st}=\min\ 
    \sum_{m \in \calU^{uv}_{\ell st}} \sum_{p \in \calP_{m}} d_{mp} w_{mp};\ \st\   \sum_{m \in \calU^{uv}_{\ell st}} \sum_{p \in \calP_{m}} D_{ps} w_{mp} = \nu;\ \text{Equations~\eqref{eq:PPw2f}--\eqref{eq:PPdomain}}
\end{align*}

Proposition~\ref{pricingproblem} shows that we can solve one pricing problem for each \textit{load differential} and every pair of checkpoints. This result reduces the number of pricing problem by a factor $\calO(\max_{\ell\in\calL}C_\ell)$, while retaining the finite convergence and exactness of the column generation scheme.
\begin{proposition}\label{pricingproblem} 
    $Z^{u,v,\nu}_{\ell st} - \Delta\psi^{u,v,\nu}_{\ell st}$ is the minimum reduced cost across all arc-based variables between checkpoints $u$ and $v$ with load differential $\nu$, for all $(\ell,t) \in \mathcal{L} \times \mathcal{T}_{\ell},\  s \in \mathcal{S}$.
\end{proposition}

\subsubsection*{Label setting.}

The pricing problem is a resource-constrained shortest path problem. We design a label-setting algorithm exploiting the directed and acyclic structure of $(\calU^{uv}_{\ell st}, \calH^{uv}_{\ell st})$.

\textit{State definition.} Let $(m^{\sigma}, \mathbb{P}^{\sigma})$ denote a state, where $m^{\sigma}$ tracks the ``current’' node, and $\mathbb{P}^{\sigma}$ tracks the set of served passengers $p \in \calP$ each with pickup node $\rho_p$. We track the reduced cost $G(m^{\sigma}, \mathbb{P}^{\sigma})$.

\textit{Initial state:} $(m^0=m,\mathbb{P}^0=\emptyset)$, where $m$ is such that $k_{m}=u$ and $t_{m}=T_{\ell t}(u)$; $G(m^0,P^0)=0$.

\textit{State transitions.} For each arc $(m,q)\in\calH^{uv}_{\ell st}$ and each passenger combination $\mathbb{P}_{m} \subseteq \calP_{m}$, the state is updated to $(q,\mathbb{P}^{\sigma} \cup \mathbb{P}_{m})$. For each new passenger $p \in \mathbb{P}_{m} \setminus\{ \mathbb{P}^{\sigma}\}$, the pickup point is set to $\rho_p = m$. For existing passengers $p \in \mathbb{P}_{m} \cap  \mathbb{P}^{\sigma}$, we update the pickup node to be $\rho_p=m$ if $d_{mp} < d_{\rho_p,p}$. This transition is admissible if the vehicle has enough capacity, i.e., if $\sum_{p \in \mathbb{P}^{\sigma}\cup\mathbb{P}_{m}} D_{ps}\le C_\ell$.

\textit{Reward function.}
$G(m^{\sigma},\mathbb{P}^{\sigma})=\sum_{p \in  \mathbb{P}^{\sigma}} d_{\rho_p,p}$ tracks the reduced cost of a subpath up to state $\sigma$.

\textit{Dominance rule.} $\sigma^1$ dominates $\sigma^2$ if $m^{\sigma^1} = m^{\sigma^2}$, $\mathbb{P}^{\sigma^1}=\mathbb{P}^{\sigma^2}$, and $G(m^{\sigma^1},\mathbb{P}^{\sigma^1})\le G(m^{\sigma^2}, \mathbb{P}^{\sigma^2})$.

Upon termination, we extract all non-dominated states such that $m^{\sigma}=m: k_{m}= v$ and $t_{m}=T_{\ell t}(v)$. We then add to the RBSP all arcs $a \in \calA_{\ell st} \setminus \{\calA'_{\ell st}\}$ such that $k_{start(a)}=u, \ k_{end(a)}=v,\ c_{end(a)} - c_{start(a)} = \sum_{p \in \mathbb{P}^{\sigma}} D_{ps}$, with negative reduced cost $\widehat{g}_a = G(m^{\sigma},\mathbb{P}^{\sigma}) - \psi_{\ell, s, t, start(a)} + \psi_{\ell, s, t, end(a)}<0$.

By design, the dominance rule yields the subpath of minimum reduced cost for each passenger combination---hence, for each load differential. Thus, we apply the label-setting algorithm for each pair of checkpoints $u,v\in\calI_\ell$, but do not duplicate it for each load differential. The number of checkpoint pairs grows linearly with $|\calI_\ell|$ because subpaths can skip up to $K\in\{0,1\}$ checkpoint. Combined with Proposition~\ref{pricingproblem}, we obtain the following reduction on the number of pricing problems:
\begin{proposition}\label{prop:labset}
    The label-setting algorithm generates $\calO(2^{\Xi}|\calV_{\ell st}|)$ variables at a time by only solving $\calO(I_\ell)$ pricing problems, for each reference trip $(\ell,t) \in \mathcal{L} \times \mathcal{T}_{\ell}$ and scenario $s \in \mathcal{S}$.
\end{proposition}

\subsubsection*{Heuristic label-setting algorithm.}

The algorithm can lead to a weak dominance rule with almost-identical subpaths serving slightly different passenger combinations. In fact, subpaths are relatively short, so the pricing problem rarely rejects a passenger with a negative reduced cost ($d_{mp}<0$) to free up capacity for a subsequent passenger. Moreover, it can be undesirable in practice to reject a passenger at a station visited by the vehicle. We therefore propose a heuristic acceleration such that, in each node $m \in \calU_{\ell st}^{uv}$, all candidate passengers $p \in \calP_{m}$ with negative reduced cost contribution ($d_{mp}<0$) are served, as long as the vehicle does not operate at capacity. This heuristic yields an upper-bounding approximation of the pricing problem, i.e., it generates solutions with a negative reduced cost but can potentially miss other subpaths with negative reduced cost. Then, we can switch back to the full label-setting algorithm to derive a certificate of optimality. In our experiments, the heuristic results in significant speedups with high-quality solutions.
\section{Computational Assessment of the Methodology}\label{sec:comp}

We develop a real-world experimental setup in Manhattan. We use demand data from the \cite{taxi} during the morning rush (6--9 am). We define a road network and travel times using data from Google Maps, OpenStreetMap, and \cite{speeds}. Parameter values are reported in~\ref{app:parameters}. We design candidate reference lines using breadth-first search (\ref{app:refline}).

We consider a MiND-VRP setting corresponding to a shuttle service from Manhattan to LaGuardia Airport (LGA) with vehicles of capacity 10 to 20 passengers. We vary the number of candidate reference lines (5 to 100), the planning horizon (1 to 3 hours), the number of checkpoints that can be skipped ($K=0,1,2,3$), and the number of scenarios (5 to 20). We use a 15-minute discretization to schedule transit vehicles in the first stage (sets $\calT_\ell$), and a 30-second discretization in the second stage (sets $\calT^{uv}_{\ell t}$). Our problem includes up to 1,900 passenger requests, 640 candidate stops, and 100 candidate reference lines (Figure~\ref{F:ManhattanLS}), resulting in over 1 million first-stage variables, 25,000 Benders subproblems, and 200,000 pricing problems. We also develop real-world experimental setups for the MiND-DAR in~\ref{app:DARcase} and for the MiND-Tr in~\ref{app:TRcase}.

All models are solved with Gurobi v12.0 using the JuMP package in Julia \citep{dunning2017jump}. We impose a three-hour time limit for optimization. All instances and code are available online.\footnote{\hyperlink{https://github.com/martiniradi/DeviatedFixedRouteMicrotransit}{https://github.com/martiniradi/DeviatedFixedRouteMicrotransit}}

\subsection{Benefits of Subpath Modeling and Double-decomposition Algorithm}
\label{subsec:method}

Table~\ref{tab:ResPvsSPvsA} compares the four formulations defined in Section~\ref{sec:model} in terms of solution quality (normalized to the best-found solution), computational times, and number of second-stage variables. All models are solved with off-the-shelf methods, using exhaustive enumeration of segments, subpaths or paths in the network-based reformulations (with up to 1 million paths per subproblem). The compact formulation, despite its much much smaller size, does not scale to even the smallest instances. This underscores the highly challenging structure of the two-stage stochastic optimization structure with discrete recourse. Among network-based approaches, the segment-based formulation features the most limited scalability, requiring 30 million variables in the smallest instance due to the granular time-station-load discretization. The path-based formulation scales to medium instances but its performance quickly deteriorates due to the exponential growth in the number of variables. In comparison, the subpath-based formulation requires orders of magnitude fewer variables, terminates much faster, and returns a superior solution with 10 candidate lines.

\begin{table}[h!]
\renewcommand*{\arraystretch}{1.0}
	\caption{Comparison of path-based, subpath-based, segment-based, and compact MiND-VRP formulations.}
	\label{tab:ResPvsSPvsA}
	\begin{center}
		\centering\footnotesize
		\resizebox{1.\textwidth}{!}{%
		\begin{tabular}{ccccccccccccccc}
			\toprule
			\textbf{} & \textbf{} & \textbf{} & \multicolumn{3}{c}{Path-based} & \multicolumn{3}{c}{Subpath-based} & \multicolumn{3}{c}{Segment-based} & \multicolumn{3}{c}{Compact}  \\ \cmidrule(lr){4-6}\cmidrule(lr){7-9}\cmidrule(lr){10-12}\cmidrule(lr){13-15}
			$|\calL|$&  Hor. & $K$ & Sol. & CPU (s) & Arcs & Sol. & CPU (s) & Arcs & Sol. & CPU (s) & Arcs  & Sol. & CPU (s) & S2 Vars \\ 
			\bottomrule
				5 & 60 & 0 & 100 & 137s & 3.1M & 100 & 17s & 39K & 100 & 6,431s & 30.0M & | & 10,800s+ & 58K \\
			5 & 60 & 1 & | & | & | & | & | & | & | & 10,800s+ & 30.0M & | & 10,800s+ & 68K \\
			5 & 120 & 0 & 100 & 652s & 8.6M & 100 & 206s & 86K & | & | & |  & | & 10,800s+ & 121K \\
			5 & 180 & 0 & 100 & 696s & 9.6M & 100 & 238s & 110K & | & | & |  & | & 10,800s+ & 182K \\
			10 & 60 & 0 & 100.3 & 1,490s & 29.1M & 100 & 46s & 142K & | & | & |  & | & 10,800s+ & 123K \\
			\bottomrule[1pt]
		\end{tabular}
        }
	\end{center}
	\begin{tablenotes}\footnotesize
		\item ``10,800s+'': optimization timeout; ``|’’: the algorithm does not terminate due to memory limitations.
	\end{tablenotes}
\end{table}

These results demonstrate the benefits of the subpath-based representation of the second-stage problem. We provide additional results in~\ref{app:subpath} by comparing the four formulations on the capacitated vehicle routing problem with time windows in the second stage alone. These results uncover the different bottlenecks of the algorithms: the branch-and-cut structure in the compact formulation due to its weak linear relaxation, the size of the segment-based formulation, and the enumeration of arc variables in the path-based and subpath-based formulations. Whereas the subpath-based model is the most scalable one, subpath enumeration remains intractable in medium-scale instances, thus motivating our double-decomposition algorithm.

Next, Table~\ref{tab:algorithm_PMIOLO} compares Benders decomposition with subpath enumeration and our DD methodology with exact and heuristic label setting. Benders decomposition alone remains limited with subpath enumeration. In comparison, our double-decomposition algorithm achieves much stronger scalability by leveraging column generation in the Benders subproblem.

\begin{table}[h!]
\renewcommand*{\arraystretch}{1.0}
\centering
\small
\caption{Assessment of exact algorithms for solving $(\texttt{MIO}-\texttt{LO})$.}
\label{tab:algorithm_PMIOLO}
\resizebox{1.\textwidth}{!}{%
\begin{tabular}{cccccccccccccccc}
\toprule
&  &  &  \multicolumn{8}{c}{$K=0$} & \multicolumn{5}{c}{$K=1$} \\
\cmidrule(lr){4-11}\cmidrule(lr){12-16}
 &  &  & \multicolumn{3}{c}{Benders} & \multicolumn{3}{c}{DD (exact)} & \multicolumn{2}{c}{DD (acceleration)} & \multicolumn{3}{c}{DD (exact)} & \multicolumn{2}{c}{DD (acceleration)} \\ \cmidrule(lr){4-6}\cmidrule(lr){7-9}\cmidrule(lr){10-11}\cmidrule(lr){12-14}\cmidrule(lr){15-16}
$|\calL|$ & $|\calS|$ & Horizon & Sol. & Gap & CPU(s)  & Sol. & Gap & CPU(s) & Sol. & CPU(s)  & Sol. & Gap & CPU(s) & Sol. & CPU(s)  \\
 \cmidrule{1-16}
5 & 5 & 60 & 100 & \textbf{0.0\%} & 39 & 100 & \textbf{0.0\%} & 11 & 102.2 & 10 & 100 & \textbf{0.0\%} & 165 & 102.1 & 47 \\ 
 &  & 120 & 100 & \textbf{0.0\%} & 513 & 100 & \textbf{0.0\%} & 38 & 102 & 30 & 100 & \textbf{0.0\%} & 4,834 & 101.2 & 169 \\ 
 &  & 180 & 100 & \textbf{0.0\%} & 599 & 100 & \textbf{0.0\%} & 49 & 101.1 & 53 & 100 & \textbf{0.0\%} & 4,212 & 101.6 & 217 \\ 
  \cmidrule{2-16}
 & 20 & 60 & --- & --- & --- & 100 & \textbf{0.0\%} & 71 & 101.5 & 47 & 100 & \textbf{0.0\%} & 10,489 & 102.4 & 382 \\ 
 &  & 120 & --- & --- & --- & 100 & \textbf{0.0\%} & 349 & 101.3 & 161 & 100 & 4.9\% & 10,800 & 101.1 & 1,224 \\ 
 &  & 180 & --- & --- & --- & 100 & \textbf{0.0\%} & 535 & 101.1 & 219 & 100.1 & 3.6\% & 10,800 & 100 & 3,368 \\ 
 \cmidrule{1-16}
 10 & 5 & 60 & 100 & \textbf{0.0\%} & 93 & 100 & \textbf{0.0\%} & 77 & 101.4 & 61 & 100 & \textbf{0.0\%} & 10,380 & 102.1 & 617 \\ 
 &  & 120 & 100 & \textbf{0.0\%} & 836 & 100 & \textbf{0.0\%} & 239 & 101.3 & 158 & 107.3 & 21.3\% & 10,800 & 100 & 2,732 \\ 
 &  & 180 & 100 & \textbf{0.0\%} & 1,047 & 100.1 & \textbf{0.1\%} & 427 & 100.9 & 210 & 108 & 25.3\% & 10,800 & 100 & 10,800 \\ 
   \cmidrule{2-16}
 & 20 & 60 & --- & --- & --- & 100 & \textbf{0.1\%} & 803 & 101.1 & 558 & --- & --- & --- & --- & --- \\ 
 &  & 120 & --- & --- & --- & 100 & \textbf{0.0\%} & 2,354 & 101.2 & 1,622 & --- & --- & --- & --- & --- \\ 
 &  & 180 & --- & --- & --- & 100 & \textbf{0.0\%} & 5,473 & 100.8 & 4,412 & --- & --- & --- & --- & --- \\ 
 \cmidrule{1-16}
 50 & 5 & 60 & --- & --- & --- & 100 & \textbf{0.0\%} & 2,265 & 100.2 & 613 & --- & --- & --- & --- & --- \\ 
 &  & 120 & --- & --- & --- & 100 & 1.3\% & 10,800 & 100.7 & 10,800 & --- & --- & --- & --- & --- \\ 
 &  & 180 & --- & --- & --- & 100 & 1.8\% & 10,800 & 100.7 & 10,800 & --- & --- & --- & --- & --- \\ 
   \cmidrule{2-16}
 & 20 & 60 & --- & --- & --- & 100 & 0.9\% & 10,800 & 100.6 & 10,800 & --- & --- & --- & --- & --- \\ 
 &  & 120 & --- & --- & --- & 100 & 12.5\% & 10,800 & 109.2 & 10,800 & --- & --- & --- & --- & --- \\ 
 &  & 180 & --- & --- & --- & 100 & 23\% & 10,800 & 112.2 & 10,800 & --- & --- & --- & --- & --- \\ 
 \cmidrule{1-16}
 100 & 5 & 60 & --- & --- & --- & 100 & \textbf{0.0\%} & 5,665 & 100.2 & 1,887 & --- & --- & --- & --- & --- \\ 
 &  & 120 & --- & --- & --- & 100 & 0.8\% & 10,800 & 100.4 & 10,800 & --- & --- & --- & --- & --- \\ 
 &  & 180 & --- & --- & --- & 100.8 & 3.0\% & 10,800 & 100 & 10,800 & --- & --- & --- & --- & --- \\ 
   \cmidrule{2-16}
 & 20 & 60 & --- & --- & --- & 102.3 & 3.9\% & 10,800 & 100 & 10,800 & --- & --- & --- & --- & --- \\ 
 &  & 120 & --- & --- & --- & 100 & 28.7\% & 10,800 & 108.6 & 10,800 & --- & --- & --- & --- & --- \\ 
\bottomrule
\end{tabular}
}
\begin{tablenotes}\footnotesize
    \item ``|’’ indicates that the algorithm does not terminate due to memory limitations.\vspace{-3pt}
    \item Optimality gaps measure the difference between a feasible solution to~\eqref{OPT} and the partial relaxation $(\calP-\texttt{MIO}-\texttt{LO})$.\vspace{-3pt}
    \item Bolded values indicate instances that terminate within the time limit.
\end{tablenotes}
\end{table}

Specifically, when all checkpoints must be visited ($K=0$), the DD algorithm can solve instances with 50 candidate lines and a three-hour horizon, or instances with 100 candidate lines and a two-hour horizon. When vehicles can skip checkpoints ($K=1$), the longer subpaths result in exponentially larger second-stage problems; subpath enumeration fails to find feasible solutions, whereas DD can solve instances with up to 10 candidate reference lines and a three-hour horizon (we report more results on the role of $K$ in~\ref{app:sensitivity}). These improvements are driven by the small number of variables needed to guarantee convergence in column generation---namely, DD converges with up to 93\% fewer variables. Moreover, the DD algorithm yields a 0.0--0.1\% optimality gap whenever it terminates, confirming the tightness of our subpath-based formulation. Ultimately, the DD methodology provides certifiably optimal, or near-optimal solutions in large and otherwise-intractable instances of the problem where all benchmarks fail to even return a feasible solution.

Then, our label-setting acceleration further enhances the scalability of the algorithm. These benefits are stronger with $K=1$ because the stronger dominance criterion becomes more impactful with longer subpaths. In small instances, DD terminates up to 3 times faster with the pricing heuristic, while returning solutions within 3\% of the optimum. In medium instances, the pricing heuristic actually enables higher-quality solutions in faster computational times. This is because the DD algorithm with exact label-setting algorithm does not terminate within the time limit, and the heuristic label-setting scheme enables more effective convergence due to a much smaller number of subpaths (up to 52\% and 78\% fewer subpaths when $K=0$ and $K=1$, respectively). In other words, the benefits of acceleration can outweigh the slight loss of flexibility when choosing which passengers to pick up at each station. Ultimately, by combining Benders decomposition, column generation and label-setting acceleration, our algorithm can solve realistic instances with up to 100 candidate reference lines, hundreds of stations, 5 demand scenarios and a three-hour horizon (or 20 scenarios and a one-hour horizon). These correspond to large-scale network design and routing instances, with hundreds of candidate stops and thousands of passenger requests.

Similarly, despite the higher complexity of the problem, the methodology can handle realistic MiND-DAR and MiND-Tr instances with up to 10--40 candidate lines (\ref{app:DARcase} and \ref{app:TRcase}).

Finally, Table~\ref{tab:algorithms_P} reports results of DD\&ILS and UB\&DD. Recall that Benders decomposition and DD solve the partial relaxation $(\calP-\texttt{MIO}-\texttt{LO})$, while providing a valid (and tight) optimality gap for Problem~\eqref{OPT}, whereas the other two methods certifiably solve Problem~\eqref{OPT} (Theorems~\ref{thm:DDILS}--\ref{thm:UBDD}). These results identify instances where DD leaves a small optimality gap that DD\&ILS can close (e.g., 8 lines, 5 scenarios, two-hour horizon). Yet, by solving the second-stage mixed-integer problem repeatedly, both the DD\&ILS and UB\&DD algorithms are much more computationally intensive. In medium instances, DD\&ILS return the optimal solution in longer computational times than DD, whereas UB\&DD fails to return a feasible solution. In larger instances, both DD\&ILS and UB\&DD time out whereas DD still returns a certifiably optimal solution.

\begin{table}[h!]
\renewcommand*{\arraystretch}{1.0}
\centering
\footnotesize
\caption{Assessment of exact algorithms to solve Problem~\eqref{OPT}.}
\label{tab:algorithms_P}
\begin{tabular}{ccccccccccccccccc}
\toprule
 &  &  & \multicolumn{3}{c}{Benders} & \multicolumn{3}{c}{DD} & \multicolumn{3}{c}{DD\&ILS} & \multicolumn{3}{c}{UB\&DD} \\ 
 \cmidrule(lr){4-6}\cmidrule(lr){7-9}\cmidrule(lr){10-12} \cmidrule(lr){13-15} 
$|\calL|$ & $|\calS|$ & Horizon & Sol. & Gap & CPU(s) & Sol. & Gap & CPU(s) & Sol. & Gap & CPU(s) & Sol. & Gap & CPU(s) \\
\cmidrule(lr){1-15}
5 & 5 & 60 & 100.0 & \textbf{0.0\%} & 39 & 100.0 & \textbf{0.0\%} & 11 & 100.0 & \textbf{0.0\%} & 43 & 100.0 & \textbf{0.0\%} & 41 \\
 &  & 120 & 100.0 & \textbf{0.0\%} & 513 & 100.0 & \textbf{0.0\%} & 38 & 100.0 & \textbf{0.0\%} & 348 & 100.0 & \textbf{0.0\%} & 63 &  \\
\cmidrule(lr){1-15}
8 & 5 & 60 & 100.0 & \textbf{0.0\%} & 93 & 100.0 & \textbf{0.0\%} & 100 & 100.0 & \textbf{0.0\%} & 340 & 100.0 & \textbf{0.0\%} & 1,464 &  \\
 &  & 120 & 100.2 & \textbf{0.2\%} & 908 & 100.2 & \textbf{0.2\%} & 214 & 100.0 & \textbf{0.0\%} & 916 & --- & --- & 10,800 &  \\
\cmidrule(lr){1-15}
10 & 20 & 60 & --- & --- & --- & 100.0 & \textbf{0.0\%} & 803 & --- & --- & 10,800 & --- & --- & 10,800 &  \\
 &  & 120 & --- & --- & --- & 100.0 & \textbf{0.0\%} & 2,354 & --- & --- & 10,800 & --- & --- & 10,800 &  \\
\bottomrule
\end{tabular}
\end{table}

These results underscore that the structural complexity of Problem~\eqref{OPT} stems from the exponential size of the second-stage problem rather than its discreteness. Integer L-shaped cuts and the UB\&BC scheme can be critical to solve stochastic integer problems with a polynomial second-stage formulation and a weaker relaxation; in contrast, the network-based representation in Problem~\eqref{OPT} yields a very tight second-stage reformulation, thus shifting the complexity of the problem from a discrete recourse function to an exponential number of second-stage variables---and motivating our DD algorithm. As our results show, the integer L-shaped cuts and UB\&BC provide limited marginal computational benefits in this setting whereas the DD methodology is instrumental in enabling convergence and deriving high-quality solutions to Problem~\eqref{OPT}.

\subsection{Benefits of Stochastic Optimization Methodology}\label{subsec:StochVal}

Table~\ref{tab:VSS-obj} compares the stochastic optimization solution against a deterministic baseline and a clairvoyant benchmark. The MiND-VRP reduces unmet demand by 6-7\% on average, while reducing passengers' walking time by 25-35\% from the deterministic baseline, resulting in a high VSS---5-7\% on average and up to 10\%. In fact, our solution bridges 40-50\% of the gap on average between the deterministic and perfect-information benchmarks. These results highlight the benefits of our two-stage stochastic optimization model (and our DD methodology to solve it) to increase demand coverage while maintaining or even improving level of service, as compared to a deterministic model (which can be solved via off-the-shelf methods).

\begin{table}[h!]
\centering
\renewcommand{\arraystretch}{1.0}
\small
\caption{Value of Stochastic Solution (VSS) and Expected Value of Perfect Information (EVPI)}
\label{tab:VSS-obj}
\centering
\resizebox{\textwidth}{!}{%
 \begin{tabular}{cccccccccccccc}
 \toprule[1pt]
 &  &  &  &  & \multicolumn{3}{c}{Performance assessment} & \multicolumn{6}{c}{VSS breakdown} \\ \cmidrule(lr){6-8}\cmidrule(lr){9-14}
 $\calL$ & $|\calS|$&Horizon   & $K$ & Heur. & $|\frac{\text{VSS}}{Sol.}|$ & $|\frac{\text{EVPI}}{Sol.}|$ & $\frac{\text{VSS}}{(\text{VSS}+\text{EVPI})}$ & \begin{tabular}[c]{@{}c@{}}Unmet\\demand (\%)\end{tabular} & \begin{tabular}[c]{@{}c@{}}Walking\\time (\%)\end{tabular} & \begin{tabular}[c]{@{}c@{}}Waiting\\ time (\%)\end{tabular} & \begin{tabular}[c]{@{}c@{}}Earliness\\(\%)\end{tabular} & \begin{tabular}[c]{@{}c@{}}Delay\\(\%)\end{tabular} & \begin{tabular}[c]{@{}c@{}}Detour\\ (\%)\end{tabular}\\  \hline
10 & 5 & 60 & 0 & \ding{55} & 5.8 & 8.2 & 41.2 & -6 & -59.3 & 0.7 & -2.1 & 13.2 & 1 \\
 &  &  & 1 &\ding{55}  & 2.7 & 9.7 & 22 & -2.7 & 36.3 & -1.6 & -9.9 & 3.9 & 0.4 \\
 &  & 120 & 0 &  \ding{55}& 4 & 7.4 & 35.2 & -4.3 & -77 & 1.2 & 1.7 & 13 & 0.5 \\
 &  &  & 1 & \checkmark & 3 & 7.4 & 28.6 & -2.8 & -1.6 & -6.8 & 3.3 & -4.8 & -0.8 \\
 &  & 180 & 0 & \ding{55} & 9.9 & 7.4 & 57.4 & -11 & -70.6 & 7.4 & 3.4 & -1.5 & -0.4 \\
 &  &  & 1 & \checkmark & 7 & 5.7 & 55.2 & -7.5 & -17.4 & 0.4 & -5.8 & 4.1 & 0.5 \\ \cline{2-14} 
 & 20 & 60 & 0 & \ding{55} & 7.6 & 5.9 & 56.3 & -8.1 & -91.6 & 1.1 & -5.6 & 15.7 & 1 \\
 &  &  & 1 & \checkmark & 9.8 & 6.1 & 61.5 & -10.4 & -8.1 & -10.2 & -11.4 & 23.9 & -1.2 \\
 &  & 120 & 0 & \checkmark & 9.7 & 9.2 & 51.1 & -10.8 & -81.2 & 4.3 & 10.5 & 4.2 & 0.2 \\
 &  &  & 1 & \checkmark & 4.8 & 8 & 37.8 & -5 & 49.7 & 7.4 & -6.9 & 10.3 & 0 \\
 &  & 180 & 0 & \checkmark & 6.1 & 7.2 & 45.7 & -6.6 & -82.9 & 6.4 & 6.9 & 14.5 & 3.2 \\
 &  &  & 1 & \checkmark & 2.7 & 7.3 & 27.2 & -2.8 & -2.8 & 5.1 & -0.6 & 14 & 1.4 \\ \hline
50 & 5 & 60 & 0 & \ding{55} & 5.8 & 2.5 & 70.3 & -6.5 & -6.1 & 0 & -13 & 6.2 & 0.3 \\
 &  & 120 & 0 & \checkmark & 5.3 & 7 & 43.1 & -5.7 & -2 & 1.3 & 2.7 & 1.8 & 1 \\
 &  & 180 & 0 & \ding{55} & 4.3 & 10.7 & 30.5 & -4.5 & -16.7 & -2.1 & 8.3 & 4.2 & 0.8 \\ \hline
100 & 5 & 60 & 0 & \checkmark & 8.8 & 2.3 & 79.4 & -11.1 & -56.2 & -19.5 & 1.3 & -1.5 & -0.7 \\
 &  & 120 & 0 & \checkmark & 5 & 7.4 & 40.2 & -5.6 & -25.2 & -7.6 & 0.4 & 0.1 & 0.6 \\
 &  & 180 & 0 & \checkmark & 2.4 & 15.5 & 13.4 & -2.7 & -3.7 & -5.6 & 0.3 & -3.2 & 1 \\ \hline
\multicolumn{5}{l}{\textbf{Average 5 scenarios}} & \textbf{5.3} & \textbf{7.6} & \textbf{43.0} & \textbf{-5.9} & \textbf{-25.0} & \textbf{-2.7} & \textbf{-0.8} & \textbf{3.0} & \textbf{0.4} \\
\multicolumn{5}{l}{\textbf{Average 20 scenarios}} & \textbf{6.8} & \textbf{7.3} & \textbf{46.6} & \textbf{-7.3} & \textbf{-36.2} & \textbf{2.4} & \textbf{-1.2} & \textbf{13.8} & \textbf{0.8}\\\bottomrule[1pt]
 \end{tabular}%
 }
 \begin{tablenotes}
     \footnotesize 
    \vspace{-6pt}
     \item ``Heur.'': solution with heuristic (\checkmark) vs exact (\ding{55}) pricing algorithm; ``Sol.'': stochastic optimization solution.
        \vspace{-6pt}
     \item  Unmet demand is measured in number of passengers; all other components are measured per served passenger.
 \end{tablenotes}
 \vspace{-6pt}
\end{table}
\section{Practical Assessment of Deviated Fixed-route Microtransit}\label{sec:practical}

Finally, we conduct a comprehensive assessment of microtransit against fixed-route transit (a single-stage problem without second-stage deviations) and ride-sharing (on-demand system with vehicle capacities of 1, 2, and 4, described in \ref{A:rideshare}). We use the same experimental setup as in Section~\ref{sec:comp}. Recall that, since Manhattan represents a high-density region, the results can be seen as conservative estimates of the impact of microtransit in lower-density areas with fewer transit options. Again, all our insights hold in the MiND-DAR and MiND-Tr, as shown in~\ref{app:DARcase} and~\ref{app:TRcase}.

\subsection{Value of Microtransit Flexibility}
\label{subsec:flex}

\subsubsection*{Second-stage microtransit operations.} 
\begin{figure}[htbp!]
       \centering
       \subfloat[\footnotesize Fixed-route transit]{\label{F:Case2_T}\includegraphics[width=0.49\textwidth]{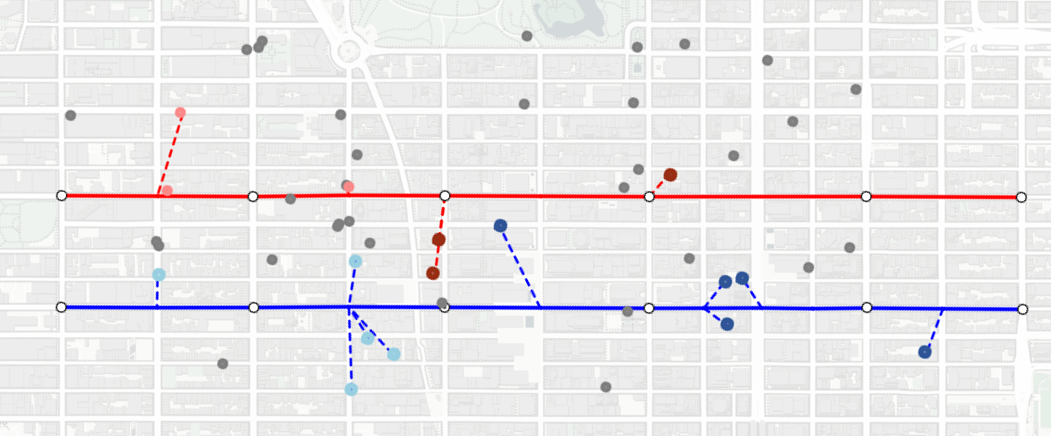}}
       \hspace{0.1cm}
       \subfloat[\footnotesize Deviated fixed-route microtransit]{\label{F:Case2_MT}\includegraphics[width=0.49\textwidth]{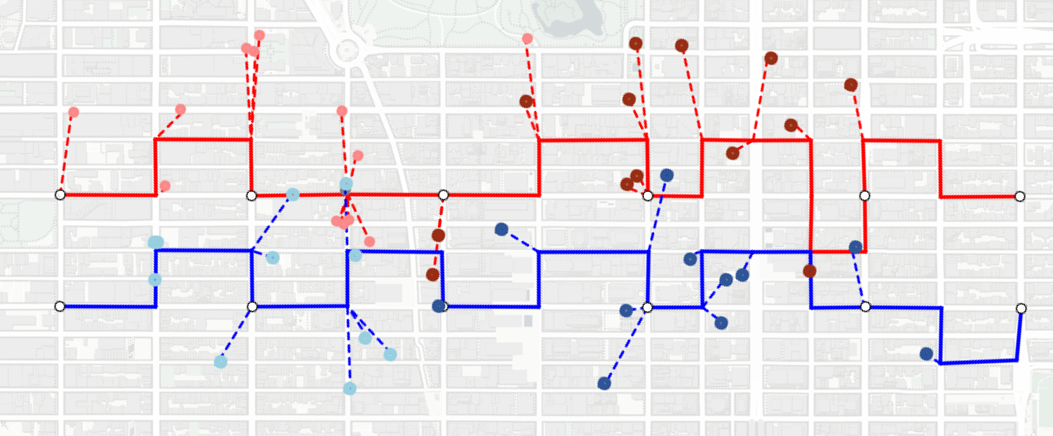}}
        \caption{Illustration of transit and microtransit operations in the MiND-DAR for two reference lines [light (resp. dark) blue/red circles: origins (resp. destinations) of passengers served by the blue/red line; grey circles: origins and destinations of unserved passengers; white circles: checkpoints].}
    \label{F:DARops}
\end{figure}

Figure~\ref{F:DARops} illustrates MiND-DAR operations along two lines in Midtown Manhattan. By design, transit follows the reference line whereas microtransit deviates from the reference line in all but one checkpoint pair. Thus, the microtransit system serves more passengers (24 versus 8), at the cost of a longer distance (8.5 vs. 5.5 km). Still, the higher vehicle loads leads to a smaller distance per passenger (356 vs. 699 meters), which can translate into lower costs for the operator, lower fares for passengers, and a smaller environmental footprint.

Table~\ref{T:mtLOS} reports average performance and level of service in the second stage of the MiND-VRP, with different vehicle capacities and extents of flexibility ($\Delta=600$ vs. $\Delta=1,200$ meters; $K=0$ vs. $K=1$). In these experiments, we use the same reference trips for transit and microtransit. On average, microtransit can add 1--4 passengers per vehicle, while reducing walking times by 50\% and wait times by 2 minutes. These benefits come with a small increase in detours (+2\%) and an increase in distance traveled (+15--25\%). Still, due to the large increase in utilization, distance per passenger is reduced by up to 500 meters, or 23\%. These benefits become stronger with larger vehicles. Interestingly, even when microtransit vehicles are constrained to stay close to the reference lines (low deviation) and to visit all checkpoints ($K=0$), the microtransit system can significantly improve coverage (0.5 to 3 extra passengers per vehicle, on average) and level of service (reduction in walking times by 40 seconds and in waiting times by 1 minute). In other words, even limited extents of demand-responsiveness can achieve significant performance improvements through stronger demand consolidation, higher level of service, and a smaller cost per passenger.

\begin{table}[htbp!]
    \centering
       \footnotesize 
    \caption{Average operating performance and level of service for fixed-route transit and microtransit.}
    \label{T:mtLOS}
        \begin{tabular}{clllcccccccc}
         \toprule[1pt]
&\multicolumn{3}{c}{Operating model}&\multicolumn{4}{c}{Average operating performance}&\multicolumn{4}{c}{Average level of service}\\\cmidrule(lr){2-4}\cmidrule(lr){5-8}\cmidrule(lr){9-12}
  Cap. & Mode & Dev. & Skip? & \#pass./vehicle & Util. & Dist. & Dist./pass. & Walk & Wait & Detour & Delay \\\hline
10	& Transit & --- & ---	        & 8.05	& 80.50\%	& 14.78	& 2.66	& 2.19	& 6.77	& 152.06\%	& -0.87 \\
    & Microtransit & Low & $K = 0$	& 8.61	& 86.12\%	& 16.48	& 2.53	& 1.48	& 5.69	& 154.36\%	& -0.47 \\
    & Microtransit & High & $K = 0$	& 8.63	& 86.32\%	& 16.57	& 2.51	& 1.44	& 5.58	& 153.62\%	& -0.44 \\
    & Microtransit & Low & $K = 1$	& 8.78	& 87.81\%	& 17.08	& 2.54	& 1.17	& 4.51	& 156.50\%	& -0.28 \\
    & Microtransit & High & $K = 1$	& 8.99	& 89.87\%	& 17.39	& 2.38	& 1.03	& 4.33	& 153.72\%	& -0.28 \\\hline
15	& Transit & --- & ---	        & 10.72	& 71.47\%	& 15.06	& 2.40	& 2.27	& 6.88	& 150.73\%	& -1.29 \\
    & Microtransit & Low & $K = 0$	& 12.20	& 81.32\%	& 17.17	& 2.17	& 1.57	& 5.90	& 151.82\%	& -0.49 \\
    & Microtransit & High & $K = 0$	& 12.29	& 81.96\%	& 17.34	& 2.14	& 1.50	& 5.74	& 151.20\%	& -0.46 \\
    & Microtransit & Low & $K = 1$	& 12.56	& 83.70\%	& 17.78	& 2.15	& 1.31	& 4.83	& 154.83\%	& -0.26 \\
    & Microtransit & High & $K = 1$	& 12.89	& 85.95\%	& 18.15	& 1.97	& 1.15	& 4.63	& 151.74\%	& -0.31 \\\hline
20	& Transit & --- & ---	        & 12.24	& 61.21\%	& 15.16	& 2.34	& 2.30	& 6.94	& 150.38\%	& -1.84 \\
    & Microtransit & Low & $K = 0$	& 15.28	& 76.42\%	& 17.52	& 2.02	& 1.69	& 6.21	& 150.77\%	& -0.52 \\
    & Microtransit & High & $K = 0$	& 15.46	& 77.32\%	& 17.72	& 1.98	& 1.62	& 6.04	& 150.08\%	& -0.50 \\
    & Microtransit & Low & $K = 1$	& 15.90	& 79.49\%	& 18.13	& 1.96	& 1.48	& 5.25	& 153.51\%	& -0.28 \\
    & Microtransit & High & $K = 1$	& 16.16	& 80.78\%	& 18.57	& 1.81	& 1.43	& 5.39	& 151.17\%	& -0.33 \\\bottomrule[1pt]
    \end{tabular}
    \begin{tablenotes}
        \footnotesize
        \vspace{-6pt}
        \item ``Cap.'' -- Capacity; ``Pass.'' -- Passenger; ``Util.'' -- Utilization; ``Dist.'' -- Distance; ``Dev.'' -- deviation.
        \vspace{-6pt}
        \item Units: distance, distance per passenger -- kilometers; walk, wait, delay/earliness -- minutes.
        \vspace{-6pt}
        \item Parameters: two-hour horizon; 10 weekday scenarios, maximum walk: 7 minutes, maximum wait: 10 minutes.
    \end{tablenotes}
\end{table}

\begin{figure}[htbp!]
     \centering
     \includegraphics[height=0.35\textwidth]{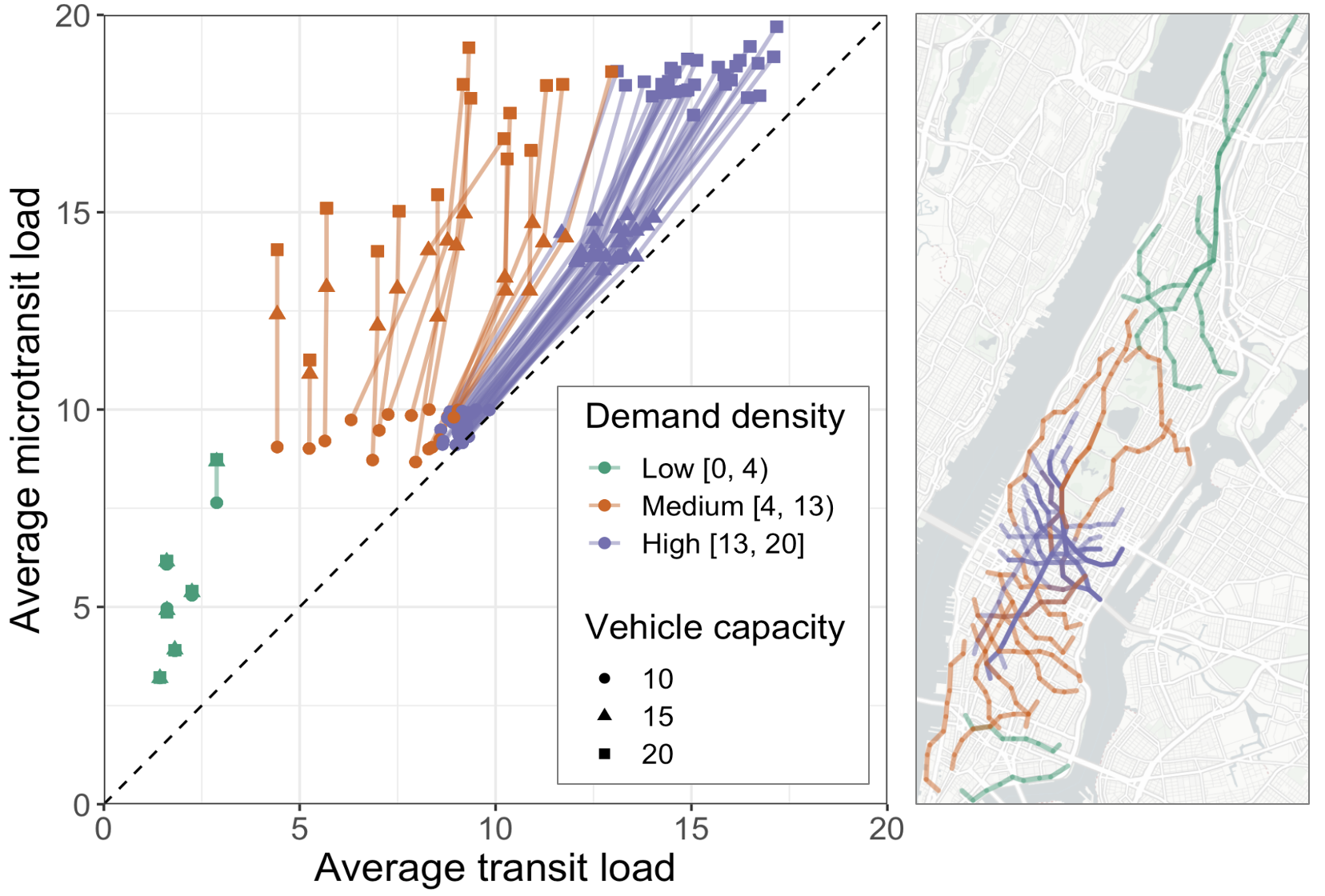}
     \caption{Value of operating flexibility ($\Delta=1,200$ m., $K=0$). Low (resp. medium, high) density: lines with maximum load less than 4 (resp. 4 to 13, more than 13) passengers on average under transit.}
     \label{F:density}
     \vspace{-24pt}
 \end{figure}

Figure~\ref{F:density} plots the average vehicle load in fixed-route transit vs. deviated fixed-route microtransit, for each reference line broken down into low, medium and high density (colored lines) and for each vehicle capacity (dot shapes). All observations lie above the 45-degree line as microtransit makes use of the deviations to increase vehicle load. In low-density regions, microtransit vehicles do not operate at capacity, whereas low-occupancy ride-sharing can provide high levels of service at limited detours and delays. In high-density regions, fixed-route transit already provides high demand coverage due to high synergies across passengers, so the marginal improvements from microtransit are more limited. In-between, microtransit can increase average vehicle loads by 1--5 and the number of pickups by 5--10. These results identify a medium-density regime where deviated fixed-route microtransit can be most impactful, as population density is high enough to consolidate demand into high-occupancy vehicles but too low for fixed-route transit to be as effective.

\subsubsection*{Network design.}

\begin{figure}[htbp!]
       \subfloat[\footnotesize Fixed-route transit]{\label{F:transitNetwork}\includegraphics[height=0.5\textwidth]{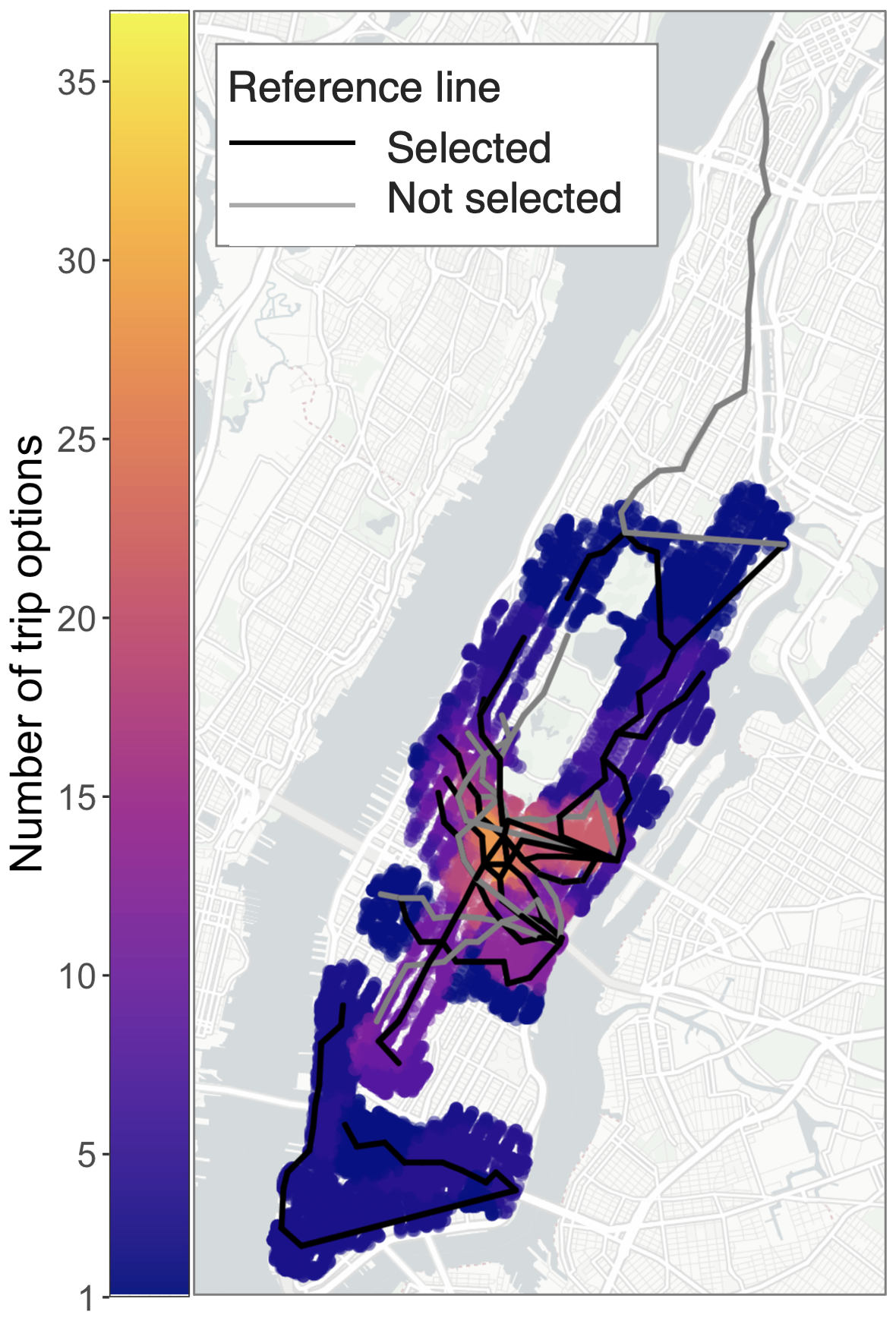}}
   \centering
   ~ ~
   \subfloat[\footnotesize Deviated fixed-route microtransit]{\label{F:microtransitNewtork}\includegraphics[height=0.5\textwidth]{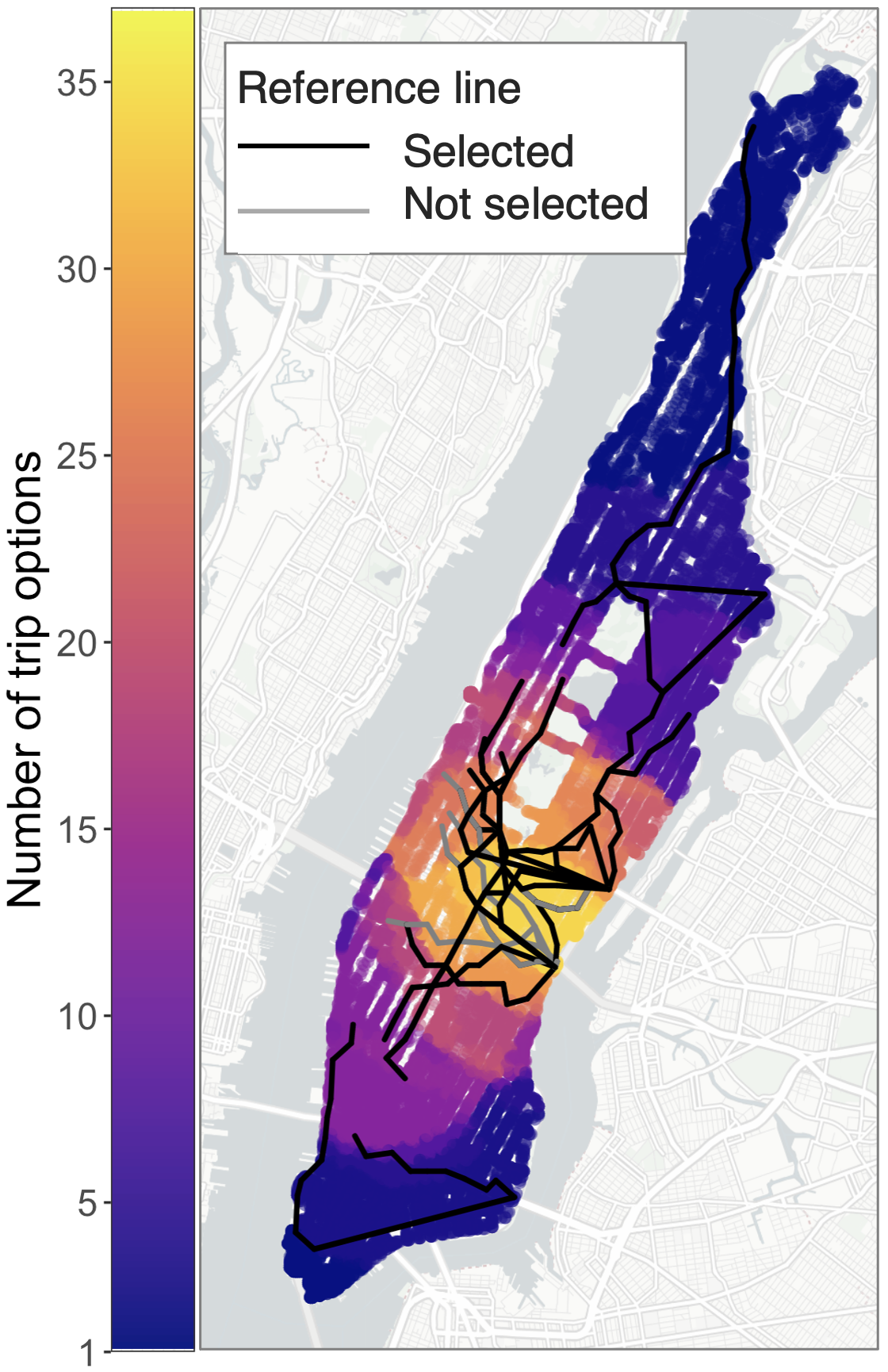}}
    \caption{Reference lines and catchment areas. Parameters: 25 candidate reference lines, 2-hour horizon, 20 vehicles with a 20-passenger capacity each. Microtransit parameters: $\Delta=1,200$ m., $K=0$.}
    \label{F:networkDesign}
    \vspace{-12pt}
\end{figure}

Figure \ref{F:networkDesign} depicts the optimized first-stage networks, with reference lines labeled as ``selected'' if at least one reference trip is selected over the planning horizon. The figure also depicts the number of trip options from each of Manhattan's 21,000 roadway intersections, defined as the number of reference trips with a candidate pickup location within a 5-minute walking radius. Note that the microtransit network expands the catchment area from fixed-route transit. Consistently with Figure~\ref{F:density}, the fixed-route transit system mostly selects lines in high-demand areas. Thanks to its operating flexibility, microtransit provides more trip options: in Midtown Manhattan for instance, the number of trip options increases with microtransit from 20--25 to 30--35; overall, the average number of trip options per intersection increases by a factor of 3 (8.31 vs. 2.61). As a result, the microtransit system reaches low-demand regions, such as Uptown Manhattan. Specifically, microtransit covers 60\% more intersections with at least one trip option (53.8\% vs. 85.4\%). By enhancing service options in high-density regions, microtransit can allocate some budget to expand its geographic reach to under-served regions, thus enhancing accessibility across the population.

\subsection{Performance Assessment}
\label{subsec:performance}

We now compare the performance of microtransit against fixed-route transit and ride-sharing. To establish an apples-to-apples comparison, we fix total seating capacity across all systems (e.g., 10 transit/microtransit vehicles of capacity 10, ride-sharing with 100/50/25 vehicles of capacity of 1/2/4), and perform an out-of-sample assessment corresponding to five new weekdays. Unlike in Table~\ref{T:mtLOS} and Figure~\ref{F:density}, we consider here the optimized network of reference lines in transit and microtransit. Table~\ref{T:los} reports average coverage, level of service, and distance traveled.

\begin{table}[htbp!]
    \footnotesize
    \caption{Average level of service of fixed-route transit, microtransit ($\Delta=1,200$ m., $K=0$), and ride-sharing.}
    \label{T:los}
    \begin{center}
    \begin{tabular}{llcccccc}
    \toprule[1pt]
    Mode & Design & Coverage  & Walk  & Wait & Detour  &Delay   & Distance\\  \hline
Transit        & 5 candidate lines       & 13.9\%     & 2.06     & 7.06       & 158.56\%     & -1.17       &   356  \\
               & 10 candidate lines      & 20.4\%     & 2.21     & 6.91       & 146.22\%     & -0.79       &   384  \\
               & 25 candidate lines      & 29.8\%     & 2.03     & 6.8        & 136.98\%     &  0.13       &   435  \\
               & 50 candidate lines      & 33.6\%     & 2.03     & 6.65       & 137.34\%     & -0.06       &   472  \\\hline
Microtransit   & 5 candidate lines       & 22.3\%     & 1.68     & 6.22       & 159.99\%     & -0.01       &   419  \\
               & 10 candidate lines      & 30.0\%     & 1.68     & 6.22       & 146.11\%     & -0.15       &   462  \\
               & 25 candidate lines      & 35.6\%     & 1.53     & 5.82       & 138.52\%     & -0.16       &   471  \\
               & 50 candidate lines      & 36.6\%     & 1.36     & 5.55       & 141.00\%     &  0.03       &   468  \\\hline
Rideshare      & Cap. 4        & 36.3\%     & 0        & 4.2        & 150.68\%     & 13.4        &  1,883  \\
               & Cap. 2        & 44.7\%     & 0        & 3.74       & 124.60\%	    &  8.17       &  3,359  \\
               & Cap. 1        & 50.5\%     & 0        & 1.79       & 100.00\%     &  1.79       &  5,671  \\
         \bottomrule[1pt]
    \end{tabular}
    \end{center}
    \begin{tablenotes}
        \footnotesize
        \vspace{-6pt}
        \item Coverage: percentage of served requests; distance in kilometers; walk, wait, delay/earliness in minutes.
    \end{tablenotes}
    \vspace{-12pt}
\end{table}

These results confirm that microtransit increases demand coverage and reduces walk and wait times from fixed-route transit, at virtually no cost in terms of detours and delays. The differences in coverage go down with more candidate lines. With 5--10 lines, the operating flexibility in microtransit is most valuable in sparse networks. With 25--50 lines, the transit network adds more reference trips in high-density regions. Since the microtransit network already achieves higher coverage with fewer lines, it allocates more lines to lower-demand regions (Section~\ref{subsec:flex}), resulting in a smaller marginal increase in coverage. At the other extreme, ride-sharing achieves high coverage with no walking (by design) and short waits, but much longer distances traveled with low-occupancy vehicles. Thus, microtransit defines a middle ground with less walk and less wait for passengers than transit, less delays than ride-sharing, and intermediate ridership and costs.

Another interesting observation stems from the comparison of microtransit to ride-pooling, both of which leverage on-demand operations to consolidate demand into multi-occupancy vehicles. On-demand door-to-door ride-pooling results in no walk and low wait times but increases detours and delays---underscoring the impact of spatiotemporal externalities, even with small-occupancy vehicles. By consolidating demand into high-capacity vehicles along reference lines, deviated fixed-route microtransit reduces distance traveled by a factor of 4 but reaches similar demand coverage and a comparable level of service---no delay, smaller detours, moderate walking times, and slightly longer wait times. These results identify deviated fixed-route microtransit as a possible pathway to provide efficient and convenient urban mobility options with high-capacity vehicles.

\begin{figure}
    \centering
   \subfloat[\footnotesize Average total distance traveled.]{\label{F:footprint}\includegraphics[height=0.3\textwidth]{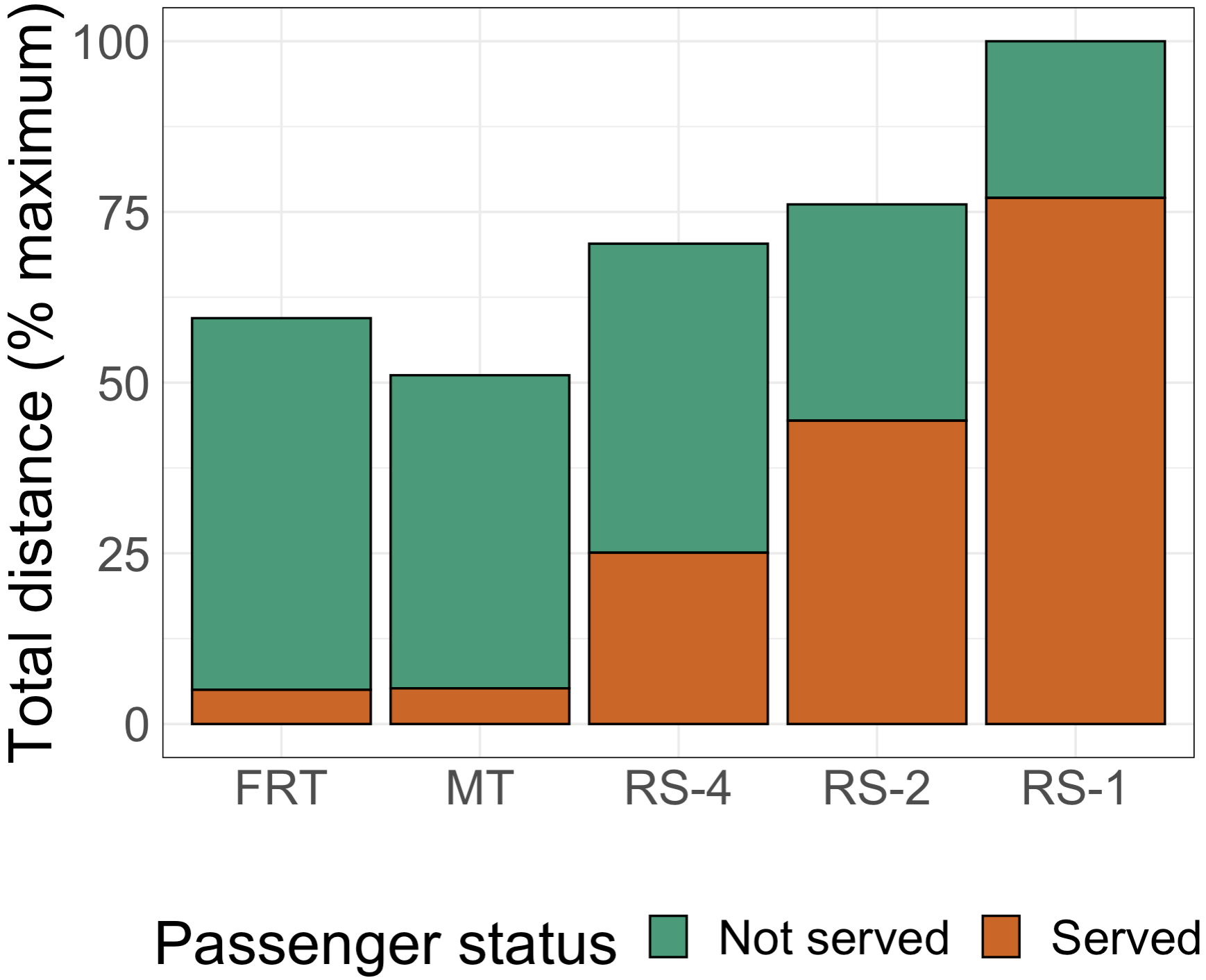}}
   ~ ~
   \subfloat[\footnotesize Average vehicle load.]{\label{F:load}\includegraphics[height=0.3\textwidth]{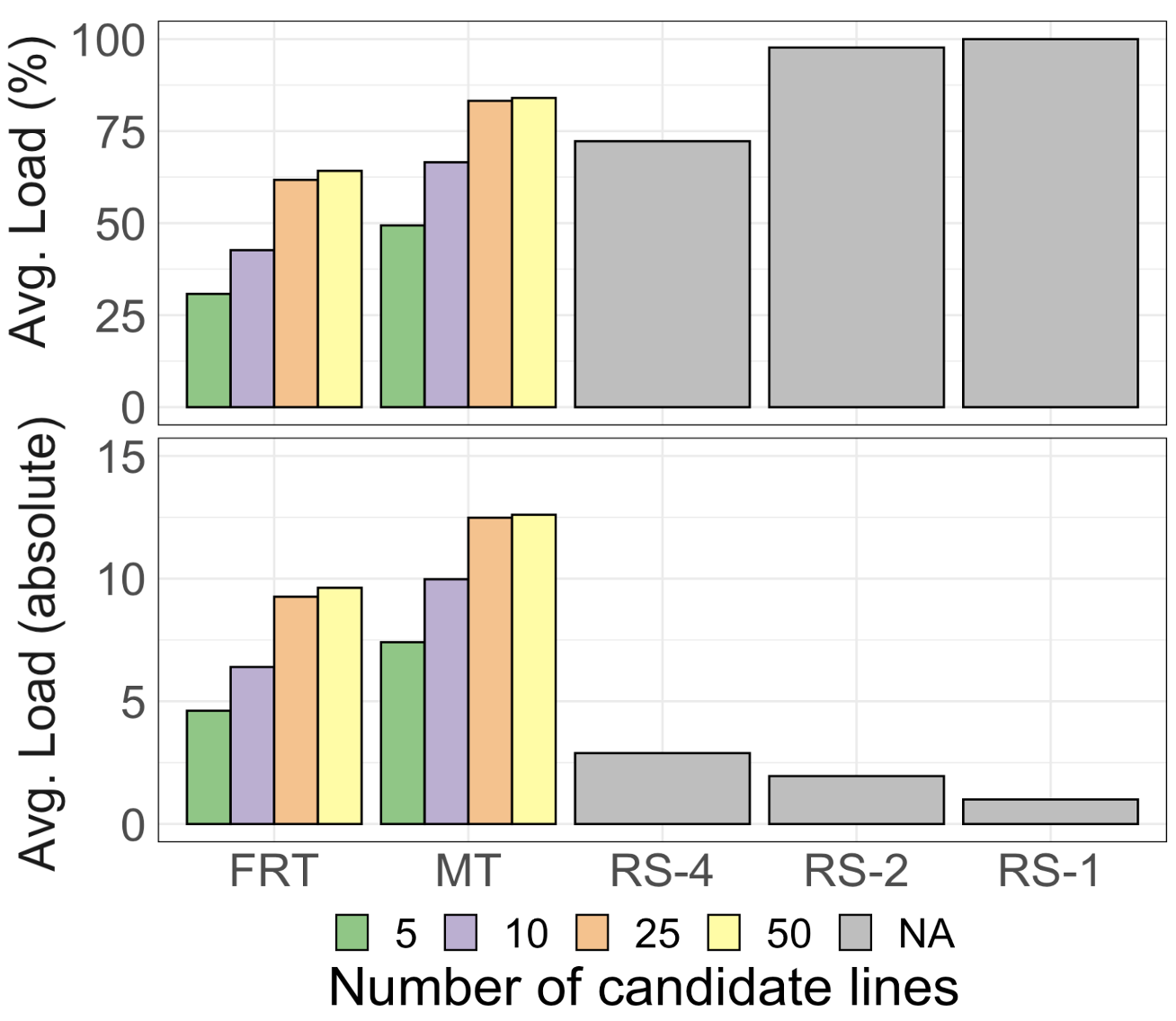}}
    \caption{System-wide assessment of fixed-route transit (FRT), microtransit ($\Delta=1,200$ m., $K=0$) (MT), and ride-sharing systems with capacities 4, 2 and 1 (RS-4, RS-2 and RS-1, respectively).}
    \label{F:coverageVsFootprint}
    \vspace{-12pt}
\end{figure}

Figure~\ref{F:coverageVsFootprint} provides an out-of-sample system-wide assessment. Figure~\ref{F:footprint} plots total distance traveled, used as a proxy of operating costs and environmental footprint; it includes both the ``internal'' distance within the system plus the ``external'' distance from single-occupancy trips (e.g., taxi trips) for all unserved passengers. Microtransit reduces total distance by 10-15\% versus fixed-route transit, by 20-30\% versus ride-pooling, and by 50\% versus single-occupancy ride-sharing. From Table~\ref{T:los}, these benefits are driven by a much smaller internal distance than ride-sharing that outweigh the impact of smaller demand coverage, and by higher demand coverage than fixed-route transit that outweigh the slightly longer internal distances. Altogether, these results identify potential operating, economic and environmental benefits from deviated fixed-route microtransit from stronger demand consolidation than ride-sharing and from higher demand coverage than fixed-route transit.
 
We establish the robustness of these results in~\ref{app:sensitivity}, and extend them in~\ref{app:DARcase} and~\ref{app:TRcase} to the MiND-DAR and MiND-Tr. For instance, in the MiND-DAR, deviated fixed-route microtransit increases demand coverage versus fixed-route transit, improves demand consolidation versus ride-sharing and ride-pooling, and reduces total distance versus all benchmarks (by 5-15\% vs. fixed-route transit, by 40-50\% vs. ride-pooling, and by over 100\% vs. single-occupancy ride-sharing).

In conclusion, deviated fixed-route microtransit can contribute to efficient, equitable, and sustainable mobility. Efficiency stems from high levels of service, low operating costs and high demand coverage enabled by reference lines and on-demand flexibility. Equity stems from a microtransit network with broader geographic reach. Sustainability stems from a smaller distance traveled per passenger enabled by high coverage and consolidation into high-capacity vehicles.

\section{Conclusion}

This paper optimizes the design and operations of a deviated fixed-route microtransit system endowed with advance planning capabilities along reference lines (as in public transit) and on-demand adjustments in response to passenger demand (as in ride-sharing). We formulated a \textit{Microtransit Network Design (MiND)} model via two-stage stochastic integer optimization. The model features a first-stage network design and service scheduling structure and a second-stage capacitated vehicle routing structure with time windows. We proposed a subpath-based representation of microtransit operations in a load-expanded network, which results in a tight second-stage relaxation but an exponential number of variables. To solve it, we have developed a double-decomposition algorithm, leveraging Benders decomposition to decompose the problem per scenario and reference trip, as well as subpath-based column generation to further decompose operations between checkpoints.

Using New York City data, results showed that the methodology scales to real-world and otherwise-intractable problems, with up to 100 candidate reference lines, hundreds of stations, thousands of requests, and 5-20 demand scenarios. Practical results suggested that even limited on-demand flexibility can provide significant operating benefits by consolidating demand into high-capacity vehicles while leveraging on-demand deviations to enhance passenger level of service and increase demand coverage. At a time where hybrid solutions are emerging to design new mobility services, this paper suggests that deviated fixed-route microtransit can contribute to efficient mobility (high demand coverage, low operating costs, high levels of service), equitable mobility (high accessibility with broad geographic reach), and sustainable mobility (low environmental footprint). Based on these results, we have been collaborating with transit operators toward the deployment of deviated fixed-route microtransit, with a pilot implementation targeted in 2025.

\section*{Acknowledgments}
The authors thank Joseph Kajon for helpful assistance with data analyses and experimentation.

\bibliographystyle{informs2014}
\bibliography{ref}

\newpage
\ECSwitch
\ECHead{A Double Decomposition Algorithm for Network Planning and Operations in Deviated Fixed-route Microtransit\\Electronic Companion}

\section{Details on Model Formulations}\label{app:segment_path_benchmarks}

\subsection{Notation Tables}
\label{app:notation}

Table \ref{T:notation} summarizes all notation for the MiND-VRP formulation.
\begin{table}[h!]
    \centering
{\footnotesize \begin{tabular}{llp{13.5cm}}
\toprule[1pt]
    \bf Component & \bf Type & \bf Description \\ \hline
     $\calN^S$ & Set & Stations: checkpoints and pickup locations \\
    $\calE$  & Set &Directed arcs in $\mathcal{N} \times \mathcal{N}$ corresponding to roadways \\
    $\calL$ & Set & Candidate reference lines\\
    $\calP$ & Set & Passenger types\\
    $\calS$ & Set & Demand scenarios\\
    $\calC_\ell$ & Set & Vehicle loads on reference line $\ell\in\calL$ \\
    $\calI_{\ell}$ & Set & Checkpoints for line $\ell \in \calL$, of cardinality $I_{\ell}$\\
    $\calI_{\ell}^{(i)}$ & Set & $i^{th}$ checkpoint in reference line $\ell \in \calL$ for $i=1, \cdots, I_{\ell}$ \\
    $\Gamma_\ell$ &Set & Subset of checkpoint pairs in $\calI_\ell \times \calI_\ell$ for line $\ell \in \calL$ that skip up to $K$ checkpoints in between \\
    $\mathcal{N}_{uv} $ & Set & Subset of nodes in $\calN^S$ representing possible stations between checkpoints $u,v \in \calI_\ell$ for each line $\ell \in \calL$ \\
    $\calT_{\ell}$ & Set & Allowable departure times of a vehicle from the beginning of line $\ell \in \calL$  \\
    $\calT^{uv}_{\ell t}$ & Set &  Time intervals between the scheduled times $T_{\ell t}(u)$ and $T_{\ell t}(v)$ for checkpoint pair $(u,v) \in \Gamma_\ell$\\ 
    $\calM_p$ & Set & Compatible trips in $\calL \times \calT_{\ell}$ for passenger type $p \in \calP$\\
    $\calR_{\ell st}$& Set& Subpaths corresponding to reference trip $(\ell, t) \in \calL\times\calT_{\ell}$ in scenario $s \in \calS$.\\
    &&Each subpath $r\in\calR_{\ell st}$ starts in $u_r\in\calI_{\ell}$ and ends in $v_r\in\calI_{\ell}$.\\
    $(\calV_{\ell st}, \calA_{\ell st})$& Graph & Load-expanded network of trip $(\ell, t) \in \calL \times \calT_{\ell}$ in scenario $s \in \calS$. \\
    &&Every trip starts at $u_{\ell st}\in\calV_{\ell st}$ and ends at $v_{\ell st}\in\calV_{\ell st}$\\
    $\calA_r$ & Set  & Arcs in $\calA_{\ell st}$ corresponding to subpath $r \in \calR_{\ell st}$ for $(\ell, t) \in \calL \times \calT_{\ell},$ $s \in \calS$ \\
    $\calA^v_{\ell st}$ & Set  & Arcs in $\calA_{\ell st}$ connecting line destination to sink node for $(\ell, t) \in \calL \times \calT_{\ell},$ $s \in \calS$\\
    $(\calU^{uv}_{\ell st}, \calH^{uv}_{\ell st})$& Graph&Time-expanded network from $(u,T_{\ell t}(u))$ to $(v,T_{\ell t}(v))$. Node $m \in \calU^{uv}_{\ell st}$ is characterized by a location-time tuple $(k_{m},t_{m})$ \\
    $\calP_{m}$ & Set & Passengers in $\calP$ that can be picked up in node $m \in \calU^{uv}_{\ell st}$ \\
    $\calP_r$ & Set & Passenger types in $\calP$ picked up by subpath $r \in \calR_{\ell st}$ for $(\ell, t) \in \calL\times \calT_{\ell}$, $s \in \calS$ \\\hline
    $K$ & Parameter & Number of consecutive checkpoints that can be skipped (0 or 1) \\
    $C_\ell$ & Parameter &  Vehicle capacity on reference line $\ell\in\calL$ \\
    $F$ & Parameter & Fleet size\\
    $h_{l}$ & Parameter & Cost to operate one trip via line $\ell \in \calL$\\
    $D_{ps}$ & Parameter & Number of passengers of type $p \in \calP$ in scenario $s \in \calS$\\
    $T_{\ell t}(n)$ & Parameter & Time at which trip $(\ell, t) \in \calL \times \calT_{\ell}$ must visit checkpoint $n \in \calI_{\ell}$ \\
    $\pi_s$ &Parameter & Probability of scenario $s \in \calS$\\
    $g_a$ &Parameter & Cost of arc $a \in \calA_{\ell st}$ for trip $(\ell, t) \in \calL \times \calT_{\ell},$ scenario $s \in \calS$ (Equation~\eqref{eq:arccost})\\
    $\Delta$ & Parameter & Maximum vehicle deviation from reference line \\
    $\Omega$ & Parameter & Maximum walking distance for passengers \\
    $\Psi$ & Parameter & Maximum waiting time for passengers \\
    $\alpha$ & Parameter & Time window radius around  passengers' requested drop-off times to build $\calM_p$ \\
    $\omega_{o,d}$ & Parameter & Walking distance between locations $o$ and $d$ \\
    $\psi_{o,d}$ & Parameter & Walking time between locations $o$ and $d$ \\
    $\tau_{rp}^{\text{walk}}$ &Parameter & Walk time of passenger $p \in \calP_r$ via subpath $r \in \calR_{\ell st}$, for $(\ell, t) \in \calL \times \calT_{\ell}$, $s \in \calS$ \\
    $\tau_{rp}^{\text{wait}}$& Parameter & Wait time of passenger $p \in \calP_r$ via subpath $r \in \calR_{\ell st}$, for $(\ell, t) \in \calL \times \calT_{\ell}$, $s \in \calS$ \\
    $\tau^{\text{travel}}_{rp}$&Parameter& In-vehicle time of passenger $p\in \calP_r$ via subpath $r \in \calR_{\ell st}$, for $(\ell, t) \in \calL \times \calT_{\ell}$, $s \in \calS$ \\
    $\tau_{\ell tp}^{\text{late}}$& Parameter & Delay of passenger type $p \in \calP$ when taking trip $(\ell, t) \in \calL \times \calT_{\ell}$ \\
    $\tau_{\ell tp}^{\text{early}}$& Parameter & Earliness of passenger type $p \in \calP$ when taking trip $(\ell, t) \in \calL \times \calT_{\ell}$ \\
    $\tau_p^{\text{dir}}$& Parameter & Direct travel time for passenger type $p \in \calP$\\
    $tt(e)$ & Parameter & Travel time corresponding to road segment $e \in \calE$ \\
    $\tau_{mp}^{\text{walk}}$  &Parameter & Walk time of passenger $p \in \calP_{m}$ when picked up at node $m \in \calU^{uv}_{\ell st}$\\
    $\tau_{mp}^{\text{wait}}$ & Parameter &Wait time of passenger $p \in \calP_{m}$ when picked up at node $m \in \calU^{uv}_{\ell st}$\\
    $\tau_{mp}^{\text{travel}}$  &Parameter & In-vehicle travel time of passenger $p \in \calP_m$ when picked up at node $m \in \calU^{uv}_{\ell st}$ \\
    $M$ & Parameter &  Reward for each passenger pickup\\
    $\lambda, \mu,\sigma,\delta $&Parameters &Penalties on passenger walk time, wait time, detour, and displacement\\
    $\kappa$ &Parameter & Target vehicle load in the first-stage network design problem\\
    \bottomrule[1pt]
\end{tabular}}
\caption{Notation for the MiND-VRP model and its decomposition.}\label{T:notation}
\end{table}

\subsection{Compact Formulation for the Second-stage Problem}
\label{app:compact}

We define the following additional parameters:
\begin{align*}
    \calN^S_{\ell} &= \text{set of stations for line $\ell$, within a distance $\Delta$ to the reference line}\\
    \calN^-_{\ell i} &= \text{set of locations that can be visited immediately after location $i$ on line $\ell$}\\
    \calN^+_{\ell i} &= \text{set of locations that can be visited immediately before location $i$ on line $\ell$}\\
    \calN^{\text{pickup}}_{\ell t p} &= \text{set of possible pickup locations for passenger $p$ on trip $(\ell, t)$, within $\Omega$ of their origin}\\
    \tau^{\text{req}}_p &= \text{requested time for passenger $p$}\\
    \tau^{\text{walk}}_{ip}&= \text{walking time from origin of passenger $p$ to pickup location $i$}\\
    \tau_{\ell t}^{\text{dropoff}} &= \text{dropoff time for reference trip $(\ell,t)$ at the destination}
\end{align*}
We define the following decision variables:
\begin{align*}
    v_{\ell t s i} &= \begin{cases}
    1 &\text{if reference trip $(\ell,t)$ visits checkpoint $i \in \calI_\ell$ in scenario $s\in\calS$,}\\
    0 &\text{otherwise.}
    \end{cases} \\
    y_{\ell t s i j} &= \begin{cases}
    1 &\text{reference trip $(\ell,t)$ travels from location $i\in\calN^S_\ell$ to location $j\in\calN^S_\ell$ in scenario $s\in\calS$,}\\
    0 &\text{otherwise.}
    \end{cases} \\
    w_{\ell t s p i} &= \begin{cases}
    1 &\text{reference trip $(\ell,t)$ picks up passenger $p\in\calP$ at location $i\in\calN^S_\ell$ in scenario $s\in\calS$,}\\
    0 &\text{otherwise.}
    \end{cases} \\
    t^{\text{stop}}_{\ell t s i} &= \text{time at which reference trip $(\ell,t)$ stops at location $i\in\calN^S_\ell$ in scenario $s\in\calS$} \\
    t^{\text{pickup}}_{\ell t s p} &= \text{time at which reference trip $(\ell,t)$ picks up passenger $p\in\calP$ in scenario $s\in\calS$} 
\end{align*}

The compact formulation of the second-stage problem in scenario $s \in \calS$ and for reference trip $(\ell,t)$ is then given as follows:
\begin{align}
    \min \quad & \sum_{p \in \calP} D_{ps} \Bigg( \lambda \sum_{i \in \calN^{\text{pickup}}_{\ell t p}} \tau^{\text{walk}}_{ip} w_{\ell t s p i} + \mu \left(t^{\text{pickup}}_{\ell t s p} -  \sum_{i \in \calN^{\text{pickup}}_{\ell t p}} (\tau^{\text{req}}_p + \tau^{\text{walk}}_{ip}) w_{\ell t s p i} \right) \nonumber \\
    & \quad + \frac{\sigma}{\tau_p^{dir}} \left( \tau_{\ell t}^{\text{dropoff}} \sum_{i \in \calN^{\text{pickup}}_{\ell t p}} w_{\ell t s p i} - t^{\text{pickup}}_{\ell t s p} \right) + \left(\delta \frac{\tau^{\text{late}}_{\ell tp}}{\tau_p^{\text{dir}}} + \frac{\delta}{2} \frac{\tau^{\text{early}}_{\ell tp}}{\tau_p^{\text{dir}}}\right)\sum_{i \in \calN^{\text{pickup}}_{\ell t p}}  w_{\ell t s p i} \nonumber\\
    &\quad -M\sum_{i \in \calN^{\text{pickup}}_{\ell t p}}  w_{\ell t s p i}\Bigg)\label{comp_obj}\\
    \text{s.t.} \quad & \sum_{j \in \calN^-_{\ell i}} y_{\ell ltsij} - \sum_{j \in \calN^+_{\ell i}} y_{\ell ltsji} = \begin{cases}
    x_{\ell t} &\text{ if $i = o_\ell$}\\
    -x_{\ell t} &\text{ if $i = d_\ell$}\\
    0 &\text{ otherwise.}
    \end{cases} \quad \forall i \in \calN^S_{\ell} \label{comp_fb}\\
    & t^{\text{stop}}_{\ell t s j} \geq t^{\text{stop}}_{\ell t s i} + tt_{ij} - M^{\text{stop}}(1 - y_{\ell t s i j}) \quad \forall i \in \calN^S_{\ell}, j \in \calN^-_{\ell i} \label{comp_tt}\\
     &\sum_{i \in \calN^{\text{pickup}}_{\ell t p}} w_{\ell t s p i} \leq z_{p \ell s t} \quad \forall p \in \calP \label{comp_wzconsist}\\
     &\sum_{p \in \calP} \sum_{i \in \calN^{\text{pickup}}_{\ell t p}} D_{ps} w_{\ell t s p i} \leq C_\ell x_{\ell t} \label{comp_capacity}\\
    & t^{\text{pickup}}_{\ell t s p} \geq t^{\text{stop}}_{\ell t s i} - M^{\text{stop}} (1 - w_{\ell t s p i} ) \quad \forall p \in \calP, i \in \calN^S_{\ell} \label{comp_tconsist}\\
    & w_{\ell t s p i} \leq \sum_{j \in  \calN^-_{\ell i}} y_{\ell t s i j} + \sum_{j \in  \calN^+_{\ell i}} y_{\ell t s j i}   \quad  \forall p \in \calP, i \in \calN^{\text{pickup}}_{\ell t p} \label{comp_ywconsist} \\
    & T_{\ell t}(i) - M^T (1 - v_{\ell t s i}) \leq t^{\text{stop}}_{\ell t s i} \leq T_{\ell t}(i) + M^T (1 - v_{\ell t s i}) \quad \forall i \in \calI_\ell \label{comp_cpstop} \\
     & t^{\text{stop}}_{\ell t s i} \leq M^{y} \cdot \left(\sum_{j \in  \calN^-_{\ell i}} y_{\ell t s i j} + \sum_{j \in  \calN^+_{\ell i}} y_{\ell t s j i}  \right) \quad \forall i \in \calN^S_{\ell}\label{comp_ty}\\
    & \sum_{j \in  \calN^-_{\ell i}} y_{\ell t s i j} + \sum_{j \in  \calN^+_{\ell i}} y_{\ell t s j i}   \geq v_{\ell t s i} \quad \forall i \in \calI_\ell \label{comp_yvconsist} \\
    & \sum_{j=i}^{i+K} v_{\ell t s \calI_\ell^{(j)}} \geq x_{\ell t} \quad \forall i = 1,\ldots, I_\ell-K \label{comp_skipcp}\\
    & \bv, \by, \bw \text{ binary} \\
    & \bt^{\text{stop}}, \bt^{\text{pickup}} \geq \boldsymbol{0} 
\end{align}

The objective \eqref{comp_obj} minimizes the weighted sum of walking time, waiting time, relative detour, and relative delay across passengers, minus the reward for passenger pickups. Constraints \eqref{comp_fb} enforce flow balance at each station. Constraints \eqref{comp_tt} coordinate travel time between locations on each trip, with a big-M parameter $M^{\text{stop}}$ activating them for pairs of consecutive locations. Constraints \eqref{comp_wzconsist} and \eqref{comp_capacity} ensure consistency with the first-stage solution, and Constraints \eqref{comp_capacity} enforce vehicle capacity. Constraints~\eqref{comp_tconsist} determine the pick-up time for each passenger as the vehicle stop time, again with big-M parameter $M^{\text{pickup}}$ activating them for appropriate passenger-location pairs. Constraints \eqref{comp_ywconsist} ensure that passengers are only assigned to pick-up locations where the vehicle stops. Constraints~\eqref{comp_cpstop} ensure that the checkpoints on selected reference trips are visited at the pre-specified time, with big-M parameter $M^{T}$ activating them only for visited checkpoints. Constraints~\eqref{comp_ty} define the stopping time only for those locations visited by the vehicles, again with a corresponding big-M parameter $M^y$. Constraints \eqref{comp_yvconsist} ensure consistency of the checkpoint stops with the routing variables, and Constraints \eqref{comp_skipcp} ensure that the vehicle skips at most $K$ consecutive checkpoints.

Note that the definition of the $t^{\text{stop}}_{\ell t s i}$ variables, along with Constraints \eqref{comp_tt}, impose the implicit condition that each reference trip visits each stop at most once in each scenario. In principle, we could relax that assumption by defining variables of the form $t^{\text{stop}}_{\ell t s i \kappa}$, where $\kappa$ counts the number of times reference trip $(\ell,t)$ stops at location $i\in\calN^S$ in scenario $s\in\calS$. The formulation would be augmented with precedence constraints accordingly. We omit these details for simplicity.

This formulation features a polynomial number of decision variables and constraints, but suffers from a weak linear relaxation due to the several big-M constraints (Equations~\eqref{comp_tt}, \eqref{comp_tconsist}, \eqref{comp_cpstop}, and \eqref{comp_ty}). Results in Section~\ref{subsec:method} and in~\ref{app:subpath} show that, in turn, it features much more limited scalability than our network-based reformulation.

\subsection{Proof of Proposition~\ref{prop:compact}}

Consider scenario $s$ and reference trip $(\ell,t)$. Recall that arc $a=(m,n)\in \calA_{\ell s t}$ corresponds to a subpath $r(a)$ defined between two checkpoints $u_r,v_r \in \calI_\ell$. Each node $m\in\calV_{\ell s t}$ corresponds to a checkpoint-load pair $(k_m, c_m)$ on the load-expanded network.

We introduce additional notation to describe the series of stops and passenger pickups for each subpath. We index the checkpoints by $\text{ind}(i)$, so that $\text{ind}(i)<\text{ind}(j)$ if checkpoint $i\in\calI_\ell$ is visited earlier than checkpoint $j\in\calI_\ell$ (if they are visited at all). The route of subpath $r$ is defined by a sequence of $H$ stops, indexed by $h=1,\cdots,H$. Each stop $h$ encodes a location-time pair $\text{Stops}(r) = \{(i^h,t^h) \text{ for }h=1,\ldots,H\}$. For ease of notation, we also define the set of \emph{pairs} of consecutive location-time pairs, $\text{Path}(r) = \{((i^h,t^h), (i^{h+1},t^{h+1})), \: \forall h = 1,\ldots, H-1\}$. Finally, let $\text{Pax}(r,i)$ be the set of passengers picked up at location $i$ on subpath $r$. Each subpath satisfies four conditions:
\begin{enumerate}[label=(\roman*)]
    \item The distance to the reference line never exceeds $\Delta$
    \item The load satisfies $\sum_{p \in \calP_r} D_{ps} \leq C_\ell$ 
    \item The travel time does not exceed $T_{\ell t}(v_r) - T_{\ell t}(u_r)$
    \item Up to $K$ checkpoints are skipped
\end{enumerate}

Let  $(\bv^*, \bw^*, \by^*, \bt^{\text{stop}*}, \bt^{\text{pickup}*})$ be an optimal solution to the compact formulation. We construct an equivalent feasible solution $\by$ for the subpath formulation with the same objective value.  Let $\calB_{ij}$ be the set of intermediate stops between visited checkpoints $i$ and $j$. Specifically, set $\calB_{ij} = \emptyset$ for $i,j \in \calI_\ell$ such that $v^*_{\ell t s i} = v^*_{\ell t s j} = 1$ and $v^*_{\ell t s k} = 0$ for all $k\in\calI_\ell$ such that $\text{ind}(i) < \text{ind}(k) < \text{ind}(j)$, i.e., if $i$ and $j$ are consecutive checkpoints in the solution. Then, construct $\calB$ iteratively by starting from $i=k_m$, then identifying $j$ such that $y^*_{\ell t s i j}=1$ and setting $\calB \gets \calB \cup \{j\}$. Continue until $j = k_n$. 

We construct $\by^*$ by setting $y^*_a = 1$ for $a=(m,n) \in \calA_{\ell s t}$ if:
\begin{enumerate*}[label=(\arabic*)]
    \item $v^*_{\ell t s i} = 1$ for $i=k_m$ and $i=k_n$;
    \item $v^*_{\ell t s i} = 0$ for all $i\in\calI_\ell$ such that $\text{ind}(k_m) < \text{ind}(i) < \text{ind}(k_n)$;
    \item $\text{Pax}(r(a),i) = \{p\in\calP : w^*_{\ell t s p i}=1\}$ for $i \in \calB_{k_m,k_n}$; and
    \item the load $c_m$ is equal to the total number of passengers picked up along subpaths prior to node $m$: $\displaystyle c_m = \sum_{a^\prime = (m^\prime, n^\prime) \in \calA_{\ell s t} \setminus \{a\} :c_{n^\prime} \leq c_m} \sum_{p \in \calP_{r(a^\prime)}} D_{ps} y_a$
\end{enumerate*}
In particular, condition (3) implies that $\calP_{r(a)} = \{p\in\calP : \sum_{i \in \calB_{k_m,k_n}} w^*_{\ell t s p i}=1\}$. Otherwise, $y_a = 0$.

We first show that for each $a$ with $y_a=1$, $r(a)$ is a valid subpath.

\textbf{Condition (i).} The distance to the reference line never exceeds $\Delta$ by construction, since we only consider $i \in \calN^S_\ell$. 

\textbf{Condition (ii).} By construction of $\calP_{r(a)}$ and \eqref{comp_capacity}, 
\begin{align*} 
\sum_{p\in\calP_{r(a)}} D_{ps} = \sum_{p\in\calP_{r(a)}} \sum_{i\in \calN^{\text{pickup}}_{\ell t p}} D_{ps} w^*_{\ell t s p i} \leq \sum_{p\in\calP} \sum_{i\in \calN^{\text{pickup}}_{\ell t p}} D_{ps} w^*_{\ell t s p i} \leq C_\ell
\end{align*}

\textbf{Condition (iii).} Note that $t^{\text{stop}*}_{\ell t s k_m} = T_{\ell t}(k_m)$ and $t^{\text{stop}*}_{\ell t s k_n} = T_{\ell t}(k_n)$ by \eqref{comp_cpstop}. By summing over \eqref{comp_tt} for all $i,j \in \calB_{k_m k_n}$ such that $y^*_{\ell s t i j} = 1$, the travel time does not exceed $T_{\ell t }(k_n) - T_{\ell t }(k_m)$:
\begin{align*}
    \sum_{i \in \calB_{k_m k_n}}  \sum_{j \in \calB_{k_m k_n}: y^*_{\ell s t i j} = 1} t^{\text{stop}*}_{\ell s t j} \geq \sum_{i \in \calB_{k_m k_n}}  \sum_{j \in \calB_{k_m k_n}: y^*_{\ell s t i j} = 1} ( t^{\text{stop}*}_{\ell s t i} + tt_{ij} ) \\
    t^{\text{stop}*}_{\ell s t k_n} \geq t^{\text{stop}*}_{\ell s t k_m} + \sum_{i \in \calB_{k_m k_n}} \sum_{j \in \calB_{k_m k_n}: y^*_{\ell s t i j} = 1} tt_{ij} \\
    T_{\ell t}(k_n) \geq T_{\ell t}(k_m) +  \sum_{i \in \calB_{k_m k_n}} \sum_{j \in \calB_{k_m k_n}: y^*_{\ell s t i j} = 1} tt_{ij}
\end{align*}

\textbf{Condition (iv).} By construction, $v^*_{\ell t s k_m} = 1$, $v^*_{\ell t s k_n} = 1$, and $v^*_{\ell t s i} = 0$ for all $i\in\calI_\ell$ such that $\text{ind}(k_m) < \text{ind}(i) < \text{ind}(k_n)$. By Constraints~\eqref{comp_skipcp}, $\sum_{i=\text{ind}(k_m)+1}^{\text{ind}(k_m)+K+1} v^*_{\ell t \calI_\ell^{(i)}} \geq 1$ and since $\text{ind}(k_n)$ is the next checkpoint index $i$ for which $v^*_{\ell t s i} = 1$, we must have $\text{ind}(k_n)\leq\text{ind}(k_m)+K+1$. Therefore, at most $K$ checkpoints are skipped between $k_m$ and $k_n$. 

Finally, we prove that $\by$ is feasible for the subpath formulation and achieves the same objective as $(\bv^*, \bw^*, \by^*, \bt^{\text{stop}*}, \bt^{\text{pickup}*})$. Denote the objective of the compact formulation as $\texttt{OPT}$:
\begin{align*}
    \texttt{OPT} &= \sum_{p\in\calP} \sum_{i \in \calN^{\text{pickup}}_{\ell t p}} D_{ps} \left(\frac{\delta\tau^{\text{late}}_{\ell t p}+\frac{\delta}{2}\tau^{\text{early}}_{\ell t p} + \sigma \tau^{\text{dropoff}}_{\ell t}}{\tau^{\text{dir}}_p}+ \lambda \tau^{\text{walk}}_{ip} + \mu( - \tau^{\text{req}}_{p} - \tau^{\text{walk}}_{ip}) -M\right) \cdot w^*_{\ell t s p i}   \\
    &\quad\quad + \sum_{p\in\calP} D_{ps} \left(\frac{-\sigma }{\tau^{\text{dir}}_p} + \mu \right) \cdot t_{\ell t s p}^{\text{pickup}*}\\
\end{align*}

By Constraints~\eqref{comp_wzconsist} and by condition (3) in the construction of subpath $a=(m,n)$:
\begin{align}
    w^*_{\ell t s p i} &= \sum_{a \in \calA_{\ell s t}: \text{Pax}(r(a),i)} y_a \quad \forall p\in\calP, i\in\calN^{\text{pickup}}_{\ell t p} \label{wy_equiv} 
\end{align}
Furthermore, $ t_{\ell t s p}^{\text{pickup}*}$ can be written as $\sum_{i \in \calN^{\text{pickup}}_{\ell t p}} \sum_{a \in \calA_{\ell s t}: p \in \text{Pax}(r(a),i)} \tau^{\text{stop}}_{i}(a) \cdot y_a$ where $\tau^{\text{stop}}_{i}(a)$ is the time arc $a$ stops at location $i$. Then,
\begin{align*}   
    \texttt{OPT} &= \sum_{p\in\calP} \sum_{i \in \calN^{\text{pickup}}_{\ell t p}} \sum_{a \in \calA_{\ell s t}: \text{Pax}(r(a),i)} D_{ps} \left(\frac{\delta\tau^{\text{late}}_{\ell t p}+\frac{\delta}{2}\tau^{\text{early}}_{\ell t p} + \sigma \tau^{\text{dropoff}}_{\ell t}}{\tau^{\text{dir}}_p}+ \lambda \tau^{\text{walk}}_{ip} + \mu( - \tau^{\text{req}}_{p} - \tau^{\text{walk}}_{ip}) -M\right) \cdot y_a   \\
    &\quad\quad + \sum_{p\in\calP} \sum_{i \in \calN^{\text{pickup}}_{\ell t p}} \sum_{a \in \calA_{\ell s t}: \text{Pax}(r(a),i)} D_{ps} \left(\frac{-\sigma \tau^{\text{stop}}_{i}(a)}{\tau^{\text{dir}}_p} + \mu \tau^{\text{stop}}_{i}(a) \right) \cdot y_a \\
    &= \sum_{a \in \calA_{\ell s t}} 
    \sum_{p\in\calP_{r(a)}} D_{ps} \left(\frac{\delta\tau^{\text{late}}_{\ell t p}+\frac{\delta}{2}\tau^{\text{early}}_{\ell t p} + \sigma(\tau^{\text{dropoff}}_{\ell t} - \tau^{\text{stop}}_{m_{ap}}(a))}{\tau^{\text{dir}}_p}+ \lambda \tau^{\text{walk}}_{m_{ap}, p} + \mu (\tau^{\text{stop}}_{m_{ap}}(a) - \tau^{\text{req}}_{p} - \tau^{\text{walk}}_{m_{ap},p}) - M\right) \cdot y_a,
\end{align*}
where $m_{ap}\in \calN^{\text{pickup}}_{\ell t p}$ denotes the stop on subpath $r(a)$ where passenger $p$ is picked up.

The objective of the subpath formulation relies on two preprocessed parameters, $\tau^{\text{wait}}_{mp}$ and $\tau^{\text{travel}}_{mp}$ (each dependent on arc $a$), which represent the wait time and in-vehicle time if passenger $p$ is picked up at location $m$, respectively. We can re-write them as $\tau^{\text{wait}}_{mp}(a) = \tau^{\text{stop}}_{m}(a) - \tau^{\text{req}}_{p} - \tau^{\text{walk}}_{mp}$ and $\tau^{\text{travel}}_{mp}(a) = \tau^{\text{dropoff}}_{\ell t} - \tau^{\text{stop}}_{m}(a)$. We get:
\begin{align*}   
    \texttt{OPT} &= \sum_{a \in \calA_{\ell s t}} \sum_{p\in\calP_{r(a)}} D_{ps} \left(\frac{\delta\tau^{\text{late}}_{\ell t p}+\frac{\delta}{2}\tau^{\text{early}}_{\ell t p} + \sigma\tau^{\text{travel}}_{m_{ap}, p}}{\tau^{\text{dir}}_p}+ \lambda \tau^{\text{walk}}_{m_{ap}, p} + \mu \tau^{\text{wait}}_{m_{ap}, p}-M\right) \cdot y_a \\
    &= \sum_{a \in \calA_{\ell s t}}  g_a y_a
\end{align*}

Conversely, let $\by^*$ be an optimal solution to the subpath formulation. We assume that, under that solution, each reference trip visits each stop at most once in each scenario. As discussed above, this assumption is embedded in the compact formulation itself. All arguments could be extended otherwise, but we omit these details for simplicity. In practice, it is highly unlikely that a reference trip would find it beneficial to stop multiple times in the same location.

We construct a feasible solution $(\bv, \bw, \by, \bt^{\text{stop}}, \bt^{\text{pickup}})$ for the compact formulation that corresponds to subpath solution $\by^*$ and has the same objective value. Specifically:
\begin{align}
    v_{\ell t s i} &= \sum_{a\in \calA_{\ell s t}| \exists t^\prime : (i,t^\prime) \in \text{Stops}(r(a))} y^*_a \label{comp:definev} \\
    y_{\ell t s i j} &= \sum_{a\in \calA_{\ell s t}| \exists t_1,t_2 : ((i,t_1),(j,t_2)) \in \text{Paths}(r(a))} y^*_a  \label{comp:definey}\\
    w_{\ell t s p i} &= \sum_{a\in \calA_{\ell s t}: p \in \text{Pax}(r(a),i)} y^*_a \label{comp:definew}\\
    t_{\ell t s i}^{\text{stop}} &= \sum_{a\in \calA_{\ell s t}| \exists t^\prime : (i,t^\prime) \in \text{Stops}(r(a))} \tau^{\text{stop}}_i(a) y^*_a \label{comp:definetstop}\\
    t_{\ell t s p}^{\text{pickup}} &= \sum_{a\in\calA_{\ell s t}: p\in\calP_{r(a)}} \tau^{\text{stop}}_{m_{ap}}(a) y^*_a \label{comp:definetpick}
\end{align}

From \eqref{comp:definew}, we have:
\begin{align}
    \sum_{i \in \calN^{\text{pickup}}_{\ell t p}} w_{\ell t s p i} &= \sum_{a \in \calA_{\ell s t}: p\in\calP_{r(a)}} y^*_a \quad \text{for $p\in\calP$} \label{eq:sum_w} \\
    \sum_{p\in\calP} \sum_{i \in \calN^{\text{pickup}}_{\ell t p}} w_{\ell t s p i} &= \sum_{a \in \calA_{\ell s t}} \sum_{p\in\calP_{r(a)}} y^*_a  \label{eq:sumsum_w}
\end{align}

We show that $(\bv, \bw, \by, \bt^{\text{stop}}, \bt^{\text{pickup}})$ satisfies all constraints of the compact formulation. 

\textbf{Constraints~\eqref{comp_fb}}. If $x_{\ell t} = 0$, then $y^*_a = 0$ for all $a \in \calA_{\ell s t}$, and $y_{\ell t s i j} = 0$ for all $i,j \in \calN^S_\ell$, satisfying the constraint. If $x_{\ell t} = 1$, then we have four cases.
\begin{enumerate}[label=(\arabic*)]
    \item  $i=\calI_\ell^{(1)}$:  By flow balance constraints \eqref{eq:2Sflow}, $\sum_{m:(n,m) \in \calA_{\ell s t}} y^*_{(u_{\ell s t },m)} = 1$, so there exists exactly one $m=(j,c)$ such that $y^*_{(u_{\ell s t },m)} = 1$. Note that $u_{\ell s t} = (\calI_\ell^{(1)}, T_{\ell t}(\calI_\ell^{(1)}))$. Thus, $y_{\ell t s i j} = 1$ and $y_{\ell t s i k} = 0$ for all $k \in \calN^-_{\ell i}, k \neq j$.
    
    \item $i=\calI_\ell^{(I_\ell)}$: By flow balance constraints \eqref{eq:2Sflow}, $\sum_{m:(m,n) \in \calA_{\ell s t}} y^*_{(m,v_{\ell s t})} = 1$, so there exists exactly one node $m^* = (j,c)$ such that $y^*_{(m^*,v_{\ell s t})} = 1$. Recall that $v_{\ell s t}$ is a dummy sink node and thus $j=\calI_\ell^{(I_\ell)}$. Again by \eqref{eq:2Sflow} for $n = (j,c)$, we have:
    \begin{align*}
        &\sum_{m:(m^*,m) \in \calA_{\ell st}} y^*_{(m^*, m)} - \sum_{m:(m,m^*) \in \calA_{\ell st}} y^*_{(m,m^*)} = 0 \\
        \implies &\sum_{m:(m,m^*) \in \calA_{\ell st}} y^*_{(m,m^*)} = 1
    \end{align*}
    We can then conclude that for exactly one node $m = (k,c^\prime)$, we have $ y^*_{(m,m^*)} = 1$. Therefore, $y_{\ell, t, s, k, \calI_\ell^{(I_\ell)}} = 1$ for exactly one $k\in \calN^+_{\calI_\ell^{(I_\ell)}}$ and the constraint is satisfied. 
    
    \item $i = \calI_\ell^{(2)}, \ldots, \calI_\ell^{(I_\ell-1)}$: Let $a = (m,n) = ((i,c_m),(k_n,c_n))$ and suppose $y^*_a=1$. By flow balance constraints \eqref{eq:2Sflow}, $\exists a^\prime = (o,m) = ((k_o,c_o), (i,c_m))$ such that $y^*_{a^\prime}=1$. 
    \begin{align*}
        y_a^* = 1 \implies ((i,t_1),(j,t_2)) \in \text{Paths}(r(a)) \implies y_{\ell t s i j} = 1 \\
        y_{a^\prime}^* = 1 \implies ((j^\prime,t_3),(i,t_4)) \in \text{Paths}(r(a^\prime)) \implies y_{\ell t s j^\prime i} = 1 
    \end{align*}
    Thus, Constraints~\eqref{comp_fb} are satisfied
    \item $i \notin \calI_\ell$: Assume $(i,t_1) \in \text{Stops}(r(a))$ for $a\in \calA_{\ell s t}$ with $y^*_a = 1$. Then, $(i,t_1)$ will appear exactly twice as the end point of an arc in $\text{Paths}(r(a))$: $((i,t_1), (j,t_2))$ for some $j$ and $((j^\prime,t_3), (i,t_1))$ for some $j^\prime$. By construction, we then have $y_{\ell t s i j} = y_{\ell t s j^\prime i} = 1$ and flow balance is satisfied. If $(i,t_1) \notin \text{Stops}(r(a))$ for any $t_1$, then $y_{\ell t s i j} = 0$ for all $j\in \calN_\ell^S$, again satisfying flow balance.
\end{enumerate}

\textbf{Constraints~\eqref{comp_tt}}. For $i\in \calN^S_\ell, j\in \calN^i_{\ell i}$ such that $y_{\ell t s i j}=0$,  Constraints~\eqref{comp_tt} is trivially satisfied. For $i,j$ with $y_{\ell t s i j}=1$, then $\exists t_1, t_2$ such that $((i,t_1), (j,t_2)) \in \text{Path}(r)$, and by construction of the subpaths, $t_2 = t_1 + tt_{ij}$. Thus, Constraints~\eqref{comp_tt} is satisfied for $i,j$.

\textbf{Constraints~\eqref{comp_wzconsist}}. By Equation~\eqref{eq:sum_w} and \eqref{eq:2Sy2z}:
\begin{align*}
    \sum_{i \in \calN_{\ell t p}^{\text{pickup}}} w_{\ell t s p i} = \sum_{a \in \calA_{\ell s t} : p \in \calP_{r(a)}} y^*_a \leq z_{p \ell s t} \quad \forall s \in \calS, p \in \calP, (\ell,t) \in \calM_p
\end{align*}

\textbf{Constraints~\eqref{comp_capacity}}.  If $x_{\ell t} = 0$, then $y^*_a = 0$ for all $a \in \calA_{\ell s t}$, and $w_{\ell t s p i} = 0$ for all $p\in \calP, i \in \calN^S_\ell$, satisfying the constraint. If $x_{\ell t} = 1$, by Equation~\eqref{eq:sumsum_w} we have:
\begin{align*}
   \sum_{p\in\calP} \sum_{i \in \calN^{\text{pickup}}_{\ell t p}} D_{ps} w_{\ell t s p i} &= \sum_{a \in \calA_{\ell s t}} \sum_{p\in\calP_{r(a)}} D_{ps} y^*_a  \leq C_\ell 
\end{align*}
The last inequality comes from the maximum load of a trip, which is given by $c_m$ (and $c_m \leq C_\ell$) for $m$ such that $a=(m,v_{\ell s t}) \in \calA_{\ell s t}^v$ and $y^*_a = 1$. By construction of the load-expanded network, $c_m$ is equal to the sum of the load differentials of all subpaths before it, namely $\sum_{a \in \calA_{\ell s t}} \sum_{p\in\calP_{r(a)}} D_{ps} y^*_a$.

\textbf{Constraints~\eqref{comp_tconsist}}. If $w_{\ell t s p i} = 0$ for $p\in\calP, i \in \calN_\ell^S$, then Constraints~\eqref{comp_tconsist} is trivially satisfied. If $w_{\ell t s p i} = 1$, then $\exists a \in \calA_{\ell s t}$ such that $p \in \text{Pax}(r(a),i)$ and $y^*_a = 1$. If $t^{\text{stop}}_{\ell t s i } = 0$, then Constraints~\eqref{comp_tconsist} is clearly satisfied. If $t^{\text{stop}}_{\ell t s i } = t^\prime > 0$, then by construction $(i,t^\prime) \in \text{Stops}(r(a^\prime))$ for some $a^\prime\in\calA_{\ell s t}$. Since each station is visited at most once, $a = a^\prime$ and $\exists i\in \calN_\ell^+, a\in\calA_{\ell s t}$ such that $(i,t^\prime) \in \text{Stops}(r(a))$, $p \in \text{Pax}(r(a),i)$, and $y^*_a = 1$. Thus, $t^{\text{pickup}}_{\ell t s i } = t^\prime$, satisfying Constraints~\eqref{comp_tconsist}.

\textbf{Constraints~\eqref{comp_ywconsist}}. If $w_{\ell t s p i} = 0$, the constraints are trivially satisfied. If $w_{\ell t s p i} = 1$, by construction, $\exists  a \in \calA_{\ell s t}$ such that $p \in \text{Pax}(r(a),i)$ and $y^*_a=1$. Therefore, $i \in \text{Stops}(r(a))$ and either $((i,t_1), (j,t_2)) \in \text{Paths}(r(a))$ or $((j,t_1), (i,t_2)) \in \text{Paths}(r(a))$ for some $j \in \calN^S_\ell$ and $t_1,t_2\geq 0$. We then conclude that Constraints~\eqref{comp_ywconsist} are satisfied:
\begin{align*}
    \sum_{j \in  \calN^-_{\ell i}} y_{\ell t s i j} + \sum_{j \in  \calN^+_{\ell i}} y_{\ell t s j i} \geq 1 = w_{\ell t s p i}
\end{align*} 

\textbf{Constraints~\eqref{comp_cpstop}}. For all $i \in \calI_\ell$ such that $v_{\ell t s i} = 0$, \eqref{comp_cpstop} is trivially satisfied. For $i \in \calI_\ell$ such that $v_{\ell t s i} = 1$, $\exists a\in \calA_{\ell s t}, t^\prime \in \mathbb{R}$ such that $(i,t^\prime) \in \text{Stops}(r(a))$ and $y_a^* = 1$ by \eqref{comp:definev}. Note that $\tau^{\text{stop}}_i(a) = T_{\ell t}(i)$ for $i\in\calI_\ell$. Thus, we have $t^{\text{stop}}_{\ell t s i} = T_{\ell t}(i)$, satisfying \eqref{comp_cpstop}.  

\textbf{Constraints~\eqref{comp_ty}}. If $t^{\text{stop}}_{\ell t s i} = 0$, the constraint is trivially satisfied. If $t^{\text{stop}}_{\ell t s i} > 0$ for $i\neq v_{r(a)}$, then $\exists j,t_1,t_2$ such that $((i,t_1),(j,t_2)) \in \text{Paths}(r(a))$ and $y_a^*=1$, thus $y_{\ell t s i j} = 1$. If $t^{\text{stop}}_{\ell t s i} > 0$ for $i = v_{r(a)}$, then $\exists j,t_1,t_2$ such that $((j,t_1),(i,t_2)) \in \text{Paths}(r(a))$ and $y_a^*=1$, thus $y_{\ell t s j i} = 1$. Thus, the constraint is satisfied.

\textbf{Constraints~\eqref{comp_yvconsist}}. If $v_{\ell t s i} = 0$, then the constraint is trivially satisfied. If $v_{\ell t s i} = 1$ for $i\in \calI_\ell$, then $\exists t^\prime$ such that $(i,t^\prime) \in \text{Stops}(r(a))$ and $y_a^*=1$. Because $i\in\calI_\ell$, checkpoint $i$ must be either the first or last stop for subpath $r(a)$ for $a=(m,n)$ (i.e. $k_m=i$ or $k_n=i$). If $k_m=i$, then $\exists j,t_1,t_2$ such that $((i,t_1),(j,t_2) \in \text{Path}(r(a))$ and by construction $y_{\ell t s i j} =1$ for $j\in\calN^-_{\ell i}$. If $k_n=i$, then $\exists j,t_1,t_2$ such that $((j,t_1),(i,t_2) \in \text{Path}(r(a))$ and by construction $y_{\ell t s j i} =1$ for $j\in\calN^+_{\ell i}$. Thus, the constraint is satisfied.

\textbf{Constraints~\eqref{comp_skipcp}}.
With a slight abuse of notation, we use $i\in\{1,\cdots,I_\ell-K\}$ to denote the index of a given checkpoint as well as the checkpoint itself. Due to flow balance constraints and the construction of the subpaths, all checkpoints in $\calI_\ell$ are either endpoints of a subpath $r(a)$ with $y^*_a = 1$ or are skipped by some subpath $r(a^\prime)$ with $y^*_{a^\prime} = 1$. In the first case, let $a=(m,n) \in \calA_{\ell s t}$ with either $k_m=i$ or $k_n=i$, and $y^*_a =1$. Then, $v_{\ell s t i}=1$ by construction and Constraints~\eqref{comp_skipcp} are satisfied. In the second case, there exist separate checkpoints $i^-,i^+\in\calI_\ell$, such that $i^- < i < i^+$ and $y^*_a=1$ with $a=(m,n)$ and $k_m=i^-$ and $k_n=i^+$. By property (iv), subpath $r(a)$ skips at most $K$ checkpoints between $i^-$ and $i^+$. Thus, $i^+\leq i^- + K+1$, hence $i^+\leq i + K$. This implies the constraint:
\begin{align*}
    \sum_{i^\prime = i}^{i + K} v_{\ell t s \calI_\ell^{(i^\prime)}} \geq 1
\end{align*}

\textbf{Objective~\eqref{comp_obj}}. Denote the objective of the compact formulation as $\texttt{OPT}$. Using Equation~\eqref{eq:sumsum_w} and the definition of $t^{\text{pickup}}_{\ell t s p} $ \eqref{comp:definetpick}, we have:
\begin{align*}
    \texttt{OPT} &= \sum_{a \in \calA_{\ell s t}} \sum_{p\in\calP_{r(a)}} D_{ps} \left(\frac{\delta\tau^{\text{late}}_{\ell t p}+\frac{\delta}{2}\tau^{\text{early}}_{\ell t p} + \sigma\tau^{\text{travel}}_{m_{ap}, p}}{\tau^{\text{dir}}_p}+ \lambda \tau^{\text{walk}}_{m_{ap}, p} + \mu \tau^{\text{wait}}_{m_{ap}, p}-M\right) \cdot y^*_a \\
     &= \sum_{a \in \calA_{\ell s t}} 
    \sum_{p\in\calP_{r(a)}} D_{ps} \frac{\delta\tau^{\text{late}}_{\ell t p}+\frac{\delta}{2}\tau^{\text{early}}_{\ell t p} + \sigma(\tau^{\text{dropoff}}_{\ell t} - \tau^{\text{stop}}_{m_{ap}}(a))}{\tau^{\text{dir}}_p}\cdot y^*_a\\
    &+ \sum_{a \in \calA_{\ell s t}} 
    \sum_{p\in\calP_{r(a)}}\left(\lambda \tau^{\text{walk}}_{m_{ap}, p} + \mu (\tau^{\text{stop}}_{m_{ap}}(a) - \tau^{\text{req}}_{p} - \tau^{\text{walk}}_{m_{ap},p}) - M\right) \cdot y^*_a \\
    &= \sum_{p\in\calP} \sum_{i \in \calN^{\text{pickup}}_{\ell t p}} D_{ps} \left(\frac{\delta\tau^{\text{late}}_{\ell t p}+\frac{\delta}{2}\tau^{\text{early}}_{\ell t p} + \sigma \tau^{\text{dropoff}}_{\ell t}}{\tau^{\text{dir}}_p}+ \lambda \tau^{\text{walk}}_{ip} + \mu( - \tau^{\text{req}}_{p} - \tau^{\text{walk}}_{ip}) -M\right) \cdot w_{\ell t s p i}   \\
    &\quad\quad + \sum_{p\in\calP} D_{ps} \left(\frac{-\sigma t_{\ell t s p}^{\text{pickup}}}{\tau^{\text{dir}}_p} + \mu t_{\ell t s p}^{\text{pickup}} \right) \\
    &= \sum_{p \in \calP} D_{ps} \Bigg( \lambda \sum_{i \in \calN^{\text{pickup}}_{\ell t p}} \tau^{\text{walk}}_{ip} w_{\ell t s p i} + \mu \left(t^{\text{pickup}}_{\ell t s p} -  \sum_{i \in \calN^{\text{pickup}}_{\ell t p}} (\tau^{\text{req}}_p + \tau^{\text{walk}}_{ip}) w_{\ell t s p i} \right)  \\
    & \quad + \frac{\sigma}{\tau_p^{\text{dir}}} \left( \tau_{\ell t}^{\text{dropoff}} \sum_{i \in \calN^{\text{pickup}}_{\ell t p}} w_{\ell t s p i} - t^{\text{pickup}}_{\ell t s p} \right) + \left(\delta \frac{\tau^{\text{late}}_{\ell tp}}{\tau_p^{dir}} + \frac{\delta}{2} \frac{\tau^{early}_{\ell tp}}{\tau_p^{\text{dir}}}\right)\sum_{i \in \calN^{\text{pickup}}_{\ell t p}}  w_{\ell t s p i} \\
    &\quad -M\sum_{i \in \calN^{\text{pickup}}_{\ell t p}}  w_{\ell t s p i}\Bigg)
\end{align*}

Thus $(\bv, \bw, \by, \bt^{\text{stop}}, \bt^{\text{pickup}})$ is a feasible solution to the compact formulation and achieves the same objective value as $\by^*$ in the subpath formulation. This completes the proof.
\hfill\Halmos

\subsection{Segment-based benchmark for second-stage problem}\label{app:segment}

\begin{table}[h!]
    \centering
    \small
\begin{tabular}{llp{11.5cm}}
\toprule[1pt]
    \bf Component & \bf Type & \bf Description \\ \hline
    $\overline\calE_{\ell st}$ & Set & Load-augmented road segments $e$ associated with $road(e) \in \calE$\\
    $\calT^S$ & Set & Set of time periods during the planning horizon \\ 
    $\calP_{e}$ & Set & Passengers picked up on segment $e\in \overline\calE_{\ell st}$ \\
    $(\overline\calV_{\ell st}, \overline\calA_{\ell st})$ & Graph &  Time-load-expanded road network of trip $(\ell, t) \in \calL \times \calT_{\ell}$ in scenario $s \in \calS$ \\
    $\overline\calA_{e}$ & Set & Arcs in $\overline\calA_{\ell st}$ corresponding to segment $e \in \overline\calE_{\ell st}$ for $(\ell, t) \in \calL \times \calT_{\ell},$ $s \in \calS$ \\
    $\overline\calA_{\ell st}^{idle}$ & Set  & Arcs in $\overline\calA_{\ell st}$ representing an idling vehicle \\
    $\overline\calA_{\ell st}^v$ & Set  & Arcs in $\overline\calA_{\ell st}$ connecting the line's destination to the dummy sink node \\\hline
    $\tau_{ep}^{\text{walk}}$ &Parameter & Walk time of passenger $p \in \calP_{e}$ via segment $e \in \overline\calE_{\ell st}$, $(\ell, t) \in \calL \times \calT_{\ell}$, $s \in \calS$ \\
    $\tau_{ep}^{\text{wait}}$& Parameter & Wait time of passenger $p \in \calP_{e}$ via segment $e \in \overline\calE_{\ell st}$, $(\ell, t) \in \calL \times \calT_{\ell}$, $s \in \calS$ \\
    $\tau^{\text{travel}}_{ep}$&Parameter& In-vehicle time of passenger $p\in \calP_r$ via segment $e \in \overline\calE_{\ell st}$, $(\ell, t) \in \calL \times \calT_{\ell}$, $s \in \calS$ \\
    $\overline{g}_a$ & Parameter & Cost of arc $a \in \overline\calA_{\ell st}$ on trip $(\ell, t) \in \calL \times \calT_{\ell}$ in scenario $s \in \calS$  \\
    \bottomrule[1pt]
\end{tabular}
\caption{Additional inputs of the segment-based formulation.}\label{T:notation_segment}
\end{table}

Throughout the section, we fix first-stage decisions $\bx$ and $\bz$, as well as scenario $s \in \calS$. The time horizon is discretized into $T_S+1$ intervals in the set $\calT^S = \{0,1,\cdots, T_S \}$, from the departure of the first trip ($t=0$) to the arrival of the last trip ($t=T_S$).

To capture time and capacity constraints without relying on big-$M$ constraints---therefore retaining a tight second-stage formulation---we build a time-load-expanded network $(\overline\calV_{\ell st}, \overline\calA_{\ell st})$. A dummy sink node $v_{\ell st}$ represents the end of a trip. Each other node $n\in\overline\calV_{\ell st}$ is associated with a tuple $(k_n, c_n, t_n)$, so that node $n$ represents a vehicle's arrival to station $k_n\in\calN^S$ at time $t_n\in \calT^S$ with $c_n\in\calC$ passengers. The source node is denoted by $u_{\ell st}:=(\calI_{\ell}^{(1)},0, T_{\ell t}(\calI_{\ell}^{(1)}))$. We decompose the arc set $\overline\calA_{\ell st} \subset \overline\calV_{\ell st}\times\overline\calV_{\ell st}$ into traveling arcs, idling arcs, and terminating arcs, by writing $\overline\calA_{\ell st} = \bigcup_{e \in \overline\calE_{\ell st}} \overline\calA_e \cup \overline\calA_{\ell st}^{idle} \cup \overline\calA^v_{\ell st}$.

To characterize traveling arcs, we denote by $\overline\calE_{\ell st}$ the set of possible roadways and passenger pickups. Specifically, each segment $e\in\overline\calE_{\ell st}$ is associated with a raodway $road(e) \in \calE$ and a set of passengers $\calP_e$ who are picked up. We define traveling arcs by duplicating $e\in\overline\calE_{\ell st}$ for all load pairs that correspond to the passenger pickups, and all time pairs that correspond to the travel time:
\begin{align}
    \overline\calA_e = \bigg\{(n,m) \in \overline\calV_{lst} \times \overline\calV_{lst} \, : \, &(k_n, k_m) = road(e),\nonumber\\&c_m - c_n = \sum_{p \in \calP_e} D_{ps},\nonumber\\&t_m-t_n=tt(road(e)) \bigg\} \qquad \forall e \in \overline\calE_{\ell st}
\end{align}
Next, each idling arc in $\overline\calA_{\ell st}^{idle}$ connects nodes corresponding to two consecutive time intervals at the same station:
\begin{equation}\label{eq:SegmentAH}
\overline\calA_{\ell st}^{idle} = \{(n,m) \in \overline\calV_{\ell st} \times \overline\calV_{\ell st} \, : \, k_n = k_m,\ c_n = c_m ,\ t_m - t_n = 1 \}.
\end{equation}
Finally, each terminating arc in $\overline\calA_{\ell st}^v$ connects the line's destination to the dummy sink node:
\begin{equation}\label{eq:segmentAT}
\overline\calA_{\ell st}^v = \{(n,m) \in \overline\calV_{\ell st} \times \overline\calV_{\ell st} \, : \, k_n = \calI^{\text{end}}, m = v^S_{\ell st} \}.
\end{equation}
Again, we can prune the time-load-expanded network by excluding disconnected nodes and all incident arcs. We define a segment-based cost $\overline{g}_a$ for each $a \in \overline\calA_{\ell st}$ analogously to Equation~\eqref{eq:arccost} to capture passenger walking times, waiting times, and relative arrival delays:
 \begin{equation}
     \overline{g}_a = \begin{cases}
         \sum_{p \in \calP_{e}}D_{ps} \left(\lambda \tau_{ep}^{\text{walk}} + \mu \tau_{ep}^{\text{wait}} + \sigma \frac{\tau_{ep}^{\text{travel}}}{\tau_p^{\text{dir}}} + \delta \frac{\tau_{\ell tp}^{\text{late}}}{\tau_p^{\text{dir}}} + \frac{\delta}{2} \frac{\tau_{\ell tp}^{\text{early}}}{\tau_p^{\text{dir}}} - M\right) &\text{if }e \in \overline\calE_{\ell st}, a \in\overline\calA_e\\
         0 &\text{if }a \in \overline\calA_{\ell st}^{idle}\cup \overline\calA_{\ell st}^v.
     \end{cases}
 \end{equation}

We define decision variables to select arcs in the time-load-expanded segment network:
\begin{align}
\xi_a = \begin{cases}
    1 &\text{if arc $a$ is selected, for $(\ell, t) \in \calL \times \calT_{\ell}, s \in \calS, a \in \overline\calA_{\ell st}$,} \\
    0 &\text{otherwise.}
\end{cases}
\end{align}

Recall that $\Gamma_{\ell} \subset \calI_\ell \times \calI_\ell$ denotes the set of checkpoint pairs with up to $K$ skipped checkpoints:
$$\Gamma_\ell = \left\{(\calI_\ell^{(i)},\calI_\ell^{(j)}) \in \calI_\ell\times\calI_\ell \, : \, 1 \leq i < j \leq I_\ell, j - i \leq K + 1 \right\}, \qquad \forall \ell \in \calL $$
We define additional decision variables to select the set of checkpoint pairs that are visited:
\begin{equation*}
    \beta_{uv} = \begin{cases}
     1 &\text{if checkpoints }(u,v) \in \Gamma_\ell \text{ are visited in sequence, and intermediate checkpoints are not visited,}    \\
     0 &\text{otherwise.}
    \end{cases}
\end{equation*}
Recall that $\calN^S_{uv}$ denotes the set of stations that can be visited between checkpoints $u$ and $v$, and $\calT^{uv}_{\ell t}$ denotes the valid arrival times. We link the $\beta_{uv}$ decisions with the $\xi_a$ decisions, so that the vehicle route abides by the deviation limits imposed by the reference schedule. Altogether, the segment-based formulation exhibits a double flow structure---flow from checkpoint to checkpoint along the reference line, and flow from station to station between checkpoints---with linking constraints to ensure the consistency of these two sets of decisions.

The second-stage segment-based formulation is given as follows for scenario $s\in\calS$. 
\begin{align}
    \min \quad & \sum_{(\ell, t) \in \calL \times \calT_{\ell}} \sum_{a \in \overline\calA_{\ell st}} \overline{g}_a \xi_a \label{edge:obj}\\
    \text{s.t. } &\sum_{j:(i,j)  \in \overline\calA_{\ell st} } \xi_{(i,j)} - \sum_{j:(j,i) \in \overline\calA_{\ell st}} \xi_{(j,i)}= \begin{cases}
    x_{lt} &\text{if } i = u_{\ell st}, \\
    -x_{lt} &\text{if } i = v_{\ell st}, \\
    0 &\text{otherwise,}
    \end{cases}\ \forall (\ell, t) \in \calL \times \calT_{\ell}, \forall i \in \overline\calV_{\ell st}\label{edge:flow} \\
    &\sum_{e\in\overline\calE_{\ell st}}\sum_{a \in \overline\calA_e \, : \, p \in \calP_{e}} \xi_a \leq z_{plt}  \qquad \forall p \in \calP,\ \forall(\ell, t) \in \calM_p \label{edge:deactivate}\\
     &\sum_{v \, : (u,v) \in \Gamma_{\ell}} \beta_{uv} - \sum_{v \, : (v,u) \in \Gamma_{\ell}} \beta_{vu} = \begin{cases}
        x_{\ell t} &\text{if }u = \calI_{\ell}^{(1)} \\
        -x_{\ell t} &\text{if }u = \calI^{\text{end}} \\
        0 &\text{otherwise}
    \end{cases}, \quad  \forall (\ell, t) \in \calL \times \calT_{\ell}, \forall u \in \calI_{\ell}\label{eq:flowBalanceCheckpoints}\\
    &   \sum_{\substack{(i, j) \in \overline\calA_{\ell st} \, : \\ k_j = v,\, t_j = T_{\ell t}(v)}} \xi_{(i,j)} \geq \sum_{w \in \calI_\ell \, : \, (w, v) \in \Gamma_\ell}\beta_{wv} \qquad \forall v \in \calI_{\ell} \backslash \calI^{(1)}_\ell\label{eq:visitCheckpoint}\\
    &    \xi_{(n,m)} \leq \sum_{\substack{(u, v) \in \Gamma_\ell \, : \\ k_n, k_m \in \calN^S_{uv}, \\ t_n, t_m \in \calT_{\ell t}^{uv}}} \beta_{uv} \qquad \forall (\ell, t) \in \calL \times \calT_\ell, \forall (n, m) \in \overline\calA_{\ell st}\label{eq:checkpointArcs} \\
    & \xi_a \in \{0,1\} \qquad \forall (\ell, t) \in \calL \times \calT_{\ell}, a \in \overline\calA_{\ell st}\label{edge:domain} \\
    &\beta_{uv} \in \{0,1\} \qquad \forall (\ell, t) \in \calL \times \calT_{\ell}, \forall (u,v) \in \Gamma_\ell\label{edge:domain2}
\end{align}
Equations~\eqref{edge:obj}--\eqref{edge:deactivate} are analogous to Equations~\eqref{obj}--\eqref{eq:2Sy2z}. Constraint \eqref{eq:flowBalanceCheckpoints} ensures that the vehicle does not skip more than $K$ checkpoints in a row by selecting checkpoint pairs that form a valid path along the reference line. Equations~\eqref{eq:visitCheckpoint} and~\eqref{eq:checkpointArcs} serve as the linking constraints, ensuring that selected checkpoints are visited at the time specified by the reference schedule, and that the vehicle visits any intermediate locations with the correct chronology. In other words, we can only select a segment if (i) its endpoints correspond to stations in $\calN^S_{uv}$ between selected checkpoints, and (ii) its visit times fall within the reference schedule window defined by $T_{\ell t}(u)$ and $T_{\ell t}(v).$ Constraints~\eqref{edge:domain}--\eqref{edge:domain2} apply the binary requirements to the decision variables.

\subsection{Path-based formulation for second-stage problem}\label{app:path}

\begin{table}[h!]
    \centering
    \small
\begin{tabular}{llp{11.5cm}}
\toprule[1pt]
    \bf Component & \bf Type & \bf Description \\ \hline
    $\calQ_{\ell st}$ & Set & Valid paths for reference trip $(\ell, t) \in \calL \times \calT_{\ell}$ and scenario $s \in \calS$ \\
    $\calP_q$ & Set & Passenger pickup set corresponding to each path $q \in \calQ_{\ell st}$ \\\hline 
    $\tau_{qp}^{\text{walk}}$ &Parameter & Walk time of passenger $p \in \calP_r$ via path $q \in \calQ_{\ell st}$, for $(\ell, t) \in \calL \times \calT_{\ell}$, $s \in \calS$ \\
    $\tau_{qp}^{\text{wait}}$& Parameter & Wait time of passenger $p \in \calP_r$ via path $q \in \calQ_{\ell st}$, for $(\ell, t) \in \calL \times \calT_{\ell}$, $s \in \calS$ \\
    $\tau_{qp}^{\text{travel}}$& Parameter & In-vehicle time of passenger $p \in \calP_r$ via path $q \in \calQ_{\ell st}$, for $(\ell, t) \in \calL \times \calT_{\ell}$, $s \in \calS$ \\
    $g_q^Q$ & Parameter & Cost of path $q \in \calQ_{\ell st}$ on trip $(\ell, t) \in \calL \times \calT_{\ell}$ in scenario $s \in \calS$  \\
    \bottomrule[1pt]
\end{tabular}
\caption{Additional inputs of the path-based formulation.}\label{T:notation_path}
\end{table}

Throughout the section, we fix first-stage decisions $\bx$ and $\bz$, as well as scenario $s \in \calS$.

Let $\calQ_{\ell st}$ denote the set of all valid paths to reference trip $(\ell, t) \in \calL \times \calT_{\ell}$ and scenario $s \in \calS$. Each path $q \in \calQ_{\ell st}$ corresponds to a sequence of road segments that starts at the beginning of the line, end at its destination, satisfies flow balance in between, skips at most $K$ checkpoints in a row, does not pick up more than $C_\ell$ passengers, and satisfies the reference schedule at the checkpoints. For each $q\in\calQ_{\ell st},$ we store the passenger pickups in $\calP_q \subset \calP.$ By definition, $\sum_{p \in \calP_q} D_{ps} \leq C_\ell.$ The cost $g_q^Q$ of each path is defined analogously to Equation~\eqref{eq:arccost} to capture passenger level of service:
\begin{equation}
    g_q^{Q} = \sum_{p \in \calP_q} D_{ps} \left(\lambda \tau_{qp}^{\text{walk}} + \mu \tau_{qp}^{\text{wait}} + \sigma \frac{\tau_{qp}^{\text{travel}}}{\tau_p^{\text{dir}}} + \delta \frac{\tau_{\ell tp}^{\text{late}}}{\tau_p^{\text{dir}}} + \frac{\delta}{2} \frac{\tau_{\ell tp}^{\text{early}}}{\tau_p^{\text{dir}}} - M\right), \quad \forall (\ell, t) \in \calL \times \calT_{\ell}, q \in \calQ_{\ell st}.
\end{equation}

We define the following decision variables:
\begin{align}
\zeta_q = \begin{cases}
    1 &\text{if path $q$ is selected, for $(\ell, t) \in \calL \times \calT_{\ell}, s \in \calS, q \in \calQ_{\ell st}$,} \\
    0 &\text{otherwise.}
\end{cases}
\end{align}

The path-based formulation is given as follows for scenario $s \in \calS$. 
\begin{align}
    \min \quad & \sum_{(\ell, t) \in \calL \times \calT_{\ell}} \sum_{q \in \calQ_{\ell st}} g_q^{Q} \zeta_q \label{path:obj}\\
    \text{s.t. } &\sum_{q \in \calQ_{\ell st}} \zeta_q = x_{lt} \qquad \forall (\ell, t) \in \calL \times \calT_{\ell} \label{path:flow}\\
    &\sum_{q \in \calQ_{\ell st} \, : \, p \in \calP_q } \zeta_q \leq z_{plt} \qquad\forall p \in \calP, \ \forall(\ell, t) \in \calM_p \label{path:deactivate}\\
    &\zeta_q \in \{0,1\} \qquad\forall (\ell, t) \in \calL \times \calT_{\ell}, q \in \calQ_{\ell st}\label{path:domain}
\end{align}
Equations~\eqref{path:obj} is analogous to Equation~\eqref{obj}. Constraints~\eqref{path:flow} ensure that exactly one path is selected for each selected reference trip. Constraints~\eqref{path:deactivate} ensure that selected paths only serve passengers that have been assigned to that trip, analogously to Equation~\eqref{eq:2Sy2z}.

\subsection{Proof of Proposition \ref{subpath_form}}\label{A:subpath_structure}

Throughout this proof, we fix the first-stage decisions $\mathbf{x}$, $\mathbf{z}$. We consider a fixed reference trip $(\ell, t) \in \calL \times \calT_\ell$ as well as a fixed scenario $s\in \calS$.

\subsubsection*{Equivalence of the path-based and subpath-based formulations.}\

{\it Constructing a load-expanded subpath solution from a path solution.} Let us consider a feasible solution $\widehat{\mathbf{\zeta}}$ to the path-based formulation (Equations \eqref{path:obj}--\eqref{path:domain}) and build a feasible solution to the subpath-based formulation with the same objective value.

Assume that $x_{\ell t} = 1,$ and let $q \in \calQ_{\ell st}$ be the selected path with $\widehat\zeta_q = 1$ (which exists by Equation~\eqref{path:flow}). By definition, the path corresponds to a sequence of road segments that starts at the beginning of line $\ell$, ends at its destination, picks up at most $C_\ell$ passengers, visits checkpoints without skipping more than $K$ in a row, and arrives at each checkpoint at the scheduled times. With a slight abuse of notation, let $\calI_\ell^q := \{\nu_1, \cdots, \nu_{Q}\} \subseteq \calI_\ell$ identify the ordered set of $Q$ checkpoints visited by path $q$. Similarly, we decompose path $q$ into an ordered sequence of $Q-1$ subpaths $\calR_q := \{r_1, \cdots, r_{Q-1}\}$. The subpaths in $\calR_q$ partition the served passengers $\calP_q$ on path $q$, so that $\calP_q = \bigcup_{r \in \calR_q} \calP_r$. Each subpath $r_i \in \calR_q$ induces a unique arc $a_i := (n,m) \in \calA_{\ell st}$ in the load-expanded network, such that
\begin{enumerate*}[label=(\roman*)]
    \item the arc corresponds to the subpath: $r(a_i) = r_i$;
    \item the loads are consistent with pickups: $c_n = 0$ if $i = 1,$ and $c_n = \sum_{j=1}^{i-1} |\calP_{r_j}|$ otherwise, and $c_m = \sum_{j=1}^i |\calP_{r_j}|.$
\end{enumerate*}
Let us collect these load-expanded subpath arcs into the set $\calA_q := \{a_1, \cdots, a_{Q-1}\} \subset \calA_{\ell st}.$ 

We can construct a feasible solution to the load-expanded subpath formulation.
\begin{equation*}
    \widehat{y}_a = \begin{cases}
        1 &\text{if }a \in \bigcup_{q \in \calQ_{\ell st} \, : \, \widehat\zeta_q = 1}\calA_q, \\ 
        0 &\text{otherwise,}
    \end{cases} \qquad \forall (\ell, t) \in \calL \times \calT_\ell, \forall a \in \calA_{\ell st}.
\end{equation*}
This solution satisfies the flow balance constraints in Equation~\eqref{eq:2Sflow}:
\begin{itemize}
    \item[--] If $x_{\ell t} = 0$, no path in $\calQ_{\ell st}$ is selected, so no arc in $\calA_{\ell st}$ is selected either. Therefore,
\begin{align*}
     \sum_{m:(n,m) \in \calA_{\ell st}} 
    \widehat{y}_{(n,m)} - \sum_{m:(m,n) \in \calA_{\ell st}} 
    \widehat{y}_{(m,n)} = 0 = \begin{cases}
         x_{\ell t}&\text{if }n = u_{\ell st}, \\ 
         -x_{\ell t} &\text{if }n = v_{\ell st}, \\ 
        0 &\text{otherwise.}
    \end{cases}  
\end{align*}
\item[--] If $x_{\ell t} = 1$, we have, for each node $n \in \calV_{\ell st}$:
\begin{equation*}
\sum_{m : (n, m) \in \calA_{\ell st}} \widehat{y}_{(n,m)} - \sum_{m :  (m, n) \in \calA_{\ell st}} \widehat{y}_{(m, n)} = \begin{cases}
    \widehat{y}_{a_1} - 0 = 1 &\text{if }n = u_{\ell st}, \\
    0 - \widehat{y}_{a_{Q-1}} = -1 & \text{if }n = v_{\ell st}, \\ 
    \widehat{y}_{a_{i}} - \widehat{y}_{a_{i-1}} = 1 - 1 = 0 &\text{if } k_n = \nu_i \in \calI_\ell^q \setminus \{\nu_1, \nu_{Q}\}, \\
    &c_n = \sum_{j=1}^{i-1} |\calP_{r_j}|, \\ 
    0-0=0 &\text{otherwise.}
\end{cases}
\end{equation*}

\end{itemize}

The solution also satisfies the passenger linking constraints in Equations~\eqref{eq:2Sy2z}.
\begin{itemize}
    \item[--]  Consider a passenger request $p \in \calP \setminus \bigcup_{q \in \calQ_{\ell st} \, : \, \widehat\zeta_q = 1} \calP_q$. The set of pickups on the selected paths $\bigcup_{q \in \calQ_{\ell st} \, : \, \widehat\zeta_q = 1} \calP_q$ induces the set of pickups on the selected subpaths, so that $p \in \calP \setminus \bigcup_{q \in \calQ_{\ell st} \, : \, \widehat\zeta_q = 1} \bigcup_{r \in \calR_q} \calP_r.$ Thus, $\widehat{y}_a = 0$ for each arc $a \in \calA_{\ell st}$ with $p \in \calP_{r(a)},$ and
    $$\sum_{a \in \calA_{\ell st} \, : \, p \in \calP_{r(a)}} \widehat{y}_a = 0 \leq z_{\ell pst}. $$
    \item[--] Consider passenger request $p \in \calP$ served by some path $q' \in \calQ_{\ell st}$, so that $\widehat\zeta_{q'} = 1$ and $p \in \calP_{q'}.$ Each pickup set $\calP_{q'}$ has been partitioned into pickup subsets at the subpath level, so there exists $r' \in \calR_{q'}$ with $p \in \calP_{r'}$. This subpath has been mapped to a unique arc $a_{r'} \in \calA_{q'}$, so that
    $$\sum_{a \in \calA_{\ell st}  :  p \in \calP_{r(a)}} \widehat{y}_a = \widehat{y}_{a_{r'}} = 1 = \widehat\zeta_{q'} = \sum_{q \in \calQ_{\ell st} : p \in \calP_{q'}} \widehat\zeta_q \leq z_{\ell pst} $$
    where the first three equalities come from the construction of paths and subpaths, the fourth equality stems from the fact that passenger $p$ can be picked up by one subpath, and the final one stems from Equation~\eqref{path:deactivate}.
\end{itemize}

Therefore, the solution $\widehat\by$ is feasible in the subpath-based formulation. We now show that it achieves the same objective value as $\widehat{\mathbf{\zeta}}$:
\begin{align}
    \sum_{a \in \calA_{\ell st}} g_a \widehat{y}_a &= \sum_{a \in \calA_{\ell st} : \widehat{y}_a = 1} g_a  \nonumber\\
    &= \sum_{q \in \calQ_{\ell st} : \widehat\zeta_q = 1} \sum_{a \in \calA_q} g_a  \nonumber\\
    &= \sum_{q \in \calQ_{\ell st} : \widehat\zeta_q = 1}  \sum_{r_i \in \calR_q} g_{a_i}  \nonumber\\
    &= \sum_{q \in \calQ_{\ell st} : \widehat\zeta_q = 1}\sum_{r \in \calR_q} \sum_{p \in \calP_r}D_{ps} \left(\lambda \tau_{rp}^{\text{walk}} + \mu \tau_{rp}^{\text{wait}} + \sigma \frac{\tau_{rp}^{\text{travel}}}{\tau_p^{\text{dir}}} + \delta \frac{\tau_{\ell tp}^{\text{late}}}{\tau_p^{\text{dir}}} + \frac{\delta}{2} \frac{\tau_{\ell tp}^{\text{early}}}{\tau_p^{\text{dir}}} - M\right)  \nonumber\\
    &= \sum_{q \in \calQ_{\ell st} : \widehat\zeta_q = 1}\sum_{p \in \calP_q} D_{ps} \left(\lambda \tau_{qp}^{\text{walk}} + \mu \tau_{qp}^{\text{wait}} + \sigma \frac{\tau_{qp}^{\text{travel}}}{\tau_p^{\text{dir}}} + \delta \frac{\tau_{\ell tp}^{\text{late}}}{\tau_p^{\text{dir}}} + \frac{\delta}{2} \frac{\tau_{\ell tp}^{\text{early}}}{\tau_p^{\text{dir}}} - M\right)  \nonumber\\
    &= \sum_{q \in \calQ_{\ell st} : \widehat\zeta_q = 1} g_q^Q \nonumber\\
    &= \sum_{q \in \calQ_{\ell st}}g_q^Q \widehat\zeta_q \label{objComputation}
\end{align}
The first two equalities come from the construction of paths and subpaths; the third equality leverages the uniqueness of the load-expanded subpath arc induced by the subpath sequence; the fourth equality is due to the definition of a load-expanded subpath arc cost; the fifth is due to the partition of $\calP_q = \bigcup_{r \in \calR_q} \calP_r$; and the last two equalities stem from the definition of path costs $g_q^Q$.

In conclusion, any path solution can be mapped into a feasible subpath solution with the same objective value. Therefore, the subpath-based formulation achieves an objective that is at most equal to an optimum of the path-based formulation.

{\it Constructing a path solution from a load-expanded subpath solution.} Let us consider a feasible solution $\widehat\by$ to the subpath-based formulation (Equations~\eqref{eq:2Sflow}--\eqref{domain}), and build a feasible solution $\widehat{\mathbf{\zeta}}$ to the path-based formulation (Equations~\eqref{path:obj}--\eqref{path:domain}) with the same objective value. 

Assume that $x_{\ell t} = 1$. We leverage Equation~\eqref{eq:2Sflow} to construct a path from $u_{\ell st}$ to $v_{\ell st}$ in the load-expanded subpath network $(\calV_{\ell st}, \calA_{\ell st})$. Beginning from the source, we select the unique arc $a_1 \in \calA_{\ell st}$ incident with $u_{\ell st}$ for which $\widehat{y}_{a_1} = 1,$ proceeding sequentially along the directed network until $v_{\ell st}$ is reached and $Q-1$ arcs are retrieved. A unique outgoing arc is guaranteed at every intermediate node by Equation~\eqref{eq:2Sflow}.

Each arc $a_i \in \calA_q := \{a_1, \cdots, a_{Q-1}\}$ corresponds to a subpath $r_i := r(a_i) \in \calR_{\ell st}$ and a passenger pickup set $\calP_{r_i}.$ The sequence of subpaths $\calR_q := \{r_1, \cdots, r_{Q-1}\}$ defines a path $q$ from $u_{\ell st}$ to $v_{\ell st}$ (by Equation~\eqref{eq:2Sflow}), skipping at most $K$ checkpoints in a row (by definition of the subpaths $r_i \in \calR_{\ell st}$), cohering with the scheduled arrival times associated with trip $(\ell, t)$ (again by definition of $r_i$), obeying the vehicle's capacity (by definition of the node set in the load-expanded network $\calV_{\ell st}$), and picking up the passengers in $\calP_q := \bigcup_{i=1}^{Q-1} \calP_{r_i}$ (which is unique due to Equations~\eqref{eq:2Sy2z}). Thus, $\widehat\by$ defines a unique and valid path in $\calQ_{\ell st}$ for reference trip $(\ell,t)$ if $x_{\ell t} = 1$.

Let us collect all such paths in the set $\calQ(\widehat\by).$ We construct solution $\widehat{\mathbf{\zeta}}$ from $\widehat\by$:
    \begin{equation*}
        \widehat\zeta_q = \begin{cases}
            1 &\text{if }q \in \calQ(\widehat\by), \\
            0 &\text{otherwise,}
        \end{cases} \qquad \forall (\ell, t) \in \calL \times \calT_\ell, \forall q \in \calQ_{\ell st}.
    \end{equation*}
    
By construction, the solution satisfies Equations \eqref{path:flow}. If $x_{\ell t} = 1,$ we constructed a single path based on the subpath solution. If $x_{\ell t} = 0,$ there was no path to construct, as no arcs were selected from $u_{\ell st}$ to $v_{\ell st}$ by Equation~\eqref{eq:2Sflow}. Therefore:
    \begin{align*}
        \sum_{q \in \calQ_{\ell st}}\widehat\zeta_q = \sum_{q \in \calQ_{\ell st}} \mathbbm{1}(q \in \calQ(\widehat\by)) = x_{\ell t}, \quad \forall (\ell, t) \in \calL \times \calT_\ell.
    \end{align*}

The solution also satisfies Equation~\eqref{path:deactivate}. For each reference trip $(\ell, t) \in \calL \times \calT_\ell$ and some passenger $p \in \calP,$ we obtain, from the construction of the path solution and Equation~\eqref{eq:2Sy2z}:
    \begin{align*}
        \sum_{q \in \calQ_{\ell st} : p \in \calP_q} \widehat\zeta_q &= \sum_{q \in \calQ(\widehat\by)} \mathbbm{1}(p \in \calP_q) \\
        &= \sum_{q \in \calQ(\widehat\by)}  \sum_{r \in \calR_q} \mathbbm{1}\left(p \in \calP_{r}\right) \\
        &= \sum_{q \in \calQ(\widehat\by)}  \sum_{a \in \calA_q} \mathbbm{1}\left(p \in \calP_{r(a)}\right) \\ 
        &= \sum_{a \in \calA_{\ell st} \, : \, \widehat{y}_a = 1} \mathbbm{1}(p \in \calP_{r(a)}) \\
        &= \sum_{a \in \calA_{\ell st} \, : \, p \in \calP_{r(a)}} \widehat{y}_a \\
        &\leq z_{\ell pst}
    \end{align*}
    
Finally, the solutions $\widehat{\mathbf{\zeta}}$ and $\widehat\by$ achieve the same objective values, which can be shown similarly to Equations~\eqref{objComputation}. Therefore, any subpath solution can be mapped into a feasible path solution with the same objective value, and the path-based formulation achieves an objective that is at most equal to the optimum of the subpath-based formulation. This concludes the proof of equivalence of the path-based and subpath-based formulations.

\textit{Equivalence of path-based and subpath-based relaxations.} The arguments employed in this proof do not require the integrality of the path solution $\widehat{\bzeta}$ and of the subpath solution $\widehat{\by}$. By following the same steps as above, we can map any non-integral path solution $\widehat{\bzeta}$ into a feasible subpath solution with the same objective value, as follows:
\begin{equation*}
    \widehat{y}_a = \sum_{q \in \calQ_{\ell st} : a \in \calA_q} \widehat\zeta_q \qquad \forall (\ell, t) \in \calL \times \calT_\ell, \forall a \in \calA_{\ell st}.
\end{equation*}
Similarly, we can map any non-integral subpath solution $\widehat{\by}$ into a feasible path solution with the same objective value. Alternatively, we can observe that the path-based formulation is a Dantzig-Wolfe reformulation of the subpath-based formulation where Equations~\eqref{eq:2Sflow} are convexified into Equations~\eqref{path:flow}. Since Equations~\eqref{eq:2Sflow} already form an integral polyhedron, both formulations contain the same convex hull. This proves that the path-based and subpath-based formulations define the same linear relaxations.

\subsubsection*{Equivalence of the segment-based and subpath-based formulations.}\

{\it Constructing a time-load-expanded segment solution from a load-expanded subpath solution.} Let us consider a feasible solution $\widehat{\mathbf{y}}$ to the subpath-based formulation (Equations \eqref{obj}--\eqref{domain}) and build a feasible solution to the segment-based formulation with the same objective value.

Assume that $x_{\ell t} = 1,$ and let $a \in \calA_{\ell st}$ be a selected subpath-based arc with $\widehat{y}_a = 1$ (which exists by Equation~\eqref{eq:2Sflow}). By definition, the subpath-based arc corresponds to the load expansion of a subpath that traverses a sequence of road segments starting at checkpoint $u \in \calI_\ell$ at time $T_{\ell t}(u)$, ending at checkpoint $v \in \calI_\ell$ at time $T_{\ell t}(v)$, skipping up to $K$ checkpoints in-between (i.e., $(u,v) \in \Gamma_\ell$), carrying $c_{start(a)}$ passengers in $u$, and carrying $c_{end(a)}$ passengers in $v$. Let us store the stations visited by subpath $r$ in an ordered set $\calN^S_r := \{\nu_1, \cdots, \nu_{N}\} \subseteq \calN^S, $ where $\nu_1=u$, $\nu_N=v$, and $\nu_2,\cdots,\nu_{N-1}$ denote intermediate stations. Similarly, we decompose subpath $r$ into a sequence of $N-1$ segments $\calE_r := \{e_1, \cdots, e_{N-1}\}$, where segment $e_i$ connects stations $\nu_i$ and $\nu_{i+1}$ with travel time $tt_{e_i}$ (potentially with idling time). The segments in $\calE_r$ partition the passengers in $\calP_r$: $\calP_r = \bigcup_{e \in \calE_r} \calP_e.$

To obtain the corresponding segment solution, we need to specify an appropriate time discretization. Due to the adherence to the reference schedule, the discretization in the segment-based formulation does not introduce errors as long as all viable subpaths are feasible in that formulation. We show that there exists a discrete time unit for which this is the case, in the following lemma.

\begin{lemma}\label{lem:disc}
    Assume that the elapsed time between the scheduled arrival times at the checkpoints along the reference line are strictly larger than the travel times of the corresponding subpaths. Then, there exists a discrete time unit such that, in the corresponding time-expanded network, all feasible subpaths have an estimated travel time that is less than the elapsed time between the corresponding checkpoints' scheduled arrival times.
\end{lemma}

Let $\Delta_{uv} := T_{\ell t}(v) - T_{\ell t}(u)$ denote the travel time between checkpoints $u$ and $v$, determined by the scheduled arrival times at both checkpoints, and let us denote the travel time of subpath $r$ by
$\Delta_r:= \sum_{e\in\calE_r}tt_e$. Due to the maximum deviation from the reference line, the number of passenger pickups, and the upper bound on passengers' walking distance, the set of potential subpaths $\calR_{\ell st}^{uv}$ between checkpoints $u$ and $v$ is finite. For convenience, let us denote this subset by
$$\calR_{\ell st}^{uv} := \left\{r \in \calR_{\ell st} : u_r = u, v_r = v\right\}.$$ 

By assumption, all subpaths $r\in\calR_{\ell st}$ satisfy $\Delta_r \leq T_{\ell t}(v_r) - T_{\ell t}(u_r)$, so that $\Delta_r<\Delta_{uv}$ for each $r \in \calR_{\ell st}^{uv}$.
We define the discrete time unit between checkpoints $u$ and $v$ as:
\begin{equation}\label{eq:rhouv}
    \rho_{uv}=\min_{r \in \calR_{\ell st}^{uv}}\frac{\Delta_{uv} - \Delta_r}{|\calE_r|}>0.
\end{equation}

Without loss of generality, we assume that $\rho_{uv}$ is rational; otherwise, we can define it as the largest rational number bounded from above by the minimum given in Equation~\eqref{eq:rhouv}. We define the universal discrete time unit as
\begin{equation}\label{eq:gcdTimeUnit}
\rho = \text{GCD}\left(\{\rho_{uv} : (u, v) \in \Gamma_\ell \} \right),
\end{equation}
where GCD denotes the greatest common divisor. By construction, for each $(u,v) \in \Gamma_\ell$, there exists $R_{uv}\in\Z_+$ such that $\rho_{uv}=R_{uv}\rho\geq\rho$.

In the segment-based formulation, travel times are rounded up to the nearest discrete time step on each segment. The estimated travel time on each segment $e\in\calE_{r}$, denoted by $\overline{\Delta}_e$, is therefore
$$\overline{\Delta}_e=\left\lceil\frac{tt_e}{\rho}\right\rceil\cdot\rho.$$

The travel time estimate of subpath $r \in \calR_{\ell st}^{uv}$ in the segment-based formulation, denoted by $\overline{\Delta}_r$, is then given by:
$$\overline{\Delta}_r=\sum_{e\in\calE_r}\overline{\Delta}_e=\sum_{e\in\calE_r}\left\lceil\frac{tt_e}{\rho}\right\rceil\rho.$$

We make use of the following property:
\begin{equation*}
    \left\lceil\frac{tt_e}{\rho}\right\rceil\cdot\rho\leq\left\lceil\frac{tt_e}{\rho_{uv}}\right\rceil\cdot\rho_{uv}\leq tt_e+\rho_{uv}.
\end{equation*}
The first inequality stems from the fact that $\lceil a/R_{uv}\rceil\leq\lceil a\rceil/R_{uv}$ for any $a>0$. The second inequality follows from the definition of the ceiling function. Thus, we obtain:
\begin{equation}\label{eq:Deltabar}
    \overline{\Delta}_r\leq \sum_{e\in\calE_r}(tt_e+\rho_{uv})=\Delta_r+\sum_{e\in\calE_r}\rho_{uv}\leq\Delta_r+\sum_{e\in\calE_r}\left(\frac{\Delta_{uv} - \Delta_r}{|\calE_r|}\right)=\Delta_{uv},\ \forall r\in\calR_{\ell st}^{uv}.
\end{equation}

This completes the proof of the lemma.\hfill\Halmos

Lemma~\ref{lem:disc} shows that there exists a discrete time unit for which all feasible subpaths in the subpath-based formulation are also feasible in the segment-based formulation. With this discretization, each segment $e_i \in \calE_r$ induces a unique arc $\overline{a}_i := (n,m) \in \overline\calA_{\ell st}$ in the time-load-expanded network, such that:
\begin{enumerate*}[label=(\roman*)]
    \item the arc corresponds to the segment: $e(\overline{a}_i) = e_i$;
    \item the capacities are consistent with pickups: $c_n = c_{start(a)}$ if $i = 1,$ and $c_n = c_{start(a)} + \sum_{j=1}^{i-1} |\calP_{e_j}|$ otherwise, and $c_m = c_{start(a)} + \sum_{j=1}^i |\calP_{e_j}|$; and
    \item the time is consistent with travel times: $t_n = T_{\ell t}(u)$ if $i = 1,$ and $t_n = T_{\ell t}(u) + \rho_{uv} \cdot \sum_{j=1}^{i-1} \lceil \frac{tt_{e_j}}{\rho_{uv}} \rceil$ otherwise, and $t_m = T_{\ell t}(u) + \rho_{uv} \cdot \sum_{j=1}^i \lceil \frac{tt_{e_j}}{\rho_{uv}} \rceil.$
\end{enumerate*}
Let us collect these time-load-expanded segment arcs into the set $\overline\calA_a := \{\overline{a}_1, \cdots, \overline{a}_{N-1}\} \subset \overline\calA_{\ell st}.$ 

With these arcs, we construct a feasible solution to the time-load-expanded segment formulation.
\begin{align*}
    \widehat{\xi}_{\overline{a}} &= \begin{cases}
        1 &\text{if }\overline{a} \in \bigcup_{a \in \calA_{\ell st} \, : \, \widehat{y}_a = 1}\overline\calA_a, \\ 
        0 &\text{otherwise,}
    \end{cases}&& \forall (\ell, t) \in \calL \times \calT_\ell, \forall \overline{a} \in \overline\calA_{\ell st}\\
    \widehat{\beta}_{uv} &= \begin{cases}
        1 &\text{if there exists } a \in \calA_{\ell st} \, : \, \widehat{y}_a = 1, u_{r(a)}=u, v_{r(a)}=v, \\ 
        0 &\text{otherwise,}
    \end{cases}&& \forall (\ell, t) \in \calL \times \calT_\ell, \forall (u,v) \in \Gamma_{\ell}.
\end{align*}

This solution satisfies the flow balance constraints in Equation~\eqref{edge:flow}. 
\begin{itemize}
    \item[--] If $x_{\ell t} = 0$,  then no subpath arc in $\calA_{\ell st}$ is selected, so no arc in $\overline\calA_{\ell st}$ is selected either. Therefore
\begin{align*}
    \sum_{j:(i,j)  \in \overline\calA_{\ell st} } \widehat\xi_{(i,j)} - \sum_{j:(j,i) \in \overline\calA_{\ell st}} \widehat\xi_{(j,i)}= 0=\begin{cases}
    x_{lt} &\text{if } i = u_{\ell st}, \\
    -x_{lt} &\text{if } i = v_{\ell st}, \\
    0 &\text{otherwise.}
    \end{cases}
\end{align*}
\item[--] Suppose that $x_{\ell t} = 1$. Recall the sequence of subpath arcs $a \in \calA_{\ell st}$ such that $\widehat y_a = 1,$ which exist by Equation~\eqref{eq:2Sflow}. We collect the corresponding segment arcs from $\overline\calA_a$ in order into set $\overline\calA_{\text{all}} = \bigcup_{a \in \calA_{\ell st}: \widehat y_a = 1} \overline\calA_a := \{\overline{a}_1,\cdots,\overline{a}_{M-1} \} $ and corresponding time-load-expanded nodes $\overline\calV^{\text{all}}  = \{n_1, \cdots, n_M\}$, where $n_1=u_{\ell st}$ and $n_M=u_{\ell st}$ and $n_2,\cdots,n_{M-1}$ refer to intermediate nodes. We have, for each node $n \in \overline\calV_{\ell st}:$
\begin{equation*}
\sum_{m:(n,m)  \in \overline\calA_{\ell st} } \widehat\xi_{(n,m)} - \sum_{m:(m,n) \in \overline\calA_{\ell st}} \widehat\xi_{(m,n)}= \begin{cases}
    \widehat\xi_{\overline{a}_{1}} - 0 = 1 = x_{\ell t} &\text{if } n = u_{\ell st}, \\
    0 - \widehat\xi_{\overline{a}_{M}} = -1 = -x_{lt} &\text{if } n = v_{\ell st}, \\
    \widehat\xi_{\overline{a}_{i}} - \widehat\xi_{\overline{a}_{i-1}} = 0 &\text{if }n = n_i \in \overline\calV^{\text{all}} \setminus \{n_1, n_M\}, \\
    0-0=0 &\text{otherwise.}
    \end{cases}
\end{equation*}
\end{itemize}

The solution also satisfies the passenger linking constraints in Equations~\eqref{edge:deactivate}.
\begin{itemize}
    \item[--]  Consider a passenger request $p \in \calP\setminus \bigcup_{a \in \calA_{\ell st} \, : \, \widehat y_a = 1} \calP_{r(a)}$. The set of pickups on the selected subpaths $\bigcup_{a \in \calA_{\ell st} \, : \, \widehat y_a = 1} \calP_a$ induces the set of pickups on the selected segments, so that $p \in \calP \setminus \bigcup_{a \in \calA_{\ell st} \, : \, \widehat y_a = 1} \bigcup_{e \in \calE_{r(a)}} \calP_e.$ Thus, $\widehat{\xi}_{\overline{a}} = 0$ for each arc $\overline{a} \in \overline\calA_{\ell st}$ with $p \in \calP_{e(\overline{a})},$ and
    $$\sum_{\overline{a} \in \overline\calA_{\ell st} \, : \, p \in \calP_{e(\overline{a})}} \widehat{\xi}_{\overline{a}} = 0 \leq z_{\ell pst}. $$
    \item[--] Consider passenger request $p \in \calP$ served by some subpath arc $a' \in \calA_{\ell st}$, so that $\widehat y_{a'} = 1$ and $p \in \calP_{r(a')}.$ Each pickup set $\calP_{r(a')}$ has been partitioned into pickup subsets at the segment level, so that there exists some $e \in \calE_{r(a')}$ with $p \in \calP_e$. This segment has been mapped to a unique arc $\overline{a}' \in \calA_{a'}$. Using Equation~\eqref{eq:2Sy2z}, we obtain:
    $$\sum_{\overline{a} \in \overline\calA_{\ell st}  :  p \in \calP_{e(\overline{a})}} \widehat{\xi}_{\overline{a}} = \widehat{\xi}_{\overline{a}'} = 1 = \widehat y_{a'} = \sum_{a \in \calA_{\ell st} : p \in \calP_{r(a)}} \widehat y_a \leq z_{\ell pst}. $$
\end{itemize}

Next, the solution satisfies the flow balance between checkpoints in Equations~\eqref{eq:flowBalanceCheckpoints}. 
\begin{itemize}
    \item[--] If $x_{\ell t} = 0$, then no subpath arcs in $\calA_{\ell st}$ are selected, so $\widehat\beta_{uv}=0$ for all $(u,v)\in\Gamma_\ell$. Therefore, for each checkpoint $u \in \calI_{\ell}$, we have:
\begin{align*}
    \sum_{v \, : (u,v) \in \Gamma_{\ell}} \widehat\beta_{uv} - \sum_{v \, : (v,u) \in \Gamma_{\ell}} \widehat\beta_{vu} = 0 = \begin{cases}
        x_{\ell t} &\text{if }u = \calI_{\ell}^{(1)} \\
        -x_{\ell t} &\text{if }u = \calI^{\text{end}} \\
        0 &\text{otherwise}
    \end{cases}
\end{align*}
\item[--] If $x_{\ell t} = 1$, then we identify the sequence of subpath arcs $a \in \calA_{\ell st}$ such that $\widehat y_a = 1,$ which exists and defines the unique sequence of checkpoints per Equation~\eqref{eq:2Sflow}. With a slight abuse of notation, this sequence is denoted by $\calI_{\ell t} := \{\omega_1:= \calI_\ell^{(1)},\cdots,\omega_{O}:= \calI_\ell^{(I_\ell)} \}$. We obtain the flow balance constraints for each checkpoint $u \in \calI_\ell$:
\begin{equation*}
\sum_{v \, : (u,v) \in \Gamma_{\ell}} \widehat\beta_{uv} - \sum_{v \, : (v,u) \in \Gamma_{\ell}} \widehat\beta_{vu} = \begin{cases}
        \widehat\beta_{\omega_1,\omega_2} - 0 = 1 = x_{\ell t} &\text{if }u =\omega_1 \\
        0 - \widehat\beta_{\omega_{O-1},\omega_O} = -1 = -x_{\ell t} &\text{if }u =\omega_O\\
        \widehat\beta_{\omega_i,\omega_{i+1}} - \widehat\beta_{\omega_{i-1},\omega_i} = 1 - 1 = 0 &\text{if }u  \in \calI_{\ell t} \backslash \{ \omega_1, \omega_{O}\} \\
        0 &\text{otherwise}
    \end{cases}
\end{equation*}
\end{itemize}

Next, the solution satisfies the checkpoint visit constraints given in Equation~\eqref{eq:visitCheckpoint}:
\begin{itemize}
    \item[--] If $x_{\ell t} = 0$, then $\widehat\beta_{uv}=0$ for all $(u,v)\in\Gamma_\ell$, so the equation is trivially satisfied.
    \item[--] If $x_{\ell t} = 1,$ then we enumerate the set of visited checkpoints by the subpaths in $\calI_{\ell t}$ using Equations~\eqref{eq:2Sflow}. We consider a checkpoint $v\in\calI_{\ell}\setminus\calI_\ell^{(1)}$. If $\sum_{w \in \calI_\ell \, : \, (w, v) \in \Gamma_\ell}\beta_{wv}=1$, then there exists a subpath-based arc $a \in \calA_{\ell st}$ such that $\widehat{y}_a = 1, u_{r(a)}=w,$ and $ v_{r(a)}=v$, which terminates in $v$. Per construction of the segment-based arcs, there exist segments $\overline{a}_1,\cdots,\overline{a}_{N-1}\in\overline\calA_a$ such that $\widehat{\xi}_{\overline{a}_1}=\cdots=\widehat{\xi}_{\overline{a}_{N-1}}=1$, corresponding to segments $e_1,\cdots,e_{N-1}$. Then, $$\sum_{\substack{(i, j) \in \overline\calA_{\ell st} \, : \\ k_j = v,\, t_j = T_{\ell t}(v)}} \widehat\xi_{(i,j)}=1,$$ and the constraint is satisfied. The constraint is trivially satisfied if $\sum_{w \in \calI_\ell \, : \, (w, v) \in \Gamma_\ell}\beta_{wv}=0$.
\end{itemize}

Finally, the solution satisfies the checkpoint sequencing constraints given in Equation~\eqref{eq:checkpointArcs}, by construction of the $\widehat{\bxi}$ and $\widehat{\bbeta}$ variables. Indeed, $\widehat{\beta}_{uv}=1$ whenever there exists an arc $\overline{a}\in \bigcup_{a \in \calA_{\ell st} \, : \, \widehat{y}_a = 1}\overline\calA_a$ between checkpoints $u$ and $v$ and between times $T_{\ell t}(u)$ and $T_{\ell t}(v)$ such that $\widehat{\xi}_{\overline{a}}=1$.

Next, the solution $\widehat{\bxi}$ achieves the same objective value as $\widehat\by$:
\begin{align}
    \sum_{\overline{a} \in \overline\calA_{\ell st}} \overline{g}_{\overline{a}} \widehat{\xi}_{\overline{a}} &= \sum_{\overline{a} \in \overline\calA_{\ell st} : \widehat{\xi}_{\overline{a}} = 1} \overline{g}_{\overline{a}}  \nonumber\\
    &= \sum_{a \in \calA_{\ell st} : \widehat y_a = 1} \sum_{\overline{a} \in \calA_{\overline{a}}} \overline{g}_{\overline{a}}  \nonumber\\
    &= \sum_{a \in \calA_{\ell st} : \widehat y_a = 1}\sum_{e \in \calE_r} \sum_{p \in \calP_e}D_{ps} \left(\lambda \tau_{ep}^{\text{walk}} + \mu \tau_{ep}^{\text{wait}} + \sigma \frac{\tau_{ep}^{\text{travel}}}{\tau_p^{\text{dir}}} + \delta \frac{\tau_{\ell tp}^{\text{late}}}{\tau_p^{\text{dir}}} + \frac{\delta}{2} \frac{\tau_{\ell tp}^{\text{early}}}{\tau_p^{\text{dir}}} - M\right)  \nonumber\\
    &= \sum_{a \in \calA_{\ell st} : \widehat y_a = 1}\sum_{p \in \calP_r} D_{ps} \left(\lambda \tau_{rp}^{\text{walk}} + \mu \tau_{rp}^{\text{wait}} + \sigma \frac{\tau_{rp}^{\text{travel}}}{\tau_p^{\text{dir}}} + \delta \frac{\tau_{\ell tp}^{\text{late}}}{\tau_p^{\text{dir}}} + \frac{\delta}{2} \frac{\tau_{\ell tp}^{\text{early}}}{\tau_p^{\text{dir}}} - M\right)  \nonumber\\
    &= \sum_{a \in \calA_{\ell st} : \widehat y_a = 1} g_a \nonumber\\
    &= \sum_{a \in \calA_{\ell st}}g_a \widehat y_a \label{objComputationSeg}
\end{align}

In conclusion, any subpath solution can be mapped into a feasible segment solution with the same objective value. Therefore, the segment-based formulation achieves an objective that is at most equal to the optimum of the subpath-based formulation.

{\it Constructing a subpath solution from a time-load-expanded segment solution.}

Suppose that $\widehat{\bxi}$ is a feasible solution to the segment-based formulation (Equations~\eqref{edge:obj}--\eqref{edge:domain2}). Assume that $x_{\ell t} = 1.$ We leverage Equations~\eqref{edge:flow} to construct a subpath between checkpoints $u$ and $v$ and between times $T_{\ell t}(u),T_{\ell t}(v)$. Starting from the source checkpoint $u$, we select the arc $\overline{a}_1 \in \overline\calA_{\ell st}$ incident with $u_{\ell st}$ for which $\widehat{\xi}_{a} = 1,$ proceeding sequentially along the directed network until reaching the node corresponding to checkpoint $v$ at time $T_{\ell t}(v)$. An outgoing arc is guaranteed at every intermediate node by Equation~\eqref{edge:flow}, and boundary conditions at the checkpoints are guaranteed by Equation~\eqref{eq:visitCheckpoint}.

Each arc $\overline{a}_i \in \overline\calA_a := \{\overline{a}_1, \cdots, \overline{a}_{N-1}\}$ corresponds to a segment $e_i := e(\overline{a}_i) \in \calE_{\ell st}$ and a passenger pickup set $\calP_{e_i}.$ The sequence of segments $\calE_r := \{e_1, \cdots, e_{N-1}\}$ defines a subpath $r$ from $u_r$ to $v_r$, skipping at most $K$ checkpoints in a row (by Equations~\eqref{eq:flowBalanceCheckpoints} and the definition of checkpoint pairs $\Gamma_\ell$), adhering to the scheduled arrival times at the checkpoints (defined by Equations~\eqref{eq:visitCheckpoint}), obeying the vehicle's capacity (by definition of the node set $\overline\calV_{\ell st}$ in the time-load-expanded network), and picking up the passengers in $\calP_r := \bigcup_{i=1}^{N-1} \calP_{e_i}$ (who are unique due to Equations~\eqref{edge:deactivate}). Thus, we obtain a unique and valid subpath-based arc in $\calA_{\ell st}$, induced by $\widehat{\bxi}$. 

Let us collect all such subpath arcs in the set $\calA(\widehat{\bxi}).$ 
We use $\calA(\widehat{\bxi})$ to construct solution $\widehat{\by}$ from $\widehat{\bxi}$:
    \begin{equation*}
        \widehat y_a = \begin{cases}
            1 &\text{if }a \in \calA(\widehat{\bxi}), \\
            0 &\text{otherwise,}
        \end{cases} \qquad \forall (\ell, t) \in \calL \times \calT_\ell, \forall a \in \calA_{\ell st}.
    \end{equation*}
    
By construction, the solution satisfies Equations \eqref{eq:2Sflow}. If $x_{\ell t} = 1$, we constructed a unique subpath-based solution for each pair $(u,v)\in\Gamma_{\ell}$. If $x_{\ell t} = 0,$ there was no subpath to construct. Therefore:
    \begin{align*}
        \sum_{a \in \calA_{\ell st}}\widehat y_a = \sum_{a \in \calA_{\ell st}} \mathbbm{1}(a \in \calA(\widehat{\bxi})) = x_{\ell t}, \quad \forall (\ell, t) \in \calL \times \calT_\ell.
    \end{align*}

The solution also satisfies Equation~\eqref{eq:2Sy2z}. For passenger $p \in \calP,$ we obtain, from the construction of the subpath solution and Equation~\eqref{edge:flow}:
    \begin{align*}
        \sum_{a \in \calA_{\ell st} : p \in \calP_{r(a)}} \widehat y_a 
        &= \sum_{a \in \calA(\widehat{\bxi})} \mathbbm{1}(p \in \calP_{r(a)}) \\
        &= \sum_{a \in \calA(\widehat{\bxi})}  \sum_{e \in \calE_{r(a)}} \mathbbm{1}\left(p \in \calP_{e}\right) \\
        &= \sum_{a \in \calA(\widehat{\bxi})}  \sum_{\overline{a} \in \overline\calA_a} \mathbbm{1}\left(p \in \calP_{e(\overline{a})}\right) \\ 
        &= \sum_{\overline{a} \in \overline\calA_{\ell st} \, : \, \widehat{\xi}_{\overline{a}} = 1} \mathbbm{1}(p \in \calP_{e(\overline{a})}) \\
        &= \sum_{\overline{a} \in \overline\calA_{\ell st} \, : \, p \in \calP_{e(\overline{a})}} \widehat{\xi}_{\overline{a}}\\
        &\leq z_{\ell pst}
    \end{align*}

Finally, the solutions $\widehat{\bxi}$ and $\widehat\by$ achieve the same objective values, which can be shown similarly to Equation~\eqref{objComputationSeg}. Therefore, any segment solution can be mapped into a feasible subpath solution with the same objective value, and the subpath-based formulation achieves an objective that is at most equal to the optimum of the segment-based formulation. This concludes the proof of equivalence of the subpath-based and segment-based formulations.

\textit{Proof that the subpath-based relaxation is at least as strong as the segment-based relaxation.}

The subpath-based formulation is a Dantzig-Wolfe reformulation of the segment-based formulation. Alternatively, we can observe that the arguments to map a segment-based solution into a subpath-based solution do not require the integrality of the subpath solution $\widehat{\by}$. By following the same steps as above, we can map any non-integral subpath solution $\widehat{\by}$ into a feasible segment solution with the same objective value, as follows:
\begin{align*}
    \widehat{\xi}_{\overline{a}} &= \sum_{a \in \calA_{\ell st} : \overline{a} \in \overline\calA_a} \widehat{y}_a \qquad \forall (\ell, t) \in \calL \times \calT_\ell, \forall \overline{a} \in \overline\calA_{\ell st} \\
    \widehat{\beta}_{uv} &= \sum_{a \in \calA_{\ell st} : u_{r(a)}=u, v_{r(a)}=u} \widehat{y}_a \qquad \forall (\ell, t) \in \calL \times \calT_\ell, \forall u,v \in \Gamma_\ell.
\end{align*}

However, a non-integral segment solution cannot be mapped directly to a subpath solution. We demonstrate this claim with an example with three checkpoints (A, B and C). Figure \ref{fig:LRsegEx} shows a non-integral segment-based solution---all shown segments have a flow of 0.5. The solution satisfies the flow balance constraints from station to station given in Equation~\eqref{edge:flow}, the flow balance constraints from checkpoint to checkpoint given in Equation~\eqref{eq:flowBalanceCheckpoints}, as well as the consistency constraints between checkpoint-checkpoint flows and station-station flows given in Equations~\eqref{eq:visitCheckpoint}--\eqref{eq:checkpointArcs}. In this solution $\widehat\bbeta$ values are $\widehat{\beta}_{(A,B)} = \widehat{\beta}_{(B,C)} = \widehat{\beta}_{(A,C)} = 0.5$, so that the flows are split between a subpath from Checkpoint A to Checkpoint B, a subpath from Checkpoint B to Checkpoint C, and a subpath from Checkpoint A to Checkpoint C, each with a flow of 0.5. The critical observation is that the solution leverages the segments (shown in solid lines) that fall outside the spatial scope of the deviations between Checkpoints A and B as part of the subpath connecting Checkpoints A and B. Specifically, there exists a segment $(n,m)\in \calA_{\ell st}$ with $k_a \in \calN^S_{13} \setminus \calN^S_{12}$ or $k_b \in \calN^S_{13} \setminus \calN^S_{12}$. This solution belongs to the polyhedron defined by the segment-based formulation, because $\xi_{(n,m)}=0.5 \leq \beta_{(A,C)} = 0.5$. However, the resulting subpath is infeasible because it connects Checkpoints A and B without adhering to the maximum deviation $\Delta$. This proves that the subpath-based relaxation is at least as strong as the segment-based relaxation.
\hfill\Halmos

\begin{figure}[h!]
    \centering
    \begin{tikzpicture}[scale=0.85,transform shape]
            \node[] at (0,0){};
            \node[] at (15,0){};
            \node[] at (4,7){$0$};
            \node[] at (5,7){$t_1$};
            \node[] at (6,7){$t_2$};
            \node[] at (7,7){$t_3$};
            \node[] at (8,7){$t_4$};
            \node[] at (9,7){$t_5$};
            \node[] at (10,7){$t_6$};
            \node[] at (11,7){$t_7$};
            \node[] at (12,7){$t_8$};
            \node[] at (13,7){$t_9$};
            \node[] at (14,7){$t_{10}$};
            \node[] at (15,7){$T$};
            \node[] at (0,6){Checkpoint A};
            \node[] at (0,5){Station 1};
            \node[] at (0,4){Station 2};
            \node[] at (0,3){Checkpoint B};
            \node[] at (0,2){Station 3};
            \node[] at (0,1){Station 4};
            \node[] at (0,0){Checkpoint C};
            \foreach \i in {0,...,6}{\foreach \t in {1,...,10}{
                    \filldraw[black] (4+\t,\i) circle (2pt);
            }}
            \foreach \i in {0,...,6}{\foreach \j in {0,...,6}{\foreach \t in {1,...,9}{
                    \draw[draw=black,dotted] (4+\t,\i) to (5+\t,\j);
            }}}
            \foreach \j in {0,...,6}{
                \draw[draw=black,dotted] (4,6) to (5,\j);
            }
            \foreach \i in {0,...,6}{
                    \draw[draw=black,dotted] (14,\i) to (15,0);
            }
            \node[ultra thick,draw=black,fill=black,minimum size=10pt] (v1_b) at (10,3){};
            \node[ultra thick,draw=black,fill=black,minimum size=10pt] (O) at (4,6){};
            \node[ultra thick,draw=black,fill=black,minimum size=10pt] (D) at (15,0){};
            \draw [ultra thick, draw=black, fill=myblue, opacity=0.2]
            (3.5,-0.5) rectangle (15.5,6.5);
            \draw [ultra thick, draw=black, fill=myred, opacity=0.3]
            (3.75,2.75) rectangle (10.25,6.25);
            \draw [ultra thick, draw=black, fill=mygreen, opacity=0.3]
            (9.75,-0.25) rectangle (15.25,3.25);
            \draw[->,ultra thick,dotted,draw=myred] (O) to (5,4);
            \draw[->,ultra thick,draw=myred] (5,4) to (7,1);
            \draw[->,ultra thick,draw=myred] (7,1) -- (8,2) node[midway, sloped, above] {$(n,m)$};
            \draw[->,ultra thick,draw=myred] (8,2) to (9,4);
            \draw[->,ultra thick,dotted,draw=myred] (9,4) to (v1_b);
            \draw[->,ultra thick,dotted,draw=mygreen] (v1_b) to (11,2);
            \draw[->,ultra thick,dotted,draw=mygreen] (11,2) to (13,1);
            \draw[->,ultra thick,dotted,draw=mygreen] (13,1) to (D);
            \draw[->,ultra thick,dotted,draw=myblue] (O) to (5,5);
            \draw[->,ultra thick,dotted,draw=myblue] (5,5) to (6,3);
            \draw[->,ultra thick,dotted,draw=myblue] (6,3) to (7,4);
            \draw[->,ultra thick,dotted,draw=myblue] (7,4) to (9,3);
            \draw[->,ultra thick,dotted,draw=myblue] (9,3) to (10,2);
            \draw[->,ultra thick,dotted,draw=myblue] (10,2) to (12,3);
            \draw[->,ultra thick,dotted,draw=myblue] (12,3) to (13,4);
            \draw[->,ultra thick,dotted,draw=myblue] (13,4) to (14,1);
            \draw[->,ultra thick,dotted,draw=myblue] (14,1) to (D);
        \end{tikzpicture}
        \caption{Example of a non-integral segment solution that cannot be mapped to a subpath solution. For simplicity, the load dimension of the time-load-expanded network is omitted. The black squares encode the reference schedule at the checkpoints. The red (resp. green, blue) area represents the stations that can be reached between Checkpoints A and B (resp. between Checkpoints B and C, between Checkpoints A and C). All thick segments are associated with a flow of 0.5. The solid segments are outside the allowable region in the subpath-based formulation.}
    \label{fig:LRsegEx}
    \vspace{-12pt}
\end{figure}

\subsection{Proof of Proposition~\ref{form_complexity}}

Let $D_{\text{min}}$ be the minimum distance between pickup locations (a constant dictated by the station set $\calN^S$) and let $\Pi$ denote the maximum distance between any checkpoint pair (a constant dictated by the candidate reference lines and the value of $K=0$ vs. $K=1$). Recall that $\Delta$ denotes the maximum deviation from the reference line for microtransit vehicles, and $\Omega$ denotes the maximum walking distance. Therefore, the rectangular service area associated with a checkpoint pair has side lengths $\Pi + 2(\Delta + \Omega)$ and $2 (\Delta + \Omega)$. The maximum number of stations between checkpoints is $\Xi = \lfloor\frac{(\Pi + 2(\Delta + \Omega))}{D_{\text{min}}}\rfloor \cdot \lfloor\frac{2\cdot (\Delta + \Omega)}{D_{\text{min}}}\rfloor$. In our case study, $\Xi$ is significantly less than $|\calN^S| = 640$.

\subsubsection*{Segment-based model.}\

The number of variables scales in $\calO(T_S \cdot C_{\ell}^2 \cdot I_\ell\cdot \Xi^2)$. 

\begin{itemize}
    \item[--] The $\bbeta$ variables are indexed over the set of valid directed checkpoint pairs $\Gamma_\ell.$ Since $K$ is a small constant, the number of valid checkpoint pairs scales linearly with the number of checkpoints, so that $\calO(|\Gamma_\ell|) = \calO\left(\sum_{i=1}^{I_\ell - (K+1)} (K+1)\right) = \calO(I_\ell).$
    \item[--] The $\bxi$ variables are indexed over the set $\overline{\calA}_{\ell st},$ i.e., the arcs in the time-load-expanded network. Arcs define connections between consecutive stations between a checkpoint pair (which scale in $\calO(I_\ell \cdot \Xi^2)$) for each time period, and they can correspond to any vehicle load pair, so the number of $\bxi$ variables scales in $\calO(\overline{\calA}_{\ell st}) =  \calO(T_S \cdot C_{\ell}^2 \cdot I_\ell\cdot \Xi^2)$.
\end{itemize}

The number of constraints scales in $\calO(|\calP| + T_S \cdot C_{\ell} \cdot |\calN^S| + T_S\cdot C_{\ell}^2 \cdot I_\ell\cdot \Xi^2)$. 

\begin{itemize}
    \item[--] Equations~\eqref{edge:flow}:  There is one constraint per node in the time-load-expanded network. There is one node per combination of time periods in $\calT^S$, vehicle loads in $\calC_{\ell},$ and stations in $\calN^S$, so that there are $\calO(T_S \cdot C_{\ell} \cdot |\calN^S|)$ flow balance constraints. 
    \item[--] Equations~\eqref{edge:deactivate}: The passenger linking constraints scale with $\calO(\sum_{p \in \calP}|\calM_p|)$. The cardinality of each set $\calM_p$ is bounded by a small constant, so there are $\calO(|\calP|)$ linking constraints.
    \item[--] Equations~\eqref{eq:flowBalanceCheckpoints}--\eqref{eq:visitCheckpoint}: There are $\calO(I_{\ell})$ flow balance constraints for checkpoint-to-checkpoint flows and $\calO(I_\ell)$ schedule adherence constraints. 
    \item[--] Equations~\eqref{eq:checkpointArcs}: 
    There is one constraint per arc in the time-load-expanded network to ensure consistency between the station-to-station and checkpoint-to-checkpoint flows, which grows in $\calO(\overline\calA_{\ell st}) = \calO(C_{\ell}^2 \cdot T_S \cdot I_\ell \cdot \Xi^2)$, as previously established.
\end{itemize}

The complexity of Equations~\eqref{eq:flowBalanceCheckpoints}--\eqref{eq:visitCheckpoint} is dominated by that of Equations~\eqref{edge:flow} and~\eqref{eq:checkpointArcs}. The result follows.

\subsubsection*{Subpath-based model.}\ 
The number of variables scales in $\calO(I_\ell \cdot C_{\ell} \cdot 2^\Xi)$. In particular, $\by$ scales with $\calO(|\calA_{\ell st}|),$ the number of arcs in the load-expanded subpath network. By definition:
\begin{align*}
    |\calA_{\ell st}| &= |\calA^v_{\ell st}| + \sum_{r \in \calR_{\ell st}} |\calA_r|.
\end{align*}
As for $|\calA^v_{\ell st}|,$ there are $C_\ell + 1$ arcs connecting the last stop to the sink node (one per vehicle load). Turning to $\calA_r$, a subpath $r \in \calR_{\ell st}$ is the shortest path to serve the corresponding passenger set $\calP_{r}$. Thus, the number of subpath variables is proportional to the number of possible sets $\calP_{r}.$ The number of different passengers that can be picked up at each station is bounded by a small constant, so we use the number of stations as a proportional proxy for the number of passengers that can be picked up. There are up to $\binom{\Xi}{c}$ station combinations that pick up $c$ passengers, each of which can be replicated $C_{\ell} - c + 1$ times in the arc set (corresponding to initial loads $0,1,\cdots,C_\ell-c$). Therefore, the number of subpaths is
    $$\sum_{c=0}^{C_{\ell}} \binom{\Xi}{c} \cdot (C_{\ell} - c + 1) \leq (C_{\ell}+1) \cdot  \sum_{c=0}^{C_{\ell}} \binom{\Xi}{c}.$$
When $\Xi \leq C_{\ell},$ the binomial sum above is equal to $2^\Xi.$ When $\Xi > C_{\ell}$, it is equal to
$$\sum_{c=0}^{C_{\ell}}\binom{\Xi}{c} = 2^\Xi - \sum_{c = C_{\ell} + 1}^\Xi \binom{\Xi}{c} \leq 2^\Xi.  $$
Therefore, $\calO(|\calA_{\ell st}|) = \calO(2^{\Xi} \cdot C_{\ell} \cdot I_\ell)$.
    
The number of constraints scales in $\calO(|\calP| + C_{\ell} \cdot I_\ell)$.
\begin{itemize}
    \item[--] Equations~\eqref{eq:2Sflow}: There are $\calO(\calV_{\ell st}) = \calO(C_{\ell} \cdot I_\ell)$ flow balance constraints, one per node in the load-expanded subpath network.
    \item[--] Equations~\eqref{eq:2Sy2z}: There are $\calO(|\calP|)$ linking constraints, one per passenger that can be picked up by reference trip $(\ell, t)$.
\end{itemize}

\subsubsection*{Path-based model.}
The number of variables scales in $\calO(2^{\Xi\cdot I_{\ell}})$. The $\mathbf{\zeta}$ variables are indexed over the path set $\calQ_{\ell st}.$ Each path $q\in\calQ_{\ell st}$ can be decomposed into a sequence of subpaths in $\calR_{\ell st}$ by partitioning the path-based passenger set $\calP_q$ into subpath-based passenger sets $\calP_r$, that is $\calP_q = \bigcup_{r \in \calR_q} \calP_r.$ Recall that there are $O(2^\Xi)$ possible subpaths between each checkpoint pair, and there are $\calO(I_{\ell})$ checkpoint pairs, so we obtain $\calO(2^{\Xi\cdot I_{\ell}})$ overall paths.
    
The number of constraints scales in $\calO\left(|\calP|\right)$, i.e., with the number of linking constraints (Equations~\eqref{path:deactivate}). The model also comprises a single partitioning constraint (Equation~\eqref{path:flow}), which does not affect the constraint complexity.
\hfill\Halmos
\section{Details on Solution Algorithm}\label{app:alg}

\subsection{Proof of Theorem \ref{thm:exact}}

First, we show that the inner column generation procedure converges to an optimal solution of $\overline{\texttt{SP}}(\bx)$ for a given first stage solution $\bx$. The restricted subproblem $\texttt{RSP}(\bx,\calA_{sj}^\prime)$ has optimal value $\varphi^\prime(\bx,\calA_{sj}^\prime) \geq \varphi(\bx)$ because $\calA_{sj}^\prime \subset \calA_{sj}$ for all $s\in\calS,j\in\calJ$. Let $\by^\prime$ be the primal solution to $\texttt{RSP}(\bx,\calA_{sj}^\prime)$ corresponding to dual solution $(\boldsymbol{\psi}, \bgamma)$. At each column generation iteration and for each $s\in\calS,j\in\calJ$, if $\texttt{RC} < 0$, we add the corresponding arc to $\calA_{sj}^\prime$. If $\texttt{RC} \geq 0$ for all $s\in\calS,j\in\calJ$, then: 
\begin{align*}
     \bar{c}_{(m,n)} = g_{(m,n)} - (\psi_{sjm} - \psi_{sjn}) - \bgamma^\top_{sj} \bff_{(m,n),s,j} \geq 0 \quad \forall (m,n)\in\calA_{sj}
\end{align*}
We construct a feasible solution $\by$ to $\overline{\texttt{SP}}(\bx)$ by letting $y_a = y_a^\prime$ for all $a\in\calA_{sj}^\prime$ and $y_a = 0$ for all $a\in\calA_{sj}\setminus\calA_{sj}^\prime$. Then, $\by$ is optimal for $\overline{\texttt{SP}}(\bx)$ by the non-negativity of the reduced costs. Thus, $\varphi_{sj}^\prime(\bx,\calA_{sj}^\prime) = \overline{\varphi}_{sj}(\bx)$. There are finitely many arcs in $\calA_{sj}$ and at each iteration with $\texttt{RC} < 0$ for some $s\in\calS,j\in\calJ$, at least one new arc is added to $\calA_{sj}$. Thus, the algorithm converges to an optimal solution of $\overline{\texttt{SP}}(\bx)$ in finitely many iterations.

Then, we show that the outer Benders decomposition scheme converges to an optimal solution of the partial relaxation $(\texttt{MIO}-\texttt{LO})$. Suppose that, at iteration $t$ in the outer Benders decomposition loop, we solve $\texttt{MP}(\calU^t,\calV^t)$ and obtain solution $(\bx^t, \btheta^t)$. The value $\bc^\top \bx^t + \sum_{s\in\calS} \sum_{j\in\calJ} \pi_s \theta_{sj}^t$ satisfies a subset of Benders cuts, and thus provides a valid lower bound for $(\texttt{MIO}-\texttt{LO})$. We then solve the relaxed subproblem $\overline{\texttt{SP}}(\bx^t)$ via column generation (proved to converge above). 
\begin{itemize}
    \item[--] If $\overline{\texttt{SP}}(\bx^t)$ is unbounded, then $\exists s \in\calS,j \in\calJ$ and a corresponding direction of unboundedness $(\boldsymbol{\psi}^v_{sj},\bgamma^v_{sj})$ which is dual feasible (i.e. $(\boldsymbol{\psi}^v_{sj},\bgamma^v_{sj}) \in \calP_{sj}$) and $ \sum_{n\in\calN_{sj}}\sum_{k\in\calK_j}b_{nsk}x^t_k\psi^v_{sjn}+\bh_{sj}^\top\bgamma^v_{sj} > 0$. This generates a valid feasibility cut \eqref{bendersfeascut}, which is clearly violated by $(\bx^t, \btheta^t)$. We add the corresponding cut to the master problem, eliminating the incumbent solution.
    \item[--] If $\overline{\texttt{SP}}(\bx^t)$ admits a feasible solution, let $(\boldsymbol{\psi}^u,\bgamma^u)$ be the optimal dual solution. This generates a valid optimality cut \eqref{bendersoptcut}. If $(\bx^t, \btheta^t)$ violates the cut (i.e. there exists $s\in\calS,j\in\calJ$ such that $\varphi_{sj}(\bx^t) > \theta^t_{sj}$), then we add the corresponding cut to the master problem, eliminating the incumbent solution. We also obtain a valid upper bound for $(\texttt{MIO}-\texttt{LO})$ from $\bc^\top \bx^t + \sum_{s\in\calS} \sum_{j\in\calJ} \pi_s \varphi_{sj}(\bx^t)$. If no cuts are violated, then the solution is optimal for $(\texttt{MIO}-\texttt{LO})$. 
\end{itemize}

Each Benders iteration either generates a new cut to exclude the incumbent $(\bx^t,\btheta^t)$ or identifies that the incumbent satisfies all cuts. Because the set of cuts is finite, the algorithm terminates in finitely many iterations.

Upon convergence, we restore integrality requirements in the subproblem by finding a feasible integer solution to $\texttt{SP}(\bx^*)$, with optimal value given by $\varPhi^\prime(\bx^*)$. The result is a feasible integer solution to Problem~\eqref{OPT} and thus $\bc^\top \bx^* + \varPhi^\prime(\bx^*)$ provides an upper bound. The optimal solution to the partial relaxation $(\texttt{MIO}-\texttt{LO})$ provides a valid lower bound for Problem~\eqref{OPT}. These together constitute a valid optimality gap for Problem~\eqref{OPT} upon convergence of the DD algorithm.
\hfill\Halmos

\subsection{Application of the double-decomposition algorithm to the MiND-VRP}
\label{app:VRPalg}

We provide details on the double decomposition algorithm applied to the MiND-VRP partial relaxation. Note that the problem has relatively complete resource, since it is always feasible to merely follow the reference line without deviations. Therefore, the methodology involves optimality cuts but no feasibility cuts. 

In the MiND-VRP, the Benders subproblem is decomposable across reference trips $(\ell,t) \in \mathcal{L} \times \mathcal{T}_{\ell}$ and scenarios $s \in \mathcal{S}$. It is given by the following Benders subproblem, for each first-stage decision $(\bx,\bz)$:
\begin{align}
    \overline\varphi_{sj}(\bx)= \min_{\by\geq\bo} &\quad \sum_{a \in \mathcal{A}_{\ell st}} g_a y_a \label{BSP:obj}     \qquad\st\quad\text{Equations~\eqref{eq:2Sflow}-\eqref{eq:2Sy2z}}
\end{align}

Let $\boldsymbol{\psi}$ and $\bgamma$ respectively denote the dual variables corresponding to Equations~\eqref{eq:2Sflow} and~\eqref{eq:2Sy2z}, respectively. The dual Benders subproblem is then formulated as follows:
\begin{align}
    \max \quad &x_{\ell t} \cdot (\psi_{\ell s t u_{\ell st}} - \psi_{\ell s t v_{\ell st}}) - 
    \sum_{p \in \calP \, : \,(\ell,t) \in \mathcal{M}_p}
    z_{\ell pst} \cdot \gamma_{\ell s t p}  \\
    \text{s.t. } \quad &\psi_{\ell s t n} - \psi_{\ell s t m} -  \sum_{p \in \calP_a } \gamma_{\ell s t p}  \leq g_a && \forall  a=(n,m) \in \mathcal{A}_{\ell st} \label{Benders:dualArc} \\
    &\psi_{\ell s t i} \in \mathbb{R} && \forall i \in \mathcal{V}_{\ell st}  \\
    &\gamma_{\ell s t p} \geq 0 && \forall p \in \calP \, : \,(\ell ,t) \in \mathcal{M}_p
\end{align}

Let $\Lambda_{\ell st}$ store the extreme points of the dual second-stage polyhedron, each corresponding to a second-stage solution $(\boldsymbol{\psi},\bgamma)$ for reference trip $(\ell,t) \in \mathcal{L} \times \mathcal{T}_{\ell}$ and scenario $s \in \mathcal{S}$. Let $\mathbf{\Lambda} = (\Lambda_{\ell st})_{(\ell, t) \in \calL \times \calT_{\ell}, s \in \calS}$ store all extreme points. The MiND-VRP partial relaxation optimizes network design and passenger assignments subject to a piece-wise linear recourse approximation:
{\small\begin{align}
    \text{BMP}(\mathbf{\Lambda}) \ \ \min  \quad &\sum_{\ell \in \calL} \sum_{t \in \calT_\ell} \left( h_\ell x_{\ell t} + \sum_{s \in \mathcal{S}} \pi_s  \theta_{\ell st}\right) \label{BD:obj}\\
    \text{s.t.} \quad & \text{Equations~\eqref{eq:budget}--\eqref{eq:load2}} \\
    \quad&\theta_{\ell st}\geq x_{\ell t} \cdot (\psi_{\ell s t  u_{\ell st}} - \psi_{\ell s t v_{\ell st}}) - 
    \sum_{p \in \calP \, : \,(\ell,t) \in \mathcal{M}_p}
    z_{\ell pst} \cdot \gamma_{\ell s t p},\ \  \forall (\ell, t) \in \calL \times \calT_{\ell},\ s\in\calS,\ (\boldsymbol{\psi}, \bgamma) \in  \Lambda_{\ell st} \label{BD:genCut}\\
    \quad&\bx, \bz \text{ binary} \label{BD:domain}
\end{align}}

The restricted Benders subproblem simply solves the Benders subproblem with a subset of subpath-based arcs by $\calA'_{\ell st}\subseteq\mathcal{A}_{\ell st}$. It is formulated as follows:
\begin{align}
    \varphi'_{s,j}(\bx,\calA'_{sj})=\min_{\by\geq\bo} &\quad \sum_{a \in \calA'_{\ell st}} g_a y_a \label{RBSP:obj} \\
    \st&\quad\sum_{m:(n,m) \in \calA'_{\ell st}} 
    y_{(n,m)} - \sum_{m:(m, n) \in \calA'_{\ell st}} 
    y_{(m, n)} = \begin{cases}
    x_{\ell t} &\text{if } n = u_{\ell st} \\
    -x_{\ell t} &\text{if } n = v_{\ell st} \\
    0 &\text{ otherwise}
    \end{cases} \ \forall n \in \calV_{\ell st}\label{RBSP:2Sflow} \\
    &\sum_{a \in \calA'_{\ell st} \, : \, p \in \calP_{r(a)}} y_a \leq z_{\ell pst} \quad \forall p \in \calP: (\ell,t) \in \calM_p \label{RBSP:2Sy2z}
\end{align}

\subsection{Single-tree implementation of the double-decomposition algorithm}
\label{app:singletree}

Figure~\ref{fig:singletree_DD} visualizes the DD methodology in a single-tree Benders decomposition implementation. This implementation solves Problem~$(\texttt{MIO}-\texttt{LO})$ using a single branch-and-cut tree, by solving the master problem relaxation (continuous first- and second-stage variables) at each node, and solving the Benders subproblem whenever an integral solution is found. In other words, multi-tree implementation, shown in Figure~\ref{fig:multitree_DD}, relies on a cutting-plane version of Benders decomposition, whereas single-tree implementation adds Benders cuts via lazy constraints in the branch-and-cut tree.

\begin{figure}[h!]
    \centering
    \begin{tikzpicture}[scale=0.4, every node/.style={scale=0.75},every text node part/.style={align=center}]
        \tikzset{test/.style={diamond,fill=mygray!20,draw,aspect=2,minimum width = 50pt, minimum height=30pt,text badly centered}}
        \tikzset{reg/.style={rectangle,draw, minimum width = 60pt, minimum height=20pt}}
        \node[test] (termtest) at (0,18.5){$LB=UB$};
        \node[rectangle,draw] (term) at (9,18.5){Terminate};

        \node[reg] (LO) at (-10.5,13){Choose node, solve $\overline{\texttt{MP}}(\calU^0,\calV^0)$\\Solution $(\bx,\by)$, objective $MP$};
        
        \node[test] (bad) at (-10.5,8.5){Infeasible\\or $MP\geq UB$};
        \node[reg,anchor=north,rectangle,draw] (pruneB) at (-10.5,4.5){Prune};
        \draw[->,thick] (bad.south) -- node[right,midway]{YES} (pruneB.north);
        
        \node[test] (integer) at (-2.5,8.5){$\bx\in\calX^{\text{MIP}}$};
        \node[reg,fill=myred!20] (BSP) at (6,8.5){Solve $\texttt{RSP}(\bx,\calA'_{sj})$};
        \node[reg,fill=myred!20] (PP) at (6,4){Solve $\texttt{PP}(\bpsi,\bgamma)$};
        \node[test] (RCneg) at (6,0){$RC<0$};
        \node[test] (violation) at (15,8.5){Violated cuts};
        
        \draw[->,thick,color=myred] (BSP.south) -- (PP.north);
        \draw[->,thick,color=myred] (PP.south) -- (RCneg);
        \draw[->,thick,color=myred] (RCneg.west) -- node[above,midway]{YES} +(-2,0) |- (BSP.south west);
        \draw[->,thick,color=myred] (RCneg.east) -- node[above,midway]{NO} +(2,0) |- (violation.west);
        
        \node[reg,anchor=north,rectangle,draw] (pruneF) at (15,4.5){Prune; $UB\gets MP$};
        \draw[->,thick] (violation.south) -- node[right,very near start]{NO} (pruneF.north);
        
        \node[reg,draw] (Benders) at (2,13){Add Benders cuts\\Expand $\calU^0,\calV^0$};
        
        \node[test] (cut) at (-2.5,0){MIO cuts};
        \node[reg] (branch) at (-2.5,-4){Branch};
        \node[reg] (bound) at (-10.5,-4){Update lower bound};
        
        \draw[->,thick] (termtest.east) -- node[above,midway]{YES}(term);
        \draw[->,thick] (termtest.south) -- node[right,midway]{NO} +(0,-1.5) -- +(-(10.5,-1.5) -- (LO.north);
        \draw[->,thick] (LO.south) -- (bad.north);
        \draw[->,thick] (bad.east) -- node[above,midway]{NO}(integer.west);
        \draw[->,thick] (integer.east) -- node[above,midway]{YES} (BSP.west);
        \draw[->,thick] (violation.north) -- node[right,very near start]{YES} (15,13) -- (Benders.east);
        \draw[->,thick] (Benders.west) -- (LO.east);
        
        \draw[->,thick] (integer.south) -- node[right,very near start]{NO} (cut.north);
        \draw[->,thick] (cut.west) -- node[above,very near start]{YES} +(-12.5,0) |- (LO.west);
        \draw[->,thick] (cut.south) -- node[right,very near start]{NO} (branch.north);
        \draw[->,thick] (branch) -- (bound);
        \draw[->,thick] (bound.west) -- +(-4.5,0) |- (termtest.west);

        \draw[->,dashed,thick] (integer.north) -- +(0,1.5) node[right,midway]{NO} -- +(7,1.5) -- +(7,-0.5);

        \draw[color=myred,dotted,ultra thick] (.5,-1.5) rectangle (11,10.5);
    \end{tikzpicture}
    \caption{Single-tree DD algorithm to solve Problem~$(\texttt{MIO}-\texttt{LO})$.}
    \label{fig:singletree_DD}
\end{figure}

\subsection{Details on the DD\&ILS algorithm}
\label{app:DDILS}

Recall that $\underline{\varPhi}_{sj}$ denotes a global lower bound of the second-stage cost $\varphi_{sj}(\bx)$. We index the first-stage solutions visited throughout the algorithm by $\{\widehat\bx^{w}:w\in\calW^0\}$. The master problem, now referred to as $\texttt{MP-ILS}(\calU^0,\calW^0)$, comprises Benders optimality cut (it omits feasibility cuts due to the relatively complete recourse assumption) and integer L-shaped cuts, which express that $\theta_{sj}\geq\varphi_{sj}(\widehat\bx^{w})$ whenever $\bx=\widehat\bx^w$ and that $\theta_{sj}\geq\underline{\varPhi}_{sj}$ otherwise:
\begin{align}
    \min\quad   &   \sum_{j \in \calJ}\sum_{k\in\calK_j}c_kx_k+\sum_{s\in\calS}\sum_{j\in\calJ}\pi_s\theta_{sj}\label{MPILS}\tag{$\texttt{MP-ILS}(\calU^0,\calW^0)$}\\
    \st\quad    &   \bA\bx\geq\bb\nonumber\\
                &   \theta_{sj}\geq\sum_{n\in\calN_{sj}}\sum_{k\in\calK_j}b_{nsk}x_k\psi^u_{sjn}+\bh_{sj}^\top\bgamma^u_{sj},\ \forall s\in\calS,\ \forall j\in\calJ,\ \forall u\in\calU_{sj}^0\nonumber\\
                &   \theta_{sj}\geq\underline{\varPhi}_{sj}+(\varphi_{sj}(\widehat\bx^{w})-\underline{\varPhi}_{sj})\left(1-\sum_{k\in\calK_j:\widehat{x}^{w}_k=1}(1-x_k)-\sum_{k\in\calK_j:\widehat{x}^{w}_k=0}x_k\right),\ \forall s\in\calS,j\in\calJ,w\in\calW^0\nonumber\\
                &   \bx\in\calX^{\text{MIO}}\nonumber
\end{align}

The DD\&ILS algorithm is shown in Figure~\ref{fig:DD_ILS}, with our DD methodology highlighted in red and the integer L-shaped cuts in purple. The algorithm relies on the multi-tree implementation of the DD algorithm (Figure~\ref{fig:multitree_DD}). Upon convergence, the DD algorithm returns an optimal solution to Problem~$(\texttt{MIO}-\texttt{LO})$. We then restore integrality requirements in the second-stage problem and solve $\texttt{SP}(\bx)$. If the optimality gap lies within some tolerance, the algorithm terminates with an optimality guarantee for Problem~\eqref{OPT}. Otherwise, the algorithm generates new integer L-shaped cuts. These cuts make the incumbent solutions $\{\theta_{sj}:s\in\calS,j\in\calJ\}$ infeasible and the algorithm proceeds to the master problem. The overall DD\&ILS algorithm therefore converges to an exact solution to Problem~\eqref{OPT} in a finite number of iterations, due to the binary first-stage structure.

\begin{figure}[h!]
    \centering
    \begin{tikzpicture}[scale=0.4, every node/.style={scale=0.75},every text node part/.style={align=center}]
        \tikzset{test/.style={diamond,fill=mygray!20,draw,aspect=2,minimum width = 50pt, minimum height=30pt,text badly centered}}
        \tikzset{reg/.style={rectangle,draw, minimum width = 60pt, minimum height=20pt}}

        \node[reg] (LO) at (-10.5,13){Solve $\texttt{MP-ILS}(\calU^0,\calW^0)$\\Solution $(\bx,\by)$, objective $MP$};
        
        \node[test] (bad) at (-10.5,8.5){Infeasible};
        \node[reg,anchor=north,rectangle,draw] (terminateF) at (-10.5,4.5){Terminate\\Problem infeasible};
        \draw[->,thick] (bad.south) -- node[right,midway]{YES} (terminateF.north);
        
        \node[reg,fill=myred!20] (BSP) at (2,8.5){Solve $\texttt{RSP}(\bx,\calA'_{sj})$};
        \node[reg,fill=myred!20] (PP) at (2,4){Solve $\texttt{PP}(\bpsi,\bgamma)$};
        \node[test] (RCneg) at (2,0){$RC<0$};
        \node[test] (violation) at (12,8.5){Violated cuts};
        
        \draw[->,thick,color=myred] (BSP.south) -- (PP.north);
        \draw[->,thick,color=myred] (PP.south) -- (RCneg);
        \draw[->,thick,color=myred] (RCneg.west) -- node[above,midway]{YES} +(-2,0) |- (BSP.south west);
        \draw[->,thick,color=myred] (RCneg.east) -- node[above,midway]{NO} +(2,0) |- (violation.west);

        \node[reg,fill=DarkPurple!20] (SPMIO) at (12,4.5){Solve $\texttt{SP}(\bx)$};
        \node[test] (intgap) at (12,0.5){Integrality gap};
        \node[reg,anchor=north,rectangle,draw,fill=DarkPurple!20] (terminateO) at (12,-3.5){Terminate\\Return optimal solution};
        \draw[->,thick,color=DarkPurple] (violation.south) -- node[right,very near start]{NO} (SPMIO.north);
        \draw[->,thick,color=DarkPurple] (SPMIO.south) -- (intgap.north);
        \draw[->,thick,color=DarkPurple] (intgap.south) -- node[right,midway]{NO} (terminateO.north);
        
        \node[reg,draw,fill=DarkPurple!20] (ILS) at (2,-8){Add integer L-shaped cuts\\Expand $\calW^0$};
                
        \draw[->,thick,color=DarkPurple] (intgap.east) -- node[above,midway]{YES} +(2,0) |- (ILS.east);
        \draw[->,thick,color=DarkPurple] (ILS.west) -- +(-14,0) |- (LO.west);
                
        \node[reg,draw] (Benders) at (2,13){Add Benders cuts\\Expand $\calU^0$};
        
        \draw[->,thick] (LO.south) -- (bad.north);
        \draw[->,thick] (bad.east) -- node[above,midway]{NO}(BSP.west);
        \draw[->,thick] (violation.north) -- node[right,very near start]{YES} (12,13) -- (Benders.east);
        \draw[->,thick] (Benders.west) -- (LO.east);

        \draw[color=myred,dotted,ultra thick] (-3.5,-2) rectangle (7.5,10);
    \end{tikzpicture}
    \caption{DD\&ILS methodology to solve Problem~\eqref{OPT}.}
    \label{fig:DD_ILS}
\end{figure}

\subsubsection*{Proof of Theorem \ref{thm:DDILS}}

Consider master problem solution $\left(\widehat{\bx}^t, \widehat{\btheta}^t\right)$ at Benders iteration $t$. Suppose that the double decomposition algorithm has converged, meaning that $\widehat{\theta}_{sj}^t \geq \overline\varphi_{sj}(\widehat{\bx}^t)$ for all $s\in\calS$ and $j\in\calJ$; in fact, $\theta_{sj}^t = \overline\varphi_{sj}(\widehat{\bx}^t)$ by optimality (Theorem \ref{thm:exact}). We solve $\texttt{SP}(\widehat{\bx}^t)$ and obtain integer optimal value $\varphi_{sj}(\widehat{\bx}^t) \geq \overline\varphi_{sj}(\widehat{\bx}^t)$ for all $s\in\calS$ and $j\in\calJ$. Recall that $\texttt{SP}(\widehat{\bx}^t)$ is feasible due to the assumption of relatively complete recourse. The integrality gap is defined as the difference in optimal value between the feasible solution to Problem~\eqref{OPT} and the optimal solution to the partial relaxation, that is:
\begin{align*}
    \texttt{GAP} &= \left( \bc^\top \widehat{\bx}^t + \sum_{s\in\calS} \pi_s \sum_{j\in\calJ} \varphi_{sj}(\widehat{\bx}^t) \right) - \left( \bc^\top \widehat{\bx}^t + \sum_{s\in\calS} \pi_s \sum_{j\in\calJ} \widehat{\theta}^t_{sj} \right) \\
    &= \sum_{s\in\calS} \pi_s \sum_{j\in\calJ} (\varphi_{sj}(\widehat{\bx}^t)-\widehat{\theta}^t_{sj})
\end{align*}
If $\texttt{GAP} \geq \varepsilon$ for a small tolerance $\varepsilon>0$, then $\widehat{\theta}^t_{sj} < \varphi_{sj}(\widehat{\bx}^t)$ for some $s\in\calS, j\in \calJ$. We add the integer L-shaped cuts given in Equation~\eqref{ils:cut} to the master problem, which evaluate to the valid inequality, $\theta_{sj}\geq\underline{\varPhi}_{sj}$ for all $\bx \neq \widehat{\bx}^t$. Meanwhile, it eliminates the incumbent solution $\left(\widehat{\bx}^t, \widehat{\btheta}^t\right)$ because it imposes that $\theta_{sj}\geq\varphi_{sj}(\widehat\bx^{t})$ when $\bx=\widehat\bx^{t}$:
$$\theta_{sj}\geq\underline{\varPhi}_{sj}+(\varphi_{sj}(\widehat\bx^{t})-\underline{\varPhi}_{sj})\Bigg(1-\underbrace{\sum_{k\in\calK_j:\widehat{x}^{t}_k=1}(1-x_k)-\sum_{k\in\calK_j:\widehat{x}^{t}_k=0}x_k}_{\text{=0 for $\bx=\widehat\bx^{t}$}}\Bigg) \quad \forall s\in\calS,j\in\calJ$$

If $\texttt{GAP} < \varepsilon$, then the algorithm has converged and $(\widehat{\bx}^t, \btheta^t)$ is an optimal solution of Problem~\eqref{OPT}. Since $\bx\in \{0,1\}^{n_Z}$, there are finitely many first stage solutions $\bx$ and thus finitely many cuts of the form \eqref{ils:cut}. Therefore, the algorithm converges in a finite number of iterations to an optimal solution of Problem~\eqref{OPT}.
\hfill\Halmos

\subsection{Details on the UB\&DD algorithm}
\label{app:UBDD}

At each node, the UB\&BC algorithm solves the master problem relaxation $\overline{\texttt{MP}}(\calU^0,\calV^0)$ and stores the first-stage solution $\widehat\bx$. By design, $\overline{\texttt{MP}}(\calU^0,\calV^0)$ incorporates a subset of Benders cuts and relaxes the integrality constraints. Therefore, if $\texttt{MP}(\calU^0,\calV^0)$ is infeasible or returns a solution that is larger than an upper bound of Problem~\eqref{OPT}, the branch can be pruned. If the solution does not satisfy first-stage integrality requirements, the algorithm proceeds to its branch-and-cut elements. Otherwise, it solves the second-stage relaxation to add Benders cuts or prove that none exist. These elements mirror the single-tree implementation of Benders decomposition (Figure~\ref{fig:singletree_DD}); by themselves, they would solve the partial relaxation with mixed-integer first-stage variables and continuous second-stage variables (that is, Problem $(\texttt{MIO}-\texttt{LO})$).

The difference occurs when the solution $\widehat\bx$ of $\overline{\texttt{MP}}(\calU^0,\calV^0)$ satisfies first-stage integrality requirements and all Benders cuts at a given node. By design, $\widehat{\bx}$ solves the full relaxation with continuous first- and second-stage variables at that node (because no Benders cut is violated), hence the partial relaxation with mixed-integer first-stage variables and continuous second-stage variables (because it satisfies the integrality constraints). For instance, at the root node, $\widehat\bx$ solves $(\texttt{MIO}-\texttt{LO})$. However, $\widehat\bx$ may not solve the full problem at the node. Unlike traditional single-tree implementations of Benders decomposition, the UB\&BC algorithm proceeds with branching rather than pruning. It first updates a lower bound corresponding to the incumbent first-stage solution $\widehat\bx$:
\begin{equation}\label{eq:LB}
    LB(\widehat{\bx})=\bc^\top\widehat\bx+\sum_{s\in\calS}\sum_{j\in\calJ}\pi_s\overline\varphi_{sj}(\widehat\bx)\leq\texttt{OPT}(\widehat\bx)
\end{equation}
If the lower bound is higher than the upper bound, the branch can be pruned; otherwise, the algorithm stores the first-stage solution solves the mixed-integer second-stage problem---with an exact or heuristic algorithm---and updates the upper bound as needed:
\begin{equation}\label{eq:UB}
    UB\gets\min\left\{UB,\texttt{OPT}(\widehat\bx)\right\}\geq\texttt{OPT}
\end{equation}

Upon termination, the algorithm returns the solution $\widehat\bx$ with the smallest value of $\texttt{OPT}(\widehat\bx)$.

The UB\&DD algorithm, shown in Figure~\ref{fig:UBDD}, embeds our DD methodology (highlighted in red) into the UB\&BC methodology from \cite{maheo2024unified} (in blue).
\begin{figure}[h!]
    \centering
    \begin{tikzpicture}[scale=0.4, every node/.style={scale=0.75},every text node part/.style={align=center}]
        \tikzset{test/.style={diamond,fill=mygray!20,draw,aspect=2,minimum width = 50pt, minimum height=30pt,text badly centered}}
        \tikzset{reg/.style={rectangle,thick,draw, minimum width = 60pt, minimum height=20pt}}
        
        \node[test] (termtest) at (0,18.5){Empty queue};
        \node[reg,draw] (term) at (9,18.5){Terminate};

        \node[reg] (LO) at (-10.5,13){Choose node, solve $\overline{\texttt{MP}}(\calU^0,\calV^0)$\\Solution $(\bx,\by)$, objective $MP$};
        
        \node[test] (bad) at (-10.5,8.5){Infeasible\\or $MP\geq UB$};
        \node[reg,anchor=north,rectangle,draw] (pruneB) at (-10.5,4.5){Prune};
        \draw[->,thick] (bad.south) -- node[right,midway]{YES} (pruneB.north);
        
        \node[test] (integer) at (-2.5,8.5){$\bx\in\calX^{\text{MIO}}$};
        \node[reg,fill=myred!20] (BSP) at (6,8.5){Solve $\texttt{RSP}(\bx,\calA'_{sj})$};
        \node[reg,fill=myred!20] (PP) at (6,4){Solve $\texttt{PP}(\bpsi,\bgamma)$};
        \node[test] (RCneg) at (6,0){$RC<0$};
        \node[test] (violation) at (15,8.5){Violated cuts};
        
        \draw[->,thick,color=myred] (BSP.south) -- (PP.north);
        \draw[->,thick,color=myred] (PP.south) -- (RCneg);
        \draw[->,thick,color=myred] (RCneg.west) -- node[above,midway]{YES} +(-2,0) |- (BSP.south west);
        \draw[->,thick,color=myred] (RCneg.east) -- node[above,midway]{NO} +(2,0) |- (violation.west);
        
        \node[reg,fill=myblue!20] (LB) at (15,4){Update lower bound\\(Equation~\eqref{eq:LB})};
        \node[test] (good2) at (15,0){$LB(\bx)\geq UB$};
        \node[reg,rectangle,draw,fill=myblue!20] (pruneB2) at (22.5,0){Prune};
        \node[reg,fill=myblue!20] (SPMIO) at (15,-4){Solve $\texttt{SP}(\bx)$};
        \node[reg,fill=myblue!20] (UB) at (15,-8){Update upper bound\\(Equation~\eqref{eq:UB})};
        
        \draw[->,thick,color=myblue] (violation.south) -- node[right,very near start]{NO} (LB.north);
        \draw[->,thick,color=myblue] (LB.south) -- (good2.north);
        \draw[->,thick,color=myblue] (good2.east) -- node[above,midway]{YES} (pruneB2.west);
        \draw[->,thick,color=myblue] (good2.south) -- node[right,midway]{NO} (SPMIO.north);
        \draw[->,thick,color=myblue] (SPMIO.south) -- (UB.north);
        
        \node[reg,draw] (Benders) at (1.5,13){Add Benders cuts\\Expand $\calU^0,\calV^0$};

        \node[test] (cut) at (-2.5,0){MIO cuts};
        \node[reg] (branch) at (-2.5,-10){Branch};
        
        \draw[->,thick] (termtest.east) -- node[above,midway]{YES}(term);
        \draw[->,thick] (termtest.south) -- node[right,midway]{NO} +(0,-1.5) -- +(-(10.5,-1.5) -- (LO.north);
        \draw[->,thick] (LO.south) -- (bad.north);
        \draw[->,thick] (bad.east) -- node[above,midway]{NO}(integer.west);
        \draw[->,thick] (integer.east) -- node[above,midway]{YES} (BSP.west);
        \draw[->,thick] (violation.north) -- node[right,very near end]{YES} (15,13) -- (Benders.east);
        \draw[->,thick] (Benders.west) -- (LO.east);
        
        \draw[->,thick] (integer.south) -- node[right,very near start]{NO} (cut.north);
        \draw[->,thick] (cut.west) -- node[above,very near start]{YES} +(-12.5,0) |- (LO.west);
        \draw[->,thick] (cut.south) -- node[right,very near start]{NO} (branch.north);
        \draw[->,color=myblue,thick] (UB.south) |- (branch.east);
        \draw[->,thick] (branch.west) -- +(-14,0) |- (termtest.west);

        \draw[->,dashed,thick] (integer.north) -- +(0,1.5) node[right,midway]{NO} -- +(8.5,1.5) -- +(8.5,-0.5);

        \draw[color=myred,dotted,ultra thick] (.5,-1.5) rectangle (10.2,10.5);
        \draw[color=myblue,dotted,ultra thick] (10.8,-10.5) rectangle (25,10.5);

    \end{tikzpicture}
    \caption{UB\&DD methodology to solve Problem~\eqref{OPT}.}
    \label{fig:UBDD}
\end{figure}

The final question involves solving the second-stage subproblem $\texttt{SP}(\bx)$ at each node where the UB\&BC component (blue part in Figure~\ref{fig:UBDD}) is visited. At that point, the algorithm has already solved the second-stage relaxation $\overline{\texttt{SP}}(\bx)$ via column generation. We can then derive a feasible second-stage solution by solving $\texttt{RSP}(\bx,\calA'_{sj})$ upon restoring integrality requirements. This is a common heuristic in column generation, which consists of solving a mixed-integer restricted master problem (i.e., the restricted Benders subproblem in our DD methodology) upon termination. The alternative would be to embed the column generation algorithm into a branch-and-price algorithm in each relevant node, which would come with huge computational costs. In fact, the column-generation heuristic is consistent with the approach from \cite{maheo2024unified}, who also solve $\texttt{SP}(\widehat\bx)$ via a heuristic in the branch-and-bound tree. Following \cite{maheo2024unified}, we then add a post-processing procedure to determine the optimal solution among all visited solutions $\widehat\bx$.

\subsubsection*{Proof of Theorem \ref{thm:UBDD}}

\cite{maheo2024unified} showed that UB\&BC converges in a finite number of iterations to an optimal solution of the two-stage stochastic mixed-integer optimization problem. In the UB\&BC scheme, nodes cannot be pruned by integrality, and so branching continues until all first-stage solution nodes have been explored or pruned (by bound or infeasibility). Thus, the algorithm stores a set of mixed-integer first-stage solutions, and evaluates them in a post-processing phase where integrality requirements are restored in the second-stage problem. As shown by \cite{maheo2024unified}, an optimal first-stage solution $\bx^*$ for Problem~\eqref{OPT} will be stored in the pool because:
\begin{enumerate}[label=(\roman*)]
    \item The node of an optimal first-stage solution $\bx^*$ will be visited. Assume that no optimal solution is in the pool, hence the upper bound satisfies $UB>\texttt{OPT}$. Any parent node of $\bx^*$ solves a relaxation, so it will not be pruned by infeasibility; moreover, it yield a solution at most equal to $\bc^\top \bx^* + \overline{\varPhi}(\bx^*) \leq \bc^\top \bx^* + \varPhi(\bx^*)=\texttt{OPT}<UB$, so it will not be pruned by bound.
    \item In the node itself, $\bx^*$ cannot be eliminated by Benders cuts \eqref{bendersoptcut}-\eqref{bendersfeascut} because the cuts are valid for the relaxation $(\texttt{MIO}-\texttt{LO})$, and thus for Problem~\eqref{OPT}.
\end{enumerate}
Combining that result with the exactness of our DD methodology for the partial relaxation $(\texttt{MIO}-\texttt{LO})$ from Theorem \ref{thm:exact}, we obtain that UB\&DD converges to an exact solution to Problem~\eqref{OPT}.

\subsection{Proof of Proposition~\ref{pricingproblem}}

Fix a reference trip $(\ell,t)\in\calL\times\calT_\ell$ and a scenario $s\in\calS$.

Let us consider a pair of checkpoints $(u,v)\in\Gamma_\ell$ and two load values $c_1\leq c_2\in\calC_\ell$. Let us define the load differential as $\nu=c_2-c_1$. By construction, the load component of the reduced cost satisfies:
\begin{equation}\label{proof:load}
    \Delta\psi^{u,v,\nu}_{\ell st}\geq\psi_{\ell,s,t,(u,c_1)} - \psi_{\ell,s,t,(v,c_2)}
\end{equation}

Consider a solution $\bff^*$, $\bw^*$, $\bxi^*$ of the pricing problem $\text{PP}_{\ell st}^{u,v,c_1,c_2}$. With a slight abuse of notation, we also refer to its optimal value as $\text{PP}_{\ell st}^{u,v,c_1,c_2}$. By construction, the solution $\bff^*$, $\bw^*$ defines a feasible solution to the problem defining $Z_{\ell st}^{u,v,\nu}$. Indeed, the load differential satisfies
    $$\sum_{m \in \calU^{uv}_{\ell st}} \sum_{p \in \calP_{m}} D_{ps} w^*_{mp} = \sum_{(m,q) \in \calH^{uv}_{\ell st}\ :\ f_{mq}=1}(\xi^*_q-\xi^*_m) = \xi^*_{(v,T_{\ell t}(v))}-\xi^*_{(u,T_{\ell t}(u))}=c_{(v,c_2)}-c_{(u,c_1)}=\nu,$$
where the first equality is induced by Equations~\eqref{eq:PPdiff}--\eqref{eq:PPdiff2}, the second equality is induced by telescoping the sum from Equation~\eqref{eq:PPflow}, the third equality is induced by Equation~\eqref{eq:PPse}, and the last inequality is by assumption.

Therefore, the routing component of the reduced cost expression satisfies:
\begin{equation}\label{proof:subpath}
    Z^{u,v,\nu}_{\ell st}\leq
    \sum_{m \in \calU^{uv}_{\ell st}} \sum_{p \in \calP_{m}} d_{mp} w^*_{mp}
\end{equation}

From Equations~\eqref{proof:load} and~\eqref{proof:subpath}, we obtain:
$$Z^{u,v,\nu}_{\ell st}-\Delta\psi^{u,v,\nu}_{\ell st}\leq\sum_{m \in \calU^{uv}_{\ell st}} \sum_{p \in \calP_{m}} d_{mp} w^*_{mp}+\psi_{\ell,s,t,(v,c_2)} - \psi_{\ell,s,t,(u,c_1)}=\text{PP}_{\ell st}^{u,v,c_1,c_2}$$

By taking the minimum over all arcs with a load differential $\nu$, we obtain:
\begin{equation}\label{proof:Prop3firstpart}
Z^{u,v,\nu}_{\ell st}-\Delta\psi^{u,v,\nu}_{\ell st}\leq\min_{c_1,c_2\in\calC_\ell:c_2-c_1=\nu}\text{PP}_{\ell st}^{u,v,c_1,c_2},\ \forall (u,v)\in \Gamma_\ell
\end{equation}

Vice versa, let us consider two checkpoints $(u,v)\in\Gamma_\ell$ and a load differential $\nu\in\calC_\ell$. Consider an arc $a^*\in\calA_{\ell st}$ that maximizes the load component of the reduced cost and a solution $\bff^*$, $\bw^*$ that minimizes the routing component for that load differential. Specifically, the arc $a^*\in\calA_{\ell st}$ defines a subpath that starts in checkpoint $u = k_{start(a^*)}\in\calI_\ell$ at time $T_{\ell t}(u)$ with vehicle load $c_{start(a^*)}$, that ends in checkpoint $v = k_{end(a^*)} \in \calI_\ell$ at time $T_{\ell t}(v)$ with load $c_{end(a^*)}=c_{start(a^*)}+\nu$, and that satisfies
$$\psi_{\ell,s,t,start(a^*)}-\psi_{\ell,s,t,end(a^*)}=\Delta\psi^{u,v,\nu}_{\ell st}$$

The solution $\bff^*$, $\bw^*$ satisfies Equations~\eqref{eq:PPw2f}--\eqref{eq:PPdomain} by construction. We then define a load variable $\xi_m$, keeping track of the load at node $m \in \calU^{u,v}_{\ell st}$. We initialize it with:
$$\xi_{(u,T_{\ell t}(u))}=c_{start(a^*)}.$$
Following solution $\bff^*$, $\bw^*$, we increase $\xi_m$ by $\sum_{p \in \calP_m} D_{ps} w^*_{mp}$ if we traverse $(m,q) \in \calH^{u,v}_{\ell st}$:
$$ \xi_q -\xi_m = \sum_{p \in \calP_m} D_{ps} w^*_{mp}, \qquad \forall (m,q) \in \calH^{u,v}_{\ell st}:f^*_{mq}=1.$$
The variables $\xi_m$ satisfy Equations~\eqref{eq:PPdiff}--\eqref{eq:PPdiff2} and \eqref{eq:PPdomain} by construction. By combining it with Equations~\eqref{eq:PPflow}, and telescoping the sum, we obtain:
\begin{align*}
    \xi_{(v,T_{\ell t}(v))}&=\sum_{m\in\calU^{uv}_{\ell st}:f^*_{m,(v,T_{\ell t}(v))}=1}\left(\xi_m+\sum_{p \in \calP_m} D_{ps} w^*_{mp}\right)\\
    &= \xi_{(u,T_{\ell t}(u))}+\sum_{(m,q)\in\calH^{uv}_{\ell st}:f^*_{mq}=1}\sum_{p \in \calP_m} D_{ps} w^*_{mp}\\
    &= c_{start(a^*)}+\nu\\
    &= c_{end(a^*)},
\end{align*}
where the third equality comes from the initialization $\xi_{(u,T_{\ell t}(u))}=c_{start(a^*)}$ and the constraint $\sum_{m \in \calU^{uv}_{\ell st}} \sum_{p \in \calP_{m}} D_{ps} w_{mp} = \nu$, and the last equality follows from the construction of $a^*\in\calA_{\ell st}$. Therefore, the variables $\xi_m$ also satisfy Equations~\eqref{eq:PPse}.

Therefore, solution $\bff^*$, $\bw^*$, $\xi$ defines a feasible solution for the pricing problem $\text{PP}_{\ell s t}^{u,v,c_{start(a^*)},c_{end(a^*)}}$, and we have:
$$Z^{u,v,\nu}_{\ell st}-\Delta\psi^{u,v,\nu}_{\ell st}=\sum_{p \in \calP_m} D_{ps} w^*_{mp}+\psi_{\ell,s,t,end(a^*)}-\psi_{\ell,s,t,start(a^*)}$$
Since, by construction, $c_{end(a^*)}-c_{start(a^*)}=\nu$, we obtain:
$$Z^{u,v,\nu}_{\ell st}-\Delta\psi^{u,v,\nu}_{\ell st}\geq\min_{c_1,c_2\in\calC_\ell:c_2-c_1=\nu}\text{PP}_{\ell st}^{u,v,c_1,c_2}.$$

This completes the proof that $Z^{u,v,\nu}_{\ell st}-\Delta\psi^{u,v,\nu}_{\ell st}$ is equal to the minimum reduced cost across all variables with load differential $\nu$:
\begin{equation*}
    Z^{u,v,\nu}_{\ell st}-\Delta\psi^{u,v,\nu}_{\ell st} = \min_{c_1,c_2\in\calC_\ell:c_2-c_1=\nu}\text{PP}_{\ell st}^{u,v,c_1,c_2}
\end{equation*}
This completes the proof.
\hfill\Halmos

\subsection{Proof of Remark \ref{R:ppRHS}}

Suppose that Equation~\eqref{eq:PPflow} is replaced with the following constraints in the pricing problem.
\begin{align}
    & \sum_{q \, :\, (m, q) \in \calH_{\ell st}^{uv}} f_{mq} - \sum_{q \, : \, (q, m) \in \calH_{\ell st}^{uv}} f_{qm} = \begin{cases}
        x_{\ell t} &\text{if } m = (u, T_{\ell t}(u)),\\
        x_{\ell t} &\text{if } m = (v, T_{\ell t}(v)),\\
        0 &\text{otherwise. } \\
    \end{cases} \qquad \forall m \in \calU_{\ell st}^{uv}\label{eq:wrongFlow}
\end{align}
With Equation~\eqref{eq:wrongFlow}, the resulting optimal solution to the pricing problem could not be used to construct a subpath, which by definition is a sequence of arcs connecting checkpoints $u$ and $v$.

Suppose toward a contradiction that Equation~\eqref{eq:PPOneP} is replaced with the following constraints in the pricing problem.
\begin{align}
    &\sum_{m \in \calU_{\ell st}^{uv} \, : \, p \in \calP_m} w_{mp} \leq z_{\ell pst} \qquad \forall p \in \calP \, : \, (\ell, t) \in \calM_p \label{eq:wrongLink}
\end{align}
Consider the subset of passengers $p \in \calP_{\ell st}$ for which $z_{\ell pst} = 0.$ Then the pricing problem only constructs arcs over the following set:
$$\calA_{\ell st}(\bz) := \{a \in \calA_{\ell st} \, : \, z_{p \ell st} = 1, \forall p \in \calP_{r(a)}\}. $$
As a result, the optimal dual solution $(\boldsymbol{\psi}, \bgamma)$ to the corresponding RMP would have unknown feasibility to the following constraints:
\begin{equation}
    \psi_{\ell s t n} - \psi_{\ell s t m} - \sum_{p \in \calP_a} \gamma_{\ell s t p} \leq g_{(n,m)} \qquad \forall (n,m)\in \calA_{\ell st} \setminus \calA_{\ell st}(\bz).
\end{equation}
Thus, the solution $(\boldsymbol{\psi}, \bgamma)$ is not necessarily in $\Lambda_{\ell st}$, and the corresponding optimality cut would be violated in the Benders decomposition algorithm.
\hfill\Halmos
\section{Experimental Setup}
\label{app:setup}

\begin{figure}[htbp!]
    \centering
   \subfloat[Demand and stations.]{\label{F:pickup}\includegraphics[height=0.4\textwidth]{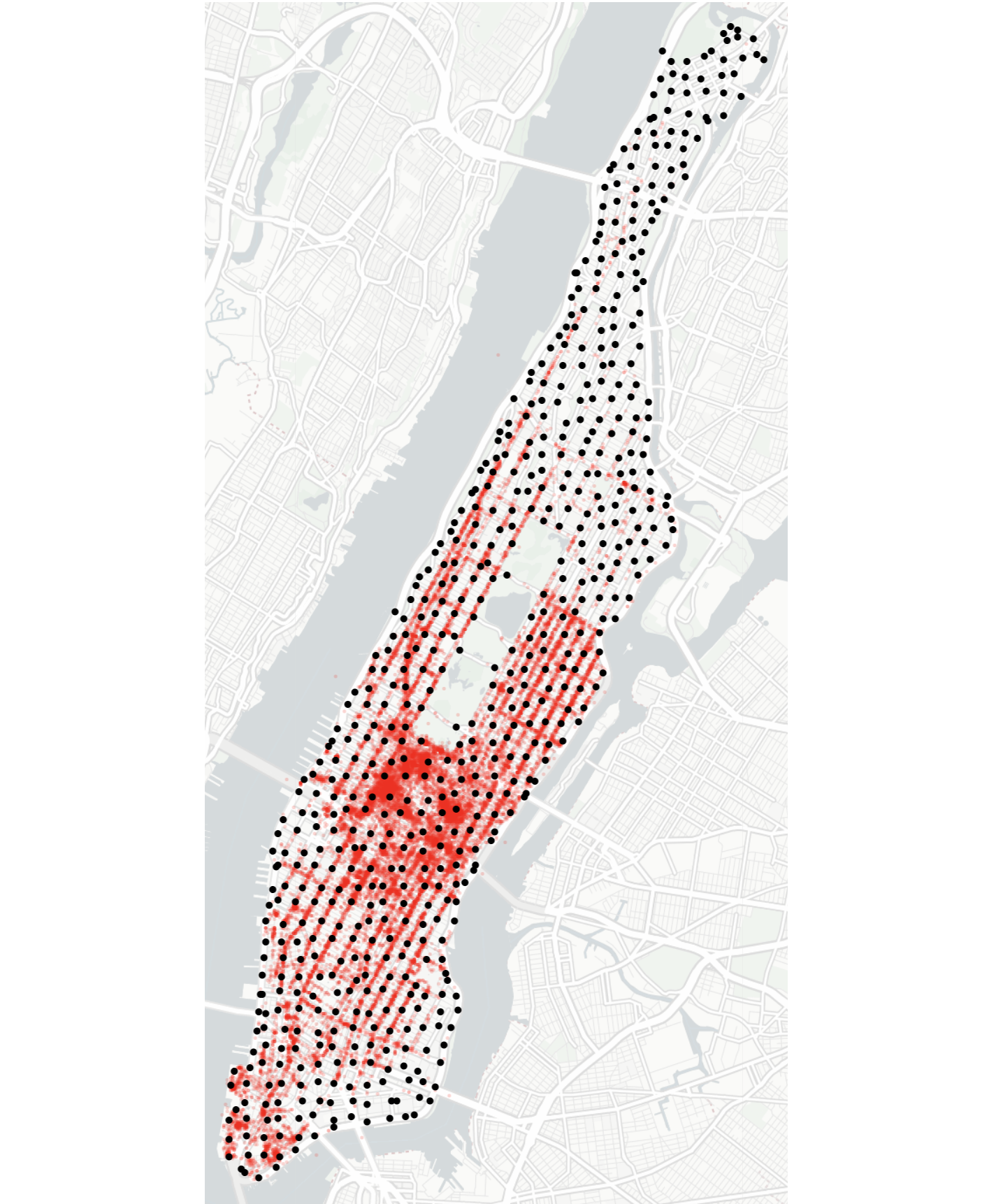}}
   \centering
   ~
   \subfloat[Reference lines.]{\label{F:lines}\includegraphics[height=0.4\textwidth]{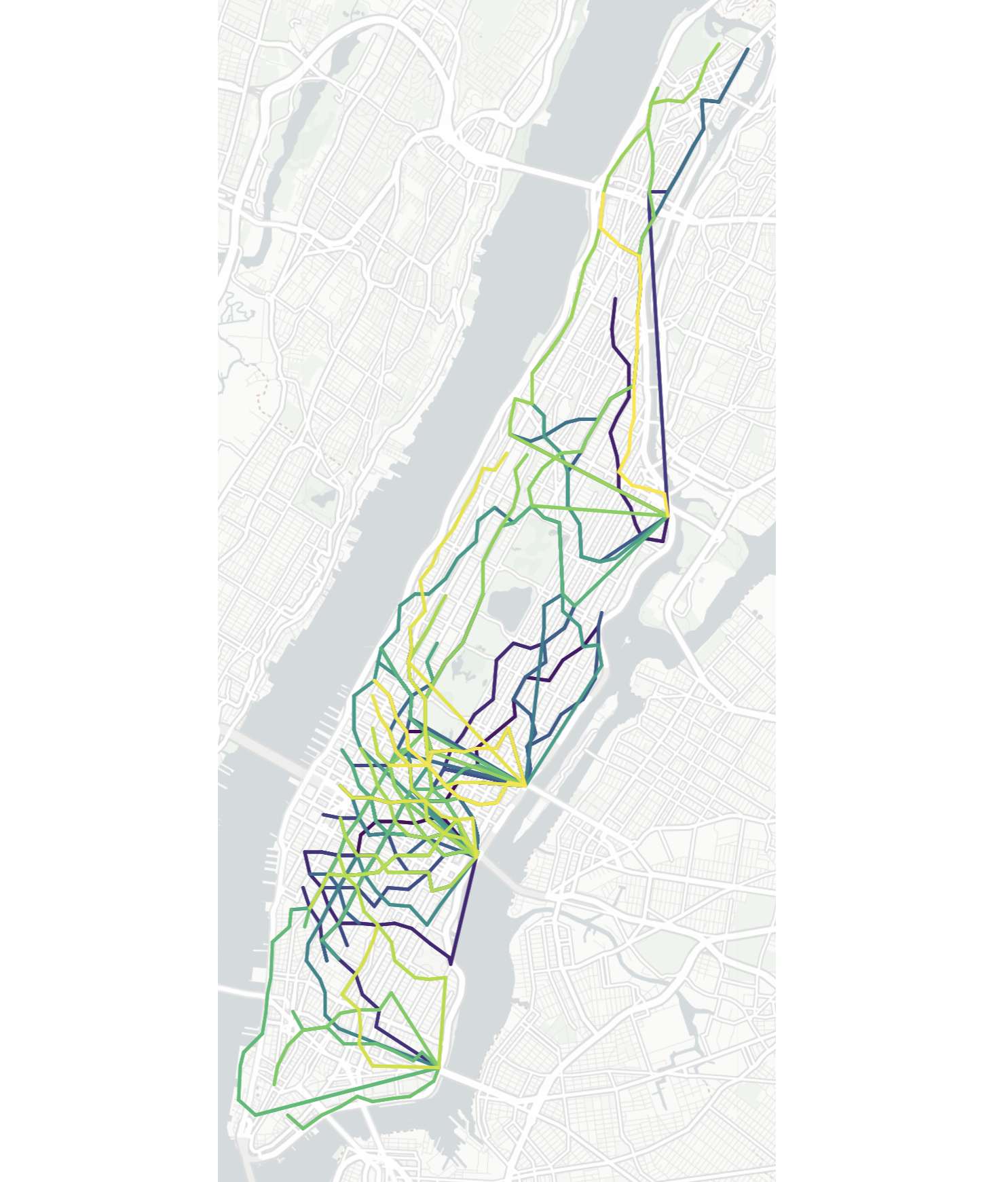}}
    \qquad \caption{Visualization of MiND-VRP inputs in Manhattan.}
    \label{F:ManhattanLS}
    \vspace{-12pt}
\end{figure}

In this appendix, we provide details on the generation of the model inputs (\ref{app:parameters}); in particular, we present a breadth-first search tree approach to define candidate reference lines (\ref{app:refline}). Figure~\ref{F:ManhattanLS} illustrates these inputs. We also detail our ride-sharing benchmarks (\ref{A:rideshare}).

\subsection{Model Inputs}\label{app:parameters}

We developed a real-world experimental setup in Manhattan, using data from the \cite{taxi}. We filtered trips to the airports during the morning rush (6--9 am), leading to up to 1,900 passenger request per instance (shown in Figure~\ref{F:pickup} of the paper). We defined a road network and travel times using data from Google Maps, OpenStreetMap, and \cite{speeds}. We considered pickup stations 300 meters apart, leading to 640 stations (also shown in Figure~\ref{F:pickup} of the paper). We assumed passengers could originate from any of the approximately 20,000 roadway intersections in Manhattan, and that they would walk from their origin to the closest station. We obtained the mapping and routing inputs from the {\tt fastest\_route} functionality in the OpenStreetMapX package in Julia \citep{osmx}. We calibrated travel time estimates to heavy Manhattan traffic using speed data from \cite{speeds}. We computed average speeds during the morning rush for each roadway type present in our Manhattan map (primary, secondary, tertiary, unclassified) and used these average speeds as input to the travel time estimation function, overriding default speeds provided by OpenStreetMapX.

Recall that our MiND-VRP experiments model a shuttle service from Manhattan to LaGuardia Airport with vehicles of capacity 10 to 20 passengers. Every trip leaves Manhattan and heads directly toward LaGuardia Airport via four possible exits: the Queensboro Bridge, the Williamsburg Bridge, the Kennedy Bridge, and the Midtown Tunnel. Travel times from each exit to LaGuardia were obtained via Google Maps estimates during the morning rush.

Table \ref{tab:parametersetting} reports the parameter values used in our computational experiments (Section~\ref{sec:comp}), and practical experiments for the MiND-VRP (Section~\ref{sec:practical}) and MiND-DAR (Appendix \ref{app:DARcase}).

\begin{table}[htbp!]
\caption{Details on input calibration for computational and practical analyses.}
\label{tab:parametersetting}
\centering
\begin{center}
\resizebox{1.\textwidth}{!}{%
\begin{tabular}{c|c|c|c|c}
\toprule
    Model component & Section \ref{sec:comp} value(s) & Section \ref{sec:practical} value(s) & \ref{app:DARcase} value(s) & \ref{app:TRcase} value(s)\\ \hline
    $\Omega$ & 210 meters & 420 meters & 250 meters & 210 meters\\  
    $\Psi$    & 10 minutes & 10 minutes  & 10 minutes & 10 minutes\\  
    $\Delta$  & 600 meters & 600 or 1,200 meters & 300 meters & 600 meters \\
    $\alpha$  & 5 minutes & 10 minutes  & 10 minutes & 10 minutes\\
    $C_\ell$  & 10 people & 10, 15, or 20 people  & 5, 10 or 20 people & 10 people \\
    $\calT_S$  & 30 seconds & 30 seconds & 30 seconds & 30 seconds\\  
    $\mu$  & 1 & 1 & 1 & 1 \\  
    $\lambda$  & 1 & 1 & 1 & 1 \\ 
    $\sigma$  & 1 & 1 & 1 & 1 \\
    $\delta$  & 1 & 1 & 1 & 1 \\
    $\kappa$  & 1 & 1 & 1 & 1 \\
    $M$  & 10,000 & 10,000 & 10,000 & 10,000 \\  
    $h_\ell$ & $T_{\ell t}(\calI_\ell^{(I_\ell)}) - T_{\ell t}(\calI_\ell^{(1)})$ & $T_{\ell t}(\calI_\ell^{(I_\ell)}) - T_{\ell t}(\calI_\ell^{(1)})$  & $T_{\ell t}(\calI_\ell^{(I_\ell)}) - T_{\ell t}(\calI_\ell^{(1)})$ & $T_{\ell t}(\calI_\ell^{(I_\ell)}) - T_{\ell t}(\calI_\ell^{(1)})$ \\
    $F$  & $|\calL|$ vehicles & 10 or 20 vehicles & 5 or 10 vehicles & $|\calL|$ vehicles \\
    $\calT_\ell$  & 15 minute intervals & 15 minute intervals & 15 minute intervals & 15 minute intervals \\
    $T_{\ell t}(\calI_\ell^{(i+1)}) - T_{\ell t}(\calI_\ell^{(i)})$  & 120\% of direct  & 120\% of direct & 110\% of direct &  120\% of direct \\ \bottomrule
\end{tabular}
}
\end{center}
\begin{tablenotes}\footnotesize
    \vspace{-6pt}
    \item $T_{\ell t}(\calI_\ell^{(i+1)}) - T_{\ell t}(\calI_\ell^{(i)})$: buffer time between arrival times at consecutive checkpoints $\calI_\ell^{(i)}$ and $\calI_\ell^{(i+1)}$.
    \vspace{-3pt}
    \item $\calT_\ell$: the frequency set is populated with departure times at evenly spaced intervals across the demand horizon.
    \vspace{-3pt}
    \item $\calT_S$: Time elapsed between consecutive discrete time units (between $t$ and $t+1$) in the discretized set $\calT^S$.
\end{tablenotes}
\vspace{-12pt}
\end{table}

\subsection{Reference Line Generation}\label{app:refline}

We describe the process of generating the set $\calL$ of candidate reference lines (shown in Figure~\ref{F:lines} of the paper). The procedure proceeds in three steps: (i) generating a comprehensive routing graph over Manhattan; (ii) using breadth-first search (BFS) trees to generate a very large set of candidate reference lines; and (iii) clustering and filtering to obtain a small but representative final set of candidate reference lines. We describe each step in detail below.

Note that our procedure to construct and optimize reference lines relies on a training set of demand data. This process avoids any bias moving from design to evaluation.

\subsubsection*{Manhattan routing graph.} 
\begin{figure}[htbp!]
    \centering
   \subfloat[Candidate checkpoints.]{\label{F:stopGridCC}\includegraphics[width=0.45\textwidth]{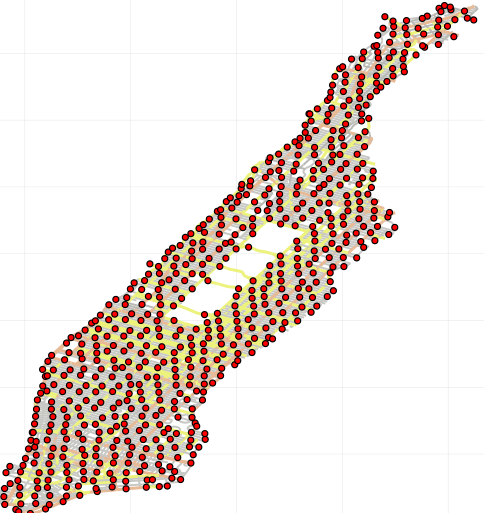}}
   \centering
   ~ ~
   \subfloat[Breadth-first search tree.]{\label{F:stopGridBFS}\includegraphics[width=0.45\textwidth]{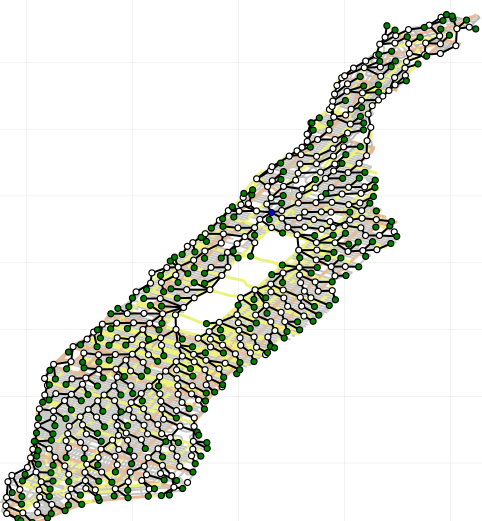}}
    \qquad \caption{Candidate checkpoints and BFS tree (blue: root node; green: leaves; white: intermediate nodes)}
    \label{F:stopGrid}
\end{figure}

We build a node set using discrete locations in Manhattan by generating a grid of GPS coordinates spanning Manhattan that were each 300 meters apart, and snapping each node to the closest road intersection. The outcome of this process is a list of candidate checkpoints $\calN$, shown in Figure \ref{F:stopGridCC}. We then build an edge set over this routing graph by connecting each node to its six closest neighbors according to their Euclidean distance. We used OpenStreetMapX to remove any edges that were impossible for a vehicle to traverse.

\subsubsection*{BFS trees.} To generate a large set of reference line candidates, we build BFS trees over the routing graph. Specifically, we let each node be the root of a BFS tree over the routing network (see Figure~\ref{F:stopGridBFS}). We then build reference line candidates over each BFS tree, by constructing node sequences from the root node to each leaf. Ultimately, we obtain tens of thousands of distinct candidate reference lines, across all BFS trees.

\subsubsection*{Clustering and filtering.} We first filter out many candidate lines that are illogical (e.g., indirect lines, very short or very long lines). We developed several metrics of line quality to systematically filter out low-quality options:
\begin{itemize}
    \item[--] {\it Minimum number of checkpoints.} Each line must visit a minimum of 10 stations.
    \item[--] {\it Low average and maximum detour.} For each checkpoint, we compute the relative detour as the ratio of the travel time from the checkpoint to the destination (LaGuardia) with the reference line and the corresponding direct travel time. The average detour across all checkpoints should not exceed 200\%, and the maximum detour should not exceed 250\%.
    \item[--] {\it Limited wrong-way travel.} To measure travel in the ``wrong direction,'' we measure the percentage of a reference line's checkpoints that are farther away from LaGuardia than their immediate predecessors.
    \item[--] {\it Demand coverage.} We assigned a popularity score to each checkpoint based on the frequency of trip requests with pickup locations close to that stop---in a training dataset. We filter out lines with a low average popularity score across its checkpoints.
\end{itemize}

Then, we remove redundancy over overlapping candidate lines, which is especially present among lines constructed from the same BFS tree. We measure the dissimilarity of two candidate lines as:
    $$\text{dissim}_{k\ell} = 1 - \frac{|\calI_k \cap \calI_\ell|}{\min\{I_k, I_\ell\}}.$$
When $\text{dissim}_{k\ell}=0$, lines $k$ and $\ell$ share as many stops as possible and are therefore substitutable. We collect these substitutable pairs into an undirected graph, and define an updated set of candidate lines $\calL'$ by computing a minimum vertex cover over that graph.

At this point, we are left with approximately 3,000 candidate lines in $\calL'$. In order to retain a tractable set of candidate lines in the optimization model, we cluster them into 100 representative and high-quality options. Specifically, we formulate a bi-objective clustering model to maximize medoid quality and diversity. Let $y_\ell \in \{0,1\}$ indicate whether line $l \in \calL'$ is selected in the final set $\calL$, and $x_{k\ell} \in \{0,1\}$ indicate whether line $k \in \calL'$ is assigned to the cluster with medoid line $l \in \calL'.$ We define a line-dependent parameter $q_\ell$ penalizing undesirable line characteristics based on the aforementioned metrics.

The clustering model maximizes line quality and minimizes the total dissimilarity among the line mapping (Equation~\eqref{eq:clusteringOBJ}), subject to partitioning constraints (Equation~\eqref{eq:clusteringPARTITION}), consistency constraints (Equation~\eqref{eq:clusteringCONSISTENCY}) and budget constraints (Equation~\eqref{eq:clusteringBUDGET}). We define the final reference line set as $\calL := \{l \in \calL' \, : \, y_\ell = 1\}.$
    \begin{align}
        \min \quad &\sum_{\ell \in \calL'} q_\ell y_\ell + \lambda \sum_{k \in \calL'}\sum_{\ell \in \calL'} \text{dissim}_{k\ell} x_{k\ell} \label{eq:clusteringOBJ}\\ 
        \text{s.t. } \quad &\sum_{k \in \calL'} x_{k\ell} = 1 &&\forall l \in \calL' \label{eq:clusteringPARTITION} \\ 
        &x_{k\ell} \leq |\calL'| y_\ell &&\forall k,l \in \calL'\label{eq:clusteringCONSISTENCY} \\ 
        &\sum_{\ell \in \calL'} y_\ell = 100 \label{eq:clusteringBUDGET} \\ 
        &\bx \in \{0,1\}^{\calL' \times \calL'} \\ 
        &\by \in \{0,1\}^{\calL'}
    \end{align}

We constructed three candidate line sets with 100 lines each by scaling the aforementioned quality measures with the following parameter settings:
\begin{align*}
    q_\ell^{\text{cluster}} &= 0 &&\forall \ell \in \calL' \\ 
    q_\ell^{\text{direct}} &= \frac{1}{3} \cdot \left(\text{maxDetour}_{\ell} + \text{meanDetour}_{\ell} + \text{wrongWay}_\ell \right) &&\forall \ell \in \calL;\\ 
    q_\ell^{\text{popular}} &= \text{popularity}_\ell &&\forall \ell \in \calL'
\end{align*}
Throughout the manuscript, we use $\bq^{\text{popular}}$ as the default to focus on the demand coverage objective, except for Section \ref{subsec:flex} on microtransit network design, in which we consider the line sets corresponding to all three quality measures.

\subsection{Ride-sharing Benchmark}\label{A:rideshare}

We build our ride-sharing benchmark using the cluster-then-route heuristic from \cite{by2021}, originally built to generate paratransit itineraries with up to 4 passengers per vehicle. Their approach was itself based on the maximum weighted matching over a shareability network from \cite{santi2014a}. To extend the approach from two- to four-passenger trips, \cite{by2021} first created a set of passenger pairs and then approximated the shareability network over two-passenger trips. We adopt a similar approach except that, instead of requiring all requests to be served, we maximize the number of served requests and then minimize travel times.

\subsubsection*{Single-occupancy ride-sharing.}
With single-occupancy vehicles, the clustering step is unnecessary. We simply apply the routing step from \cite{by2021} over the request set.

\subsubsection*{Two-occupancy ride-sharing.}
We build a pair-wise shareability network that encodes the pairs of trips that can share a vehicle. Let $t_i$ denote the requested pickup time of request $i$, $T_i$ the direct travel time of request $i$, and $tt(x, y)$ the travel time from location $x$ to location $y$.
\begin{itemize}
    \item[--] If $t_j \leq t_i + T_i + \Psi$, then trip $j$ can be picked up before trip $i$ is dropped off;
    \item[--] if $t_i \leq t_j + T_j + \Psi$, then trip $i$ can be picked up before trip $j$ is dropped off; and
    \item[--] otherwise, trips $i$ and $j$ cannot be shared.
\end{itemize}
Then we determine whether there exists pickup times for trips $i$ and $j$ (in that order) such that no request is picked up early and each pickup is within $\Psi$ of their requested times. The following conditions must hold, where $x$ denotes the pickup time of trip $i$:
\begin{align*}
&t_i \leq x \leq t_i + \Psi &&\text{Request $i$ has tolerable wait time}\\
&t_j \leq x + tt(o_i, o_j) \leq t_j + \Psi &&\text{Request $j$ has tolerable wait time}
\end{align*}
which reduces to finding some $x$ such that:
$$x \in [\max\{t_i, t_j - tt(o_i, o_j)\}, \min\{t_i + \Psi, t_j + \Psi - tt(o_i, o_j)\}]. $$
The two requests can also share a vehicle if the symmetric problem holds, corresponding to the instance where trip $j$ is picked up first:
$$x \in [\max\{t_j, t_i - tt(o_j, o_i)\}, \min\{t_j + \Psi, t_i + \Psi - tt(o_j, o_i)\}] $$
\cite{by2021} impose a maximum delay limit, but we remove this restriction to enable more ride-pooling.
Finally, we determine the travel time associated with each version of the trip.
\begin{align*}
    c_{i \to j} &= tt(o_i, o_j) + T_j\\
    c_{j \to i} &= tt(o_j, o_i) + T_i
\end{align*}
If $c_{i \to j} \leq T_i + T_j$ or $c_{j \to i} \leq T_j + T_i$, then the shared trip is more efficient than serving the two requests separately. If both are efficient, then we select the best option.

The shared trips satisfying the above conditions are added to the VSN with cost $T_i + T_j - \min\{c_{i \to j}, c_{j \to i}\}$ to reflect the cost savings of pooling the requests. We solve a maximum weighted matching problem to pair requests into capacity-2 trips, with some requests potentially still served in isolation if they are not matched to any other request. We first maximize the number of served requests, and then we minimize the total travel time, subject to the fleet size limit. 

\subsubsection*{Four-occupancy ride-sharing.}

We build a new shareability network that combines trip pairs from the pair-wise shareability network. For the MiND-VRP, we solve a simple vehicle routing problem for each candidate set of four trips to find the best sequence of stops within that set, while ensuring that no one is picked up earlier than their requested times and that none of passengers' wait times exceeds limit $\Psi$. For the MiND-DAR, we solve a simple dial-a-ride problem for each candidate set of four trips, which also includes precedence constraints so that each pickup occurs before the corresponding dropoff. We note that the optimal pooling configuration of two request pairs could potentially be to serve all four requests together, or to pool only a subset of these requests and serve the remaining requests separately. We proceed as in the two-occupancy case, solving a maximum weighted matching problem over the VSN to determine final trips, and then performing an identical itinerary generation procedure to the one described previously.
\section{Extension to the dial-a-ride setting (MiND-DAR)}
\label{app:DAR}

\subsection{Modeling extension}
\label{subsec:MiND-DAR}

In the dial-a-ride setting, each passenger request $p\in\calP$ is associated with an origin $o(p)$ and a destination $d(p)$. The first-stage formulation remains unchanged, except that the set $\calM_p$ is re-defined as the set of reference lines that cover both the origin and the destination of request $p\in\calP$. In the second stage, we define the sets $\calP_r^+$ and $\calP_r^-$ ($\calP_r= \calP_r^+ \cup \calP_r^-)$ as the passenger requests that are picked up and dropped off, respectively, by subpath $r \in \calR_{\ell st}$ for $(\ell, t) \in \calL \times \calT_\ell, s \in \calS$.

Level of service involves similar measures of passenger disutility. A subpath $r \in \calR_{\ell st}$ is associated with a walking cost both for pickups (from the origin of passenger $p\in\calP_r^+$ to the pickup location) and for dropoffs (from the dropoff location to the destination of passenger $p\in\calP_r^-$); with a waiting cost for pickups; and with a delay cost for dropoffs. To capture detour costs, we denote by $T_{r(a)p}^+$ (resp. $T_{r(a)p}^-$) the pickup (resp. dropoff) time of passenger $p$ on arc $a \in \calA_{\ell st}$ such that $p \in \calP_{r(a)}^+$ (resp. $p \in \calP_{r(a)}^-$). The arc costs $g_a$ are re-derived as follows.
\begin{align}\label{eq:arccostDAR}
    g_a^{DAR} = \begin{cases}
        \sum_{p \in \calP_{r(a)}^+} D_{ps} \left(\lambda \tau_{r(a)p}^{walk} + \mu \tau_{r(a)p}^{wait} - \sigma \frac{T^+_{r(a)p}}{\tau^{dir}_p} - M \right) + \\
        \qquad \qquad\sum_{p \in \calP_{r(a)}^-} D_{ps} \left(\lambda \tau_{r(a)p}^{walk} + \sigma \frac{T^-_{r(a)p}}{\tau_{p}^{dir}} +\delta \frac{\tau^{late}_{\ell r(a)p}}{\tau_p^{dir}} + \frac{\delta}{2} \frac{\tau^{early}_{\ell r(a)p}}{\tau_p^{dir}}\right) &\forall a \in \bigcup_{r \in \calR_{\ell st}} \calA_r, \\
        0 &\forall a \in \calA_{\ell st}^v.
    \end{cases}
\end{align}

The MiND-DAR is then formulated as follows. The only difference with the MiND-VRP is the additional constraint ensuring that a passenger who is picked up needs to be dropped off (Equation~\eqref{eq:DAR}). Note that the precedence constraint is captured by the set $\calM_p$ and therefore does not need to be enforced explicitly in the MiND-DAR formulation.
\begin{align}
    \min   \quad&\sum_{\ell \in \calL} \sum_{t \in \calT_\ell} h_\ell x_{\ell t} +  \sum_{s \in \calS} \pi_s \left(\sum_{\ell \in \calL }\sum_{t \in \calT_{\ell}}\sum_{a \in \calA_{\ell st}} g_a^{DAR} y_a\right)\label{DARobj}\\
    \st\quad & \text{First-stage constraints: Equations~\eqref{eq:budget}--\eqref{eq:load2}}\\
    & \text{Second-stage constraints: Equations~\eqref{eq:2Sflow}--\eqref{eq:2Sy2z}}\\
    & \sum_{a \in \calA_{\ell st} \, : \, p \in \calP^+_{r(a)}} y_a - \sum_{a \in \calA_{\ell st} \, : \, p \in \calP^-_{r(a)}} y_a = 0 \quad \forall s \in \calS, p \in \calP, (\ell,t) \in \calM_p\label{eq:DAR}\\
    &\bx, \by, \bz \text{ binary} \label{DARdomain}
\end{align}

\subsection{Algorithmic extension}
\label{app:DARalg}

\paragraph{Benders decomposition.}

For a reference trip $(\ell, t)$ and a scenario $s$, let $\zeta_{\ell s t p}$ denote the dual variable associated to the new consistency constraint between pickup and drop-off decisions (Equation~\eqref{eq:DAR}). The Benders dual subproblem becomes:
\begin{align}
    \max \quad &x_{\ell t} \cdot (\psi_{\ell,s,t,u_{\ell st}} - \psi_{\ell,s,t,v_{\ell st}}) - 
    \sum_{p \in \calP \, : \,(\ell,t) \in \mathcal{M}_p}
    z_{p\ell st} \cdot \gamma_{\ell s t p}  \\
    \text{s.t. } \quad &\psi_{\ell,s,t,n} - \psi_{\ell,s,t,m} -  \sum_{p \in \calP^+_{r(a)} } (\gamma_{\ell s t p}  - \zeta_{\ell s t p}) - \sum_{p \in \calP^-_{r(a)}} \zeta_{\ell s t p} \leq g^{DAR}_a  && \forall  a=(n,m) \in \mathcal{A}_{\ell st} \label{DARBenders:dualArc} \\
    &\psi_{\ell,s,t,i} \in \mathbb{R} && \forall i \in \mathcal{V}_{\ell st}  \\
    &\gamma_{\ell s t p} \geq 0 && \forall p \in \calP \, : \,(\ell ,t) \in \mathcal{M}_p \\
    &\zeta_{\ell s t p} \in \mathbb{R} && \forall p \in \calP \, : \,(\ell ,t) \in \mathcal{M}_p 
\end{align}

Note that the new dual variables do not appear in the dual objective function, so the Benders optimality cut remains unchanged.

\paragraph{Column generation.}

The restricted Benders subproblem is still obtained by restricting the decisions to a subset of arc-based variables in $\calA'_{\ell st}$:
\begin{align}
    \text{RBSP}(\calA'_{\ell st}, \bx, \bz) \quad\min_{\by\geq\bo} &\quad \sum_{a \in \calA'_{\ell st}} g_a^{DAR} y_a \label{RMPDAR:obj} \\
    \st\quad&\sum_{m:(n,m) \in \calA'_{\ell st}} 
    y_{(n,m)} - \sum_{m:(m, n) \in \calA'_{\ell st}} 
    y_{(m, n)} = \begin{cases}
    x_{\ell t} &\text{if } n = u_{\ell st}, \\
    -x_{\ell t} &\text{if } n = v_{\ell st}, \\
    0 &\text{ otherwise,}
    \end{cases} \ \forall n \in \calV_{\ell st}\label{RMPDAR:2Sflow} \\
    & \sum_{a \in \calA'_{\ell st} \, : \, p \in \calP^+_{r(a)}} y_a - \sum_{a \in \calA'_{\ell st} \, : \, p \in \calP^-_{r(a)}} y_a = 0 \quad \forall p \in \calP, (\ell,t) \in \calM_p\label{RMPDAR:DAR}\\
    &\sum_{a \in \calA'_{\ell st} \, : \, p \in \calP_{r(a)}} y_a \leq z_{p\ell st} \quad \forall p \in \calP: (\ell,t) \in \calM_p \label{RMPDAR:2Sy2z}
\end{align}

In the pricing problem, we split the level-of-service parameter $d_{mp}$ into $d^+_{mp}$ and $d^-_{mp}$, corresponding to the level-of-service components associated with pickups and dropoffs, respectively. Following Section~\ref{subsec:MiND-DAR}, we denote by $\calP^+_m$ (resp. $\calP^-_m$) the set of passengers that can be picked up (resp. dropped off) and by $T^{+}_{mp}$ (resp. $T^{-}_{mp}$) the pickup time (resp. dropoff time) of passenger $p\in\calP^+_m$ (resp. $p\in\calP^-_m$). We then define:
\begin{align*}
    d^+_{mp}&= D_{ps} \left( \lambda \tau^{\text{walk}}_{mp} + \mu \tau^{\text{wait}}_{mp}-\sigma\frac{T^{+}_{mp}}{\tau^{\text{dir}}_p} - M\right) + \gamma_{\ell s t p} - \zeta_{\ell s t p},&& \quad \forall m \in \calU^{uv}_{\ell st}, p \in\calP^+_{m}\\
    d^-_{mp}&= D_{ps} \left(\frac{\delta\tau^{\text{late}}_{mp}+\frac{\delta}{2}\tau^{\text{early}}_{mp} + \sigma T^{-}_{mp}}{\tau^{\text{dir}}_p} + \lambda \tau^{\text{walk}}_{mp} \right) + \zeta_{\ell s t p},&& \quad \forall m \in \calU^{uv}_{\ell st}, p \in\calP^-_{m}
\end{align*}

Similarly, we define the following decision variables to split pickups and dropoffs:
\begin{align*}
    f_{mq}&=\begin{cases}1&\text{if arc $(m,q) \in \calH^{uv}_{\ell st}$ is traversed in the time-expanded road segment network,}\\0&\text{otherwise.}\end{cases}\\
    w^+_{mp}&=\begin{cases}1&\text{if passenger $p \in \calP^+_{m}$ is picked up at node $m \in \calU^{uv}_{\ell st}$,}\\0&\text{otherwise.}\end{cases}\\
    w^-_{mp}&=\begin{cases}1&\text{if passenger $p \in \calP^-_{m}$ is dropped off at node $m \in \calU^{uv}_{\ell st}$,}\\0&\text{otherwise.}\end{cases}\\
    \xi_m&=\ \text{vehicle load in node $m \in \calU^{uv}_{\ell st}$}
\end{align*}

The pricing problem is then formulated as follows. Equation~\eqref{eq:PPobj+} minimizes the reduced cost. Constraints~\eqref{eq:PPseDAR}--\eqref{eq:PPdiff2DAR} define the load at each node based on the pickups and dropoffs. Constraints~\eqref{eq:PPw2f+} and~\eqref{eq:PPw2f-} ensure that a passenger can only be picked up or dropped off in a node that is visited. Constraints~\eqref{eq:PPOneP+} and~\eqref{eq:PPOneP-} guarantee that a passenger is picked up and dropped off at most once, respectively. Constraints~\eqref{eq:PPflow} apply flow balance in the time-expanded road segment network. The remaining constraints enforce binary requirements.
\begin{align}
    \min  \quad & 
    \sum_{m \in \calU^{uv}_{\ell st}} \left(\sum_{p \in \calP^+_{m}} d^+_{mp} w^+_{mp} +\sum_{p \in \calP^-_{m}} d^-_{mp} w^-_{mp}\right) + \psi_{\ell,s,t,end(a)}- \psi_{\ell,s,t,start(a)} \label{eq:PPobj+}\\
    \st\quad
    &\xi_{(u,T_{\ell t}(u))}=c_{(u,c_1)},\ \xi_{(v,T_{\ell t}(v))}=c_{(v,c_2)} \label{eq:PPseDAR}\\
    & \xi_q-\xi_m \leq \left(\sum_{p \in \calP^+_{m}} D_{ps} w^+_{mp} - \sum_{p \in \calP^-_{m}} D_{ps} w^-_{mp}\right) + C_\ell (1 - f_{mq}),\quad\forall(m,q)\in\calH^{uv}_{\ell st} \label{eq:PPdiffDAR}\\
    & \xi_q-\xi_m \geq \left(\sum_{p \in \calP^+_{m}} D_{ps} w^+_{mp} - \sum_{p \in \calP^-_{m}} D_{ps} w^-_{mp}\right) - C_\ell (1 - f_{mq}),\quad\forall(m,q)\in\calH^{uv}_{\ell st} \label{eq:PPdiff2DAR}\\
    &w^+_{mp} \leq \sum_{q:(m,q) \in \calH^{uv}_{\ell st}} f_{mq} \quad \forall m \in \calU^{uv}_{\ell st},\ \forall p \in \calP^+_{m} \label{eq:PPw2f+}\\
    &w^-_{mp} \leq \sum_{q:(m,q) \in \calH^{uv}_{\ell st}} f_{mq} \quad \forall m \in \calU^{uv}_{\ell st},\ \forall p \in \calP^-_{m} \label{eq:PPw2f-}\\
    &\sum_{m \in \calU^{uv}_{\ell st} \, : \, p \in \calP^+_m} w^+_{mp} \leq 1 \quad \forall p \in \calP \label{eq:PPOneP+} \\
    &\sum_{m \in \calU^{uv}_{\ell st} \, : \, p \in \calP^-_m} w^-_{mp} \leq 1 \quad \forall p \in \calP \label{eq:PPOneP-} \\
    &\sum_{q:(m,q) \in \calH^{uv}_{\ell st}} f_{mq} - \sum_{q:(q,m) \in \calH^{uv}_{\ell st}} f_{qm} = \begin{cases}
    1 &\text{if } m = (u,T_{\ell t}(u)), \\
    -1 &\text{if } m = (v,T_{\ell t}(v)), \\
    0 &\text{ otherwise.}
    \end{cases}\quad \forall m \in \calU^{uv}_{\ell st} \label{eq:PPflow+}\\
    &f_{mq} \in \{0,1\} \quad \forall (m,q) \in \calH^{uv}_{\ell st} \label{eq:PPf+}\\
    &w^+_{mp} \in \{0,1\} \quad \forall m \in \calU^{uv}_{\ell st}, p \in \calP^+_{m} \label{eq:PPw+}\\
    &w^-_{mp} \in \{0,1\} \quad \forall m \in \calU^{uv}_{\ell st}, p \in \calP^-_{m} \label{eq:PPw-}
\end{align}

\paragraph{Label setting algorithm.}

To distinguish pickups and dropoffs, we extend the label-setting algorithm from a two-dimensional to a three-dimensional state space. Dropoffs are treated the same way as pickups; for instance, the state transition includes checking all passenger combinations for pickups and all passenger combinations for dropoffs. This extension has two major implications that increase the computational requirements in the pricing problem. First, the dominance rule requires the dominating state to have the same set of pickups \textit{and} the same set of dropoffs as the dominated state. Second, the set of load differential needs to be extended from $\{0,1,\cdots,C_\ell\}$ to $\{-C_\ell,\cdots,-1,0,1,\cdots,C_\ell\}$. Nonetheless, our results show that our methodology scales to meaningful practical instances of the MiND-DAR model in Manhattan, with up to 10 candidate lines, hundreds of stations, thousands of passenger requests and 5 demand scenarios---resulting in over 60,000 first-stage variables and 700 second-stage problems.

\subsection{Experimental results}
\label{app:DARcase}

We construct a case study setting in Midtown Manhattan, with 10 candidate lines traveling West to East from the 11th to the 1st avenue along every other street between 36th and 54th. Each line contains a checkpoint at every other avenue, and each street-avenue intersection defines a station---leading to a total of 168 stations. We calibrate demand inputs by collecting all West-to-East requests in Midtown Manhattan during the morning rush from 6 to 9 am, amounting to over 3,000 passenger requests. We set up one-hour, two-hour and three-hour instances (from 6 to 7 am, 6 to 8 am, and 6 to 9 am, respectively). For each one, we run the deviated fixed-route microtransit as well as the fixed-line transit benchmark and ride-sharing benchmarks with single-occupancy, two-occupancy and four-occupancy vehicles (see~\ref{A:rideshare}). We consider five demand scenarios. Again, for apples-to-apples comparison, we group results by total seating capacity (e.g., 10 transit/microtransit vehicles of capacity 10, ride-sharing with 100/50/25 vehicles of capacity of 1/2/4), and perform an out-of-sample assessment corresponding to five new weekdays.

We evaluate the system-wide performance of all optimized transportation modes in Table~\ref{T:losDAR}, broken down into level of service (demand coverage and average walking time, waiting time, delay and detour), vehicle utilization (passengers served divided by vehicle capacity), and distance traveled (internal distance for served passengers plus external distance for unserved passengers). These results confirm and extend all takeaways from the MiND-VRP (Table~\ref{T:los} and Figure~\ref{F:coverageVsFootprint}).
 
\begin{table}[h!]
\renewcommand{\arraystretch}{1.0}
\caption{Average performance of fixed-route transit, microtransit, and ride-sharing in a dial-a-ride setting.}
    \label{T:losDAR}
\resizebox{\textwidth}{!}{%
\begin{tabular}{lll rrrrrrrrrr}
\toprule[1pt]
\multicolumn{3}{c}{Setting}   &
  \multicolumn{5}{c}{Average level of service} &
  \multicolumn{2}{c}{Vehicle utilization} &
  \multicolumn{3}{c}{Distance traveled (km)} \\ \cmidrule(lr){1-3}\cmidrule(lr){4-8}\cmidrule(lr){9-10}\cmidrule(lr){11-13}
   Horizon & Capacity & Mode & Coverage & Walk & Wait & Delay & Detour & Absolute & Relative & Internal & External & Total \\ \hline
1 hour  & 50  & Transit & 6.7\%   & 1.68 & 3.57 & 1.45 & 200\% & 3.04             & 40.5\%             & 60            & 881         & 941          \\
    &     & Microtransit         & 16.5\%  & 1.92 & 3.71 & 1.70 & 144\% & 6.30             & 81.4\%             & 94            & 798         & 893          \\
    &     & RS Cap. 4    & 23.8\%  & 0.00 & 3.63 & 6.74 & 183\% & 3.81             & 95.3\%             & 293           & 1,082       & 1,375        \\
    &     & RS Cap. 2    & 36.6\%  & 0.00 & 3.56 & 4.30 & 120\% & 1.95             & 97.3\%             & 567           & 942         & 1,509        \\
    &     & RS Cap. 1    & 53.3\%  & 0.00 & 2.38 & 2.38 & 100\% & 1.00             & 100.0\%            & 1,040         & 754         & 1,793        \\ \cline{2-13} 
    & 100 & Transit   & 7.6\%   & 1.68 & 3.51 & 1.64 & 201\% & 3.18             & 22.7\%             & 63            & 881         & 944          \\
    &     & Microtransit         & 22.8\%  & 2.04 & 3.73 & 1.68 & 154\% & 8.29             & 57.4\%             & 96            & 798         & 894          \\
    &     & RS Cap. 4    & 42.4\%  & 0.00 & 3.48 & 6.92 & 190\% & 3.81             & 95.3\%             & 541           & 833         & 1,374        \\
    &     & RS Cap. 2    & 63.3\%  & 0.00 & 3.71 & 4.57 & 122\% & 1.94             & 97.2\%             & 1,057         & 569         & 1,626        \\
    &     & RS Cap. 1    & 85.2\%  & 0.00 & 2.47 & 2.47 & 100\% & 1.00             & 100.0\%            & 1,876         & 268         & 2,144        \\ \cline{2-13} 
    & 200 & Transit   & 8.7\%   & 1.64 & 3.54 & 2.07 & 199\% & 2.74             & 13.7\%             & 84            & 906         & 990          \\
    &     & Microtransit         & 26.8\%  & 2.05 & 3.80 & 1.65 & 156\% & 7.21             & 36.1\%             & 127           & 755         & 882          \\
    &     & RS Cap. 4    & 70.7\%  & 0.00 & 3.55 & 7.24 & 193\% & 3.86             & 96.6\%             & 965           & 421         & 1,386        \\
    &     & RS Cap. 2    & 94.9\%  & 0.00 & 4.01 & 5.05 & 125\% & 1.95             & 97.3\%             & 1,708         & 76          & 1,784        \\
    &     & RS Cap. 1    & 100.0\% & 0.00 & 2.66 & 2.66 & 100\% & 1.00             & 100.0\%            & 2,131         & 0           & 2,131        \\ \hline
2 hours & 50  & Transit   & 6.9\%   & 1.87 & 3.76 & 1.05 & 195\% & 3.13             & 52.1\%             & 153           & 2,403       & 2,557        \\
    &     & Microtransit         & 15.2\%  & 1.83 & 3.72 & 1.69 & 162\% & 6.62             & 101.5\%            & 220           & 2,196       & 2,416        \\
    &     & RS Cap. 4    & 19.2\%  & 0.00 & 3.88 & 6.77 & 182\% & 3.82             & 95.4\%             & 592           & 3,039       & 3,631        \\
    &     & RS Cap. 2    & 29.7\%  & 0.00 & 3.69 & 4.38 & 120\% & 1.96             & 98.1\%             & 1,154         & 2,770       & 3,924        \\
    &     & RS Cap. 1    & 43.9\%  & 0.00 & 2.46 & 2.46 & 100\% & 1.00             & 100.0\%            & 2,184         & 2,337       & 4,521        \\ \cline{2-13} 
    & 100 & Transit   & 5.9\%   & 0.94 & 3.77 & 1.32 & 196\% & 3.63             & 22.2\%             & 120           & 2,428       & 2,548        \\
    &     & Microtransit         & 19.9\%  & 2.12 & 3.95 & 1.65 & 155\% & 11.05            & 66.8\%             & 175           & 2,071       & 2,246        \\
    &     & RS Cap. 4    & 33.9\%  & 0.00 & 3.84 & 7.11 & 189\% & 3.84             & 96.0\%             & 1,113         & 2,535       & 3,649        \\
    &     & RS Cap. 2    & 53.0\%  & 0.00 & 3.78 & 4.61 & 123\% & 1.96             & 98.2\%             & 2,252         & 1,948       & 4,200        \\
    &     & RS Cap. 1    & 74.3\%  & 0.00 & 2.53 & 2.53 & 100\% & 1.00             & 100.0\%            & 4,207         & 1,214       & 5,421        \\ \cline{2-13} 
    & 200 & Transit   & 8.5\%   & 0.91 & 3.70 & 1.67 & 198\% & 3.23             & 16.2\%             & 187           & 2,331       & 2,517        \\
    &     & Microtransit         & 27.3\%  & 2.09 & 3.88 & 1.65 & 152\% & 9.42             & 47.1\%             & 258           & 1,848       & 2,106        \\
    &     & RS Cap. 4    & 59.1\%  & 0.00 & 3.80 & 7.41 & 194\% & 3.85             & 96.3\%             & 2,114         & 1,611       & 3,724        \\
    &     & RS Cap. 2    & 86.1\%  & 0.00 & 4.03 & 5.07 & 126\% & 1.96             & 97.9\%             & 4,117         & 635         & 4,752        \\
    &     & RS Cap. 1    & 99.5\%  & 0.00 & 2.76 & 2.76 & 100\% & 1.00             & 100.0\%            & 6,253         & 30          & 6,283        \\ \hline
3 hours & 50  & Transit   & 6.7\%   & 1.88 & 3.77 & 1.13 & 198\% & 4.53  & 57.0\%  & 192    & 3,895 & 4,087  \\
    &     & Microtransit         & 13.5\%  & 1.74 & 3.69 & 1.74 & 154\% & 8.25  & 104.0\% & 268    & 3,600 & 3,868  \\
    &     & RS Cap. 4    & 16.8\%  & 0.00 & 3.91 & 6.78 & 182\% & 3.80  & 95.0\%  & 887    & 5,444 & 6,331  \\
    &     & RS Cap. 2    & 26.9\%  & 0.00 & 3.80 & 4.47 & 120\% & 1.97  & 98.3\%  & 1,750  & 4,996 & 6,746  \\
    &     & RS Cap. 1    & 39.7\%  & 0.00 & 2.50 & 2.50 & 100\% & 1.00  & 100.0\% & 3,348  & 4,314 & 7,662  \\ \cline{2-13} 
    & 100 & Transit   & 7.6\%   & 1.96 & 3.69 & 0.85 & 201\% & 5.09  & 30.5\%  & 174    & 3,866 & 4,040  \\
    &     & Microtransit         & 20.7\%  & 1.99 & 3.83 & 1.74 & 132\% & 14.01 & 80.0\%  & 229    & 3,298 & 3,527  \\
    &     & RS Cap. 4    & 30.0\%  & 0.00 & 3.87 & 7.09 & 189\% & 3.81  & 95.2\%  & 1,691  & 4,670 & 6,361  \\
    &     & RS Cap. 2    & 48.2\%  & 0.00 & 3.89 & 4.69 & 122\% & 1.96  & 98.0\%  & 3,451  & 3,753 & 7,204  \\
    &     & RS Cap. 1    & 68.4\%  & 0.00 & 2.56 & 2.56 & 100\% & 1.00  & 100.0\% & 6,558  & 2,520 & 9,077  \\ \cline{2-13} 
    & 200 & Transit   & 8.4\%   & 1.85 & 3.83 & 1.93 & 195\% & 3.66  & 18.3\%  & 262    & 3,678 & 3,940  \\
    &     & Microtransit         & 26.1\%  & 2.13 & 3.99 & 1.69 & 127\% & 10.37 & 51.8\%  & 365    & 2,991 & 3,356  \\
    &     & RS Cap. 4    & 53.0\%  & 0.00 & 3.82 & 7.41 & 195\% & 3.84  & 95.9\%  & 3,264  & 3,235 & 6,499 \\
    &     & RS Cap. 2    & 80.3\%  & 0.00 & 4.03 & 5.05 & 126\% & 1.96  & 97.9\%  & 6,549  & 1,561 & 8,110  \\
    &     & RS Cap. 1    & 97.8\%  & 0.00 & 2.78 & 2.78 & 100\% & 1.00  & 100.0\% & 10,751 & 226   & 10,977 \\\bottomrule[1pt]
\end{tabular}%
}
\begin{tablenotes}
        \footnotesize
        \vspace{-6pt}
        \item Walk, wait, delay and detour are averaged across all passengers. Walk, wait, and delay are in minutes.
    \end{tablenotes}
\end{table}

Note, first, the benefits of on-demand flexibility versus fixed-line transit: by leveraging on-demand deviations, microtransit enables significant increases in demand coverage. Specifically, microtransit serves 2 to 3 times more passengers; in the three-hour case for example, this increase translates into an improvement in vehicle utilization from 30\% to 80\% on average with medium system capacity and from 18\% to 52\% with high system capacity. Unlike in the MiND-VRP, higher demand coverage comes with a slight increase in passenger walking and waiting, primarily due to an adverse selection effect---by serving passengers with pickup or drop-off locations further away from the reference lines, for example. Nonetheless, level of service remains comparable to fixed-line transit, with walking and waiting times around 2--3 minutes on average.

Next, results underscore the impact of demand consolidation: by relying on higher-capacity vehicles along reference lines, microtransit serves fewer passengers but travels much shorter distances than ride-sharing systems. As expected, ride-sharing results in higher demand coverage with no walking and short wait times. On the other hand, ride-sharing induces longer delays because of on-demand dispatches. Four-occupancy ride-pooling can also result in higher detours than microtransit, due to the negative externalities of door-to-door transportation---even with small-occupancy vehicles---and the comparative benefits of line regularization in microtransit. Moreover, the microtransit system travels much smaller (internal) distances by using higher-capacity vehicles.

At the aggregate level, microtransit induces strong system-wide improvements against all benchmarks. As compared to fixed-line transit, on-demand deviations increase distance traveled but this effect is more than compensated by the increase in demand coverage---leading to a decrease in  distance per passenger by a factor of 1.4 to 2.3. As compared to ride-sharing, microtransit decreases distance traveled by a much higher factor than the corresponding loss in demand coverage, leading to a smaller distance per passenger by a factor of 4--11  (resp. 3--6) as compared to single-occupancy ride-sharing (resp. four-occupancy ride-pooling). When accounting for the ``external'' distance from single-occupancy trips for all unserved passengers (assuming for instance that all unserved passengers take a taxi to their destination), microtransit reduces total distance from fixed-line transit by 5\%, 13\% and 15\% in the three-hour case with small, medium and high system capacity, respectively; it reduces total distance from four-occupancy ride-pooling by 39\%, 45\% and 48\%; and it reduces total distance from single-occupancy ride-sharing by 98\%, 157\% and 227\%.

These results confirm the potential of deviated fixed-route microtransit to improve demand coverage as compared to fixed-line transit---thanks to demand-responsive operations---and to improve demand consolidation as compared to ride-sharing---thanks to high-occupancy vehicles. These combined effects can induce strong reductions in distance traveled per passenger, which can ultimately contribute to creating more effective and more affordable mobility options and to mitigating the environmental footprint of urban mobility.

\section{Extension to incorporate transfers (MiND-Tr)}
\label{app:Tr}

\subsection{Modeling extension}
\label{subsec:MiND-Tr}

In the MiND-VRP, all passengers share the same destination at the end of each reference line (or the same origin at the start). Transfers are of little use in this setting, since all lines can drop off all passengers at their destination. We define the model with transfers, referred to as MiND-Tr, in a routing setting with a set $\calD$ of destinations. We assume that each line passes through one transfer point; that all passenger are picked up ahead of the transfer stations; and that their destinations are after the transfer stations. In other words, all passengers request transportation from an ``origin zone'' to a ``destination zone''; first-leg trips can pick up passengers near their origin and drop them off at a transfer station; and second-leg trips will pick them up at the transfer station and drop them off at their destination. In our experiments, we consider a similar use case as in the MiND-VRP, except that passengers' destinations are split between two airports (e.g., JFK and LGA). This setting preserves the key structural components of the second-stage model.

We define a set $\calF$ of transfer stations. Each reference trip $(\ell,t)$ transits through transfer station $f_\ell\in\calF$ at time $T^{\text{tr}}_{\ell t}$, and ends in destination $d_\ell \in \calD$. Each passenger $p\in\calP$ is bound for destination $d(p) \in \calD$. Passengers can be picked up on a line $\ell$ bound for destination $d(p)$ ($d_\ell=d(p)$); alternatively, they can be picked up by line $\ell$ and then transfer at transfer point $f_\ell \in \calF$ to another line $\ell'$ bound for destination $d(p)$ ($d_{\ell'}=d(p)$). We impose a maximum transfer time $W^{\text{tr}}$. We ensure that $f_\ell\in\calI_\ell$, i.e. the transfer station is included in the list of checkpoints; even when we allow to skip checkpoints ($K>0$), the vehicle must visit the transfer station to allow transfers.

Recall that, in the core MiND-VRP, $\calM_p \subseteq \calL \times \calT_{\ell}$ denotes the subset of reference trips that can serve request $p\in\calP$ within a tolerance $\alpha$ of their requested drop-off time. We extend this definition into a set of first-leg trips $\calM_p^1$ (from the passenger's pickup location to a transfer point) and second-leg trips $\calM_p^2$ (from a transfer point to the destination). Let $\calM_p^2(\ell, t)$ denote the set of second-leg trips $(\ell^\prime, t^\prime)$ that can serve passenger $p\in\calP$ within a tolerance $\alpha$ of their requested drop-off time, i.e., $\left|T_{\ell t}(\calI_{\ell}^{(I_{\ell})})-t^{\text{req}}_p\right|\leq\alpha$; and that are compatible with first-leg trip $(\ell, t)$, i.e., (i) reference lines $\ell$ and $\ell'$ transit through the same transfer station, i.e., $f_\ell=f_{\ell'}$; and (ii) reference trip $(\ell, t)$ arrives to the transfer point before reference trip $(\ell', t')$ departs and within the maximum transfer time, i.e., $T_{\ell t}^{\text{tr}} \leq T_{\ell^\prime t^\prime}^{\text{tr}}\leq T_{\ell t}^{\text{tr}} + W^{\text{tr}}$. If a reference trip $(\ell,t)$ can pick up and drop off passenger $p$ without a transfer (i.e., if $d_\ell=d(p)$), then the reference trip is thus included in both $\calM_p^1$ and $\calM_p^2(\ell, t)$.

The first-stage problem still selects reference trips (via binary variables $x_{\ell t}$) and assigns passengers to reference trips. The latter decisions are disaggregated into first-leg (pickup) and second-leg (dropoff) assignments, with binary variables $z^1_{\ell p s t}$ and $z^2_{\ell p s t}$ for passenger $p\in\calP$ and scenario $s\in\calS$. Each passenger will be assigned to first-leg and second-leg components (pickup and dropoff).

In the second stage, we define subpaths between checkpoints for the first leg of each line. The second leg travels directly between the line's transfer station and its destination, since all passengers on the line are then bound for the same destination; in particular, each line does not need to stop or deviate between the transfer station and the destination. Thus, we only define a set of pick-up subpaths $\calR^1_{\ell s t}$ for reference trip $(\ell,t)$ in scenario $s\in\calS$, between checkpoints. We denote by $\calP^1_{r(a)}\subseteq\calP$ the set of passengers picked up by subpath $r(a)\in\calR^1_{\ell s t}$. The parameters $g^1_a$ and $g^2_{\ell p s t}$ capture coverage, walk, wait, and in-vehicle time in the first leg, and lateness/earliness in the second leg.
\begin{align}\label{eq:arccost}
	g^{1}_{a} = &\begin{cases}
		\displaystyle \sum_{p \in \calP^1_{r(a)}}D_{ps} \left(\lambda \tau^{walk}_{r(a)p} + \mu \tau^{wait}_{r(a)p} + \sigma \frac{\tau_{r(a)p}^{travel}}{\tau_p^{dir}} - M  \right) \quad \forall a \in \bigcup_{r \in \calR_{\ell st}}\calA_r,\\
		0 \hspace{230pt} \forall a  \in\calA^v_{\ell st}.
	\end{cases}\\
	g^{2}_{\ell p s t} = &\displaystyle D_{ps} \left( \delta \frac{\tau^{late}_{\ell tp}}{\tau_p^{dir}} + \frac{\delta}{2} \frac{\tau^{early}_{\ell tp}}{\tau_p^{dir}} \right) \hspace{130pt} \forall s \in \calS, p \in \calP, (\ell,t) \in \calM_p^2
\end{align}

The MiND-Tr is formulated as follows.  We redefine binary variable $y^1_a$ to indicate whether a microtransit vehicle traverses subpath-based arc $a \in \calA_{\ell st}$ in the first-leg trip; and we introduce a new binary variable $y^2_{\ell p s t}$ to assign passenger $p$ in scenario $s$ to reference trip $(\ell,t)$ from the transfer station to the destination (this is a simpler representation than subpath-based arcs, exploiting the fact that the vehicle does not deviate after the transfer point). Equation~\eqref{trans:obj} minimizes expected costs. Constraints~\eqref{trans:assign} allow each passenger to be picked up at most once, and constraints~\eqref{trans:pudoconsist} ensure consistency between first- and second-leg trips. Constraints~\eqref{trans:load1}-\eqref{trans:load2} enforce the target load factor on the first-and second-leg trips. Constraints~\eqref{trans:pickup} and \eqref{trans:dropoff} ensure consistent assignments between the first stage and the second stage. Constraints~\eqref{trans:transfer} state that picked up passengers are dropped off, either by the same trip or by a compatible trip via a transfer. Finally, Constraints~\eqref{trans:capacity} enforce post-transfer capacity constraints for the trips.
\begin{align}
	\min  \quad&\sum_{\ell \in \calL} \sum_{t \in \calT_\ell} \left(h_\ell x_{\ell t} + \sum_{s \in \calS} \pi_s \sum_{a \in \calA_{\ell st}} g^1_a y^1_a + \sum_{s \in \calS} \pi_s \sum_{p \in \calP: (\ell,t) \in \calM_p^2} g^2_{\ell p st} y^2_{\ell p s t}\right)\label{trans:obj}\\
	\st\quad
	& \text{Fleet budget constraints (1)} \nonumber\\
	&\sum_{(\ell, t) \in \calM^1_p} z^1_{\ell pst} \leq 1,\qquad \forall p \in \calP, \forall s \in \calS \label{trans:assign}\\
	& z^1_{\ell_1 pst_1} = \sum_{(\ell_2, t_2) \in \calM^2_p(\ell_1 t_1)} z^2_{\ell_2 pst_2},\qquad \forall p \in \calP, (\ell_1, t_1) \in \calM^1_p, \forall s \in \calS \label{trans:pudoconsist}\\
        &(1-\kappa) C_\ell x_{\ell t} \leq \sum_{p \in \calP \, : \, (\ell, t) \in \calM^1_p} D_{ps}z^1_{\ell pst} \leq (1+\kappa)C_\ell x_{\ell t} \quad \forall (\ell, t) \in \calL \times \calT_\ell, \forall s \in \calS\label{trans:load1} \\
        &(1-\kappa) C_\ell x_{\ell t} \leq \sum_{p \in \calP \, : \, (\ell, t) \in \calM^2_p} D_{ps}z^2_{\ell pst} \leq (1+\kappa)C_\ell x_{\ell t} \quad \forall (\ell, t) \in \calL \times \calT_\ell, \forall s \in \calS\label{trans:load2} \\
	& \text{Subpath flow balance constraints (9)} \nonumber\\
	&\sum_{a \in \calA_{\ell st} \, : \, p \in \calP^1_{r(a)}} y^1_a \leq z^1_{\ell pst} \quad \forall
	s \in \calS, p \in \calP, (\ell,t) \in \calM^1_p \label{trans:pickup}\\
	& y^2_{\ell p s t} \leq z^2_{\ell pst} \quad \forall
	s \in \calS, p \in \calP, (\ell,t) \in \calM_p^2 \label{trans:dropoff}\\
	& \sum_{a \in \calA_{\ell_1 s t_1}:p\in \calP^1_{r(a)}}y^1_a = \sum_{(\ell_2, t_2) \in \calM_p^2(\ell_1,t_1)} y^2_{\ell_2 p s t_2} \quad \forall s \in \calS, p \in \calP, (\ell_1,t_1) \in \calM^1_p,  \label{trans:transfer} \\
	& \sum_{p \in \calP: (\ell,t) \in \calM_p^2} y^2_{\ell p s t} \leq C_\ell x_{\ell t} \quad \forall s \in \calS, \ell \in \calL, t \in \calT_\ell \label{trans:capacity} \\
	&\bx, \by, \bz^1, \bz^2, \bf \text{ binary} \label{trans:domain}
\end{align}~

\subsection{Algorithmic extension}
\label{app:TRalg}

\paragraph{Benders decomposition.}

Due to Equation~\eqref{trans:transfer}, the Benders subproblem is no longer separable by reference trip $(\ell,t) \in \calL \times \calT_\ell$; indeed, operations on several lines need to be coordinated so that, after transfers, the load of the vehicles will not exceed their capacity.

Let $\boldsymbol{\psi}$ denote the dual variable associated to the flow balance constraints, $\bgamma^1$ the one associated to the pickup trip linking constraints (Equation~\eqref{trans:pickup}), $\bgamma^2$ the one associated to the dropoff trip linking constraints (Equation~\eqref{trans:dropoff}), $\boldsymbol{\eta}$ the one associated to the transfer consistency constraints (Equation~\eqref{trans:transfer}), and $\bnu$ the dual variable associated to the final leg capacity constraints (Equation~\eqref{trans:capacity}). The Benders dual subproblem for scenario $s$ becomes:
\begin{align}
    \max \quad & \sum_{\ell \in \calL} \sum_{t \in \calT_\ell} \left(x_{\ell t} \cdot (\psi_{\ell s t u_{\ell s t}} -\psi_{\ell s t v_{\ell s t}}) - C_\ell \nu_{\ell s t} \right)\\
	&\qquad - \sum_{p \in \calP} \left(\sum_{(\ell, t) \in \calM^1_p} z^1_{\ell p s t} \cdot \gamma^1_{\ell s t p}
	+ \sum_{(\ell,t) \in \calM_p^2} z^2_{\ell p s t} \cdot \gamma^2_{\ell s t p} \right) \\
	\text{s.t. } \quad &\psi_{\ell s t n} - \psi_{\ell s t m} -  \sum_{p \in \calP^1_{r(a)} } (\gamma^1_{\ell s t p}  - \eta_{\ell s t p}) \leq g^{1}_a  \quad \forall  \ell \in \calL, t\in \calT_\ell, a=(n,m) \in \mathcal{A}_{\ell st}\\
	& - \gamma^2_{\ell s t p}  - \sum_{(\ell_2, t_2) \in \calM_p^2(\ell, t)} \eta_{\ell_2 s t_2 p} - \nu_{\ell st} \leq g^{2}_{\ell p s t} \quad \forall p \in \calP, (\ell,t) \in \calM_p^2  \\
	&\bgamma^1,\bgamma^2,\boldsymbol{\eta},\boldsymbol{\nu}\geq\bo
\end{align}

The corresponding Benders optimality cut becomes:
\begin{align}
	\theta_s \geq 
	 \sum_{\ell \in \calL} \sum_{t \in \calT_\ell} \left(x_{\ell t} \cdot (\psi_{\ell s t u_{\ell s t}} -\psi_{\ell s t v_{\ell s t}}) - C_\ell \nu_{\ell st} \right)
	- \sum_{p \in \calP}  \left(\sum_{(\ell, t) \in \calM^1_p} z^1_{\ell p s t} \cdot \gamma^1_{\ell s t p}
	+ \sum_{(\ell,t) \in \calM_p^2} z^2_{\ell p s t} \cdot \gamma^2_{\ell s t p} \right)
\end{align}

\paragraph{Column generation.}

The restricted Benders subproblem for each scenario $s$ is still obtained by restricting the decisions to a subset of arc-based variables in $\calA'_{\ell st}$ for each reference trip $(\ell, t)$:
\begin{align}
    \min_{\by^1, \by^2 \geq\bo} &\quad \sum_{a \in \calA'_{\ell st}} g_a^{1} y^1_a  + \sum_{p \in \calP: (\ell,t) \in \calM_p^2} g^{2}_{\ell p s t} y^2_{\ell p s t}  \label{RMPTr:obj} \\
    \st\quad&\sum_{m:(n,m) \in \calA'_{\ell st}} 
    y^1_{(n,m)} - \sum_{m:(m, n) \in \calA'_{\ell st}} 
    y^1_{(m, n)} = \begin{cases}
    x_{\ell t} &\text{if } n = u_{\ell st}, \\
    -x_{\ell t} &\text{if } n = v_{\ell st}, \\
    0 &\text{ otherwise,}
    \end{cases} \nonumber \\
    & \hspace{230pt} \forall \ell \in \calL, t \in \calT_\ell, n \in \calV_{\ell st}\label{RMPTr:2Sflow} \\
    &\sum_{a \in \calA^\prime_{\ell st} \, : \, p \in \calP^1_{r(a)}} y^1_a \leq z^1_{\ell pst} \quad \forall
	p \in \calP, (\ell,t) \in \calM^1_p \label{RMPTr:pickup}\\
	& y^2_{\ell p s t} \leq z^2_{\ell pst} \quad \forall
	p \in \calP, (\ell,t) \in \calM_p^2 \label{RMPTr:dropoff}\\
	&\sum_{a \in \calA^\prime_{\ell_1 st_1} \, : \, p \in \calP^1_{r(a)}} y^1_a = \sum_{(\ell_2, t_2) \in \calM_p^2(\ell_1,t_1)} y^2_{\ell_2 p s t_2} \quad \forall  p \in \calP, (\ell_1,t_1) \in \calM^1_p,  \label{RMPTr:transfer} \\
	& \sum_{p \in \calP: (\ell, t) \in \calM^2_p} y^2_{\ell p s t} \leq C_\ell x_{\ell t} \quad \forall \ell \in \calL, t \in \calT_\ell \label{RMPTr:capacity}  
\end{align}

The pricing problem seeks a subpath-based arc of minimum reduced cost:
$$\min_{\ell \in \calL, t\in \calT_\ell}\min_{a=(n,m) \in \mathcal{A}_{\ell st}}\left(g^{1}_a+\psi_{\ell s t m}-\psi_{\ell s t n}+\sum_{p \in \calP^1_{r(a)} } (\gamma^1_{\ell s t p}- \eta_{\ell s t p}) \right)$$

In the pricing problem, we define the level-of-service parameter $d^1_{mp}$ for passenger $p$ in location $m$ only over the level-of-service components associated with the first-leg pickups, as the second-leg dropoff components are contained in $g^2_{\ell p s t}$. We denote by $\calP^1_m$ the set of passengers that can be picked up in node $m\in\calU^{uv}_{\ell st}$, and by $T^1_{mp}$ the pickup time of passenger $p\in\calP^1_m$. We then define:
\begin{align*}
	d^1_{mp}&= D_{ps} \left( \lambda \tau^{\text{walk}}_{mp} + \mu \tau^{\text{wait}}_{mp}-\sigma\frac{T^1_{mp}}{\tau^{\text{dir}}_p} - M\right) + \gamma^1_{\ell s t p} - \eta_{\ell s t p},&& \quad \forall \ell \in \calL, t \in \calT_\ell, m \in \calU^{uv}_{\ell st}, p \in\calP^1_{m}
\end{align*}
 
Because the transfers do not impact subpaths before the transfer point, the pricing problem is identical to the one in the MiND-VRP (Equations~\eqref{eq:PPobj}--\eqref{eq:PPdomain}) apart from dual information in the objective function. Recall the MiND-VRP pricing problem has a binary variable $w_{mp}$ to indicate whether passenger $p$ is picked up at location $m$. The objective for the MiND-Tr pricing problem for trip $(\ell, t)$ in scenario $s\in\calS$ between checkpoints $u$ and $v$ is given as:
\begin{align}
    \min  \quad  &\sum_{m \in \calU^{uv}_{\ell st}} \sum_{p \in \calP_{m}} d^1_{mp} w_{mp}  + \psi_{\ell, s, t, end(a)} - \psi_{\ell, s, t, start(a)}\label{transPP:obj}
\end{align}

\paragraph{Label setting algorithm.}

The structure of the pricing problem remains unchanged, so the label-setting algorithm is identical to the one in the MiND-VRP.

\subsection{Experimental results}
\label{app:TRcase}

We adapt the MiND-VRP case study to include passenger requests to LaGuardia Airport (LGA) and John F. Kennedy International Airport (JFK) between 6 and 9 am. Passenger origins are located in Manhattan and lines travel to the airport via one of four bridges. We set up instances with 10, 20, and 40 candidate lines, half of which are bound for LGA and half of which are bound for JFK. We locate the transfer stations at the last stop before each bridge.

\begin{figure}[h]
	\centering
	\includegraphics[width=10cm]{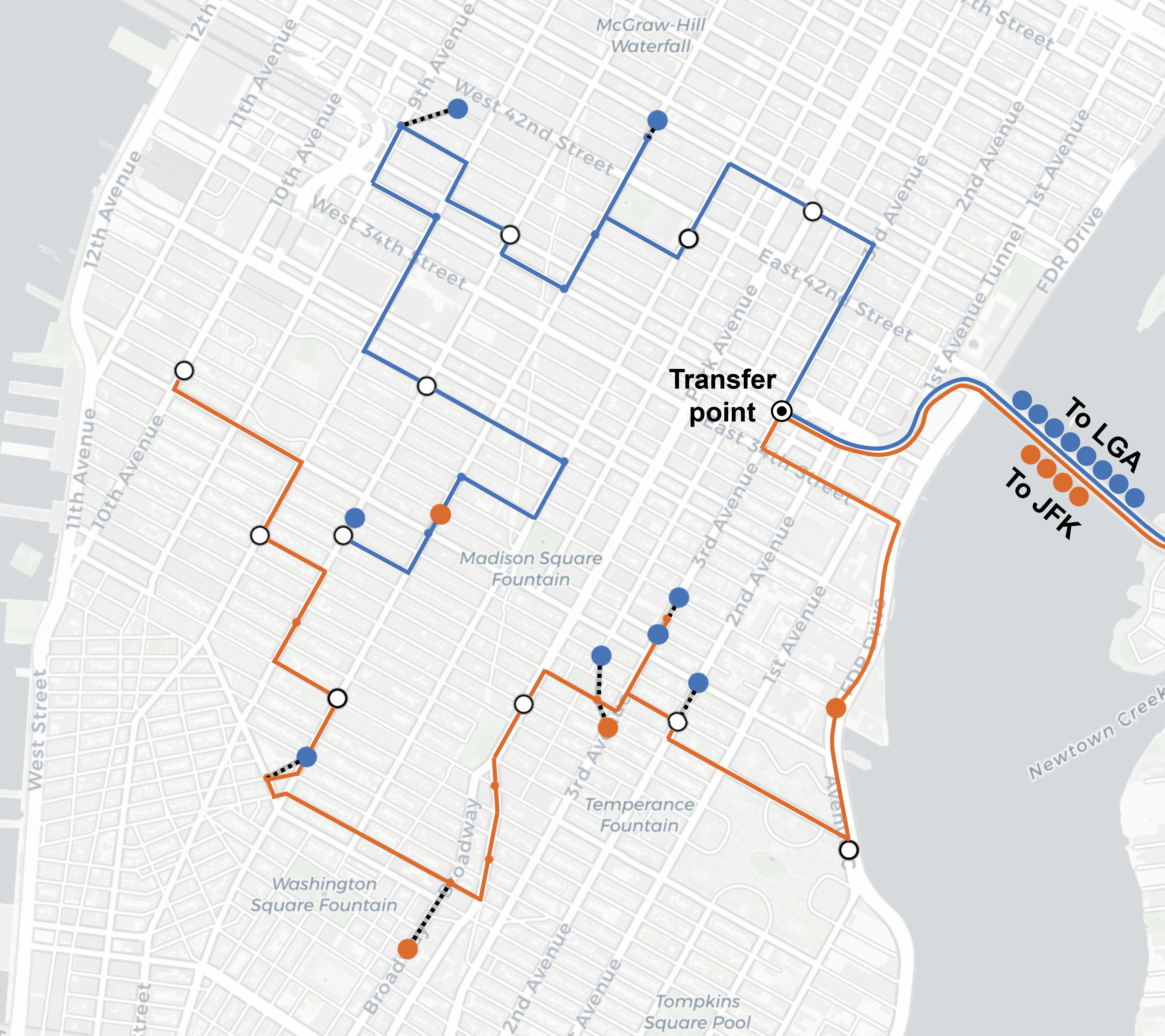} 
	\caption{A line bound for LaGuardia (blue) and a line bound for JFK (orange) with a transfer location just before the Queensboro bridge. Passengers are shown as circles, shaded to match the color of their destination airport, with their walking path shown as black dotted lines. } \label{fig:transfersol}
\end{figure}

Figure~\ref{fig:transfersol} shows a sample solution for two candidate lines. The line bound for LGA (in blue) picks up four passengers, with three headed to LGA (blue circles) and one to JFK (orange circle). The line bound for JFK picks up eight passengers, five headed to LGA and three headed to JFK. The two lines meet at a transfer point before the Queens-Midtown Tunnel where passengers can transfer to the vehicle bound for their destination.

Table~\ref{tab:transfers} compares the deviated fixed-route microtransit solution and the fixed-line transit benchmark, each with and without transfers. As expected, transfers increase problem complexity. The MiND-Tr scales up to 10 candidate lines and 5 scenarios but leaves a larger optimality gap in larger instances. This mainly comes from the fact that the Benders subproblem is no longer separable by reference trip, resulting in longer computational requirements and denser Benders cuts. Nonetheless, transfers can improve solution performance even in larger instances that the MiND-Tr does not solve to optimality, by aggregating demand across destinations. Specifically, 31-45\% of microtransit passengers transfer between lines. Transfers increase coverage by 0-3\% for transit and 1-3\% for microtransit. Accommodating these additional passengers comes at a minor cost in terms of wait times, detours, and delays, as well as a more significant 1-4 minute increase in walking times. Nonetheless, the net benefit of transfers is positive, resulting in 3--7\% cost reductions.

\begin{table}[h!]
\renewcommand*{\arraystretch}{1.0}
\centering
\small
\caption{Comparison of fixed-route transit vs. microtransit with and without transfers}
\label{tab:transfers}
\resizebox{1.\textwidth}{!}{%
    \begin{tabular}{llcrccccccccccc}
        \toprule
        &&&&& \multicolumn{5}{c}{Average level of service} & \multicolumn{2}{c}{Vehicle utilization} & \multicolumn{3}{c}{Distance traveled (km)}\\ 
        \cmidrule(lr){6-10} \cmidrule(lr){11-12} \cmidrule(lr){13-15} 
        $|\calL|$ & Mode & Sol. & Gap & Tr. & Coverage & Walk & Wait & Detour & Delay & Absolute & Relative & Internal & External & Total\\ 
        \midrule
        10 & Transit & 163.5 & \textbf{0.0} & --- & 10.8\% & 0 & 3.13 & 132\% & 0.04 & 1.48 & 38.9\% & 59 & 4,180 & 4,239 \\
         & Transit, transfers & 163.5 & \textbf{0.0} & 16\% & 10.8\% & 0 & 3.08 & 137\% & 0.09 & 1.48 & 38.9\% & 56 & 4,180 & 4,236 \\
         & MT & 104.4 & \textbf{0.0} & --- & 26.5\% & 3.08 & 2.99 & 134\% & 0.03 & 3.64 & 95.8\% & 75 & 4,031 & 4,106 \\
         & MT, transfers & 100 & \textbf{0.0} & 41\% & 27.7\% & 7.37 & 3.27 & 136\% & 0.1 & 3.8 & 100\% & 76 & 4,023 & 4,099 \\
        \midrule
        20 & Transit & 157.3 & \textbf{0.0} & --- & 15.7\% & 0 & 2.18 & 127\% & 0 & 2.2 & 44.4\% & 122 & 4,985 & 5,107 \\
         & Transit, transfers & 152.3 & \textbf{0.0} & 34\% & 17.4\% & 0 & 2.18 & 127\% & 0.01 & 2.44 & 49.2\% & 113 & 4,970 & 5,083 \\
         & MT & 107.9 & \textbf{0.0} & --- & 32.7\% & 1.13 & 2.25 & 125\% & 0.01 & 4.58 & 92.3\% & 153 & 4,832 & 4,985 \\
         & MT, transfers & 100 & 14.5 & 31\% & 35.4\% & 3.77 & 2.32 & 127\% & 0.02 & 4.96 & 100\% & 154 & 4,805 & 4,959 \\
        \midrule
        40 & Transit & 147.1 & \textbf{0.0} & --- & 24.4\% & 1.15 & 2.2 & 120\% & 0 & 2.43 & 54.4\% & 236 & 3,626 & 3,862 \\
         & Transit, transfers & 140.2 & 3.1 & 44\% & 27.5\% & 1.54 & 2.2 & 124\% & 0.05 & 2.73 & 61.1\% & 242 & 3,552 & 3,794 \\
         & MT & 103.1 & 1.7 & --- & 43.7\% & 1.29 & 2.08 & 117\% & -0.01 & 4.34 & 97.1\% & 296 & 3,120 & 3,416 \\
         & MT, transfers & 100 & 27.9 & 45\% & 45\% & 2.5 & 2.16 & 120\% & 0.04 & 4.47 & 100\% & 301 & 3,093 & 3,394 \\
        \midrule
        \bottomrule
    \end{tabular}
}\begin{tablenotes}\footnotesize
    \item ``Sol.'': total expected cost, normalized to the microtransit solution with transfers.
    \item Walk, wait, delay and detour are averaged across all passengers. Walk, wait, and delay are in minutes.
\end{tablenotes}
\end{table}

Most importantly, these new results strengthen our main takeaways. Indeed, the microtransit system still provides improvements over the fixed-route transit benchmark in the presence of transfers, consistent with the MiND-VRP and MiND-DAR results. In fact, the performance of the microtransit system \textit{without transfers} also outperforms the one of the transit benchmark \textit{with transfers}. In other words, the benefits of transfers to the transit system do not cover the significant service gap between transit and microtransit. This can be explained, in part, by the spatiotemporal coordination required across reference trips to take advantage of transfers, which restricts the model's flexibility when selecting reference trips. In comparison, these results reinforce the benefits of the routing flexibility from the microtransit system: microtransit provides a 15-17\% increase in coverage and a 44-59\% decrease in total costs without transfers; and then, transfers enable further cost reductions of 3--7\%. These benefits translate into a 3--10\% reduction in the ``total distance'' metric in the microtransit solution as compared to transit, with or without transfer, leading to reductions in operating costs per passenger and the environmental footprint of mobility systems.
\section{Additional results}
\label{app:results}

\subsection{Additional results on the benefits of the subpath-based model}
\label{app:subpath}

Recall that our results in Section~\ref{subsec:method} showed the benefits of our subpath-based formulation in view of the overall two-stage stochastic integer optimization problem, as compared to the compact, segment-based and path-based benchmarks. In this appendix, we compare the four formulations for the second-stage problem alone---that is, on the capacitated vehicle routing problem with time windows. We consider instances with 5, 8, and 10 candidate lines, 5 scenarios, and 1-, 2-, and 3-hours horizons. We fix the optimal first-stage decisions, and then solve the corresponding subproblems using the direct implementation of each formulation. By design, the subpath-based and path-based models are solved via exhaustive enumeration. Table~\ref{tab:secondstage} reports pre-processing times for arc enumeration, solution times, and optimality gaps.

\begin{table}[h!]
\renewcommand*{\arraystretch}{1.0}
\centering
\footnotesize
\caption{Comparison of the four models for the second-stage problem. }
\label{tab:secondstage}
\resizebox{1.\textwidth}{!}{%
\begin{tabular}{cccccccccccccc}
\toprule
 &   & \multicolumn{3}{c}{Compact} & \multicolumn{3}{c}{Segment} & \multicolumn{3}{c}{Path} & \multicolumn{3}{c}{Subpath} \\ 
 \cmidrule(lr){3-5}\cmidrule(lr){6-8}\cmidrule(lr){9-11} \cmidrule(lr){12-14} 
&  &   Enum. & Solve &  & Enum. & Solve &  & Enum. & Solve &  & Enum. & Solve &  \\ 
$|\calL|$ & Hor.  & CPU(s) & CPU(s) & Gap & CPU (s) & CPU(s) & Gap & CPU (s) & CPU(s) & Gap & CPU (s) & CPU(s) & Gap \\ 
 \midrule
5 & 60 & --- & $>$600 & 3.8 & 68 & 1.05 & 0.0 & 76 & 0.05 & 0.0 & 19 & $<$0.01 & 0.0 \\
 & 120 & --- & $>$600 & 5.7 & 283 & 1.03 & 0.0 & 539 & 0.05 & 0.0 & 252 & $<$0.01 & 0.0 \\
 & 180 & --- & $>$600 & --- & 565 & 1.01 & 0.0 & 602 & 0.05 & 0.0 & 296 & $<$0.01 & 0.0 \\
   \midrule 
8 & 60 & --- & $>$600 & 6.5 & 160 & 1.10 & 0.0 & 738 & 0.07 & 0.0 & 37 & $<$0.01 & 0.0 \\
 & 120 & --- & $>$600 & --- & 616 & 1.04 & 0.0 & 2,735 & 0.11 & 0.0 & 381 & $<$0.01 & 0.0 \\
 & 180 & --- & $>$600 & --- & 1,323 & 1.08 & 0.0 & 3,027 & 0.08 & 0.0 & 495 & $<$0.01 & 0.0 \\
  \midrule 
10 & 60 & --- & $>$600 & 6.1 & 228 & 0.97 & 0.0 & 919 & 0.07 & 0.0 & 47 & $<$0.01 & 0.0 \\
 & 120 & --- & $>$600 & --- & 1,004 & 1.10 & 0.0 & 4,072 & 0.12 & 0.0 & 481 & $<$0.01 & 0.0 \\
 & 180 & --- & $>$600 & --- & 2,232 & 1.09 & 0.0 & 10,800+ & --- & --- & 555 & $<$0.01 & 0.0 \\
\bottomrule
\end{tabular}
}
\begin{tablenotes}
        \item \scriptsize Time limit: 3 hours for enumeration, 10 minutes per subproblem.
        \item \scriptsize Enum CPU: time to enumerate the decision variables (segments, subpaths, paths), Solve CPU: average time to solve single subproblem, Gap: average optimality gap at the time limit and  ``---'' if no feasible solution was found.
    \end{tablenotes}
\end{table}

Even in the smallest instances, the compact formulation takes over than ten minutes to solve for each subproblem. This stems from the weak relaxation leading to limited scalability of off-the-shelf branch-and-cut algorithms. This shows that, in our MiND problem, even the second-stage is challenging to solve by itself. All network-based formulations are comparatively more efficient, terminating in fractions of a second. Still, the subpath-based formulation solves orders of magnitude faster than the segment-based one. Although the segment-based model retains a polynomial size, its slightly weaker relaxation and, most importantly, its large number of variables in the dense time-station-load network create significant computational complexities. In the context of the full MiND problem, these differences in solution times have a significant impact because hundreds of subproblems are solved at each Benders iteration. Finally, the subpath-based formulation involves much fewer arcs than the path-based one, thereby resulting in significant speedups in pre-processing.

The path-based model could also be solved via column generation---like our subpath-based model. However, we have found that subpath-based column generation also converges orders of magnitude faster than path-based column generation (0.03 second vs. 4.1 seconds on average per pricing problem). In other words, the reductions in the number of variables and in pre-processing times identified in Table~\ref{tab:secondstage} lead to similar reductions in computational times for the subpath-based vs. path-based column generation. This stems from the fact that the label-setting algorithm is applied between checkpoints in subpath-based column generation, as opposed to from the beginning to the end of a reference trip in path-based column generation---leading to an exponential decrease in the number of paths and passenger combinations.

These results underscore that the different formulations face different bottlenecks: the branch-and-cut structure in the compact formulation due to its weak linear relaxation, the linear relaxation in the segment-based formulation due to its large size, and the generation of arc variables in the path-based and subpath-based formulations either through exhaustive enumeration at the pre-processing stage or through column generation. In larger instances, subpath or path enumeration becomes intractable, which motivated our double-decomposition algorithm.

\subsection{Sensitivity analyses}
\label{app:sensitivity}

\subsubsection*{Sensitivity to the number of skipped checkpoints.}

We perform additional sensitivity analyses to study the impact of $K$ (that is, the number of consecutive checkpoints that can be skipped) on computational times and solution quality. Results are reported in Table~\ref{tab:senK}. We consider here a small instance with 5 candidate lines, 5 scenarios, and a one-hour horizon, as larger instances become intractable with $K=3$. As expected, the computational times increase significantly as $K$ becomes larger because of the exponential growth in the number of subpaths. In this setting, the algorithm can solve all instances with $K=0$ within a minute or so, but can take up to 10-25 minutes with $K=1$, over an hour with $K=2$, and reaches the 3-hour time limit with $K=3$. Even with larger values of $K$, the majority of the subpaths still visit consecutive checkpoints (even without the non-linear penalties associated with longer subpaths discussed in the first point above). This is primarily driven by the penalties on in-vehicle time and arrival delays, which discourage longer subpaths. In turn, a value of $K=1$ achieves most of the benefits obtained with the larger values of $K$ (notably, by increasing coverage by 30--40\%). In comparison, the marginal benefits with $K=2$ become smaller---although some higher-coverage solutions are obtained with $K=3$. Altogether, these results strengthen one of our practical takeaways that a limited extent of flexibility can provide strong benefits (Table~\ref{T:mtLOS}).

\begin{table}[h!]
\renewcommand*{\arraystretch}{1.0}
\centering
\small
\caption{Results for MiND-VRP with varying $K$ and $\Delta$.}
\label{tab:senK}
\resizebox{1.\textwidth}{!}{%
\begin{tabular}{ccccccccccccc}
\toprule
 &  &  & & \multicolumn{5}{c}{Passenger metrics} & \multicolumn{4}{c}{Subpaths skipping $k$ checkpoints} \\
\cmidrule(lr){5-9} \cmidrule(lr){10-13} 
K & $\Delta$ & Solution & CPU (s) & Pickups & Walk & Wait & Detour & Delay & $k=0$ & $k=1$ & $k=2$ & $k=3$\\
\midrule
0 & 600 & 138.5 & 14 & 135 & 3.11 & 2.68 & 140\% & 0.01 & 100\% & --- & --- & --- \\
 & 1200 & 136.6 & 48 & 139 & 3.02 & 2.76 & 141\% & 0.02 & 100\% & --- & --- & --- \\
 & 1800 & 136.6 & 51 & 139 & 3.02 & 2.77 & 142\% & 0.02 & 100\% & --- & --- & --- \\
 & 2400 & 136.6 & 86 & 139 & 3.02 & 2.76 & 142\% & 0.02 & 100\% & --- & --- & --- \\
\midrule
1 & 600 & 117.9 & 67 & 178 & 1.57 & 1.77 & 137\% & 0.02  & 62\% & 38\% & --- & --- \\
 & 1200 & 111.0 & 376 & 193 & 4.26 & 2.04 & 137\% & 0.02 & 56\% & 44\% & --- & --- \\
 & 1800 & 110.9 & 819 & 193 & 3.48 & 1.92 & 137\% & 0.02  & 58\% & 42\% & --- & --- \\
 & 2400 & 111.0 & 1,544 & 193 & 3.48 & 2.10 & 137\% & 0.02  & 57\% & 43\% & --- & --- \\
\midrule
2 & 600 & 116.0 & 157 & 182 & 1.54 & 1.68 & 138\% & 0.02 & 62\% & 22\% & 14\% & --- \\
 & 1200 & 107.5 & 1,223 & 200 & 1.21 & 1.70 & 137\% & 0.02 & 67\% & 16\% & 17\% & --- \\
 & 1800 & 106.1 & 3,970 & 203 & 1.19 & 1.67 & 137\% & 0.02 & 63\% & 22\% & 14\% & --- \\
 & 2400 & 105.7 & 7,724 & 204 & 1.19 & 1.84 & 137\% & 0.02 & 65\% & 20\% & 15\% & --- \\
\midrule
3 & 600 & 113.1 & 260 & 188 & 2.98 & 1.53 & 139\% & 0.01 & 67\% & 8\% & 6\% & 19\% \\
 & 1200 & 101.9 & 2,843 & 212 & 1.89 & 1.71 & 139\% & 0.02 & 67\% & 7\% & 7\% & 19\% \\
 & 1800 & 100.0 & 9,370 & 217 & 2.49 & 2.10 & 149\% & 0.02 & 70\% & 3\% & 3\% & 23\% \\
 & 2400 & 116.0 & 10,800+ & 184 & 1.52 & 2.78 & 156\% & 0.02 & 69\% & 12\% & 3\% & 15\% \\
\bottomrule
\end{tabular}
}
\end{table}

\subsubsection*{Sensitivity to the specification of the objective function.}

Recall that the model formulation trades off several underlying goals: operating costs, demand coverage, and passenger level of service (itself comprising walking times, waiting times, in-vehicle travel times, and delay at destination). Table~\ref{tab:multiobj} reports sensitivity analysis results to explore the trade-off between these components, both for microtransit and fixed-route transit. The first eight rows vary the four level-of-service weights together, ranging from a setting where demand coverage is prioritized (small $\lambda$, $\mu$, $\sigma$, and $\delta$) to a setting where passenger level of service is prioritized (large $\lambda$, $\mu$, $\sigma$, and $\delta$). The next rows vary one weight parameter at a time, reflecting different balances between level-of-service components. We consider very large variations in these weight parameters to analyze outcomes where one objective component is strongly prioritized over the others.

\begin{table}[h!]
\renewcommand*{\arraystretch}{1.0}
\centering
\small
\caption{Sensitivity to weight parameters (50 candidate lines, 5 scenarios, 2 hours).}
\label{tab:multiobj}
\resizebox{1.\textwidth}{!}{%
\begin{tabular}{cccccccccccccccc}
\toprule
First stage & \multicolumn{4}{c}{Second Stage} &  &  &  & \multicolumn{5}{c}{Avg. level of service} & \multicolumn{3}{c}{Distance (km)} \\ 
 \cmidrule(lr){1-1} \cmidrule(lr){2-5} \cmidrule(lr){9-13} \cmidrule(lr){14-16}
Line cost & $\lambda$ & $\mu$ & $\sigma$ & $\delta$ & Mode & Lines & Trips & Coverage & Walk & Wait & Detour & Delay & Internal & External & Total \\ 
\midrule
1 & \textbf{0.1} & \textbf{0.1} & \textbf{0.1} & \textbf{0.1} & MT & 21 & 42 & 36.6\% & 1.56 & 5.75 & 127\% & 0.01 & 673 & 12,061 & 12,734 \\ 
 &  &  &  &  & Transit & 18 & 44 & 35.9\% & 2.16 & 6.88 & 121\% & 0.00 & 622 & 12,178 & 12,800 \\ 
 \cmidrule(lr){2-16}
 & 1 & 1 & 1 & 1 & MT & 18 & 42 & 36.6\% & 1.54 & 5.74 & 128\% & 0.01 & 681 & 12,039 & 12,720 \\ 
 &  &  &  &  & Transit & 19 & 44 & 35.9\% & 2.11 & 6.82 & 121\% & 0.01 & 623 & 12,170 & 12,793 \\ 
 \cmidrule(lr){2-16}
 & \textbf{10} & \textbf{10} & \textbf{10} & \textbf{10} & MT & 18 & 42 & 36.5\% & 1.30 & 5.23 & 128\% & 0.01 & 687 & 12,049 & 12,736 \\ 
 &  &  &  &  & Transit & 19 & 44 & 35.7\% & 2.14 & 6.73 & 122\% & 0.01 & 626 & 12,188 & 12,814 \\ 
 \cmidrule(lr){2-16}
 & \textbf{100} & \textbf{100} & \textbf{100} & \textbf{100} & MT & 11 & 37 & 10.2\% & 0.00 & 0.86 & 135\% & 0.00 & 512 & 15,518 & 16,030 \\ 
 &  &  &  &  & Transit & 6 & 31 & 2.8\% & 0.00 & 0.74 & 142\% & -0.01 & 196 & 16,476 & 16,672 \\ 
 \midrule
1 & \textbf{10} & 1 & 1 & 1 & MT & 20 & 43 & 36.8\% & 0.67 & 5.28 & 126\% & 0.01 & 688 & 12,041 & 12,729 \\ 
 &  &  &  &  & Transit & 20 & 44 & 35.7\% & 1.78 & 6.64 & 123\% & 0.01 & 622 & 12,196 & 12,818 \\ 
  \cmidrule(lr){2-16}
 & \textbf{100} & 1 & 1 & 1 & MT & 21 & 43 & 30.3\% & 0.00 & 4.09 & 127\% & 0.01 & 664 & 12,917 & 13,581 \\ 
 &  &  &  &  & Transit & 16 & 42 & 22.4\% & 0.00 & 4.77 & 117\% & 0.00 & 572 & 13,909 & 14,481 \\ 
 \midrule
1 & 1 & \textbf{10} & 1 & 1 & MT & 16 & 42 & 37.1\% & 2.38 & 6.77 & 130\% & 0.01 & 681 & 11,973 & 12,654 \\ 
 &  &  &  &  & Transit & 19 & 44 & 35.7\% & 2.41 & 7.19 & 123\% & 0.00 & 625 & 12,201 & 12,826 \\ 
  \cmidrule(lr){2-16}
 & 1 & \textbf{100} & 1 & 1 & MT & 16 & 42 & 29.9\% & 3.17 & 7.45 & 120\% & 0.01 & 641 & 12,972 & 13,613 \\ 
 &  &  &  &  & Transit & 16 & 42 & 25.0\% & 3.20 & 7.75 & 121\% & 0.01 & 577 & 13,588 & 14,165 \\ 
 \midrule
1 & 1 & 1 & \textbf{10} & 1 & MT & 22 & 42 & 36.0\% & 1.47 & 5.59 & 127\% & 0.01 & 678 & 12,149 & 12,827 \\ 
 &  &  &  &  & Transit & 18 & 44 & 35.7\% & 2.16 & 6.87 & 121\% & 0.00 & 618 & 12,191 & 12,809 \\ 
  \cmidrule(lr){2-16}
 & 1 & 1 & \textbf{100} & 1 & MT & 17 & 41 & 36.3\% & 1.57 & 5.72 & 128\% & 0.01 & 671 & 12,101 & 12,772 \\ 
 &  &  &  &  & Transit & 20 & 44 & 35.9\% & 2.15 & 6.85 & 120\% & 0.01 & 619 & 12,172 & 12,791 \\ 
 \midrule
1 & 1 & 1 & 1 & \textbf{10} & MT & 19 & 42 & 36.2\% & 1.50 & 5.61 & 128\% & 0.01 & 687 & 12,109 & 12,796 \\ 
 &  &  &  &  & Transit & 20 & 44 & 35.8\% & 2.16 & 6.89 & 121\% & 0.00 & 624 & 12,181 & 12,805 \\ 
 \cmidrule(lr){2-16}
 & 1 & 1 & 1 & \textbf{100} & MT & 20 & 42 & 36.4\% & 1.53 & 5.70 & 128\% & 0.01 & 679 & 12,084 & 12,763 \\ 
 &  &  &  &  & Transit & 20 & 44 & 35.8\% & 2.11 & 6.81 & 122\% & 0.00 & 624 & 12,172 & 12,796 \\ 
\bottomrule
\end{tabular}
}
\end{table}

The results are highly robust to the choice of the underlying weight parameters. The general structure of the optimal solution remains stable across virtually all instances. Exceptions arise when the weight of the level-of-service costs (especially, of waiting costs) is multiplied by a factor of 100, in which case the model prioritizes passengers that can be served at the requested time at the expense of not covering other passengers. Otherwise, the solution remains mostly unchanged; in particular, the number of pickups varies by less than 1\%, with similar detours and delays, and variations in walking and waiting times on the order of 0--2 minutes on average.

Our main takeaways are therefore robust to the objective specification. In particular, the benefits of microtransit over fixed-route transit hold across all parameter values. In most settings, the microtransit solution increases demand coverage and reduces walk and wait times as compared to fixed-route transit, at limited costs in terms of detours and distance traveled. In turn, the microtransit solution retains efficiency benefits (higher coverage and higher level of service) as well as sustainability benefits (lower total distance, after accounting for the external distance of non-covered passengers). At the extreme where passenger level is heavily prioritized and where, in particular, walking is heavily penalized, both solutions involve virtually no walking for passengers; but then, the flexibility of microtransit enables large increases in demand coverage as compared to the transit solution. Altogether, the finding that microtransit can result in win-win outcomes does not depend on idiosyncratic choices of the weight parameters but holds across a range of trade-off choices between demand coverage and passenger level of service.

\end{document}